\documentclass[12pt]{amsart}
\usepackage{latexsym}
\usepackage{indentfirst}
\usepackage{amssymb,amsfonts,amsmath,amsthm}
\usepackage{graphicx}
\usepackage{mathrsfs}
\usepackage{array}
\usepackage{color}
\usepackage{epstopdf}
\usepackage[all]{xy}
\usepackage{verbatim}
\usepackage{url}
\usepackage{tikz-cd}
\usepackage{relsize}

\usepackage[linktocpage=true, hidelinks, colorlinks, linkcolor=blue, citecolor=red, bookmarksopen=true, urlcolor=brown]{hyperref}

\usepackage[a4paper,bindingoffset=0.2in,           left=0.5in,right=0.51in,top=0.6in,bottom=0.5in,            footskip=.25in]{geometry}

\usepackage[bbgreekl]{mathbbol}
\DeclareSymbolFontAlphabet{\mathbb}{AMSb} 
\DeclareSymbolFontAlphabet{\mathbbl}{bbold}

\newtheorem{thm}{Theorem}[section]
\newtheorem{prop}[thm]{Proposition}
\newtheorem{lem}[thm]{Lemma}

\newtheorem{cor}[thm]{Corollary}

\numberwithin{equation}{section}

\newtheorem{theorem}[thm]{Theorem}

\theoremstyle{definition}
\newtheorem{definition}[thm]{Definition}
\newtheorem{hypothesis}[thm]{Hypothesis}
\newtheorem{construction}[thm]{Construction}

\newtheorem{convention}[thm]{Convention}

\newtheorem{dfn}[thm]{Definition}
\newtheorem{exam}[thm]{Example}

\newtheorem{rmk}[thm]{Remark}

\begin{document}

\newcommand{\fraka}{{\mathfrak a}}
\newcommand{\frakb}{{\mathfrak b}}
\newcommand{\frakc}{{\mathfrak c}}
\newcommand{\frakd}{{\mathfrak d}}
\newcommand{\frake}{{\mathfrak e}}
\newcommand{\frakf}{{\mathfrak f}}
\newcommand{\frakg}{{\mathfrak g}}
\newcommand{\frakh}{{\mathfrak h}}
\newcommand{\fraki}{{\mathfrak i}}
\newcommand{\frakj}{{\mathfrak j}}
\newcommand{\frakk}{{\mathfrak k}}
\newcommand{\frakl}{{\mathfrak l}}
\newcommand{\frakm}{{\mathfrak m}}
\newcommand{\frakn}{{\mathfrak n}}
\newcommand{\frako}{{\mathfrak o}}
\newcommand{\frakp}{{\mathfrak p}}
\newcommand{\frakq}{{\mathfrak q}}
\newcommand{\frakr}{{\mathfrak r}}
\newcommand{\fraks}{{\mathfrak s}}
\newcommand{\frakt}{{\mathfrak t}}
\newcommand{\fraku}{{\mathfrak u}}
\newcommand{\frakv}{{\mathfrak v}}
\newcommand{\frakw}{{\mathfrak w}}
\newcommand{\frakx}{{\mathfrak x}}
\newcommand{\fraky}{{\mathfrak y}}
\newcommand{\frakz}{{\mathfrak z}}

\newcommand{\frakA}{{\mathfrak A}}
\newcommand{\frakB}{{\mathfrak B}}
\newcommand{\frakC}{{\mathfrak C}}
\newcommand{\frakD}{{\mathfrak D}}
\newcommand{\frakE}{{\mathfrak E}}
\newcommand{\frakF}{{\mathfrak F}}
\newcommand{\frakG}{{\mathfrak G}}
\newcommand{\frakH}{{\mathfrak H}}
\newcommand{\frakI}{{\mathfrak I}}
\newcommand{\frakJ}{{\mathfrak J}}
\newcommand{\frakK}{{\mathfrak K}}
\newcommand{\frakL}{{\mathfrak L}}
\newcommand{\frakM}{{\mathfrak M}}
\newcommand{\frakN}{{\mathfrak N}}
\newcommand{\frakO}{{\mathfrak O}}
\newcommand{\frakP}{{\mathfrak P}}
\newcommand{\frakQ}{{\mathfrak Q}}
\newcommand{\frakR}{{\mathfrak R}}
\newcommand{\frakS}{{\mathfrak S}}
\newcommand{\frakT}{{\mathfrak T}}
\newcommand{\frakU}{{\mathfrak U}}
\newcommand{\frakV}{{\mathfrak V}}
\newcommand{\frakW}{{\mathfrak W}}
\newcommand{\frakX}{{\mathfrak X}}
\newcommand{\frakY}{{\mathfrak Y}}
\newcommand{\frakZ}{{\mathfrak Z}}

\newcommand{\bA}{{\mathbb A}}
\newcommand{\bB}{{\mathbb B}}
\newcommand{\bC}{{\mathbb C}}
\newcommand{\bD}{{\mathbb D}}
\newcommand{\bE}{{\mathbb E}}
\newcommand{\bF}{{\mathbb F}}
\newcommand{\bG}{{\mathbb G}}
\newcommand{\bH}{{\mathbb H}}
\newcommand{\bI}{{\mathbb I}}
\newcommand{\bJ}{{\mathbb J}}
\newcommand{\bK}{{\mathbb K}}
\newcommand{\bL}{{\mathbb L}}
\newcommand{\bM}{{\mathbb M}}
\newcommand{\bN}{{\mathbb N}}
\newcommand{\bO}{{\mathbb O}}
\newcommand{\bP}{{\mathbb P}}
\newcommand{\bQ}{{\mathbb Q}}
\newcommand{\bR}{{\mathbb R}}
\newcommand{\bS}{{\mathbb S}}
\newcommand{\bT}{{\mathbb T}}
\newcommand{\bU}{{\mathbb U}}
\newcommand{\bV}{{\mathbb V}}
\newcommand{\bW}{{\mathbb W}}
\newcommand{\bX}{{\mathbb X}}
\newcommand{\bY}{{\mathbb Y}}
\newcommand{\bZ}{{\mathbb Z}}

\newcommand{\bfA}{{\mathbf A}}
\newcommand{\bfB}{{\mathbf B}}
\newcommand{\bfC}{{\mathbf C}}
\newcommand{\bfD}{{\mathbf D}}
\newcommand{\bfE}{{\mathbf E}}
\newcommand{\bfF}{{\mathbf F}}
\newcommand{\bfG}{{\mathbf G}}
\newcommand{\bfH}{{\mathbf H}}
\newcommand{\bfI}{{\mathbf I}}
\newcommand{\bfJ}{{\mathbf J}}
\newcommand{\bfK}{{\mathbf K}}
\newcommand{\bfL}{{\mathbf L}}
\newcommand{\bfM}{{\mathbf M}}
\newcommand{\bfN}{{\mathbf N}}
\newcommand{\bfO}{{\mathbf O}}
\newcommand{\bfP}{{\mathbf P}}
\newcommand{\bfQ}{{\mathbf Q}}
\newcommand{\bfR}{{\mathbf R}}
\newcommand{\bfS}{{\mathbf S}}
\newcommand{\bfT}{{\mathbf T}}
\newcommand{\bfU}{{\mathbf U}}
\newcommand{\bfV}{{\mathbf V}}
\newcommand{\bfW}{{\mathbf W}}
\newcommand{\bfX}{{\mathbf X}}
\newcommand{\bfY}{{\mathbf Y}}
\newcommand{\bfZ}{{\mathbf Z}}

\newcommand{\calA}{{\mathcal A}}
\newcommand{\calB}{{\mathcal B}}
\newcommand{\calC}{{\mathcal C}}
\newcommand{\calD}{{\mathcal D}}
\newcommand{\calE}{{\mathcal E}}
\newcommand{\calF}{{\mathcal F}}
\newcommand{\calG}{{\mathcal G}}
\newcommand{\calH}{{\mathcal H}}
\newcommand{\calI}{{\mathcal I}}
\newcommand{\calJ}{{\mathcal J}}
\newcommand{\calK}{{\mathcal K}}
\newcommand{\calL}{{\mathcal L}}
\newcommand{\calM}{{\mathcal M}}
\newcommand{\calN}{{\mathcal N}}
\newcommand{\calO}{{\mathcal O}}
\newcommand{\calP}{{\mathcal P}}
\newcommand{\calQ}{{\mathcal Q}}
\newcommand{\calR}{{\mathcal R}}
\newcommand{\calS}{{\mathcal S}}
\newcommand{\calT}{{\mathcal T}}
\newcommand{\calU}{{\mathcal U}}
\newcommand{\calV}{{\mathcal V}}
\newcommand{\calW}{{\mathcal W}}
\newcommand{\calX}{{\mathcal X}}
\newcommand{\calY}{{\mathcal Y}}
\newcommand{\calZ}{{\mathcal Z}}

\newcommand{\rA}{{\mathrm A}}
\newcommand{\rB}{{\mathrm B}}
\newcommand{\rC}{{\mathrm C}}
\newcommand{\rD}{{\mathrm D}}
\newcommand{\rd}{{\mathrm d}}
\newcommand{\rE}{{\mathrm E}}
\newcommand{\rF}{{\mathrm F}}
\newcommand{\rG}{{\mathrm G}}
\newcommand{\rH}{{\mathrm H}}
\newcommand{\rI}{{\mathrm I}}
\newcommand{\rJ}{{\mathrm J}}
\newcommand{\rK}{{\mathrm K}}
\newcommand{\rL}{{\mathrm L}}
\newcommand{\rM}{{\mathrm M}}
\newcommand{\rN}{{\mathrm N}}
\newcommand{\rO}{{\mathrm O}}
\newcommand{\rP}{{\mathrm P}}
\newcommand{\rQ}{{\mathrm Q}}
\newcommand{\rR}{{\mathrm R}}
\newcommand{\rS}{{\mathrm S}}
\newcommand{\rT}{{\mathrm T}}
\newcommand{\rU}{{\mathrm U}}
\newcommand{\rV}{{\mathrm V}}
\newcommand{\rW}{{\mathrm W}}
\newcommand{\rX}{{\mathrm X}}
\newcommand{\rY}{{\mathrm Y}}
\newcommand{\rZ}{{\mathrm Z}}

\newcommand{\Zp}{{\bZ_p}}
\newcommand{\Zl}{{\bZ_l}}
\newcommand{\Qp}{{\bQ_p}}
\newcommand{\Ql}{{\bQ_l}}
\newcommand{\Cp}{{\bC_p}}
\newcommand{\Fp}{{\bF_p}}
\newcommand{\Fl}{{\bF_l}}

\newcommand{\Ainf}{{\mathbf{A}_{\mathrm{inf}}}}
\newcommand{\BdRp}{{\mathbf{B}_{\mathrm{dR}}^+}}
\newcommand{\BdR}{{\mathbf{B}_{\mathrm{dR}}}}

\newcommand{\OXp}{{\widehat \calO_X^+}}
\newcommand{\OX}{{\widehat \calO_X}}
\newcommand{\AAinf}{{\bA_{\mathrm{inf}}}}            
\newcommand{\BBinf}{{\bB_{\mathrm{inf}}}}            
\newcommand{\BBbd}{{\bB^{\mathrm{bd}}}}              
\newcommand{\BBdRp}{{\bB_{\mathrm{dR}}^+}}           
\newcommand{\BBdR}{{\bB_{\mathrm{dR}}}}              
\newcommand{\OAinf}{{\calO\bA_{\mathrm{inf}}}}         
\newcommand{\OBinf}{{\calO\bB_{\mathrm{inf}}}}         
\newcommand{\OBdRp}{{\calO\bB_{\mathrm{dR}}^+}}        
\newcommand{\OBdR}{{\calO\bB_{\mathrm{dR}}}}           
\newcommand{\OC}{{\calO\bC}}                           

\newcommand{\Aut}{{\mathrm{Aut}}}
\newcommand{\Bun}{{\mathrm{Bun}}}           
\newcommand{\Cech}{{\check{H}}}             
\newcommand{\Coker}{{\mathrm{Coker}}}       
\newcommand{\colim}{{\mathrm{colim}}}       
\newcommand{\dlog}{{\mathrm{dlog}}}         
\newcommand{\End}{{\mathrm{End}}}           
\newcommand{\Exp}{{\mathrm{Exp}}}           
\newcommand{\Ext}{{\mathrm{Ext}}}           
\newcommand{\Fil}{{\mathrm{Fil}}}           
\newcommand{\Fitt}{{\mathrm{Fitt}}}         
\newcommand{\Frac}{{\mathrm{Frac}}}         
\newcommand{\Frob}{{\mathrm{Frob}}}         
\newcommand{\Gal}{{\mathrm{Gal}}}           
\newcommand{\Gr}{{\mathrm{Gr}}}             
\newcommand{\Hom}{{\mathrm{Hom}}}           
\newcommand{\HIG}{{\mathrm{HIG}}}           
\newcommand{\id}{{\mathrm{id}}}             
\newcommand{\Ima}{{\mathrm{Im}}}            
\newcommand{\Isom}{{\mathrm{Isom}}}         
\newcommand{\Ker}{{\mathrm{Ker}}}           
\newcommand{\Lie}{{\mathrm{Lie}}}           
\newcommand{\Log}{{\mathrm{Log}}}           
\newcommand{\Mat}{{\mathrm{Mat}}}           
\newcommand{\MIC}{{\mathrm{MIC}}}           
\newcommand{\Mod}{{\mathrm{Mod}}}           
\newcommand{\Perf}{{\mathrm{Perf}}}         

\newcommand{\Perfd}{{\mathrm{Perfd}}}

\newcommand{\pr}{{\mathrm{pr}}}             
\newcommand{\Proj}{\mathrm{Proj}}        
\newcommand{\Rep}{{\mathrm{Rep}}}           
\newcommand{\Res}{{\mathrm{Res}}}           
\newcommand{\RGamma}{{\mathrm{\rR\Gamma}}}    
\newcommand{\rk}{{\mathrm{rk}}}             
\newcommand{\Rlim}{{\mathrm{R}\underleftarrow{\lim}}} 
\newcommand{\sgn}{{\mathrm{sgn}}}           
\newcommand{\Sh}{{\mathrm{Shv}}}             
\newcommand{\Sht}{{\mathrm{Sht}}}           
\newcommand{\Spa}{{\mathrm{Spa}}}           
\newcommand{\Spf}{{\mathrm{Spf}}}           
\newcommand{\Spec}{{\mathrm{Spec}}}         
\newcommand{\Strat}{{\mathrm{Strat}}}       
\newcommand{\Sym}{{\mathrm{Sym}}}           
\newcommand{\Tor}{{\mathrm{{Tor}}}}         
\newcommand{\Tot}{{\mathrm{Tot}}}           
\newcommand{\Vect}{{\mathrm{Vect}}}         


\newcommand{\GL}{{\mathrm{GL}}}             
\newcommand{\SL}{{\mathrm{SL}}}             


\newcommand{\aff}{{\mathrm{aff}}}          
\newcommand{\an}{{\mathrm{an}}}             
\newcommand{\can}{{\mathrm{can}}}           
\newcommand{\cl}{{\mathrm{cl}}}             
\newcommand{\cofib}{{\mathrm{cofib}}}        
\newcommand{\cts}{{\mathrm{cts}}}           
\newcommand{\cris}{{\mathrm{cris}}}         
\newcommand{\crys}{{\mathrm{crys}}}         
\newcommand{\cyc}{{\mathrm{cyc}}}           
\newcommand{\dR}{{\mathrm{dR}}}             
\newcommand{\et}{{\mathrm{\acute{e}t}}}    
\newcommand{\fib}{{\mathrm{fib}}}          
\newcommand{\fl}{{\mathrm{fl}}}             
\newcommand{\fppf}{{\mathrm{fppf}}}         
\newcommand{\ket}{{\mathrm{k\acute{e}t}}}    
\newcommand{\geo}{{\mathrm{geo}}}           
\newcommand{\gp}{{\mathrm{gp}}}             
\newcommand{\la}{{\mathrm{la}}}             
\newcommand{\nil}{{\mathrm{nil}}}           
\newcommand{\pd}{{\mathrm{pd}}}             
\newcommand{\perf}{\mathrm{perf}}           
\newcommand{\proet}{{\mathrm{pro\acute{e}t}}}
\newcommand{\proket}{{prok\mathrm{\acute{e}t}}}     
\newcommand{\rig}{{\mathrm{rig}}}            
\newcommand{\sm}{{\mathrm{sm}}}             
\renewcommand{\sp}{{\mathrm{sp}}}             
\newcommand{\st}{{\mathrm{st}}}             
\newcommand{\tor}{{\mathrm{tor}}}           
\newcommand{\Zar}{{\mathrm{Zar}}}           



\newcommand{\Prism}{{\mathlarger{\mathbbl{\Delta}}}} 
\newcommand{\OPrism}{{\calO_{\Prism}}}
\newcommand{\IPrism}{{\calI_{\Prism}}}

\newcommand{\ya}{{\rangle}}
\newcommand{\za}{{\langle}}

\newcommand{\smat}[1]{\left( \begin{smallmatrix} #1 \end{smallmatrix} \right)}



\newcommand{\Ctilde}{\widetilde{\bfC}}
\newcommand{\CCtilde}{\widetilde{\bC}}
\newcommand{\Cpsi}{\widetilde{\mathbf{C}}_{\tilde{\psi}}^I}
\newcommand{\CCpsi}{\widetilde{\mathbf{C}}_{\tilde{\psi},\infty}^I}
\newcommand{\Cppsi}{\widetilde{\mathbf{C}}_{\tilde{\psi}}^{I,+}}
\newcommand{\CCppsi}{\widetilde{\mathbf{C}}_{\tilde{\psi},\infty}^{I,+}}
\newcommand{\Bnpsi}{\widetilde{\bfB}_{\tilde{\psi},n}}
\newcommand{\BBpsi}{\widetilde{\mathbf{B}}_{\tilde{\psi},\infty}}
\newcommand{\Bpsi}{\widetilde{\mathbf{B}}_{\tilde{\psi}}}
\newcommand{\Bnppsi}{\widetilde{\mathbf{B}}_{\tilde{\psi},\infty,n}}
\newcommand{\Bppsi}{\widetilde{\mathbf{B}}_{\tilde{\psi},\infty}}
\newcommand{\norm}{{|\!|}}

\newcommand{\wtA}{\widetilde{A}}
\newcommand{\wtJ}{\widetilde{J}}
\newcommand{\wtS}{\widetilde{S}}
\newcommand{\wtR}{\widetilde{R}}
\newcommand{\wtX}{\widetilde{X}}
\newcommand{\wtY}{\widetilde{Y}}

\newcommand{\wtp}{\widetilde{\frakP}}
\newcommand{\wtr}{\widetilde{\calR}}
\newcommand{\wtx}{\widetilde{\frakX}}
\newcommand{\wty}{\widetilde{\frakY}}

\newcommand{\ruochuan}[1]{{\color{blue} Ruochuan: #1}}
\newcommand{\yupeng}[1]{\textcolor{red}{yupeng: #1}}
\newcommand{\xinwen}[1]{{\color{purple} Xinwen: #1}}

\title[]{A $p$-adic Riemann--Hilbert correspondence over Robba rings}


\author[]{Zekun Chen}
\address{Z.C., Academy of Mathematics and Systems Science, Chinese Academy of Sciences, Beijing, China; Morningside Center of Mathematics, Chinese Academy of Sciences, Beijing, China.}
\email{chenzekun@amss.ac.cn}

\author[]{Ruochuan Liu}
\address{R.L., Beijing International Center for Mathematical Research, Peking University, YiHeYuan Road 5, Beijing, 100190, China.}
\email{liuruochuan@bicmr.pku.edu.cn}

\author[]{Yupeng Wang}
\address{Y.W., Shanghai Center for Mathematical Sciences, Fudan University, 2005 Songhu Road, Shanghai, China.}
\email{wangyupeng@fudan.edu.cn}

\author[]{Jiahong Yu}
\address{J.Y., Academy of Mathematics and Systems Science, Chinese Academy of Sciences, Beijing, China; Morningside Center of Mathematics, Chinese Academy of Sciences, Beijing, China.}
\email{yu\_jh@amss.ac.cn}

\author[]{Xinwen Zhu}
\address{X.Z, Department of Mathematics, Building 380, Stanford, California 94305}
\email{zhuxw@stanford.edu}

\subjclass[2020]{Primary }


\begin{abstract}
    We construct overconvergent period sheaves $\mathcal{O}\widetilde{\mathbb{C}}^{I,\dagger}$ and $\mathcal{O}\mathbb{B}_{\mathrm{dR}}^{\dagger,+}$ on the $v$-site of a liftable semi-stable formal scheme $\frakX_Z$. By establishing a decompletion theorem for small representations, we prove a $p$-adic Riemann-Hilbert correspondence relating small $\widetilde{\mathbb{C}}^I$-local systems (resp. $\mathbb{B}_{\mathrm{dR}}^+$-local systems) on the generic fiber $X_Z$ to small flat connections on the associated period ringed spaces. We further analyze the local Frobenius action and show that this correspondence is compatible with the Frobenius structure. As an application, we prove a comparison theorem between the cohomology of these period sheaves and the rigid cohomology of the special fiber, refining results of Le Bras and Colmez--Nizio{\l}.
\end{abstract}

\maketitle
\setcounter{tocdepth}{1}

\tableofcontents


\section{Introduction}\label{sec:Intro}

The aim of this paper is to develop a Robba-ring refinement of the $p$-adic Riemann--Hilbert correspondence for semi-stable rigid analytic spaces. We construct overconvergent period sheaves on the $v$-site of the generic fiber of a liftable semi-stable formal scheme, prove Poincar\'e lemmas for these sheaves, establish a decompletion theorem for small representations, and deduce a Riemann--Hilbert correspondence between small local systems with Robba-type coefficients and small modules with integrable logarithmic connection. We further prove that this correspondence is compatible with Frobenius whenever a prismatic lifting is available, and we apply the theory to compare the cohomology of our overconvergent period sheaves with the rigid cohomology of the special fiber.

\subsection{Background}

Classical Riemann--Hilbert theory relates local systems to differential equations. For a smooth complex analytic variety, the Riemann--Hilbert correspondence identifies finite-dimensional complex local systems with vector bundles equipped with integrable connection.

In $p$-adic geometry, the corresponding problem is more subtle. One has to introduce Fontaine-type period rings and period sheaves. The first step toward this direction is due to Scholze \cite{Sch-Pi}. He introduces the pro-\'etale site of rigid analytic varieties together with the period sheaves $\BBdRp$ and its geometric variant $\calO\bB_{\dR}^+$. These sheaves play the role of coefficients through which $p$-adic local systems can be compared with de Rham objects and lead him to establish a Riemann--Hilbert correspondence for de Rham $\BBdRp$-local systems. The key ingredient of this framework is the Poincar\'e lemma for $\calO\bB_{\dR}^+$: The de Rham complex associated to $\calO\bB_{\dR}^+$ is a resolution of $\BBdRp$.
Building on this framework, Liu and Zhu \cite{LZ} constructed a $p$-adic Riemann--Hilbert functor for smooth rigid analytic varieties. Their construction associates to a $p$-adic local system a vector bundle with integrable connection and can be generalized to the logarithmic case \cite{DLLZ23}. Based on these works, Gao, Min and Wang established an arithmetic Riemann--Hilbert correspondence and applied it to prismatic theory \cite{GMW24}.

In the above theories, one works mainly with $\bB_{\dR}^+$-local systems, which can be viewed as vector bundles on the infinitesimal neighborhood of the infinity of the relative Fargues--Fontaine curve. A natural question is whether the $p$-adic Riemann--Hilbert correspondence can be refined from $\bB_{\dR}^+$-local systems to local systems with coefficients in Robba-type period sheaves.

There are two additional difficulties in such a refinement. First, one has to construct geometric Robba-type period sheaves with connection and prove the corresponding Poincar\'e lemma. Second, the local calculations naturally take place after passing to completed period rings on a pro-finite-\'etale cover, whereas a Robba-ring Riemann--Hilbert correspondence must be formulated over decompleted, overconvergent rings. Thus one needs a decompletion theorem showing that small representations over the completed period rings descend to the Robba period rings without losing the cohomological information needed for the correspondence.

This paper aims to answer this question in a semi-stable setting. Instead of working only with $\bB_{\dR}^+$-local systems, we construct a period sheaf $\CCtilde^I$ depending on a closed radius interval $I=[r,s]\subset(0,+\infty)$ together with its geometric variant $\calO\CCtilde^{I,\dagger}$, which plays the role of $\calO\bB_{\dR}^+$ in the previous theory. Together with the decompletion theorem for small representations, this allows us to formulate and prove a Robba-ring version of the $p$-adic Riemann--Hilbert correspondence.

We also give two applications: a computation of the pro-\'etale cohomology of affine spaces and a comparison between the resulting overconvergent period cohomology of the generic fiber of a smooth formal scheme over $\rW(\kappa)$ and rigid cohomology of the special fiber.

\subsection{Main results}

We now describe the main constructions and results of the paper.
Let $C$ be a complete algebraically closed extension of $\Qp$, and let $ Z=\Spa(A,A^+)$ be an affinoid perfectoid space over $C$. Let $\frakX_Z$ be a liftable semi-stable formal scheme over $A^+$, and let $X_Z$ denote its generic fiber. We work on the $v$-site $X_{Z,v}$. 

\subsubsection{Riemann--Hilbert correspondence over Robba rings}
For a closed interval $I=[r,s]\subset(\frac{1}{p-1},p)$ with $pr\leq s$, we construct a Robba-type coefficient sheaf $\CCtilde^I$ and a geometric period sheaf $\calO\CCtilde^{I,\dagger}$ with a connection
\[
    \nabla:
    \calO\CCtilde^{I,\dagger}
    \longrightarrow
    \calO\CCtilde^{I,\dagger}
    \otimes_{\calO_{\wtx_Z}}
    \Omega^{1,\log}_{\wtx_Z}\{-1\}.
\]
by using a fixed lifting $\wtx_Z$ of $\frakX_Z$ over $\bA_{\inf}(Z)$, where $\{-1\}$ denotes the Breuil--Kisin twist. The dagger symbol indicates overconvergence:  Locally, sections of $\calO\CCtilde^{I,\dagger}$ are not merely completed formal power series; they are overconvergent series satisfying radius estimates controlled by the interval $I$. In particular, we have the Poincar\'e Lemma holds for $\CCtilde^I$; that is, the following sequence
\[
0\to \CCtilde^I
\to \calO\CCtilde^{I,\dagger}
\xrightarrow{\nabla}
\calO\CCtilde^{I,\dagger}\otimes_{\calO_{\wtx_Z}}
\Omega^{1,\log}_{\wtx_Z}\{-1\}
\to \cdots
\]
is exact. 

Let $\wtx_Z^I$ be the base-change of the given lifting $\wtx_Z$ of $\frakX_Z$ along $\AAinf(Z)\to\CCtilde^I(Z)$.
Denote by $\rL\rS^{\dagger}(X_Z,\CCtilde^I)$ category of small $\CCtilde^I$-local systems on $X_{Z,v}$ (cf. Definition \ref{dfn:small local system}) and by $\MIC^{\dagger}(\wtx_Z^I)$ the category of small finite locally free modules with integrable connection on $\wtx_Z^I$ (cf. Definition \ref{dfn:small connection}). Then our main result is the following.

  \begin{thm}[Theorem \ref{thm:RH correspondence}]\label{thm:intro-C}
      Let $I = [r,s]\subset(\frac{1}{p-1},p)$ be a closed sub-interval and $\nu:X_{Z,v}\to \frakX_{Z,\et}$ be the morphism of sites. 
      There is an equivalence of categories
\[
    \calD^I:
    \rL\rS^{\dagger}(X_Z,\CCtilde^I)
    \xrightarrow{\sim}
    \MIC^{\dagger}(\wtx_Z^I),
\]
which is compatible with tensor products, duals, ranks, and cohomology. More precisely, if $\bM\in \rL\rS^{\dagger}(X_Z,\CCtilde^I)$ corresponds to $(\calD,\nabla_{\calD})\in \MIC^{\dagger}(\wtx_Z^I)$, then we have
\[
    \rR\nu_*(\bM)
    \simeq
    \rD\rR(\calD,\nabla_{\calD})
\]
and thus a quasi-isomorphism
\[\rR\Gamma(X_{Z,v},\bM)\simeq \rR\Gamma(\frakX_{Z,\et},\rD\rR(\calD,\nabla_{\calD})).\]
  \end{thm}
To obtain this theorem, the difficult part is to construct $\calD^I$, which can be handled with \'etale locally on $\frakX_Z$. So we can assume that $\frakX_Z=\Spf(\calR_Z)$ is a small semi-stable formal scheme with a lifting $\wtx_Z = \Spf(\widetilde{\calR}_Z)$. Then the generic fiber $X_Z$ of $\frakX_Z$ admits a pro-finite-\'etale $\Gamma$-covering $X_{\infty,Z}\to X_Z$.
Put $\CCpsi = \CCtilde^I(X_{\infty,Z})$ and then it is equipped with a continuous $\Gamma$-action and admits a $\Gamma$-stable sub-$\CCtilde^I(Z)$-algebra $\Cpsi$ that is isomorphic to $\widetilde \calR_Z\widehat \otimes_{\bA_{\inf}(Z)}\CCtilde^I(Z)$. Then one can define the categories $\Rep^{\dagger}_{\Gamma}(\Cpsi)$ and $ \Rep^{\dagger}_{\Gamma}(\CCpsi)$ of small $\Gamma$-representations over $\Cpsi$ and $\CCpsi$, respectively. It can be shown that the functor $\bM\mapsto\bM(X_{\infty,Z})$ given by taking sections on $X_{\infty,Z}$ induces an equivalence of categories
\[\rL\rS^{\dagger}(X_{Z,v},\CCtilde^I)\xrightarrow{\simeq}\Rep^{\dagger}_{\Gamma}(\CCpsi)\]
that is compatible with cohomology (cf. Lemma \ref{lem:local description of local systems}).
By a standard argument, to construct $\calD^I$, we are reduced to computing the $\Gamma$-cohomology of $(\bM\otimes\calO\CCtilde^{I,\dagger})(X_{\infty,Z})$ for a small $\CCtilde^I$-local system $\bM$. The key ingredient is the following decompletion theorem.

\begin{thm}[Corollary \ref{cor:decompletion}]\label{thm:intro-decompletion}
Let $I=[r,s]\subset(0,p)$ be a closed interval and let $\lambda>p^{\frac{1}{p-1}}$. The base-change functor induces equivalences
\[
    \Rep^{\dagger}_{\Gamma}(\Cpsi)\simeq \Rep^{\dagger}_{\Gamma}(\CCpsi).
\]
that is compatible with cohomology.
\end{thm}
Using this, we are reduced to computing the $\Gamma$-cohomology of the decompletion of $(\bM\otimes\calO\CCtilde^{I,\dagger})(X_{\infty,Z})$, which is easier to handle.

The Riemann--Hilbert correspondence \ref{thm:intro-C} is also compatible with Frobenius action. Recall that the Frobenius endomorphism on $\AAinf(Z)$ induces an isomorphism 
\[\phi:\CCtilde^I(Z)\to \CCtilde^{I/p}(Z).\]
Using this, we are able to define $\phi$-modules over $\CCtilde^I$ (cf. Definition \ref{dfn:phi-module}).
If $\frakX_Z$ admits a prismatic lifting $\wtx_Z$ over $\AAinf(Z)$ in the sense of Definition \ref{dfn:prismatic lifting}, then we get a Frobenius morphism 
\[\phi:\calO_{\wtx_Z^I}\to \calO_{\wtx_Z^{I/p}}.\] 
Using this, we are able to define $\phi$-connections over $\wtx_Z^I$ (cf. Definition \ref{dfn:phi-connection}). Denote by $\rL\rS^{\dagger,\phi}(X_Z,\CCtilde^I)$ and $ \MIC^{\dagger,\phi}(\wtx_Z^I)$ the category of small $\phi$-modules over $\CCtilde^I$ on $X_{Z,v}$ and the category of small $\phi$-connections on $\wtx_Z^I$, respectively. Then we have the following result.

\begin{theorem}[Theorem \ref{thm:prismatic RH}]\label{thm:intro-D}
Assume that $\frakX_Z$ admits a prismatic lifting $\wtx_Z$. Then the functor $\calD^I$ in Theorem \ref{thm:intro-C} induces an equivalence of categories
\[
    \rL\rS^{\dagger,\phi}(X_Z,\CCtilde^I)
    \xrightarrow{\ \simeq\ }
    \MIC^{\dagger,\phi}(\wtx_Z^I)
\]
that is compatible with Frobenius actions.
\end{theorem}
When $\frakX_Z = \Spf(\calR_Z)$ is small affine, we get an explicit description of this equivalence (cf. Proposition \ref{prop:local prismatic RH}).

\subsubsection{Riemann--Hilbert correspondence over de Rham period ring}
When $1\in I$, one gets a natural morphism
\[
    \CCtilde^I\longrightarrow \BBdRp.
\]
identifying $\BBdRp$ with the $t$-adic completion of $\CCtilde^I$, where $t$ denotes Fontaine's $p$-adic analog of $2\pi\sqrt{-1}$. After taking $t$-adic completion of $\calO\CCtilde^{I,\dagger}$, we also get an overconvergent de Rham period sheaf $\calO\bB_{\dR}^{\dagger,+}$, which carries a natural connection and satisfies an analogous Poincar\'e lemma:
\[
0\to \BBdRp
\to \calO\bB_{\dR}^{\dagger,+}
\xrightarrow{\nabla}
\calO\bB_{\dR}^{\dagger,+}\otimes_{\calO_{\wtx_Z}}
\Omega^{1,\log}_{\wtx_Z}\{-1\}
\to \cdots .
\]

Using this, we are able to prove a de Rham version of Theorem \ref{thm:intro-C} with $\BBdRp$-coefficients: For any $1\leq n\leq \infty$, denote by $\rL\rS^{\dagger}(X_Z,\BBdRp/t^n)$ category of small $\BBdRp/t^n$-local systems on $X_{Z,v}$ (cf. Definition \ref{dfn:small local system}) and by $\MIC^{\dagger}(\wtx_Z^I)$ the category of small finite locally free modules with integrable connection on $\wtX_Z/t^n$ (cf. Definition \ref{dfn:small connection}), where $\wtX_Z/t^n$ denotes the lifting of $X_Z$ over $\bB_{\dR}^+(Z)/t^n$ induced by $\wtx_Z$.

\begin{thm}[Theorem \ref{thm:RH correspondence over dR}]\label{thm:intro-DR}
      Let $1\leq n\leq \infty$. Denote by $\nu:X_{Z,v}\to \frakX_{Z,\et}$ the morphism of sites. 
            There is an equivalence of categories
\[
    \widehat \calD:
    \rL\rS^{\dagger}(X_Z,\BBdRp/t^n)
    \xrightarrow{\sim}
    \MIC^{\dagger}(\wtX_Z/t^n),
\]
which is compatible with tensor products, duals, ranks, and cohomology. More precisely, if $\bM\in  \rL\rS^{\dagger}(X_Z,\BBdRp/t^n)$ corresponds to $(\calD,\nabla_{\calD})\in \MIC^{\dagger}(\wtX_Z/t^n)$, then we have
\[
    \rR\nu_*(\bM)
    \simeq
    \rD\rR(\calD,\nabla_{\calD})
\]
and thus a quasi-isomorphism
\[\rR\Gamma(X_{Z,v},\bM)\simeq \rR\Gamma(\frakX_{Z,\et},\rD\rR(\calD,\nabla_{\calD})).\]
  \end{thm}

  Although Theorem \ref{thm:intro-DR} follows from the same argument as  \ref{thm:intro-C} by using an analogous decompletion theorem of Theorem \ref{thm:intro-decompletion}, it can also be deduced in a much easier way: Note that taking $n=1$ in Theorem \ref{thm:intro-DR}, we get the $p$-adic Simpson correspondence in \cite{Wan23} for small $v$-bundles. As the period sheaf with connection $(\calO\bB_{\dR}^{\dagger,+},\nabla)$ is $t$-complete, one can conclude by a standard d\'evissage argument.
\subsubsection{Functoriality in coefficients}
  If $J\subset I$ is another closed sub-interval, then we get a natural morphism $\CCtilde^I\to\CCtilde^J$ of Robba rings.
  Moreover, as mentioned above, when $1\in I$, there is a natural morphism $\CCtilde^I\to\bB_{\dR}^+$. By construction, Theorems \ref{thm:intro-C} and \ref{thm:intro-DR} are compatible via the base-change of coefficients.

\begin{thm}[Theorem \ref{thm:functoriality in coefficient}]\label{thm:Intro-functoriality}
      Suppose that $I = [r,s]\subset(\frac{1}{p-1},p)$ is a closed sub-interval.
      \begin{enumerate}
          \item[(1)] Let $J\subset I$ be a closed sub-interval. Then the following diagram commutes
          \begin{equation*}
              \xymatrix@C=0.5cm{
                 \rL\rS^{\dagger}(X_Z,\CCtilde^I)\ar[d]^{\calD^I}\ar[rrrr]^{-\otimes_{\CCtilde^I}\CCtilde^J}&&&&\rL\rS^{\dagger}(X_Z,\CCtilde^J)\ar[d]^{\calD^J}\\
                 \MIC^{\dagger}(\wtx_Z^I)\ar[rrrr]^{-\otimes_{\calO_{\wtx_Z^I}}\calO_{\wtx_Z^J}}&&&&\MIC^{\dagger}(\wtx_Z^J),
              }
          \end{equation*}
          where the horizontal arrows are induced by base change.

          \item[(2)] Suppose that $1\in I$ and let $1\leq n\leq \infty$. Then the functor $\calD^I$ is compatible with $\widehat \calD$ in the sense that the following diagram commutes
          \begin{equation*}
              \xymatrix@C=0.5cm{
                 \rL\rS^{\dagger}(X_Z,\CCtilde^I)\ar[d]^{\calD^I}\ar[rrrr]^{-\otimes_{\CCtilde^I}\BBdRp/t^n}&&&&\rL\rS^{\dagger}(X_Z,\BBdRp/t^n)\ar[d]^{\widehat \calD}\\
                 \MIC^{\dagger}(\wtx_Z^I)\ar[rrrr]^{-\otimes_{\calO_{\wtx_Z^I}}\calO_{\wtX_Z}/t^n}&&&&\MIC^{\dagger}(\wtX_Z/t^n),
              }
          \end{equation*}
          where the horizontal arrows are induced by base change.
      \end{enumerate}
  \end{thm}

\subsection{Applications}
We give some applications of our theory.
\subsubsection{Cohomology of affine spaces}
As the first application, we compute the cohomology of affine spaces in terms of Riemann--Hilbert correspondence, which plays an important role in \cite{LB18} and \cite{CN20}.
\begin{thm}[Theorem \ref{thm:pro-etale cohomology of Affine space}]
      Let $\bA^d_Z$ be the affine space of relative dimension $d$ over $Z$.
      \begin{enumerate}
          \item[(1)] We have $\rH^0_v(\bA^d_Z,\Qp) = \Qp$ and $\rH^0_v(\bA^d_Z,\BBdRp) = \BBdRp(Z)$.
          \item[(2)] For any $n\geq 1$, we have $\rH^n_v(\bA^d_Z,\Qp) \cong \Omega^n(\bA_Z^d)^{\rd = 0}\{-n\}\cong \rH^n_v(\bA^d_Z,\BBdRp)$.
      \end{enumerate}
      Here, $\rd:\Omega^n(\bA_Z^d)\to \Omega^{n+1}(\bA_Z^d)$ denotes the differential in the usual de Rham complex of $\bA_Z^d$:
      \[0\to \calO(\bA_Z^d)\xrightarrow{\rd}\Omega^1(\bA_Z^d)\to\cdots\to\Omega^n(\bA_Z^d)\xrightarrow{\rd}\Omega^{n+1}(\bA_Z^d)\to\cdots.\]
\end{thm}
Our proof is different from and much easier than those in the literature. By using Riemann--Hilbert correspondences \ref{thm:intro-C} and \ref{thm:intro-D}, the complex $\rR\Gamma_v(\bA_Z^d,\bP)$ for $\bP\in \{\CCtilde^I,\BBdRp\}$ can be computed easily from the usual de Rham cohomology of affine spaces. Noting that $\bA_Z^d$ admits a prismatic lifting, one can obtain $\rR\Gamma_v(\bA_Z^d,\Qp)$ by taking $\phi$-invariant of $\rR\Gamma_v(\bA_Z^d,\CCtilde^I)$.

\subsubsection{Relation to rigid cohomology}
We can also apply the theory to rigid cohomology. Let $\frakX_0$ be a quasi-compact quasi-separated smooth formal scheme over $\calO_K = \rW(\kappa)$, and let $\calX$ be its special fiber. Suppose that $\frakX_Z$ denotes the base change of $\frakX_0$ along $\calO_K\to A^+$. For $\bP\in\{\CCtilde^I,\BBdRp\}$, inspired by \cite{BMS18}, one can construct a perfect complex of $\bP(Z)$-modules
\[
    \rR\Gamma^{\dagger}(X_Z,\bP).
\]
Let $\nu:X_{Z,v}\to \frakX_{Z,\et}$ be the natural morphism of sites. When $\frakX_Z$ is proper, we define
\[\rR\Gamma^{\dagger}(X_Z,\bP) = \rR\Gamma(\frakX_Z, \rL\eta_t\rR\nu_*\bP),\]
where $\rL\eta_t$ denotes the d\'ecalage functor in literature. In general,  $\rR\Gamma^{\dagger}(X_Z,\bP)$ is defined to be the ``overconvergent'' global section of $\rL\eta_t\rR\nu_*\bP$.

As $\phi(t) = pt$, the isomorphism $\phi:\CCtilde^I\to \CCtilde^{I/p}$ induces a quasi-isomorphism of complexes
\[\rL\eta_t\rR\nu_*\CCtilde^I\to \rL\eta_t\rR\nu_*\CCtilde^{I/p}\]
and thus a quasi-isomorphism
\[\phi:\rR\Gamma^{\dagger}(X_Z,\CCtilde^I)\to \rR\Gamma^{\dagger}(X_Z,\CCtilde^{I/p}).\]
The complex $\rR\Gamma^{\dagger}(X_Z,\CCtilde^I)$ together with the above Frobenius morphism is closely related with the rigid cohomology $\rR\Gamma_{\rig}(\calX/K)$ of $\calX$:

\begin{thm}[Theorem \ref{thm:P-cohomology}]\label{thm:intro-E}
    Suppose that $I = [r,s]\subset (\frac{1}{p-1},p)$ and let $\bP\in \{\CCtilde^I,\BBdRp\}$.
    \begin{enumerate}
        \item There exists a canonical quasi-isomorphism
        \[\iota_{\bP}:\rR\Gamma_{\rig}(\calX/K)\otimes_K\bP(Z)\xrightarrow{\simeq}\rR\Gamma^{\dagger}(X_Z,\bP).\]

        \item Suppose that $s = rp$. Then the quasi-isomorphism $\iota_{\CCtilde^I}$ is $\phi$-equivariant in the sense that the following diagram commutes:
        \begin{equation*}
            \begin{tikzcd}
                {\rR\Gamma_{\rig}(\calX/K)\otimes_K\CCtilde^r(Z)} \arrow[rr, "\iota_{\CCtilde^r}"] \arrow[d,"\phi"] &  & {\rR\Gamma^{\dagger}(X_Z,\CCtilde^r)} \arrow[d,"\phi"] \\
                \rR\Gamma_{\rig}(\calX/K)\otimes_K\CCtilde^{rp}(Z) \arrow[rr, "\iota_{\CCtilde^{rp}}"] &  & {\rR\Gamma^{\dagger}(X_Z,\CCtilde^{rp})}
            \end{tikzcd}
        \end{equation*}
        where the left vertical $\phi$ is the Frobenius map on rigid cohomology $\rR\Gamma_{\rig}(\calX/K)$ (via the coefficient map $\CCtilde^r(Z)\to \CCtilde^{rp}(Z)$) and the right vertical $\phi$ is the Frobenius map constructed above.

        \item If $1\in I$, then the natural morphism $\CCtilde^I\to\BBdRp$ induces a quasi-isomorphism
        \[\rR\Gamma^{\dagger}(X_Z,\CCtilde^I)\otimes_{\CCtilde^I(Z)}\BBdRp(Z)\xrightarrow{\simeq}\rR\Gamma^{\dagger}(X_Z,\BBdRp).\]
    \end{enumerate}
\end{thm}

Recall that the (adic) Fargues--Fontaine curve $\rF\rF_Z$ relative to $Z$ is the quotient of $\Spa(\CCtilde^I(Z),\CCtilde^{I,+}(Z))$ with $I = [r,rp]$ by gluing the boundaries $\Spa(\CCtilde^r(Z),\CCtilde^{r,+}(Z))$ with $\Spa(\CCtilde^{rp}(Z),\CCtilde^{rp,+}(Z))$ via the isomorphism $\phi:\CCtilde^{rp}(Z)\to \CCtilde^r(Z)$. So one can regard $\rR\Gamma^{\dagger}(X_Z,\CCtilde^I)$ together with the quasi-isomorphism $\phi:\rR\Gamma^{\dagger}(X_Z,\CCtilde^{rp})\xrightarrow{\simeq} \rR\Gamma^{\dagger}(X_Z,\CCtilde^r)$ as a perfect complex $\rR\Gamma_{\rF\rF}(X_Z)$ of vector bundles on $\rF\rF_Z$. Then the above theorem says that the complex $\rR\Gamma_{\rF\rF}(X_Z)$ coincides with $\rR\Gamma_{\rig}(\calX/K)\otimes\calO_{\rF\rF_Z}$. Moreover, the restriction of $\rR\Gamma_{\rF\rF}(X_Z)$ to the infinitesimal neighborhood of ``infinite point'' $\infty:Z\to \rF\rF_Z$ is exactly $\rR\Gamma^\dagger(X_Z,\BBdRp)$. This settles a conjecture \cite[Conj. 6.3]{LB18-conjecture} of Le Bras in our setting.

\subsection{Comparison with previous works}\label{ssec:intro-comparison-selected}

We end this introduction by making explicit the relation between our construction and some previous works.

\begin{itemize}
    \item The comparison with $p$-adic Simpson correspondence \cite{Wan23} appears through the reduction of our overconvergent de Rham period sheaf:   The period sheaf $(\calO\bB_{\dR}^{\dagger,+},\nabla)$ with connection coincides with the period sheaf $(\calO\bC^{\dagger},\Theta)$ in \cite{Wan23} modulo $t$ (cf. Proposition \ref{prop:compare with Wang}). Thus, our $p$-adic Riemann--Hilbert correspondence above is compatible with the small $p$-adic Simpson correspondence in literature (cf. \cite{AGT,Fal05,Wan23,AHLB23}, etc.).

    \item 
    Several earlier formulations of the $p$-adic Riemann--Hilbert correspondence are naturally tied to arithmetic situations over discretely valued fields and with coefficients over $\BdRp$ (cf. \cite{Sch-Pi,LZ,DLLZ23,GMW24}, etc.), where one also keeps track of the Galois action of the ground field. In contrast, the present paper works in a geometric relative setting over an algebraically closed complete field, and more generally over an affinoid perfectoid base.

    \item Motivated by our construction, a stacky Riemann--Hilbert correspondence for Hitchin-small objects was constructed in \cite{LMNQW}. Our overconvergent de Rham period sheaf $\calO\bB_{\dR}^{\dagger,+}$ can be viewed as a sub-sheaf of the period sheaf $\calO\bB_{\dR,\pd}^+$ in \cite{LMNQW} that is compatible with the connection (cf. Proposition \ref{prop:compare with LMNQW}). So our Riemann--Hilbert correspondence (for $\bB_{\dR}^+$-local systems) is also compatible with that in literature. However, the Riemann--Hilbert correspondence with Robba-type coefficients is first constructed in this paper. 

    \item The complex $\rR\Gamma_{\rF\rF}(X_Z)$ was first constructed in \cite{LB18-conjecture} by using another Robba-type period sheaf. The construction therein relies on the concrete calculation of $\Gamma$-cohomology of the period sheaf. In \cite{Bos23}, based on the previous work of Colmez--Nizio{\l} \cite{CN20}, Bosco also constructed $\rR\Gamma_{\rF\rF}(X_Z)$ by using Hyodo--Kato cohomology in a more general setting and also confirmed \cite[Conj. 6.3]{LB18-conjecture}. Compared with the methods in the literature, our Riemann--Hilbert correspondence simplifies the local computation and sheds light on dealing with the general coefficients.

    \item  More general, Ansch\"utz, Le Bras and Mann studied the six-functor formulism on the Fargues--Fontaine curve \cite{AMLB}. The pro-\'etale cohomology of $\Qp$-local systems were also computed in a different way by Li, Nizio{\l}, Reinecke and Zavyalov \cite{LNRZ} in the proper case. But no Riemann--Hilbert correspondence for $\CCtilde^I$-local systems appeared in these works.
\end{itemize}

\subsection{Organization of the paper}

The paper is organized as follows:
In \S\ref{sec:period sheaf without connection}, we define the coefficient sheaves $\CCtilde^I$ of relative Robba rings and study their norm and basic properties. In  \S\ref{sec:ringed space}, we recall the semi-stable formal schemes $\frakX_Z$ used throughout the paper and the ringed topoi $\wtx_Z^I$ associated to a lifting $\wtx_Z$ of $\frakX_Z$. In particular, we introduce the ringed spaces $\wtx_Z^I$ over Robba rings. 
Using the lifting $\wtx_Z$, we construct the overconvergent period sheaf $\calO\CCtilde^{I,\dagger}$ with  connection and prove the associated Poincar\'e Lemma in \S\ref{sec:period sheaf with connection}. In \S\ref{sec:overconvergent dR}, we construct and study the overconvergent de Rham period sheaf $\calO\bB_{\dR}^{\dagger,+}$ with connection and compare it and its reduction with the period sheaves appearing in \cite{LMNQW} and \cite{Wan23}  respectively.
In \S\ref{sec:decompletion}, we establish the decompletion theorem for small representations over Robba rings. This is the main descent input which passes from completed local period rings to overconvergent Robba rings, and it is used in \S\ref{sec:local RH} to establish the local Riemann--Hilbert correspondence over Robba rings. Section \S\ref{sec:local RH over BdR} is devoted to proving the corresponding local result for $\BBdRp/t^n$-coefficients. In \S\ref{sec:global RH}, we globalize the local result and prove the desired Riemann--Hilbert correspondence.
In \S\ref{sec:local Frobenius action}, we construct and analyze the local Frobenius action under the assumption that a prismatic lifting exists and prove Frobenius compatibility of the Riemann--Hilbert correspondence.
In \S\ref{sec:pro-etale cohomology of Affine space}, we compute the relevant pro-\'etale cohomology of affine spaces using Riemann--Hilbert correspondence.
Finally, in \S\ref{sec:rigid cohomology}, we introduce the complex $\rR\Gamma^\dagger(X,\bP)$ for smooth rigid variety $X$ with semi-stable reduction $\frakX$ and compare it with the rigid cohomology of the special fiber of $\frakX$.
\subsection*{Acknowledgement}

Z. Chen and Y. Wang were respectively invited to give a talk on this work at Southern University of Science and Technology in June 2024 and Peking University in July 2024. They would like to thank these institutes for the hospital.
Y. Wang would like to thank Hui Gao and Hansheng Diao for their interests in this work. 

Z. Chen is partially supported by the National Key R\&D Program of China 2025YFA1018000.
R. Liu is partially supported by the National Natural Science
Foundation of China under agreement No. NSFC–12321001, and by the New Cornerstone Foundation.
Y. Wang is supported by the National Natural Science Foundation of China (Grant No. 12601024) and partially supported by the CAS Project for Young Scientists in Basic Research (Grant No. YSBR-032). X. Zhu is supported by NSF grant under DMS-2200940.

\subsection*{AI Tool}
We write Introduction \S\ref{sec:Intro} with assistance of ChatGPT (5.6-Luna) and AI tool is not used in any other places. The authors take the full responsibility of mathematics throughout this paper.
\section{Period sheaves of Robba rings}\label{sec:period sheaf without connection}


Let $K$ be a complete discrete valuation field of mixed characteristic $(0,p)$ with perfect residue field $\kappa$ and ring of integers $\calO_K$. Let $C$ be the $p$-adic completion of a fixed algebraic closure of $K$, and let $|\cdot|$ be the $p$-adic norm on $C$ normalized by $|p|=p^{-1}$. Put $\varpi:=(p,p^{1/p},\dots)\in C^{\flat}$ and $\epsilon=(1,\zeta_p,\zeta_{p^2},\dots)\in C^{\flat}$, where $\{\zeta_{p^n}\}_{n\geq 0}$ is a fixed compatible sequence of primitive $p^n$-th roots of unity in $C$. Fix an embedding $\kappa \hookrightarrow C^\flat$, which induces an embedding $\varpi^{\bQ}\subset C^{\flat}$.

Denote by $\Perfd$ the $v$-site of affinoid perfectoid spaces over $C$ in the sense of \cite{Sch-Diamond}. 
In this section, we define and study the basic properties of the relative Robba ring $\CCtilde^I$ with respect to a closed interval $I\subset(0,+\infty)$, as a sheaf on $\Perfd$, and compare it with the de Rham period sheaf $\BBdRp$.

\subsection{The period sheaf $\CCtilde^I$: Construction and properties}

For any affinoid perfectoid space $Z = \Spa(A,A^+)\in \Perfd$, denote by $Z^{\flat} = \Spa(A^{\flat},A^{\flat+})$ its tilt and by $\norm\cdot\norm$ the spectral norm on $A^{\flat}$ compatible with the norm on $C^{\flat}$ induced from the $p$-adic norm on $C$. Put $\bA_{\inf}(Z):=W(A^{\flat+})$; any $x\in\bA_{\inf}(Z)$ is of the form
\[x = \sum_{n\geq 0}[x_n]p^n\]
with each $x_n\in A^{\flat+}$. Denote by $\xi$ a generator of the kernel of the natural surjection
\[\theta:\bA_{\inf}(Z)\to A^+.\]
Common choices for $\xi$ include $p-[\varpi]$ or $\frac{\phi(u)}{u}$ with $u = [\epsilon^{1/p}]-1$.

\begin{definition}\label{dfn:norm_r}
    Keep the notation introduced above.
    \begin{enumerate}
        \item[(1)] For any $r\geq 0$, define the function
        \[|\cdot|_r:\bA_{\inf}(Z)[\frac{1}{p[\varpi]}]\to \bR_{\geq 0}\]
        by sending each $x = \sum_{n\gg 0}[x_n]p^n$ to
        \[|x|_r:=\sup_n\{\norm x_n\norm^r \cdot|p|^n\}.\]

        \item[(2)] For any closed sub-interval $I = [r,s]\subset (0,+\infty)$, define the function
        \[|\cdot|_I:\bA_{\inf}(Z)[\frac{1}{p[\varpi]}]\to \bR_{\geq 0}\]
        by sending each $x$ to 
        \[|x|_I:=\sup_{t\in I}\{|x|_t\}.\]
    \end{enumerate}
    Equivalently, for $x=\sum_{n\gg 0}[x_n]p^n$ one has $|x|_I=\sup_{t\in I}\sup_n\{\norm x_n\norm^t\cdot|p|^n\}$.
\end{definition}

\begin{exam}\label{exam:norm of u}
    \begin{enumerate}
    \item[(1)] For any unit $v\in \bA_{\inf}(Z)$, we have $|v|_r = 1$ for all $r\geq 0$. Consequently, for any generator $\xi$ of $\Ker(\theta)$, we have 
    \begin{equation*}
    |\xi|_r = |p-[\varpi]|_r =
    \left\{
    \begin{array}{cc}
       p^{-r},  & \text{if }0<r\leq 1 \\
       p^{-1},  & \text{if }r\geq 1.
    \end{array}
    \right.
    \end{equation*}

   
    \item[(2)]
    We now compute $|u|_r$ for $u = [\epsilon^{1/p}]-1$. Write $u = \sum_{n\geq 0}[u_n]p^n$ for $u_n\in \calO_{C^{\flat}}$. It remains to compute $\norm u_n\norm$ for $n\geq 0$. Consider the ghost coordinates $(w_0(T),w_1(T),\dots)$ of $[1+T]-1\in W(\bZ[T])$; then $w_n(T)\in T\bZ[T]$ has degree $\leq p^n$ and satisfies
    \begin{equation}\label{equ:ghost coordinate}
        (1+T)^{p^n}-1=w_0(T)^{p^n}+pw_1(T)^{p^{n-1}}+\cdots+p^nw_n(T)
    \end{equation}
    for all $n\geq 0$. For example, $w_0(T) = T$ and $w_1(T) = \sum_{i=1}^{p-1}\frac{\binom{p}{i}}{p}T^i$.
    
    Let $a_n\in \bZ$ be the coefficient of $T$ in $w_n(T) \pmod{T^2}$. We claim $a_n$ is a unit in $\Zp$. Indeed, this follows from the fact that for any $n\geq 1$, in $\Zp[\zeta_p]$, 
    \begin{equation}\label{equ:calculate norm of u}
        w_n(\zeta_p-1)\equiv (-1)^{1+p+\cdots+p^{n-1}}\frac{(\zeta_p-1)^{p^n}}{p^{1+p+\cdots+p^{n-1}}} \mod (\zeta_p-1)^2,
    \end{equation}
    which implies $|w_n(\zeta_p-1)| = p^{-\frac{p^n}{p-1}+(1+p+\cdots+p^{n-1})} = p^{-\frac{1}{p-1}} = |\zeta_p-1|$, and thus $|a_n|=1$.
    
    Since $u_n^{p^n} = w_n(\epsilon^{1/p}-1)$ for any $n\geq 0$, we see that $\norm u_n\norm = p^{-\frac{1}{p^n(p-1)}}$, yielding
    \begin{equation*}
        |u|_r = \sup_{n\geq 0}\{p^{-\frac{r}{p^n(p-1)}-n}\} = \left\{
          \begin{array}{cc}
              p^{-\frac{r}{p-1}}, & \text{if }0<r\leq p \\
              p^{-\frac{r}{p^n(p-1)}-n}, & \text{if }p^n\leq r\leq p^{n+1}, \forall n\geq 1.
          \end{array}
        \right.
    \end{equation*}
    \end{enumerate}
\end{exam}

\begin{lem}\label{lem:norm_r}
    Keep the notation as above.
    \begin{enumerate}
        \item[(1)] For any $r>0$, $|\cdot|_r$ defines a power-multiplicative norm on $\bA_{\inf}(Z)[\frac{1}{p[\varpi]}]$. Moreover, if the spectral norm $\norm\cdot\norm:A^{\flat}\to \bR_{\geq 0}$ is multiplicative, then so is $|\cdot|_r$.
        \item[(2)] For any closed sub-interval $I = [r,s]\subset(0,+\infty)$, we have 
        \[|\cdot|_I=\max\{|\cdot|_r,|\cdot|_s\}.\]
    \end{enumerate}
\end{lem}
\begin{proof}
    Item (1) follows from \cite[Prop. 5.1.2]{KL15} and the proof therein.

    For item (2): For any $y\in A^{\flat}\setminus\{0\}$ and any $n\in\bZ$, consider the function
    \[\begin{split}
        f_y\colon&I\to \bR_{\geq 0}\\
        & t\mapsto \norm y\norm^t\cdot|p|^n = |[y]p^n|_t.
    \end{split}\]
    It is immediate that $f_y$ is strictly increasing when $\norm y\norm >1$ (i.e., valuation is positive) and decreasing when $\norm y\norm <1$. Thus, for any $t\in I$,
    \[|[y]p^n|_t\leq \max\{|[y]p^n|_r,|[y]p^n|_s\}.\]
    For any $x = \sum_{n\gg 0}[x_n]p^n\in \bA_{\inf}(Z)[\frac{1}{p[\varpi]}]$, we have
    \[|x|_t = \sup_n\{|[x_n]p^n|_t\}\leq \max\{\sup_n\{|[x_n]p^n|_r,|[x_n]p^n|_s\}\} = \max\{|x|_r,|x|_s\}.\]
    This proves the assertion.
\end{proof}

\begin{rmk}\label{rmk:truncated norm}
    Fix a closed sub-interval $I=[r,s]\subset (0,+\infty)$ and an affinoid perfectoid space $Z = \Spa(A,A^+)\in\Perfd$. For any $n\geq 0$, we define a norm
    \[|\cdot|_{I,n}:\bA_{\inf}(Z)[\frac{1}{[\varpi]}]/p^{n+1} = W_n(S^{\flat+})[\frac{1}{\varpi}]\to \bR_{\geq 0}\]
    by letting 
    \[|x|_{I,n}:=\sup_{t\in I}\sup_{0\leq i\leq n}\{\norm x_i\norm^t \cdot |p|^i\}\]
    for any $x = \sum_{i=0}^n[x_i]p^i\in \bA_{\inf}(Z)[\frac{1}{[\varpi]}]/p^{n+1}$. By definition, for any $x\in\bA_{\inf}(Z)[\frac{1}{[\varpi]}]$, 
    \[|x|_I=\lim_{n\to +\infty}|x|_{I,n},\]
    where we implicitly identify $x$ with its image in the truncation $\bA_{\inf}(Z)[\frac{1}{[\varpi]}]/p^{n+1}$.
\end{rmk}

\begin{dfn}\label{dfn:Robba ring}

    Let $I = [r,s]\subset (0,+\infty)$ be a closed interval. For any $Z\in \Perfd$, denote by
    \[\calY_Z^I\subset \Spa(\bA_{\inf}(Z)[\frac{1}{p}],\bA_{\inf}(Z))\]
    the affinoid subspace corresponding to  the rational localization with respect to the condition $|[\varpi]|\in |p|^I$.
    Put 
    \[\CCtilde^I(Z):=\calO(\calY_Z^I)\text{ and }\CCtilde^{I,+}(Z):=\calO^+(\calY_Z^I)\]
    and then we have
    \[\calY_Z^I = \Spa(\CCtilde^I(Z),\CCtilde^{I,+}(Z)).\]
    When $I=\{r\}$, we denote $\calY_Z^I$ by
    \[\calY_Z^r = \Spa(\CCtilde^r_K(Z),\CCtilde^{r,+}_K(Z)).\]
    So we get presheaves $\CCtilde^I$ and $\CCtilde^{I,+}$ on $\Perfd$ respectively given by
    \[Z\in \Perfd\mapsto \CCtilde^I(Z) \text{ and }\CCtilde^{I,+}(Z).\]
    In what follows, we denote $\CCtilde^I(Z)$ by $\widetilde{\bfC}^I$ in the case $Z = \Spa(C,\calO_C)$, and similarly define $\widetilde{\bfC}^{I,+}, \widetilde{\bfC}^r_K, \widetilde{\bfC}^{r,+}_K$, etc.
\end{dfn}

\begin{rmk}
    The inclusion $\calY_Z^I\to \Spa(\bA_{\inf}(Z))$ canonically factors through 
    \[\Spa(\bA_{\inf}(Z))\setminus V(p[\varpi])\subset \Spa(\bA_{\inf}(Z)),\]
    the complement of the vanishing locus of $p[\varpi]$.
\end{rmk}
\begin{rmk}
    It was explicitly used in \cite{KL16} that both presheaves $\CCtilde^I$ and $\CCtilde^{I,+}$ are acturally sheaves. We will give a direct proof in the next subsection under Hypothesis \ref{hy:rationality} below (cf. Proposition \ref{prop:C^I is a sheaf}).
\end{rmk}

\begin{hypothesis}\label{hy:rationality}
    Throughout the sequel, we assume that $I = [r,s]$ is a closed sub-interval of $(0,+\infty)$ with $r,s\in\bQ$.
\end{hypothesis}

\begin{lem}\label{lem:description of C^I}
    We keep the notation introduced above. Then $\CCtilde^{I,+}(Z)$ is the $p$-adic completion of $\bA_{\inf}(Z)[\frac{[\varpi^{1/r}]}{p},\frac{p}{[\varpi^{1/s}]}]$ inside $\bA_{\inf}(Z)[\frac{1}{p[\varpi]}]$  and $\CCtilde^I(Z) = \CCtilde^{I,+}(Z)[\frac{1}{p}]$.
\end{lem}
\begin{proof}
    Note that the condition $|[\varpi]|\in|p|^I$ amounts to
    \[|p|^s \leq |[\varpi]|\leq |p|^r \quad (\iff |[\varpi^{1/r}]|\leq |p| \leq |[\varpi^{1/s}]|).\]
    The assertion follows from the definition of the rational localization $\calY_Z^I$.
\end{proof}

\begin{lem}\label{lem:second description of C^I}
    Keep the notation as above.
    \begin{enumerate}
        \item[(1)] The ring $\CCtilde^I(Z)$ is the completion of $\bA_{\inf}(Z)[\frac{1}{p[\varpi]}]$ with respect to the norm $|\cdot|_I$, and
        \[\CCtilde^{I,+}(Z) = \{x\in \CCtilde^I(Z)\mid |x|_I\leq 1\}.\]

        \item[(2)] Fix $t\in I$. For any $c\in \widetilde \bC^I$ and any $x\in \CCtilde^I(Z)$, we have
        \[|cx|_t = |c|_t\cdot|x|_t.\]

        \item[(3)] The natural Frobenius map $\phi$ on $\bA_{\inf}(Z)$ induces an isometric isomorphism
        \[\phi:(\CCtilde^I(Z),|\cdot|_I)\to (\CCtilde^{\frac{I}{p}}(Z),|\cdot|_{\frac{I}{p}}).\]
    \end{enumerate}
\end{lem}
\begin{proof}
    For item (1): By \cite[Prop. 5.1.2]{KL15}, any point $|\cdot|\in \calY_Z^I$ must be dominated by $|\cdot|_I$; that is, for any $x\in \CCtilde^I(Z)$, we must have $|x|\leq |x|_I$.
    In particular, for any $n\geq 0$, we have
    \begin{equation}\label{equ:second description of C^I}
        \{x\in \CCtilde^I(Z)\mid |x|_I\leq p^{-n}\}\subset \{x\in \CCtilde^I(Z)\mid |x|\leq p^{-n},\forall |\cdot|\in \calY^I_Z\} = p^n\CCtilde^{I,+}(Z).
    \end{equation}
    Consequently, for any formal sum $\sum_{n\geq 0}x_n$ in $\CCtilde^{I,+}(Z)$ with $\lim_{n\to +\infty}|x_n|_I = 0$, we have that $x_n \in \CCtilde^{I,+}(Z)$ and $x_n$ tends to $0$ $p$-adically as $n\to +\infty$.
    Thus, the series $\sum_{n\geq 0}x_n$ converges to a unique element $x\in \CCtilde^I(Z)$ because $\CCtilde^{I,+}(Z)$ is $p$-adically complete (Lemma \ref{lem:description of C^I}), yielding that $\CCtilde^I(Z)$ is complete with respect to $|\cdot|_I$.
    Moreover, as $|x|_I\leq 1$ for any $x\in \bA_{\inf}(Z)[\frac{[\varpi^{1/r}]}{p},\frac{p}{[\varpi^{1/s}]}]$, it follows from Lemma \ref{lem:description of C^I} that
    \[\CCtilde^{I,+}(Z)\subset\{x\in\CCtilde^I(Z)\mid |x|_I\leq 1\}.\]
    Combined with \eqref{equ:second description of C^I} (for $n=0$), this implies
    \[\CCtilde^{I,+}(Z) = \{x\in\CCtilde^I(Z)\mid |x|_I\leq 1\}.\]

    It remains to show that $\CCtilde^I(Z)$ is the completion of $\bA_{\inf}(Z)[\frac{1}{p[\varpi]}]$ with respect to $|\cdot|_I$. Since $\CCtilde^I(Z)$ is complete with respect to $|\cdot|_I$, the natural morphism
    \[\bA_{\inf}(Z)[\frac{1}{p[\varpi]}]\to\CCtilde^I(Z)\]
    induces a unique morphism $i:\widetilde A\to \CCtilde^I(Z)$, where $\widetilde A$ denotes the completion of $\bA_{\inf}(Z)[\frac{1}{p[\varpi]}]$ with respect to $|\cdot|_I$.
    Conversely, let $\widetilde A^+$ be the completion of $\bA_{\inf}(Z)[\frac{[\varpi^{1/r}]}{p},\frac{p}{[\varpi^{1/s}]}]$ with respect to $|\cdot|_I$. Since $|p^n|_I = p^{-n} \to 0$, $\widetilde A^+$ is $p$-adically complete. Therefore, the natural morphism 
    \[\bA_{\inf}(Z)[\frac{[\varpi^{1/r}]}{p},\frac{p}{[\varpi^{1/s}]}] \to \widetilde A^+[\frac{1}{p}] \to \widetilde A\]
    induces a unique morphism $j:\CCtilde^I(Z)\to \widetilde A$. It remains only to observe that $i$ and $j$ are inverses of each other.

    For item (2): We first claim that for any $a\in C^{\flat}$ and any $f\in A^{\flat}$, $\norm af\norm = \norm a\norm\cdot\norm f\norm$.
    Indeed, assuming $a\neq 0$, Lemma \ref{lem:norm_r} gives
    \[\norm f\norm \leq \norm a^{-1}\norm\cdot \norm af\norm \leq \norm a^{-1}\norm \cdot \norm a\norm \cdot \norm f\norm = \norm f\norm,\]
    so equality holds. Consequently, for any $t>0$ and $n,m\in\bZ$, $|[af]p^{n+m}|_t = |[a]p^n|_t\cdot|[f]p^m|_t$.

    Now, let $c\in \widetilde \bC^I$ and $x\in \CCtilde^I(Z)$ and fix a $t\in I$. We show $|c|_t\cdot|x|_t = |cx|_t$.
    By Item (1): We may write $c = \sum_{n\in \bZ}[c_n]p^n$ and $x = \sum_{n\in \bZ}[x_n]p^n$ such that $\lim_{n\to\infty}|[c_n]p^n|_I = \lim_{n\to\infty}|[x_n]p^n|_I=0$.
    
    First, consider the case where there is exactly one $n\in\bZ$ and one $m\in \bZ$ maximizing the norms at $t$, i.e., $|c|_t = |[c_n]p^n|_t$ and $|x|_t = |[x_m]p^m|_t$. In this case:
    \[|cx|_t = |[c_n]p^n\cdot[x_m]p^m|_t = |[c_nx_m]p^{n+m}|_t = |[c_n]p^n|_t\cdot|[x_m]p^m|_t = |c|_t\cdot|x|_t.\]
    In general, by \cite[Lem. 5.2.1]{KL15}, there exists a finite set $J\subset I$ such that for any $t\in I\setminus J$, the maximizing indices are unique. Thus, for $t\in I\setminus J$, $|cx|_t = |c|_t\cdot|x|_t$.
    Since the functions $t\mapsto |c|_t, |x|_t, |cx|_t$ are continuous on $I$, the equality holds for all $t\in I$.

    For item (3): This follows because for any $x = \sum_{n\gg0}[x_n]p^n\in \bA_{\inf}(Z)[\frac{1}{p[\varpi]}]$, $\phi(x) = \sum_{n\geq 0}[x_n^p]p^n$, yielding $|\phi(x)|_{I/p} = |x|_I$. The assertion follows from the completion property in Item (1).
\end{proof}

\begin{rmk}
    In particular, $\CCtilde^I(Z)$ is the completion of $\bA_{\inf}(Z)[\frac{1}{p[\varpi]}]$ with respect to $|\cdot|_I$, and $\CCtilde^{I,+}(Z)=\{x\in\CCtilde^I(Z)\mid |x|_I\leq 1\}$ (cf.\ Definition \ref{dfn:Robba ring}).
\end{rmk}

\begin{cor}\label{cor:dense image}
    Let $I^{\prime} = [r^{\prime},s^{\prime}]$ be a closed sub-interval of $I$. Then the natural morphism
    \[\iota:\CCtilde^I(Z)\to \CCtilde^{I^{\prime}}(Z)\]
    has dense image in the target.
\end{cor}
\begin{proof}
    For any $x\in\bA_{\inf}(Z)[\frac{1}{p[\varpi]}]$, $|x|_{I^{\prime}}\leq |x|_I$. The result follows directly from Lemma \ref{lem:second description of C^I}(1).
\end{proof}

\begin{exam}\label{exam:t in C^I}
    Let $I = [r,s]$ be as before and $u = [\epsilon^{1/p}]-1$. The following statements follow easily from Example \ref{exam:norm of u}(2).
    \begin{enumerate}
        \item[(1)] If $r>\frac{1}{p-1}$, then 
        \[t:=\log[\epsilon] = \sum_{n\geq 1}\frac{(-1)^{n-1}([\epsilon]-1)^n}{n}\]
        is well-defined and is a unit multiple of $\phi(u) = [\epsilon]-1$ in $\CCtilde^{I,+}(Z)$. Consequently, $t$ is well-defined in $\CCtilde^J(Z)$ for any closed sub-interval $J\subset(0,+\infty)$ (since $\phi(t) = pt$).
        \item[(2)] If $s<p$, then $\frac{u}{[\epsilon^{1/p}-1]}$ is a unit in $\CCtilde^{I,+}(Z)$. In particular, $u$ is invertible in $\CCtilde^I(Z)$.

        \item[(3)] We claim $\xi$ is a unit in $\Ctilde^{I}$ provided that $1\notin I=[r,s]\subset(0,+\infty)$. To see this, observe that $|\frac{[\varpi]}{p}|_h = p^{1-h}$ for any $h>0$, implying $|\frac{[\varpi]}{p}|_I = p^{1-r}<1$ when $r>1$, and $|\frac{p}{[\varpi]}|_I = p^{s-1}<1$ when $s<1$.
    Thus we obtain the following two cases:
    \begin{enumerate}
        \item[(i)] $\frac{p-[\varpi]}{p} = 1-\frac{[\varpi]}{p}$ is a unit in $\Ctilde^{I,+}$ when $I\subset (1,+\infty)$;
        \item[(ii)] $\frac{p-[\varpi]}{[\varpi]} = \frac{p}{[\varpi]}-1$ is a unit in $\Ctilde^{I,+}$ when $I\subset(0,1)$.
    \end{enumerate}
    \end{enumerate}
\end{exam}

\begin{rmk}\label{rmk:compatibility of topology}
    Since $[\varpi^{1/r}]\widetilde \bfC^{I,+}\subset p\widetilde \bfC^{I,+}\subset [\varpi^{1/s}]\widetilde \bfC^{I,+}$, for any $\widetilde \bfC^{I,+}$-algebra $\bfA$, the $p$-adic topology, $[\varpi]$-adic topology, and $(p,\xi)$-adic topology on $\bfA$ coincide.
\end{rmk}

\begin{rmk}\label{rmk:isometry on base}
    Let $I\subset(0,+\infty)$ be a closed sub-interval. Let $Z^{\prime} = \Spa(S^{\prime},S^{\prime+})\to Z = \Spa(S,S^+)$ be a morphism of affinoid perfectoid algebras in $\Perfd$, and let 
    \[\iota:\CCtilde^I(Z)\to \CCtilde^I(Z^{\prime})\]
    be the induced morphism of $\Ctilde^I$-algebras.
    \begin{enumerate}
        \item[(1)] For any $f\in \CCtilde^I(Z)$, we always have $|\iota(f)|_I\leq |f|_I$. Indeed, it suffices to deal with the case $I = [r,r]$. By items (1) and (2) of Lemma \ref{lem:second description of C^I}, we may assume $f\in\bA_{\inf}(Z)$ and reduce to the case $f = [a]$ for $a\in S^{\flat+}$. Then the claim follows from $\norm i(a)\norm\leq \norm a\norm$, where $i:S^{\flat+}\to S^{\prime \flat+}$ is the morphism associated to $Z^{\prime}\to Z$.
        \item[(2)] If $Z^{\prime}\to Z$ is a covering, then $|\iota(f)|_I = |f|_I$. Reducing to $f=[a]$ as above, we need $\norm i(a)\norm = \norm a\norm$. For any geometric point $x\in Z^{\flat}$ of rank $1$, the fiber $Z^{\prime\flat}_x$ is non-empty; thus there exists $x^{\prime}\in Z^{\prime\flat}_x$ such that $|i(a)(x')| = |a(x)|$, yielding $\norm i(a)\norm \geq \norm a\norm$.
    \end{enumerate}
\end{rmk}

Now we compare $\CCtilde^I(Z)$ with the de Rham period ring 
\[\BBdRp(Z) = \varprojlim_n\bA_{\inf}(Z)[\frac{1}{p}]/\xi^n.\]

\begin{lem}\label{lem:compare with BdR}
    Keep the notation as above and suppose that $1\in I$. Then the map $\bA_{\inf}(Z)\to \BBdRp(Z)$ induces a morphism
    \[\iota:\CCtilde^I(Z)\to \BBdRp(Z)\]
    which carries $\CCtilde^{I,+}(Z)$ onto $A^+$ modulo $\xi$. Moreover, $\iota$ induces an isomorphism 
    \[\CCtilde^I(Z)^{\wedge_{\xi}}\xrightarrow{\cong}\BBdRp(Z),\]
    where $\CCtilde^I(Z)^{\wedge_{\xi}}$ denotes the $\xi$-adic completion of $\CCtilde^I(Z)$.
\end{lem}
\begin{proof}
    As $1\in I=[r,s]$, we have $\frac{1}{s}\leq 1\leq \frac{1}{r}$. Recall that
    \[\bA_{\inf}(Z)[\frac{[\varpi^{1/r}]}{p},\frac{p}{[\varpi^{1/s}]}] = \bA_{\inf}(Z)[X,Y]/(p X-[\varpi^{1/r}],[\varpi^{1/s}]Y-p, XY-[\varpi^{\frac{1}{r}-\frac{1}{s}}]).\]
    
    To obtain $\iota$, we must show that for any $n\geq 1$, the natural map
    \[i_n:\bA_{\inf}(Z)[X,Y]/(p X-[\varpi^{1/r}],[\varpi^{1/s}]Y-p, \dots)\to\BBdRp(Z)/\xi^n\]
    sending $X \mapsto \frac{[\varpi^{1/r}]}{p}$ and $Y \mapsto \frac{p}{[\varpi^{1/s}]}$ extends uniquely to the $p$-adic completion of its source. It suffices to show that $i_n$ is bounded, i.e., it takes values in
    \[p^{-(n-1)}\bA_{\inf}(Z)/\xi^n\subset \BBdRp(Z)/\xi^n.\]
    
    Let $\xi = p-[\varpi]$. We have
    \[i_n(X) = \frac{[\varpi^{1/r}]}{p} = [\varpi^{1/r-1}](1-\frac{\xi}{p}) \quad \text{and} \quad i_n(Y) = \frac{p}{[\varpi^{1/s}]} = [\varpi^{1-1/s}](1-\frac{\xi}{p})^{-1}.\]
    For any $m\in\bZ$, $(1-\frac{\xi}{p})^m = \sum_{i=0}^{n-1}\binom{m}{i}(-1)^ip^{-i}\xi^i\in p^{-(n-1)}\bA_{\inf}(Z)/\xi^n$, yielding the boundedness.

    Since $p^{-(n-1)}\bA_{\inf}(Z)/\xi^n$ is $p$-adically complete, Lemma \ref{lem:description of C^I} implies that $i_n$ uniquely induces a morphism $\CCtilde^{I,+}(Z)\to \BBdRp(Z)/\xi^n$, and thus $\iota_n:\CCtilde^I(Z)\to \BBdRp(Z)/\xi^n$. Letting $n\to +\infty$, we get $\iota:\CCtilde^I(Z)\to\BBdRp(Z)$. Since $\BBdRp(Z)$ is $\xi$-adically complete, $\iota$ induces $i:\CCtilde^I(Z)^{\wedge_{\xi}}\to \BBdRp(Z)$.
    Conversely, the composite $\bA_{\inf}(Z)[\frac{1}{p}]\to \CCtilde^I(Z)\to\CCtilde^I(Z)^{\wedge_{\xi}}$ uniquely induces a morphism $j:\BBdRp(Z)\to \CCtilde^I(Z)^{\wedge_{\xi}}$. It remains only to observe that $j$ is the inverse of $i$.
\end{proof}
    
\subsection{The presheaf $\CCtilde^I$ is a sheaf}\label{ssec:C^I is sheaf}

This subsection is devoted to proving that the presheaves $\CCtilde^{I,+}$ and $\CCtilde^I$ respectively sending $Z\in \Perfd\mapsto \CCtilde^I(Z)$ and $\CCtilde^{I,+}(Z)$ are sheaves. 

  \begin{lem}\label{lem:Regular}
      Let $A$ be a commutative ring and $a,b\in A$ be element which form a regular sequence in $A$. For any $k,l\in\bN$, we have an exact sequence of $A[X,Y]$-modules
      \[0\to A[X,Y]\xrightarrow{w\mapsto (-gw,fw)} A[X,Y]\oplus A[X,Y]\xrightarrow{(u,v)\mapsto uf+vg}A[X,Y]\xrightarrow{X\mapsto\frac{b^{k+l}}{a},Y\mapsto\frac{a}{b^k}} A[\frac{b^{k+l}}{a},\frac{a}{b^k}]\to 0,\]
      where $f = b^kY-a$ and $g = XY-b^l$. 
  \end{lem}
  \begin{proof}
      Recall that by definition,
      \[A[\frac{b^{k+l}}{a},\frac{a}{b^k}] \cong A[X,Y]/(f, g, h) = A[X,Y]/(f,g)\]
      where $h = aX-b^{k+l}$ and the last equality follows because $h+Xf-b^kg = 0$. We have to show that the kernel $\Ker$ of the morphism
      \[A[X,Y]\oplus A[X,Y]\xrightarrow{(u,v)\mapsto uf+vg}A[X,Y]\]
      coincides with the image of the morphism
      \[A[X,Y]\xrightarrow{w\mapsto (-gw,fw)} A[X,Y]\oplus A[X,Y].\]

      Suppose that $(u:=u(X,Y)= \sum_{i\geq 0}u_i(Y)X^i,v:=v(X,Y)= \sum_{i\geq 0}v_i(Y)X^i)\in \Ker$; that is, we have
      \begin{equation*}
          \begin{split}
              u(X,Y)(b^kY-a)+v(X,Y)(XY-b^l) = & \sum_{i\geq 0}u_i(Y)X^i(b^kY-a)+\sum_{i\geq 0}v_i(Y)X^i(XY-b^l)\\
              = &\sum_{i\geq 0}\big((b^kY-a)u_i(Y)-b^lv_i(Y)+Yv_{i-1}(Y)\big)X^i\\
              = & 0,
          \end{split}
      \end{equation*}
      where we put $v_{-1}(Y) = 0$. Equivalently, we have that for any $n\geq 0$,
      \begin{equation}\label{equ:Regular}
          b^sv_n(Y) = Yv_{n-1}(Y)+(b^kY-a)u_n(Y).
      \end{equation}
      As the sequence $a,b$ is regular in $A$, the sequence $b^l,b^kY-a$ is regular in $A[Y]$. By induction on $n$, using \eqref{equ:Regular}, we see that there is a unique sequence of polynomials $\{w_n(Y)\}_{n\geq -1}$ in $A[Y]$ such that for any $n\geq 0$,
      \[v_n(Y) = (b^kY-a)w_n(Y)\]
      and that
      \[u_n(Y) = b^lw_n(Y) - Yw_{n-1}(Y),\]
      where $w_{-1}(Y) = 0$. 
      Define $w:=w(X,Y):=\sum_{n\geq 0}w_n(Y)X^n\in A[X,Y]$. Then we have that
      \[v(X,Y) = (b^rY-a)w(X,Y) = fw\]
      and that 
      \begin{equation*}
          \begin{split}
              u(X,Y) = \sum_{n\geq 0}u_n(Y)X^n 
                     = &\sum_{n\geq 0}\big(b^lw_n(Y) - Yw_{n-1}(Y)\big)X^n \\
                     = &-(XY-b^l)\sum_{n\geq 0}w_n(Y)X^n \\
                     = &-gw.
          \end{split}
      \end{equation*}
      In other words, there exists a unique $w\in A[X,Y]$ such that $(u,v) = (-gw,fw)$. This completes the proof.
  \end{proof}
  \begin{cor}\label{cor:Regular}
      For any rational numbers $s>r>0$, the following sequence of sheaves of $\bA_{\inf}$-modules on the $v$-site $\Perfd_v$ is exact
      \[0\to \AAinf[X,Y]\xrightarrow{w\mapsto (-gw,fw)} \AAinf[X,Y]\oplus \AAinf[X,Y]\xrightarrow{(u,v)\mapsto uf+vg}\AAinf[X,Y]\xrightarrow{\iota} \AAinf[\frac{[\varpi^{\frac{1}{r}}]}{p},\frac{p}{[\varpi^{\frac{1}{s}}]}]\to 0,\]
      where $f = [\varpi^{\frac{1}{s}}]Y-p$, $g = XY-[\varpi^{\frac{1}{r}-\frac{1}{s}}]$ and $\iota$ is a morphism of $\bA_{\inf}$-algebras sending $X$ and $Y$ to $\frac{[\varpi^{\frac{1}{r}}]}{p}$ and $\frac{p}{[\varpi^{\frac{1}{s}}]}$ respectively.
  \end{cor}
  \begin{proof}
      Note that for any affinoid perfectoid $U\in \Perfd$, the sequence $p,[\varpi^{\frac{1}{rs}}]$ is regular in $\bA_{\inf}(U)$. Then the result follows from Lemma \ref{lem:Regular} immediately.
  \end{proof}
  
  \begin{prop}\label{prop:C^I is a sheaf}
      Let $I=[r,s]\subset(0,+\infty)$ be a closed sub-interval with $r,s\in\bQ$. Then $\CCtilde^I$ and $\CCtilde^{I,+}$ are both sheaves. Moreover, for any $U\in \Perfd$ and for any $n\geq 1$, $\rH^n(U,\CCtilde^I) = 0$ and $\rH^n(U,\CCtilde^{I,+})$ is almost zero in the sense that it is killed by $[\varpi^h]$ for any $h\in\bQ_{>0}$.
  \end{prop}
  \begin{proof}
      Let $V\to U$ be a covering in $\Perfd$ and denote by $V^\bullet\to U$ the associated \v Cech nerve. Then the complex $\bA_{\inf}(U)\to \bA_{\inf}(V^\bullet)$ is almost acyclic.
      By Lemma \ref{lem:description of C^I} and Corollary \ref{cor:Regular}, the complex $\CCtilde^{I,+}(U)\to\CCtilde^{I,+}(V^{\bullet})$ is almost acyclic as well. Inverting $p$, we see that $\CCtilde^I(U)\to\CCtilde^I(V^{\bullet})$ is acyclic, which implies that $\CCtilde^I$ is a sheaf. By Remark \ref{rmk:isometry on base}(2), we see that 
      \[\Ker(\CCtilde^{I,+}(V^0)\to\CCtilde^{I,+}(V^1))\subset \Ker(\CCtilde^I(V^0)\to\CCtilde^I(V^1)) = \CCtilde^I(U)\]
      is exactly $\CCtilde^{I,+}(U)$. So $\CCtilde^{I,+}$ is also a sheaf. The desired cohomology vanishing follows from the \v Cech-to-derived spectral sequence.
  \end{proof}

\section{Semi-stable formal schemes, liftings, and the ringed topos $\wtx_Z^I$}\label{sec:ringed space}
In this section we fix the semi-stable geometric framework used throughout the paper. We first recall liftable semi-stable formal schemes over $A^+$ and then attach, to each chosen lifting over $\bA_{\inf}(Z)$, a ringed topos $\wtx^I_Z = (\frakX_{Z,\et},\calO_{\wtx_Z^I})$ so that $\calO_{\wtx_Z^I}$ is a sheaf of $\CCtilde^I(Z)$-algebras.

\subsection{Review of liftable semi-stable formal schemes}\label{ssec:liftable semi-stable schemes}

Let $Z = \Spa(A,A^+)\in \Perfd$ be affinoid perfectoid as before. For any $n\geq 0$, by the \emph{canonical log-structure} on $\bA_{\inf}(Z)/\xi^{n+1}$, we mean the log-structure associated to the pre-log-structure 
\[A^{\flat\times}\cap A^{\flat+}\xrightarrow{[\cdot]}\bA_{\inf}(Z)/\xi^{n+1}.\]
Here, for a pre-log structure $\iota:M\to A$, its associated log-structure is given by $\iota^a:M^a\to A$, where $M^a$ denotes the push-out of the diagram $A^\times\xleftarrow{\iota}\iota^{-1}(A^\times)\xrightarrow{\subset}M$ of morphisms of monoids, and $\iota^a$ denotes the induced map of monoids.
When $n = 0$, this coincides with the log-structure $(A^{\times}\cap A^+\hookrightarrow A^+)$ (cf. \cite[Prop. 2.2]{SW24}). Following \cite{CK19} and \cite{SW24}, we define semi-stable formal schemes over $A^+$ as follows.

\begin{dfn}\label{dfn:semi-stable formal scheme}
    Let $Z = \Spa(A,A^+)\in\Perfd$ be as before.
    \begin{enumerate}
    \item[(1)] A $p$-complete $A^+$-algebra $\calR_Z$ is called \emph{small semi-stable} of relative dimension $d$ over $A^+$ if there exists an \'etale morphism of $p$-adic formal schemes
    \[\psi: \Spf(\calR_Z)\to \Spf(A^+\za T_0,\dots,T_r,T_{r+1}^{\pm 1},\dots, T_d^{\pm 1}\ya/(T_0\cdots T_r-p^a))\]
    for some $a\in\bQ_{\geq 0}$ and $1\leq r\leq d$, where $A^+\za T_0,\dots,T_r,T_{r+1}^{\pm 1},\dots, T_d^{\pm 1}\ya/(T_0\cdots T_r-p^a)$ denotes the $p$-adic completion of
    \[A^+[ T_0,\dots,T_r,T_{r+1}^{\pm 1},\dots, T_d^{\pm 1}]/(T_0\cdots T_r-p^a).\]
    We equip $\calR_Z$ with the log-structure associated to the pre-log-structure
    \begin{equation}\label{equ:log-structure on chart-I}
        M_{r,a}(A^+):=(\oplus_{i=0}^r\bN\cdot e_i)\oplus_{\bN\cdot(e_0+\cdots+e_r)}(A^{\times}\cap A^+)\xrightarrow{(\sum_{i=0}^rn_ie_i,x)\mapsto x\prod_{i=0}^rT_i^{n_i}} \calR_Z,
    \end{equation}
    where $(\oplus_{i=0}^r\bN\cdot e_i)\oplus_{\bN\cdot(e_0+\cdots+e_r)}(A^{\times}\cap A^+)$ denotes the monoid representing the push-out of the diagram
    \[\xymatrix@C=0.5cm{
      \bN\ar[d]^{1\mapsto p^a}\ar[rrr]^{1\mapsto e_0+\cdots+e_r\qquad}&&&\bN\cdot e_0\oplus\cdots\oplus\bN\cdot e_r\\
      (A^{\times}\cap A^+).
    }\]
    In this case, we also say that $\Spf(\calR_Z)$ is \emph{small affine}.  We call such a $\psi$ (resp. $T_0,\dots,T_d$) a \emph{chart} (resp. \emph{coordinates}) on $\calR_Z$ or on $\Spf(\calR_Z)$. Note that the generic fiber $\Spa(R_Z,R_Z^+)$ of $\Spf(\calR_Z)$ is smooth over $Z$ and endowed with the associated chart $\psi$, where $R_Z^+=\calR_Z$.

    \item[(2)] A \emph{semi-stable} formal scheme $\frakX_Z$ of relative dimension $d$ over $A^+$ is a $p$-adic formal scheme $\frakX_Z$ equipped with a log-structure $\calM_{\frakX_Z}$ over $A^+$ (endowed with the canonical log-structure), which is \'etale locally of the form $\Spf(\calR_Z)$ where $\calR_Z$ is small semi-stable of relative dimension $d$ over $A^+$ with the log-structure defined above. Note that the generic fiber $X_Z$ of $\frakX_Z$ is always smooth over $Z$.
    \end{enumerate}
\end{dfn}

We record two basic classes of semi-stable formal schemes over $A^+$.

\begin{exam}\label{exam:semi-stable formal scheme}
    \begin{enumerate}
        \item[(1)] Any smooth formal scheme $\frakX_Z$ over $A^+$ with the log-structure induced from the canonical log-structure on $A^+$ (i.e., the log-structure associated to the composite $(A^{\times}\cap A^+)\hookrightarrow A^+\to \calO_{\frakX_Z}$) is always semi-stable.
        \item[(2)] Let $\frakX$ be a semi-stable formal scheme over $\calO_C$ in the sense of \cite{CK19} with the log-structure $(\calM_{\frakX} = \calO_X^{\times}\cap\calO_{\frakX}\hookrightarrow\calO_{\frakX})$ as considered in \cite[\S1.5 and \S1.6]{CK19}, where $X$ is the generic fiber of $\frakX$. For any $Z = \Spa(A,A^+)$, denote by $\frakX_Z$ the base-change of $\frakX$ along $\Spf(A^+)\to\Spf(\calO_C)$ with the log-structure $\calM_{\frakX_Z}$ induced from the fiber product. Then $\frakX_Z$ is a semi-stable formal scheme over $A^+$. 
    \end{enumerate}
\end{exam}


\begin{dfn}\label{dfn:liftablity}
    Let $Z = \Spa(A,A^+)\in\Perfd$ be as before.
    \begin{enumerate}
        \item[(1)] Let $(\frakX_Z,\calM_{\frakX_Z})$ be a semi-stable formal scheme over $A^+$. By a \emph{lifting of $\frakX_Z$ over $\bA_{\inf}(Z)$}, we mean a compatible sequence $\{(\wtx_{Z,n},\calM_{\wtx_{Z,n}})\}_{n\geq 0}$ of flat $p$-adic formal schemes $\wtx_{Z,n}$ with log-structure $\calM_{\wtx_{Z,n}}$ over $\bA_{\inf}(Z)/\xi^{n+1}$ (equipped with the canonical log-structure) satisfying the following conditions:
        \begin{enumerate}
            \item[$\bullet$] $(\wtx_{Z,0},\calM_{\wtx_{Z,0}}) = (\frakX_Z,\calM_{\frakX_Z})$.
            \item[$\bullet$] For any $n\geq m\geq 0$, we have 
            \[\calO_{\wtx_{Z,m}} = \calO_{\wtx_{Z,n}}\otimes_{\bA_{\inf}(Z)/\xi^{n+1}}\bA_{\inf}(Z)/\xi^{m+1}\]
            such that the log-structure $\calM_{\wtx_{Z,m}}$ is induced from the composite
            \[\calM_{\wtx_{Z,n}}\to\calO_{\wtx_{Z,n}}\to\calO_{\wtx_{Z,m}}.\]
        \end{enumerate}
        Define 
        \[(\calO_{\wtx_Z},\calM_{\wtx_Z}):=\varprojlim_n(\calO_{\wtx_{Z,n}},\calM_{\wtx_{Z,n}}).\]
        Then $\calO_{\wtx_Z}$ is a sheaf of flat algebras over $\bA_{\inf}(Z)$ and $\calM_{\wtx_Z}$ defines a log-structure on $\calO_{\wtx_Z}$. We denote the lifting $\{(\wtx_{Z,n},\calM_{\wtx_{Z,n}})\}_{n\geq 0}$ by $(\wtx_Z,\calM_{\wtx_Z})$, or simply as a lifting of $\frakX_Z$ over $\bA_{\inf}(Z)$. We say a semi-stable formal scheme $\frakX$ is \emph{liftable} if it admits at least one lifting $(\wtx_Z,\calM_{\wtx_Z})$ over $\bA_{\inf}(Z)$.

        \item[(2)] Let $X_Z$ be a smooth adic space over $A$. By a \emph{lifting} of $X_Z$ over $\BBdRp(Z)$, we mean a compatible sequence $\{\wtX_{Z,n}\}_{n\geq 0}$ of smooth adic spaces $\wtX_{Z,n}$ over $\BBdRp(Z)/\xi^{n+1}$ satisfying the following conditions:
        \begin{enumerate}
            \item[$\bullet$] $\wtX_{Z,0} = X_Z$.
            \item[$\bullet$] For any $n\geq m\geq 0$, we have 
            \[\calO_{\wtX_{Z,m}} = \calO_{\wtX_{Z,n}}\otimes_{\BBdRp(Z)/\xi^{n+1}}\BBdRp(Z)/\xi^{m+1}.\]
            Define 
            \[\calO_{\wtX_Z}:=\varprojlim_n\calO_{\wtX_{Z,n}}.\]
            Then $\calO_{\wtX_Z}$ is a sheaf of flat $\BBdRp(Z)$-algebras. We denote the lifting $\{\wtX_{Z,n}\}_{n\geq 0}$ by $\wtX_Z$.
        \end{enumerate}
    \end{enumerate}
\end{dfn}

It follows directly that for a liftable semi-stable formal scheme $\frakX_Z$ over $A^+$ with generic fiber $X_Z$, any lifting $\wtx_Z$ of $\frakX_Z$ over $\bA_{\inf}(Z)$ induces a lifting $\wtX_Z$ of $X_Z$ over $\BBdRp(Z)$ by letting $\wtX_{Z,n}$ be the base-change of $\wtx_{Z,n}$ along $\bA_{\inf}(Z)\to \BBdRp(Z)$. 
We shall work only with semi-stable formal schemes over $A^+$ that admit such liftings.

\begin{convention}\label{convention:tilt}
    Let $\frakX_Z$ be a semi-stable formal scheme over $A^+$ of relative dimension $d$ with a fixed lifting $\wtx_Z$ over $\bA_{\inf}(Z)$. Denote by $\frakX_Z^{\flat}$ the reduction of $\wtx_Z$ modulo $p$ and equip it with the induced log-structure. Then $\frakX_Z^{\flat}$ is a flat $\varpi$-adic formal scheme over $A^{\flat +}$ and its generic fiber $X_Z^{\flat}$ is smooth over $A^{\flat}$. 
\end{convention}

\begin{exam}\label{exam:small affine is liftable}
    Suppose $\frakX_Z = \Spf(\calR_Z)$ is small affine in the sense of Definition \ref{dfn:semi-stable formal scheme}. Then $\frakX_Z$ is liftable. Indeed, by \cite[Lem. 3.1]{SW24}, the $p$-complete log-cotangent complex (in the sense of \cite[\S 8]{Ols}) of $\calR_Z$ over $A^+$ is
    \[\widehat \rL_{(M_{r,a}(A^+)\to \calR_Z)/(A^{\times}\cap A^+\to A^+)} \simeq \Omega^{1,\log}_{\calR_Z}[0],\]
    where $\Omega^{1,\log}_{\calR_Z}$ denotes the module of continuous log-differentials of $\calR_Z$ over $A^+$, which is a finite free $\calR_Z$-module of rank $d$ and isomorphic to
    \begin{equation}
       \left((\oplus_{i=0}^r\calR_Z\cdot e_i)/\calR_Z\cdot (e_0+\cdots+e_r)\right)\oplus\left(\oplus_{j=r+1}^d\calR_Z\cdot e_j\right),
    \end{equation}
    where $e_j$ stands for $\dlog T_j$ for any $r+1\leq j\leq d$. By deformation theory, up to non-canonical isomorphisms, $\frakX_Z$ admits a unique lifting $\wtx_Z$ over $\bA_{\inf}(Z)$. More precisely, consider the log-ring 
    \[M_{r,a}(Z) = (\oplus_{i=0}^r\bN\cdot e_i)\oplus_{\bN\cdot(e_0+\cdots+e_r)}M_Z\xrightarrow{\alpha} \bA_{\inf}(Z)\za T_0,\dots,T_r,T_{r+1}^{\pm 1},\dots,T_d^{\pm 1}\ya/(T_0\cdots T_r-[\varpi^a]),\]
    where $M_Z$ denotes the canonical log-structure on $\bA_{\inf}(Z)$, the map $\alpha$ of monoids is determined by sending each $e_i$ to $T_i$, and $M_{r,a}(Z)$ denotes the monoid representing the push-out of the diagram
    \[\xymatrix@C=0.5cm{
      \bN\ar[rrrrrr]^{1\mapsto e_0+\cdots+e_r}\ar[d]_{1\mapsto[\varpi^a]}&&&&&&\bN\cdot e_0\oplus\cdots\oplus\bN\cdot e_r\\
      M_Z.
    }\]
    By the \'etaleness of the chart $\psi$, there exists a unique lifting $\wtr_Z$ of $\calR_Z$ over
    \[\bA_{\inf}(Z)\za T_0,\dots,T_r,T_{r+1}^{\pm 1},\dots,T_d^{\pm 1}\ya/(T_0\cdots T_r-[\varpi^a]),\]
    together with a unique lifting 
    \[\widetilde \psi:\bA_{\inf}(Z)\za T_0,\dots,T_r,T_{r+1}^{\pm 1},\dots,T_d^{\pm 1}\ya/(T_0\cdots T_r-[\varpi^a])\to\wtr_Z\]
    of $\psi$ such that the log-structure on $\wtr_Z$ is associated to the composite
    \[M_{r,a}(Z)\xrightarrow{\alpha} \bA_{\inf}(Z)\za T_0,\dots,T_r,T_{r+1}^{\pm 1},\dots,T_d^{\pm 1}\ya/(T_0\cdots T_r-[\varpi^a])\to\wtr_Z.\]
    Here, $\bA_{\inf}(Z)\za T_0,\dots,T_r,T_{r+1}^{\pm 1},\dots,T_d^{\pm 1}\ya/(T_0\cdots T_r-[\varpi^a])$ denotes the $(p,\xi)$-adic completion of 
    \[\bA_{\inf}(Z)[ T_0,\dots,T_r,T_{r+1}^{\pm 1},\dots,T_d^{\pm 1}]/(T_0\cdots T_r-[\varpi^a]).\]
    Then the lifting $\wtx_Z$ of $\frakX_Z$ is given by $\{\Spf(\wtr_{Z,n})\}_{n\geq 0}$ such that the log-structure on $\wtr_{Z,n}:=\wtr_Z/\xi^{n+1}$ is associated to the composite $M_{r,a}(Z)\to \wtr_Z\to\wtr_Z/\xi^{n+1}$, and the associated lifting $\wtX_Z$ of $X_Z = \Spa(R_Z,R_Z^+ = \calR_Z)$ is given by $\{\Spa(\wtR_{Z,n} := \wtr_{Z,n}[\frac{1}{p}],\wtr_{Z,n})\}_{n\geq 0}$.

    In what follows, we simply denote this lifting of $\frakX_Z$ by $\wtx_Z = \Spf(\wtr_Z)$ and call $\wtR_Z:=\varprojlim_n\wtR_{Z,n}$ the associated lifting of $R_Z$ over $\BBdRp(Z)$. Note that in this case, the $\varpi$-adic formal scheme $\frakX_Z^{\flat}$ is of the form $\frakX_Z^{\flat} = \Spf(\calR_Z^{\flat})$ for $\calR_Z^{\flat} = \wtr_Z/p$, and the morphism $\widetilde \psi$ induces a $\varpi$-completely \'etale morphism
    \[\psi^{\flat}:A^{\flat+}\za T_1,\dots,T_r,T_{r+1}^{\pm 1},\dots,T_d^{\pm 1}\ya/(T_0\cdots T_r-\varpi^a)\to\calR_Z^{\flat}.\]
    Let $M_{r,a}(A^{\flat+})$ be the monoid associated to the push-out of the diagram
    \[\xymatrix@C=0.5cm{
      \bN\ar[rrrrrr]^{1\mapsto e_0+\cdots+e_r}\ar[d]_{1\mapsto\varpi^a}&&&&&&\bN\cdot e_0\oplus\cdots\oplus\bN\cdot e_r\\
      A^{\flat\times}\cap A^{\flat+},
    }\]
    and 
    \[\alpha^{\flat}:M_{r,a}(A^{\flat+})\to A^{\flat+}\za T_1,\dots,T_r,T_{r+1}^{\pm 1},\dots,T_d^{\pm 1}\ya/(T_0\cdots T_r-\varpi^a)\]
    be the morphism of monoids induced by sending each $e_i$ to $T_i$ for any $0\leq i\leq r$. Then the log-structure on $\calR_Z^{\flat}$ is associated to the composite
    \[M_{r,a}(A^{\flat+})\xrightarrow{\alpha^{\flat}} A^{\flat+}\za T_1,\dots,T_r,T_{r+1}^{\pm 1},\dots,T_d^{\pm 1}\ya/(T_0\cdots T_r-\varpi^a)\xrightarrow{\psi^{\flat}}\calR_Z^{\flat}.\]
\end{exam}
\begin{rmk}
    Keep notations in Example \ref{exam:small affine is liftable}. We claim that $\wtr_Z$ is the only lifting of $\calR_Z$ over $\bA_{\inf}(Z)$ that is compatible with log-structure. For simplicity, we assume that 
    \[\calR_Z = A^+\za T_0,\dots,T_r,T_{r+1}^{\pm 1},\dots,T_d^{\pm 1}\ya/(T_0\cdots T_r-p^a).\]
    Since the canonical log-structure on $A^+$ is perfect in the sense of \cite[Definition 2.6]{SW24}, the only log-structure on $\bA_{\inf}(Z)$ lifting the canonical log-structure on $A^+$ along $\theta$ is also the canonical one by \cite[Lemma 2.7]{SW24}. By \cite[Lemma 3.1]{SW24} says the $p$-complete log-cotangent complex $\widehat \rL_{(M_{r,a}(Z)\to\calR_Z)/(A^{\times}\cap A^+\to A^+)}$ is discrete and isomorphic to $\Omega^{1,\log}_{\calR_Z/A^+}$, and thus finite free over $\calR_Z$. So by deformation theory, 
    \[\wtr_Z/\xi^2=\bA_{\inf}(Z)\za T_0,\dots,T_r,T_{r+1}^{\pm 1},\dots,T_d^{\pm 1}\ya/(T_0\cdots T_r-[\varpi^a])/\xi^2\]
    is the only square-zero lifting over $\bA_{\inf}(Z)/\xi^2$ as log scheme (up to non-canonical isomorphisms). By a similar argument, we see that for any $n\geq 0$, 
    \[\wtr_Z/\xi^n=\bA_{\inf}(Z)\za T_0,\dots,T_r,T_{r+1}^{\pm 1},\dots,T_d^{\pm 1}\ya/(T_0\cdots T_r-[\varpi^a])/\xi^n\]
    is the only correct lifting over $\bA_{\inf}(Z)/\xi^n$ (up to non-canonical isomorphisms).
\end{rmk}

\subsection{The ringed topos $\wtx^I$}\label{ssec:ringed space X^I}

Fix an affinoid perfectoid space $Z = \Spa(A,A^+)\in \Perfd$ and let $I=[r,s]$ be a closed sub-interval of $(0,+\infty)$ satisfying Hypothesis \ref{hy:rationality} as before. Let $\frakX_Z$ be a liftable semi-stable formal scheme over $A^+$ with generic fiber $X_Z$. We will assign a ringed topos $\wtx_Z^I$ over $\CCtilde^I(Z)$ to each lifting $\wtx_Z$ of $\frakX_Z$.

We first assume $\frakX_Z = \Spf(\calR_Z)$ is small semi-stable with a fixed lifting $\wtx_Z = \Spf(\wtr_Z)$.

\begin{dfn}\label{dfn:R^I}
    For any closed sub-interval $I\subset(0,+\infty)$, define
    \[\wtr_Z^{I,+}:=\wtr_Z\widehat \otimes_{\bA_{\inf}(Z)}\CCtilde^{I,+}(Z) \text{ and }\wtr_Z^I:=\wtr_Z^{I,+}[\frac{1}{p}],\]
    where $\widehat \otimes_{\bA_{\inf}(Z)}$ denotes the $(p,\xi)$-adically complete tensor product over $\bA_{\inf}(Z)$.
\end{dfn}

\begin{lem}\label{lem:sheaf}
    The presheaves $\calO_{\wtx_Z^{I,+}}$ and $\calO_{\wtx_Z^I}$ on $\frakX_{Z,\et}$ defined by
    \[\frakY_Z = \Spf(\calR_Z^{\prime})\mapsto \wtr_Z^{\prime I,+} \text{ and }\wtr_Z^{\prime I}, \text{ respectively,}\]
    are both sheaves. Moreover, for any affine formal scheme $\frakY_Z = \Spf(\calR_Z^{\prime})\in \frakX_{Z,\et}$ and $n\geq 0$,
    \[\rH^n(\frakY_Z,\calO_{\wtx_Z^{I,+}}) = \rH^n(\frakY_Z,\calO_{\wtx_Z^I}) = 0.\]
\end{lem}
\begin{proof}
    Let $\calR_{Z,1}\to \calR_{Z,2}$ be a morphism of $A^+$-algebras such that the associated morphism $\Spf(\calR_{Z,2})\to\Spf(\calR_{Z,1})$ of $p$-adic formal schemes over $A^+$ is a covering in $\frakX_{Z,\et}$. 
    Denote by $\wtr_{Z,1}\to\wtr_{Z,2}$ the corresponding lifting of $\calR_{Z,1}\to \calR_{Z,2}$ over $\bA_{\inf}(Z)$. Denote by $\wtr_{Z,1}\to \wtr_{Z,2}^{\bullet}$ the corresponding \v Cech nerve associated to the morphism $\wtr_{Z,1}\to\wtr_{Z,2}$. Then we have
    \[\wtr_{Z,1}\simeq \varprojlim\wtr_{Z,2}^{\bullet}.\] 
    Note that both $\wtr_{Z,1}$ and $\wtr_{Z,2}^{\bullet}$ are $(p,\xi)$-completely faithfully flat over $\bA_{\inf}(Z)$. The derived $(p,\xi)$-complete tensor products $\wtr_{Z,1}\widehat \otimes_{\bA_{\inf}(Z)}^{\rL}\CCtilde^{I,+}(Z)$ and $\wtr_{Z,2}^{\bullet}\widehat \otimes_{\bA_{\inf}(Z)}^{\rL}\CCtilde^{I,+}(Z)$ are classical. So we have
    \[\wtr_{Z,1}\widehat \otimes_{\bA_{\inf}(Z)}\CCtilde^{I,+}(Z)\simeq\varprojlim(\wtr_{Z,2}^{\bullet}\widehat \otimes_{\bA_{\inf}(Z)}\Ctilde^{I,+}(Z)).\]
    That is, the result holds true for $\calO_{\wtx^{I,+}}$. As filtered colimits always commute with finite limits, by taking naive truncations of 
    \[\varprojlim(\wtr_{Z,2}^{\bullet}\widehat \otimes_{\bA_{\inf}(Z)}\CCtilde^{I,+}(Z))\]
    and inverting $p$, we see that the result holds true for $\calO_{\wtx^I}$.
\end{proof}

Now, we define the desired ringed topos $\wtx_Z^I$ for a general liftable semi-stable $\frakX_Z$ together with a lifting $\wtx_Z$.

\begin{dfn}\label{dfn:X^I}
    Suppose that $\frakX_Z$ is a liftable semi-stable formal scheme over $A^+$ with a fixed lifting $\wtx_Z$ over $\bA_{\inf}(Z)$. Denote by $\frakX_{Z,\et}^{\Box}$ the sub-site of $\frakX_{Z,\et}$ consisting of objects $\frakY_Z \to\frakX_Z$ such that $\frakY_Z = \Spf(\calR_Z^{\prime})$ is small semi-stable. Denote by $\calO_{\wtx_Z^{I,+}}$ and $\calO_{\wtx_Z^I}$ the sheafification of the presheaves on $\frakX_{Z,\et}$ sending each $\frakY_Z = \Spf(\calR_Z^{\prime})\in \frakX_{Z,\et}^{\Box}$ to
    \[\wtr_Z^{\prime I,+} \text{ and }\wtr_Z^{\prime I},\]
    where $\wtr_Z^{\prime}$ is the lifting of $\calR_Z^{\prime}$ over $\bA_{\inf}(Z)$ induced by the given lifting $\wtx_Z$ of $\frakX_Z$ (cf. Lemma \ref{lem:sheaf}) and we denote by $\wtx_Z^{I,+}$ and $\wtx_Z^I$ the ringed topoi $(\frakX_{Z,\et},\calO_{\wtx_Z^{I,+}})$ and $(\frakX_{Z,\et},\calO_{\wtx_Z^I})$, respectively.
\end{dfn}

For any sheaf $\calE$ of locally finite free $\calO_{\wtx_Z}$-modules on $\frakX_{Z,\et}$, the base-change $\calE\otimes_{\calO_{\wtx_Z}}\calO_{\wtx_Z^{I,+}}$ (resp. $\calE\otimes_{\calO_{\wtx_Z}}\calO_{\wtx_Z^I}$ ) is a sheaf of locally finite free $\calO_{\wtx_Z^{I,+}}$-modules (resp. $\calO_{\wtx_Z^I}$-modules) on $\wtx_{Z}^{I,+}$ (resp. $\wtx_Z^I$). For example, we have 
\[\Omega^{1,\log}_{\wtx_Z^{I,+}} = \Omega^{1,\log}_{\wtx_Z}\otimes_{\calO_{\wtx_Z}}\calO_{\wtx_Z^{I,+}} \text{ (resp. $\Omega^{1,\log}_{\wtx_Z^I} = \Omega^{1,\log}_{\wtx_Z}\otimes_{\calO_{\wtx_Z}}\calO_{\wtx_Z^I}$)}.\]
The usual $(p,\xi)$-complete differential $\rd$ on $\wtx_Z$ then induces a $\CCtilde^{I,+}(Z)$-linear (resp. $\CCtilde^I(Z)$-linear) differential
\[\rd:\calO_{\wtx_Z^{I,+}}\to \Omega^{1,\log}_{\wtx_Z^{I,+}} \text{ (resp. $\rd:\calO_{\wtx_Z^I}\to \Omega^{1,\log}_{\wtx_Z^I}$)}\]
and furthermore a connection
\[\nabla:\calO_{\wtx_Z^{I,+}}\to \Omega^{1,\log}_{\wtx_Z^{I,+}}\{-1\} \text{ (resp. $\nabla:\calO_{\wtx_Z^I}\to \Omega^{1,\log}_{\wtx_Z^I}\{-1\}$)}\]
by composing $\rd$ with the inclusion 
\[\Omega^{1,\log}_{\wtx_Z^{I,+}} = \xi\cdot \Omega^{1,\log}_{\wtx_Z^{I,+}}\{-1\}\subset \Omega^{1,\log}_{\wtx_Z^{I,+}}\{-1\} \text{ (resp. $\Omega^{1,\log}_{\wtx_Z^I} = \xi\cdot \Omega^{1,\log}_{\wtx_Z^I}\{-1\}\subset \Omega^{1,\log}_{\wtx_Z^I}\{-1\}$)}.\]

\begin{rmk}\label{rmk:xi-connection}
    When $\frakX_Z = \Spf(\calR_Z)$ is small semi-stable, via the isomorphism $\Omega^{1,\log}_{\wtr_Z^{I,+}}\{-1\}\cong \oplus_{i=1}^d\wtr_Z^{I,+}\cdot e_i$ corresponding to the chart $\psi$ (together with its lifting $\widetilde \psi$, cf. Example \ref{exam:small affine is liftable}), where $e_i$ stands for $\frac{\dlog T_i}{\xi}$ for any $1\leq i\leq d$, we have by construction that $\nabla:\wtr_Z^{I,+}\to\Omega^{1,\log}_{\wtr_Z^{I,+}}\{-1\}\cong \oplus_{i=1}^d\wtr_Z^{I,+}\cdot e_i$ is of the form
    \[\nabla = \sum_{i=1}^d\xi\frac{\partial}{\partial T_i}\otimes\frac{e_i}{\xi}.\]
\end{rmk}

The assignment $(\frakX_Z,\wtx_Z)\mapsto \wtx_Z^I$ is functorial in the following sense: Let $f:\frakX^{\prime}_Z\to\frakX_Z$ be a morphism of liftable semi-stable formal schemes over $A^+$. Fix liftings $\wtx_Z^{\prime}$ and $\wtx_Z$ of $\frakX^{\prime}_Z$ and $\frakX_Z$ respectively, and suppose that $f$ admits a lifting $\widetilde f:\wtx^{\prime}_Z\to\wtx_Z$ over $\bA_{\inf}(Z)$. Then $\widetilde f$ induces a morphism of ringed topoi
\[\widetilde f^I:\wtx_Z^{\prime I}\to \wtx_Z^I.\]
Moreover, the relative log-differential $\rd:\calO_{\wtx_Z^{\prime}}\to \Omega^{1,\log}_{\wtx_Z^{\prime}/\wtx_Z}$ induces a relative differential
\[\rd:\calO_{\wtx_Z^{\prime I,+}}\to \Omega^{1,\log}_{\wtx_Z^{\prime I,+}/\wtx_Z^{I,+}} \text{ (resp. $\rd:\calO_{\wtx_Z^{\prime I}}\to \Omega^{1,\log}_{\wtx_Z^{\prime I}/\wtx_Z^I}$)}\]
and furthermore a relative connection
\[\nabla:\calO_{\wtx_Z^{\prime I,+}}\to \Omega^{1,\log}_{\wtx_Z^{\prime I,+}/\wtx_Z^{I,+}}\{-1\} \text{ (resp. $\nabla:\calO_{\wtx_Z^{\prime I}}\to \Omega^{1,\log}_{\wtx_Z^{\prime I}/\wtx_Z^I}\{-1\}$)}.\]

\subsection{The norm $|\cdot|_I$ on $\wtx_Z$}

In the following, for any affinoid perfectoid $U = \Spa(S,S^+)\in \Perfd$ and for any $n\geq 0$, we always equip $\bA_{\inf}(U)/\xi^{n+1}$ with the canonical log-structure (i.e. the log-structure associated to $S^{\flat\times}\cap S^{\flat+}\xrightarrow{[\cdot]}\bA_{\inf}(U)/\xi^{n+1}$).
From now on, we again assume $\frakX_Z = \Spf(\calR_Z)$ is small semi-stable. Then (up to non-canonical isomorphisms) the lifting $\frakX_Z$ is given by $\wtx_Z = \Spf(\wtr_Z)$ such that the log-structure is associated to the pre-log-structure $M_{r,a}(Z)\to \wtr_Z$ as described in Example \ref{exam:small affine is liftable}, and we keep notations there. We are going to construct a norm $|\cdot|_I$ on $\wtr_Z^I$ for any closed sub-interval $I=[r,s]\subset(0,+\infty)$.

\begin{construction}\label{construction:I-norm}
    Keep notations as in Example \ref{exam:small affine is liftable}. Let
    \[A^{\flat+}\za T_0^{\frac{1}{p^{\infty}}},\dots,T_r^{\frac{1}{p^{\infty}}},T_{r+1}^{\pm \frac{1}{p^{\infty}}},\dots,T_d^{\pm \frac{1}{p^{\infty}}}\ya/(T_0^{\frac{1}{p^{\infty}}}\cdots T_r^{\frac{1}{p^{\infty}}}-\varpi^{\frac{a}{p^{\infty}}})\]
    be the $\varpi$-adic completion of 
    \[\cup_{n\geq0}A^{\flat+}\za T_0^{\frac{1}{p^{n}}},\dots,T_r^{\frac{1}{p^{n}}},T_{r+1}^{\pm \frac{1}{p^{n}}},\dots,T_d^{\pm \frac{1}{p^{n}}}\ya/(T_0^{\frac{1}{p^n}}\cdots T_r^{\frac{1}{p^n}}-\varpi^{\frac{a}{p^n}}),\]
    which coincides with the perfection 
    \[\colim_{x\mapsto x^p}A^{\flat+}\za T_1,\dots,T_r,T_{r+1}^{\pm 1},\dots,T_d^{\pm 1}\ya/(T_0\cdots T_r-\varpi^a)\]
    of $A^{\flat+}\za T_1,\dots,T_r,T_{r+1}^{\pm 1},\dots,T_d^{\pm 1}\ya/(T_0\cdots T_r-\varpi^a)$. Put
    \[\widehat R_{\infty,Z}^{\flat+}:=\calR_Z^{\flat}\widehat \otimes_{\psi^{\flat},A^{\flat+}\za T_1,\dots,T_r,T_{r+1}^{\pm 1},\dots,T_d^{\pm 1}\ya/(T_0\cdots T_r-\varpi^a)}A^{\flat+}\za T_0^{\frac{1}{p^{\infty}}},\dots,T_r^{\frac{1}{p^{\infty}}},T_{r+1}^{\pm \frac{1}{p^{\infty}}},\dots,T_d^{\pm \frac{1}{p^{\infty}}}\ya/(T_0^{\frac{1}{p^{\infty}}}\cdots T_r^{\frac{1}{p^{\infty}}}-\varpi^{\frac{a}{p^{\infty}}}).\]
    As $A^{\flat+}\za T_0^{\frac{1}{p^{\infty}}},\dots,T_r^{\frac{1}{p^{\infty}}},T_{r+1}^{\pm \frac{1}{p^{\infty}}},\dots,T_d^{\pm \frac{1}{p^{\infty}}}\ya/(T_0^{\frac{1}{p^{\infty}}}\cdots T_r^{\frac{1}{p^{\infty}}}-\varpi^{\frac{a}{p^{\infty}}})$ is a perfectoid $A^{\flat+}$-algebra, so is $\widehat R_{\infty,Z}^{\flat+}$ because $\psi^{\flat}$ is \'etale. So
    \[X_{\infty,Z}^{\flat}:=\Spa(\widehat R_{\infty,Z}^{\flat}:=\widehat R_{\infty,Z}^{\flat+}[\frac{1}{\varpi}],\widehat R_{\infty,Z}^{\flat+})\]
    is affinoid perfectoid over $Z^{\flat}$. Moreover, we have that
    \[\widehat R_{\infty,Z}^{\flat+}=\big(\colim_{x\mapsto x^p}\calR_Z^{\flat}\big)^{\wedge_{\varpi}}\] 
    coincides with the $\varpi$-adic completion of the perfection of $\calR_Z^{\flat}$.
    
    By construction, $\calR_Z^{\flat}$ is a direct summand of $\widehat R_{\infty,Z}^{\flat+}$ as a sub-$\calR_Z^{\flat}$-module. Indeed, we have
    \begin{equation}\label{equ:decomposition mod p}
        \widehat R_{\infty,Z}^{\flat+} = \widehat \bigoplus_{\underline \alpha\in J_{r,d}}\calR_Z^{\flat}\cdot \underline T^{\underline \alpha},
    \end{equation}
    where $\widehat \bigoplus$ denotes the $\varpi$-adically complete direct sum, and for any 
    \begin{equation}\label{equ:index set}
        \underline \alpha = (\alpha_0,\dots,\alpha_d)\in J_{r,d}:=\{(\alpha_0,\dots,\alpha_d)\in \big(\bN[\frac{1}{p}]\cap[0,1)\big)^{d+1}\mid \alpha_0\cdots\alpha_r = 0\},
    \end{equation}
    the symbol $\underline T^{\underline \alpha}$ denotes the product $\prod_{i=0}^dT_i^{\alpha_i}\in \widehat R_{\infty,Z}^{\flat+}$.
    Using a similar argument as in the proof of \cite[Lem. 3.1]{SW24}, we see that
    \[\widehat \rL_{(\widehat R_{\infty,Z}^{\flat\times}\cap\widehat R_{\infty,Z}^{\flat+}\to \widehat R_{\infty,Z}^{\flat+})/(M_{r,a}(A^{\flat+})\to \calR_Z^{\flat})}\cong \Omega^{1,\log}_{\calR_Z^{\flat}}[0],\]
    where $\widehat \rL_{(\widehat R_{\infty,Z}^{\flat\times}\cap\widehat R_{\infty,Z}^{\flat+}\to \widehat R_{\infty,Z}^{\flat+})/(M_{r,a}(A^{\flat+})\to \calR_Z^{\flat})}$ denotes the $\varpi$-complete logarithmic cotangent complex and $\Omega^{1,\log}_{\calR_Z^{\flat}}$ denotes the module of $\varpi$-completely continuous logarithmic differentials of $\calR_Z^{\flat}$ over $A^{\flat+}$, which is a finite free $\calR_Z^{\flat}$-module of rank $d$. Thus, there is a morphism of $\bA_{\inf}(Z)$-algebras which is compatible with log-structures
    \[\widetilde \beta:\wtr_Z\to \rW(\widehat R_{\infty,Z}^{\flat+})\]
    lifting the natural inclusion $\calR_Z^{\flat}\to \widehat R_{\infty,Z}^{\flat+}$, and such a lifting is unique up to (non-canonical) isomorphisms on $\wtr_Z$.
    
    By tilting equivalence, there is a unique affinoid perfectoid space $X_{\infty,Z} = \Spa(\widehat R_{\infty},\widehat R_{\infty}^+)$ over $Z$ whose tilting is $X_{\infty,Z}^{\flat}$. More precisely, we have 
    \[\bA_{\inf}(X_{\infty,Z}) = \rW(\widehat R_{\infty,Z}^{\flat+}).\]
    By \eqref{equ:decomposition mod p}, $\widetilde \beta$ identifies $\wtr_Z$ with a direct summand of $\bA_{\inf}(X_{\infty,Z})$ as a sub-$\wtr_Z$-module such that
    \begin{equation}\label{equ:decomposition over Ainf}
        \bA_{\inf}(X_{\infty,Z}) = \widehat \bigoplus_{\underline \alpha\in J_{r,d}}\wtr_Z\cdot[\underline T^{\underline \alpha}],
    \end{equation}
    where $\widehat \bigoplus$ denotes the $(p,\xi)$-adically direct sum. After $(p,\xi)$-complete base-change along $\bA_{\inf}(Z)\to \CCtilde^{I,+}(Z)$ and inverting $p$, we get morphisms of $\CCtilde^{I,+}(Z)$-algebras
    \[\wtr_Z^{I,+}\to\CCtilde^{I,+}(X_{\infty,Z})\text{ and }\wtr_Z^I\to\CCtilde^I(X_{\infty,Z}),\]
    which will still be denoted by $\widetilde \beta$ and, respectively, induce decompositions
    \begin{equation}\label{equ:decomposition over C^Iplus}
        \CCtilde^{I,+}(X_{\infty,Z}) = \widehat \bigoplus_{\underline \alpha\in J_{r,d}}\wtr_Z^{I,+}\cdot[\underline T^{\underline \alpha}]
    \end{equation}
    and 
    \begin{equation}\label{equ:decomposition over C^I}
        \CCtilde^I(X_{\infty,Z}) = \widehat \bigoplus_{\underline \alpha\in J_{r,d}}\wtr_Z^I\cdot[\underline T^{\underline \alpha}].
    \end{equation}

    Now, for any such $\widetilde \beta:\wtr_Z^I\to\CCtilde^I(X_{\infty,Z})$, define
    \[|\cdot|_{I,\widetilde \beta}:\wtr_Z^I\to \bR_{\geq 0}\]
    by letting $|x|_{I,\widetilde \beta}:=|\widetilde \beta(x)|_I$ for any $x\in \wtr_Z$. As $\wtr_Z^{I,+} = \wtr_Z^I\cap\CCtilde^{I,+}(X_{\infty,Z})$, we see that for any $x\in \wtr_Z^I$, it belongs to $\wtr_Z^{I,+}$ if and only if $|x|_{I,\widetilde \beta}\leq 1$.
\end{construction}

We provide examples for constructing a morphism $\widetilde \beta:\wtr_Z\to\bA_{\inf}(X_{\infty,Z})$ of $\bA_{\inf}(Z)$-algebras compatible with log-structures as mentioned above.

\begin{exam}\label{exam:construct beta}
    Keep notations as in Example \ref{exam:small affine is liftable}. Suppose that there exists a homomorphism $\phi:\wtr_Z\to\wtr_Z$ lifting the absolute Frobenius endomorphism of $\calR_Z^{\flat} = \wtr_Z/p$ which is compatible with log-structures. (This amounts to equipping $\wtr_Z$ with a $\delta_{\log}$-structure in the sense of \cite{Kos20}.) Such a $\phi$ always exists.
    For example, we may first equip 
    \[\bA_{\inf}(Z)\za T_0,\dots,T_r,T_{r+1}^{\pm 1},\dots,T_d^{\pm 1}\ya/(T_0\cdots T_r-[\varpi^a])\]
    with the $\delta$-structure compatible with that on $\bA_{\inf}(Z)$ such that $\delta(T_i) = 0$ for any $0\leq i\leq d$, and then by \cite[Lem. 2.18]{BS22}, there exists a unique $\delta$-structure on $\wtr_Z$ such that the $(p,\xi)$-completely \'etale morphism $\widetilde \psi$ is compatible with $\delta$-structures. One can directly check that the Frobenius morphism on $\wtr_Z$ associated to this $\delta$-structure is the desired one. (Note that $(\wtr_Z,(\xi))$ defines a prism in the sense of \cite{BS22}).
    
    Now, let $(\colim_{\phi}\wtr_Z)^{\wedge_{(p,\xi)}}$ be the $(p,\xi)$-adic completion of the perfection of $\wtr_Z$ with respect to the given morphism $\varphi$ (cf. \cite[Lem. 3.9]{BS22}). By construction, we have 
    \[(\colim_{\phi}\wtr_Z)^{\wedge_{(p,\xi)}}/p\cong (\colim_{x\mapsto x^p}\calR_Z^{\flat})^{\wedge_{\varpi}},\]
    yielding that 
    \[(\colim_{\phi}\wtr_Z)^{\wedge_{(p,\xi)}}\cong \bA_{\inf}(X_{\infty,Z}).\]
    Then the natural morphism 
    \[\widetilde \beta:\wtr_Z\to (\colim_{\phi}\wtr_Z)^{\wedge_{(p,\xi)}}\cong\bA_{\inf}(X_{\infty,Z})\] 
    defines a lifting of the natural inclusion $\calR_Z^{\flat}\hookrightarrow(\colim_{x\mapsto x^p}\wtr_Z^{\flat})^{\wedge_{\varpi}}$.
\end{exam}

Note that one can also realize $X_{\infty,Z}$ as a pro-finite-\'etale covering of $X_Z$.

\begin{exam}\label{exam:Gamma torsor}
    Suppose that $\frakX_Z = \Spf(\calR_Z)$ is small semi-stable with a chart 
    \[\psi:\frakX_Z\to \Spf(A_0^+:=A^+\za T_0,\dots,T_r,T_{r+1}^{\pm 1},\dots,T_d^{\pm 1}\ya/(T_0\cdots T_r-p^a))\]
    and a lifting $\wtx_Z = \Spf(\wtr_Z)$, and let $X_Z = \Spa(R_Z,R_Z^+=\calR_Z)$ be its generic fiber as in Example \ref{exam:small affine is liftable} and keep notations there. For any $n\geq 0$, define
    \[A_n^+:= A^+\za T_0^{\frac{1}{p^n}},\dots,T_r^{\frac{1}{p^n}},T_{r+1}^{\pm \frac{1}{p^n}},\dots,T_d^{\pm \frac{1}{p^n}}\ya/((T_0\cdots T_r)^{\frac{1}{p^n}}-p^{\frac{a}{p^n}})\]
    and 
    \[\widehat A_{\infty}^+:=(\cup_{n\geq 0}A_n^+)^{\wedge_p} = A^+\za T_0^{\frac{1}{p^{\infty}}},\dots,T_r^{\frac{1}{p^{\infty}}},T_{r+1}^{\pm \frac{1}{p^{\infty}}},\dots,T_d^{\pm \frac{1}{p^{\infty}}}\ya/((T_0\cdots T_r)^{\frac{1}{p^{\infty}}}-p^{\frac{a}{p^{\infty}}})\]
    where $((T_0\cdots T_r)^{\frac{1}{p^{\infty}}}-p^{\frac{a}{p^{\infty}}})$ denotes the closed ideal generated by $\{(T_0\cdots T_r)^{\frac{1}{p^n}}-p^{\frac{a}{p^n}}\}_{n\geq 0}$.
    Let $X_{\infty,Z}:=\Spa(\widehat R_{\infty,Z},\widehat R_{\infty,Z}^+)$ be the base-change of $X_Z$ along the pro-\'etale $\Gamma$-torsor
    \[\Spa(\widehat A_{\infty} =\widehat A_{\infty}^+[\frac{1}{p}],\widehat A_{\infty}^+)\to \Spa(A_0 = A_0^+[\frac{1}{p}],A_0^+),\]
    where $\Gamma$ is the following group
    \begin{equation}\label{equ:gamma group-I}
        \Gamma = \{\delta_0^{n_0}\cdots\delta_r^{n_r}\delta_{r+1}^{n_{r+1}}\cdots\delta_d^{n_d}\in \bZ_p^{r+1}\mid n_0,\dots,n_d\in\bZ_p \text{ such that } n_0+\cdots+n_r = 0\}
    \end{equation}
    and for any $0\leq i,j\leq d$ and any $n\geq 0$, the $\Gamma$-action on $\widehat A_{\infty}$ is determined by
    \begin{equation}\label{equ:gamma acts on T}
        \delta_i(T_j^{\frac{1}{p^n}}) = \zeta_{p^n}^{\delta_{ij}}T_j^{\frac{1}{p^n}}.
    \end{equation}
    Here, $\delta_{ij}$ denotes the Kronecker delta. In particular, we get a $\Gamma$-equivariant decomposition
    \begin{equation}\label{equ:decomposition original}
        \widehat R_{\infty,Z}^+ = \widehat \bigoplus_{\underline \alpha\in J_{r,d}}\calR_Z\cdot \underline T^{\underline \alpha}.
    \end{equation}
    If we put $\gamma_i:=\delta_0^{-1}\delta_i$ for any $1\leq i\leq r$ and $\gamma_j = \delta_j$ for any $r+1\leq j\leq d$, then we have an isomorphism
    \begin{equation}\label{equ:gamma group-II}
        \Gamma \cong \oplus_{i=1}^d\Zp\cdot\gamma_i.
    \end{equation}
    Put $T_i^{\flat}:=(T_i,T_i^{1/p},\dots)\in \widehat R_{\infty,Z}^{\flat+}$ for any $0\leq i\leq d$, and then we have
    \[T_0^{\flat}\cdots T_r^{\flat} = \varpi^a\] 
    and the $\Gamma$-action on $\widehat R_{\infty,Z}^{\flat+}$ is determined by
    \begin{equation}\label{equ:gamma acts on T^flat}
        \delta_i((T_j^{\flat})^{\frac{1}{p^n}}) = \epsilon^{\frac{\delta_{ij}}{p^n}}(T_j^{\flat})^{\frac{1}{p^n}}
    \end{equation}
    for any $0\leq i,j\leq d$ and $n\geq 0$. Now, consider the morphism of $\bA_{\inf}(Z)$-algebras
    \[\bA_{\inf}(Z)\za T_0,\dots,T_r,T_{r+1}^{\pm 1},\dots,T_d^{\pm 1}\ya/(T_0\cdots T_r-[\varpi^a]))\to \bA_{\inf}(X_{\infty,Z})\]
    sending each $T_i$ to $[T_i^{\flat}]$ for any $0\leq i\leq d$. It follows directly that this map is compatible with log-structures, and by the \'etaleness of $\widetilde \phi$, it extends uniquely to a morphism
    \[\widetilde \beta:\wtr_Z\to \bA_{\inf}(X_{\infty,Z})\]
    which is compatible with log-structures as well. One can check that the morphism $\widetilde \beta$ coincides with that in Example \ref{exam:construct beta} in the case where the Frobenius endomorphism $\phi$ on $\wtr_Z$ maps each $T_i$ to $T_i^p$ for any $0\leq i\leq d$ (which amounts to the $\delta$-structure on $\wtr_Z$ being uniquely determined by letting $\delta(T_i) = 0$ for all $i$).
\end{exam}

More generally, we have the following result, which will be used several times in the following sections.

\begin{lem}\label{lem:map into perfectoid}
    For any affinoid perfectoid $U = \Spa(S,S^+)\in X_{Z,v}$, consider the natural map $\calR_Z\to S^+$. Then for any $0\leq i\leq r$, there exist $\underline t_i\in S^{\flat\times}\cap S^{\flat+}$ and $u_i\in \bA_{\inf}(U)^{\times}$ that satisfy the following conditions:
    \begin{enumerate}
        \item[$\bullet$] The image of $t_i:=[\underline t_i]u_i\in \bA_{\inf}(U)$ modulo $\xi$ coincides with $T_i$ for all $i$.
        \item[$\bullet$] We have $\prod_{i=0}^rt_i = [\varpi^a]$.
    \end{enumerate}
    As a consequence, the map $\calR_Z \to S^+$ of $A^+$-algebras admits a lifting $\wtr_Z\to\bA_{\inf}(U)$ over $\bA_{\inf}(Z)$ which is compatible with log-structures.
\end{lem}
\begin{proof}
    See \cite[Lem. 2.9 and Cor. 2.10]{LMNQW}.
\end{proof} 

\begin{dfn}\label{dfn:I-norm on R^I}
    Let $\calS$ be the set of all morphisms $\widetilde \beta:\wtr_Z\to\bA_{\inf}(X_{\infty,Z})$ as in Construction \ref{construction:I-norm}. Define
    \[|\cdot|_I:\wtr_Z^I\to \bR_{\geq 0}\]
    such that
    \[|\cdot|_I:=\sup_{\widetilde \beta\in \calS}|\cdot|_{I,\widetilde \beta}.\]
    It follows directly that we have 
    \[\wtr_Z^{I,+} = \{x\in \wtr_Z^I\mid |x|_I\leq 1\}.\]
\end{dfn}

\begin{lem}\label{lem:isometry embedding}
    Fix a $\widetilde \beta:\wtr_Z\to\bA_{\inf}(X_{\infty,Z})$ as in Construction \ref{construction:I-norm}. Define
    \[|\cdot|_{I,\widetilde \beta}:\CCtilde^I(X_{\infty,Z})\to\bR_{\geq 0}\]
    such that for any $x\in\CCtilde^I(X_{\infty,Z})$ of the form $x=\sum_{\underline \alpha\in J_{r,d}}\widetilde \beta(x_{\underline \alpha})[\underline T^{\underline \alpha}]$ with $x_{\underline \alpha}\in\wtr_Z^I$ for each $\underline \alpha$ via the decomposition \eqref{equ:decomposition over C^Iplus}, 
    \[|x|_{I,\widetilde \beta} = \sup_{\underline \alpha\in J_{r,d}}|x_{\alpha}|_{I,\widetilde \beta}.\]
    Then we have $|\cdot|_{I,\widetilde \beta} = |\cdot|_I$ on $\CCtilde^I(X_{\infty,Z})$.
\end{lem}
\begin{proof}
    Note that for any $\underline \alpha\in J_{r,d}$, we have $|[\underline T^{\underline \alpha}]|_I = 1$. Thus, we have
    \[|x|_I = |\sum_{\underline \alpha}\widetilde \beta(x_{\underline \alpha})[\underline T^{\underline \alpha}]|_I\leq \sup_{\underline \alpha}|\widetilde \beta(x_{\underline \alpha})|_I = \sup_{\underline \alpha}|x_{\underline \alpha}|_{I,\widetilde \beta}.\]
    
    Conversely, for any $x = \sum_{\underline \alpha\in J_{r,d}}\widetilde \beta(x_{\underline \alpha})[\underline T^{\underline \alpha}]\in \CCtilde^I(X_{\infty,Z})$, we have to show that $|x|_{I,\widetilde \beta}\leq |x|_I$.
    Noting that 
    \[|\cdot|_{I,\widetilde \beta} = \sup_{t\in I}|\cdot|_{t,\widetilde \beta} \text{ and }|\cdot|_I = \sup_{t\in I}|\cdot|_t,\]
    we may assume that $I = \{r\}$ for some $r>0$. Without loss of generality, we may assume $x\neq 0$. Note that, by Lemma \ref{lem:second description of C^I}(2), for any $h\in \bQ$ and any $x\in \CCtilde^I(X_{\infty,Z})$, we have 
    \[|[\varpi^h]x|_{r,\widetilde \beta} = |[\varpi^h]|_r|x|_{r,\widetilde \beta} \text{ and }|[\varpi^h]x|_r = |[\varpi^h]|_r|x|_r.\]
    For any $\varepsilon>0$, replacing $x$ by $[\varpi^h]x$ for some $h\in\bQ$, we may assume that $1-\varepsilon<|x|_I\leq 1$.
    As $x\in \CCtilde^{I,+}(X_{\infty,Z})$, by decomposition \eqref{equ:decomposition over C^I}, we see that $x_{\underline \alpha}\in \wtr_Z^{I,+}$ and thus $|x_{\underline \alpha}|_{I,\widetilde \beta}\leq 1$ for any $\underline \alpha\in J_{r,d}$.
    This forces that
    \[|x|_{I,\widetilde \beta}\leq 1<(1-\varepsilon)^{-1}|x|_I,\]
    yielding that
    \[|x|_{r,\widetilde \beta}\leq |x|_r\]
    as desired, by letting $\varepsilon\to 0$.
\end{proof}

\begin{prop}\label{prop:Independence of embedding}
    Fix a $\widetilde \beta:\wtr_Z\to\bA_{\inf}(X_{\infty,Z})$ as in Construction \ref{construction:I-norm}.
    \begin{enumerate}
        \item[(1)] The norm $|\cdot|_{I,\widetilde \beta}$ is complete on $\wtr_Z^I$.
        \item[(2)] For another $\widetilde \beta^{\prime}:\wtr_Z\to\bA_{\inf}(X_{\infty,Z})$, the norms $|\cdot|_{I,\widetilde \beta}$ and $|\cdot|_{I,\widetilde \beta^{\prime}}$ are same on $\wtr_Z^I$.
    \end{enumerate}
\end{prop}
\begin{proof}
    For item (1): We have to show that for any sequence $\{x_n\}_{n\geq 1}\in \wtr_Z$ such that $|x_n|_{I,\widetilde \beta}\to 0$ with respect to the norm $|\cdot|_{I,\widetilde \beta}$, the series $\sum_{n\geq 1}x_n$ converges to a unique element $x\in \wtr_Z^I$.
    Regard $\{x_n\}_{n\geq 1}$ as a sequence on $\CCtilde^I(X_{\infty,Z})$ via the morphism $\widetilde \beta$. As $\CCtilde^I(X_{\infty,Z})$ is complete with respect to $|\cdot|_I$ (cf. Lemma \ref{lem:second description of C^I}), by Lemma \ref{lem:isometry embedding}, $\sum_{n\geq 1}x_n$ uniquely converges to an element $x\in \CCtilde^I(X_{\infty,Z})$.
    It is enough to show that $x\in \wtr_Z^I$; that is, if we write $x = \sum_{\underline \alpha\in J_{r,d}}x_{\underline \alpha}[\underline T^{\underline \alpha}]$ via the decomposition \eqref{equ:decomposition over C^I}, we have to show that $x_{\underline \alpha} = 0$ for any $\underline \alpha\neq 0$.
    If this is not true, we have $x-x_0\neq 0$ and thus $|x-x_0|_{I,\widetilde \beta}> 0$.
    Then for any $N\geq 0$, we have
    \[x-\sum_{n=1}^Nx_n = (x_0-\sum_{n\geq 1}x_n)+\sum_{\underline \alpha\neq 0}x_{\underline \alpha}[\underline T^{\underline \alpha}]\]
    yielding that
    \[|x-\sum_{n=1}^Nx_n|_{I,\widetilde \beta}\geq |\sum_{\underline \alpha\neq 0}x_{\underline \alpha}[\underline T^{\underline \alpha}]|_{I,\widetilde \beta}=|x-x_0|_{I,\widetilde \beta}>0\]
    by the previous lemma. This violates the fact that $|x - \sum_{n=1}^Nx_n|_{I,\widetilde \beta}\to 0$ as $N\to+\infty$.
    So we have $x\in \wtr_Z^I$ as desired.

    For item (2): Note that 
    \[\{x\in \wtr_Z^I\mid |x|_{I,\widetilde \beta}\leq 1\} = \wtr_Z^{I,+} = \{x\in \wtr_Z^I\mid |x|_{I,\widetilde \beta^{\prime}}\leq 1\}.\]
    The assertion follows from the same argument in the second paragraph of the proof of Lemma \ref{lem:isometry embedding}.
\end{proof}

\begin{cor}\label{cor:Independence of embedding}
    For any $\widetilde \beta:\wtr_Z\to\bA_{\inf}(X_{\infty,Z})$ as in Construction \ref{construction:I-norm}, we have 
    \[|\cdot|_I = |\cdot|_{I,\widetilde \beta}:\wtr_Z^I\to \bR_{\geq 0}\]
    such that the associated morphism
    \[\widetilde \beta:(\wtr_Z^I,|\cdot|_I)\to (\CCtilde^I(X_{\infty,Z}),|\cdot|_I)\]
    is an isometry.
\end{cor}
\begin{proof}
    This follows directly from Proposition \ref{prop:Independence of embedding} and Lemma \ref{lem:isometry embedding}.
\end{proof}

We collect some properties about the norm $|\cdot|_I$ on $\wtr_Z^I$.

\begin{lem}\label{lem:property on I-norm}
    \begin{enumerate}
        \item[(1)] Equip $\wtr_Z[\frac{1}{p[\varpi]}]$ with the norm $|\cdot|_I$ induced from the restriction of $|\cdot|_I$ on $\wtr_Z^I$ along the morphism $\wtr_Z[\frac{1}{p[\varpi]}]\to\wtr_Z^I$.
        Then $\wtr_Z^I$ is the completion of $\wtr_Z[\frac{1}{p[\varpi]}]$ with respect to the norm $|\cdot|_I$.
        \item[(2)] Suppose that $I^{\prime}\subset I$ is a closed sub-interval. Then the natural morphism $\wtr_Z^I\to \wtr_Z^{I^{\prime}}$ has dense image.
        \item[(3)] Suppose that $I=\{t\}$ for some $t>0$. Then for any $c\in \Ctilde^t$ and $x\in \wtr_Z^I$, we have 
        \[|cx|_t = |c|_t\cdot|x|_t.\]

        \item[(4)] Suppose that $1\in I$. Then there is a natural morphism 
        \[\wtr_Z^I\to \wtR_Z\]
        such that the composite $\wtr_Z^I\to \wtR_Z\xrightarrow{\mod \xi}R_Z$ sends $\wtr_Z^{I,+}$ onto $R_Z^+=\calR_Z$.
        Moreover, the morphism $\wtr_Z^I\to \wtR_Z$ induces an isomorphism
        \[(\wtr_Z^I)^{\wedge_{\xi}}\xrightarrow{\cong} \wtR_Z,\]
        where $(\wtr_Z^I)^{\wedge_{\xi}}$ denotes the $\xi$-adic completion of $\wtr_Z^I$.
        \item[(5)] Suppose that there exists a homomorphism $\phi:\wtr_Z\to\wtr_Z$ lifting the absolute Frobenius endomorphism $\calR_Z^{\flat}\xrightarrow{x\mapsto x^p}\calR_Z^{\flat}$.
        Then it induces a homomorphism $\varphi:\wtr_Z^I\to \wtr_Z^{\frac{I}{p}}$ such that for any $x\in \wtr_Z^I$, we have
        \[|\phi(x)|_{\frac{I}{p}} = |x|_I.\]
    \end{enumerate}
\end{lem}
\begin{proof}
    For item (1): Choose a morphism $\widetilde \beta:\wtr_Z\to\bA_{\inf}(X_{\infty,Z})$ as in Construction \ref{construction:I-norm}. By \eqref{equ:decomposition over Ainf}, we have a decomposition
    \[\bA_{\inf}(X_{\infty,Z})[\frac{1}{p[\varpi]}] = \widehat \bigoplus_{\underline \alpha\in J_{r,d}}\wtr_Z[\frac{1}{p[\varpi]}]\cdot[\underline T^{\underline \alpha}].\]
    Then the result follows from Lemma \ref{lem:second description of C^I}(1) and Lemma \ref{lem:isometry embedding}.
    Item (2) follows from Item (1) directly.

    Item (3) follows from Lemma \ref{lem:second description of C^I}(2).
    Item (4) follows from the same argument as in the proof of Lemma \ref{lem:compare with BdR}.
    For item (5), let $\widetilde \beta:\wtr_Z\to \bA_{\inf}(X_{\infty,Z})$ be the morphism associated to the given Frobenius morphism as in Example \ref{exam:construct beta}.
    Then we have the following commutative diagram
    \begin{equation*}
        \xymatrix@C=0.5cm{
          \wtr_Z^I\ar[d]^{\phi}\ar[rrr]^{\widetilde \beta}&&&\CCtilde^I(X_{\infty,Z})\ar[d]^{\phi}\\
          \wtr_Z^{\frac{I}{p}}\ar[rrr]^{\widetilde \beta}&&&\CCtilde^{\frac{I}{p}}(X_{\infty,Z}).
    }
    \end{equation*}
    As the horizontal arrows are both isometric by Corollary \ref{cor:Independence of embedding}, the result follows from Lemma \ref{lem:second description of C^I}(3).
\end{proof}

\begin{rmk}\label{rmk:isometry for R^I}
    Suppose that $\frakX^{\prime}_Z = \Spf(\calR^{\prime}_Z)\to \frakX_Z=\Spf(\calR_Z)$ is a morphism of small semi-stable formal schemes over $A^+$ and denote by $f:\calR_Z\to \calR_Z^{\prime}$ the associated morphism of $A^+$-algebras. Suppose that $f$ admits a lifting $\widetilde f:\wtr_Z\to\wtr_Z^{\prime}$ over $\bA_{\inf}(Z)$, and consider the induced morphism
    \[\widetilde f^I:\wtr_Z^I\to\wtr_Z^{\prime I}.\] 
    By Remark \ref{rmk:isometry on base}, for any $x\in \wtr_Z^I$, we have $|\widetilde f^I(x)|_I\leq |x|_I$.
    If moreover, $\frakX^{\prime}_Z\to \frakX_Z$ is a smooth covering (which implies that $X_{\infty,Z}^{\prime}\to X_{\infty,Z}$ is a $v$-covering), then we have $|\widetilde f^I(x)|_I = |x|_I$.
\end{rmk}

\section{Period sheaves with connection over Robba rings}\label{sec:period sheaf with connection}

Fix an affinoid perfectoid space $Z = \Spa(A,A^+)\in \Perfd$. Let $\frakX_Z$ be a semi-stable formal scheme of relative dimension $d$ over $A^+$ which is liftable with a fixed lifting $\wtx_Z$ over $\bA_{\inf}(Z)$, and denote by $X_Z$ the generic fiber of $\frakX_Z$ and by $\wtX_Z$ its lifting over $\BBdRp(Z)$ associated to $\wtx_Z$ as before. In this section we construct the overconvergent period sheaf $\calO\CCtilde^{I,\dagger}$ on the $v$-site $X_{Z,v}$. It is a sheaf of $\CCtilde^I$-algebras equipped with a natural $\CCtilde^I$-linear connection
\[\nabla:\calO\CCtilde^{I,\dagger}\to \calO\CCtilde^{I,\dagger}\otimes_{\calO_{\wtx_Z}}\Omega^{1,\log}_{\wtx_Z}\{-1\}\]
such that the following sequence of sheaves of $\CCtilde^I$-modules on $X_{Z,v}$ is exact:
\[0\to \CCtilde^I\to \calO\CCtilde^{I,\dagger}\xrightarrow{\nabla}\calO\CCtilde^{I,\dagger}\otimes_{\calO_{\wtx_Z}}\Omega^{1,\log}_{\wtx_Z}\{-1\}\to\cdots\to \calO\CCtilde^{I,\dagger}\otimes_{\calO_{\wtx_Z}}\Omega^{d,\log}_{\wtx_Z}\{-d\}\to 0.\]

\subsection{Construction for small affines}\label{ssec:period sheaf-small case}
We first treat the case in which $\frakX_Z = \Spf(\calR_Z)$ is small affine, so the lifting of $\frakX_Z$ is given by $\wtx_Z = \Spf(\wtr_Z)$ as described in Example \ref{exam:small affine is liftable}. Denote by $X_Z = \Spa(R_Z,R_Z^+)$ the generic fiber of $\frakX_Z$. For any affinoid perfectoid $U = \Spa(S,S^+)\in X_{Z,v}$, denote by $\Sigma_U$ the set of all \'etale morphisms
\[i_{\frakY_Z}:\frakY_Z = \Spf(\calR_Z^{\prime})\to\frakX_Z\]
over $A^+$ such that the natural morphism $U\to X_Z$ factors through the induced morphism of generic fibers $Y_Z\to X_Z$. In other words, the morphism $\wtr_Z\to\wtr_Z^{\prime}$ is $p$-completely \'etale and the morphism $\calR_Z^+\to S^+$ factors as $\calR_Z\to\calR_Z^{\prime}\to S^+$. The lifting $\wtx_Z$ uniquely induces a lifting $\wty_Z = \Spf(\wtr_Z^{\prime})$ of $\frakY_Z$ and a lifting $\wty_Z\to\wtx_Z$ of $i_{\frakY_Z}$. Note that the log-structure on $\wtr_Z^{\prime}$ is associated to the composite $M_{r,a}(Z)\to \wtr_Z\to\wtr_Z^{\prime}$ (cf. Example \ref{exam:small affine is liftable}).

From now on, we write $M_U$ for the canonical log-structure on $\bA_{\inf}(U)$, i.e. the log-structure associated to $S^{\flat\times}\cap S^{\flat+}\xrightarrow{[\cdot]}\bA_{\inf}(U)$. Then the natural map
\[\theta_{\frakY_Z}:\wtr_Z^{\prime}\widehat \otimes_{\bA_{\inf}(Z)}\bA_{\inf}(U)\to S^+\]
is compatible with log-structures, where $\wtr_Z^{\prime}\widehat \otimes_{\bA_{\inf}(Z)}\bA_{\inf}(U)$ is equipped with the tensor product log-structure, i.e. the log-structure associated to $M_{r,a}(Z)\oplus_{M_Z}M_U$, which denotes the monoid obtained as the push-out of the diagram
\[M_{r,a}(Z)\leftarrow M_Z\rightarrow M_U.\]
Denote the $(p,\xi)$-adic completion of the exactification of $\theta_{\frakY_Z}$ in the sense of \cite[Exam. 2.11]{LMNQW} by
\begin{equation}\label{equ:exactification}
    \theta_{\frakY_Z}^{\prime}:A_{\frakY_Z}^+:=\big(\wtr_Z^{\prime}\widehat \otimes_{\bA_{\inf}(Z)}\bA_{\inf}(U)\big)[\big(\frac{T_0}{t_0}\big)^{\pm 1},\dots,\big(\frac{T_r}{t_r}\big)^{\pm 1}]^{\wedge_{(p,\xi)}}\to S^+,
\end{equation}
where $t_0,\dots,t_r$ are elements given by Lemma \ref{lem:map into perfectoid}. As $\prod_{i=0}^r\frac{T_i}{t_i} = 1$, the ideals of $A_{\frakY_Z}^+$ generated by $1-\frac{T_1}{t_1}, \dots, 1-\frac{T_d}{t_d}$ and $1-\frac{T_0}{t_0},\dots, 1-\frac{T_d}{t_d}$ are the same, where for any $r<i\leq d$, $t_i\in \bA_{\inf}(U)^{\times}$ denotes the image of $T_i$ via the morphism $\wtr_Z\to \bA_{\inf}(U)$ in Lemma \ref{lem:map into perfectoid}. Note that when $\frakX_Z$ is smooth; that is, $a = 0$ in Definition \ref{dfn:semi-stable formal scheme}(1), we have $A_{\frakY_Z}^+ = \wtr_Z^{\prime}\widehat \otimes_{\bA_{\inf}(Z)}\bA_{\inf}(U)$ and $\theta_{\frakY_Z} = \theta_{\frakY_Z}^{\prime}$.

\begin{lem}\label{lem:finite ideal}
    Keep the notation as above and fix $t_{r+1},\dots,t_d\in\bA_{\inf}(U)$ such that their reductions modulo $\xi$ are $T_{r+1},\dots,T_d$ respectively. The ideal $\Ker(\theta_{\frakY_Z}^{\prime})$ is generated by 
    \[1-\frac{T_1}{t_1}, \dots, 1-\frac{T_d}{t_d},e,\xi,\]
    where $e\in \wtr_Z^{\prime}\widehat \otimes_{\bA_{\inf}(Z)}\bA_{\inf}(U)$ satisfies the condition that
    \[e\equiv e^2 \mod (1-\frac{T_1}{t_1}, \dots, 1-\frac{T_d}{t_d}).\]
\end{lem}
\begin{proof}
    This is \cite[Lem. 3.1]{LMNQW}.
\end{proof}

Now, we equip $\wtr_Z^{\prime I}\widehat \otimes_{\CCtilde^I(Z)}\CCtilde^I(U)$ with the norm $|\cdot|_I$ associated to the tensor product of Banach-$\CCtilde^I(Z)$-algebras $(\wtr_Z^{\prime I},|\cdot|_I)$ and $(\CCtilde^I(U),|\cdot|_I)$. Then the natural map 
\[\wtr_Z^{\prime}\widehat \otimes_{\bA_{\inf}(Z)}\bA_{\inf}(U)\to \wtr_Z^{\prime I}\widehat \otimes_{\CCtilde^I(Z)}\CCtilde^I(U)\]
induces a norm $|\cdot|_I$ on $\wtr_Z^{\prime}\widehat \otimes_{\bA_{\inf}(Z)}\bA_{\inf}(U)$. Note that
\[A_{\frakY_Z}^+ = \big(\wtr_Z^{\prime}\widehat \otimes_{\bA_{\inf}(Z)}\bA_{\inf}(U)\big)\za X_0,\dots,X_r\ya/(T_0-X_0t_0,\dots,T_r-X_rt_r,X_0\cdots X_r-1),\]
where $\big(\wtr_Z^{\prime}\widehat \otimes_{\bA_{\inf}(Z)}\bA_{\inf}(U)\big)\za X_0,\dots,X_r\ya$ denotes the $(p,\xi)$-adic completion of 
\[\big(\wtr_Z^{\prime}\widehat \otimes_{\bA_{\inf}(Z)}\bA_{\inf}(U)\big)[X_0,\dots,X_r].\]
So the Gauss norm $|\cdot|_G$ on $\big(\wtr_Z^{\prime}\widehat \otimes_{\bA_{\inf}(Z)}\bA_{\inf}(U)\big)\za X_0^{\pm 1},\dots,X_r^{\pm 1}\ya$ with $|X_i|_G=1$ for any $0\leq i\leq r$ induces a quotient norm on $A_{\frakY_Z}^+$ that we still denote by $|\cdot|_I$; this norm extends uniquely to $A_{\frakY_Z}:=A_{\frakY_Z}^+[\frac{1}{p}]$ by letting $|\frac{x}{p^n}|_I:=|x|_I\cdot |p|^{-n}$ for $x\in A_{\frakY_Z}^+, n\in\bZ$.

\begin{rmk}\label{rmk:independence of chart and embedding}
    The norm $|\cdot|_I$ on $A_{\frakY_Z}^+$ is independent of 
    \begin{enumerate}
        \item[(1)] the choice of the chart on $\frakX_Z$, and
        \item[(2)] the choice of the morphism $\iota:\wtr_Z\to\bA_{\inf}(U)$ in Lemma \ref{lem:map into perfectoid}.
    \end{enumerate}
    For item (1): Suppose that 
    \[\psi^{\prime}:A^+\za T_0^{\prime},\dots,T_r^{\prime},T_{r+1}^{\prime\pm 1},\dots,T_d^{\prime \pm 1}\ya/(T_0^{\prime}\cdots T_r^{\prime}-p^a)\] 
    is another chart on $\frakX_Z$. As the log-structure on $\wtr_Z$ is also associated to the composite
    \[M_{r,a}(Z)\xrightarrow{\alpha^{\prime}}\bA_{\inf}(Z)\za T_0^{\prime},\dots,T_r^{\prime},T_{r+1}^{\prime\pm 1},\dots,T_d^{\prime \pm 1}\ya/(T_0^{\prime}\cdots T_r^{\prime}-[\varpi^a])\xrightarrow{\widetilde \psi^{\prime}}\wtr_Z,\]
    where $M_{r,a}(Z)$ is the monoid and $\widetilde \psi^{\prime}$ is the analogue of $\widetilde \psi$ appearing in Example \ref{exam:small affine is liftable} while the morphism of monoids $\alpha^{\prime}$ sends each $e_i$ to $T_i^{\prime}$, we see that there are units $z_0,\dots,z_r\in \wtr_Z$ satisfying $z_0\cdots z_r=1$ such that $T_i^{\prime} = T_iz_i$ for any $0\leq i\leq r$. 
    Indeed, as $(M_{r,a},\alpha)$ and $(M_{r,a},\alpha^\prime)$ induce the same log-structure on $\wtr_Z$, we get an isomorphism of monoids
    \[\begin{split}\beta:(\oplus_{i=0}^r\bN\cdot e_i)\oplus_{\bN\cdot (e_0+\cdots+e_r)}(M_Z/\bA_{\inf}(Z)^{\times})&\cong M_{r,a}(Z)/\alpha^{-1}(\wtr_Z^{\times})\\
    &\xrightarrow{\cong} M_{r,a}(Z)/\alpha^{\prime-1}(\wtr_Z^{\times})\cong(\oplus_{i=0}^r\bN\cdot e_i)\oplus_{\bN\cdot (e_0+\cdots+e_r)}(M_Z/\bA_{\inf}(Z)^{\times}), \end{split}\]
    which is the identity map when we restrict it to $M_Z/\bA_{\inf}(Z)^{\times}$. It suffices to show that, up to a permutation of subscript, we have $\beta(e_i) = e_i$. To do so, for any $0\leq i\leq r$, write
    \[\beta(e_i) = \left(\sum_{j=0}^rn_{ij}e_j,m_i\right)\]
    with $m_i\in M_Z/\bA_{\inf}(Z)^{\times}$ and $n_{ij}\in\bN$ such that $\prod_{j=0}^rn_{ij} = 0$. Since $\sum_{i=0}^re_i = [\varpi^a]$, we see that
    \[[\varpi^a] = \left(\sum_{j=0}^rN_je_j,m\right),\]
    where $N_j := \sum_{i=0}^rn_{ij}$ and $m:=\prod_{i=0}^rm_i$. Put $N:=\min\{N_0,\dots,N_r\}$. Then we have
    \[[\varpi^a] = \left(\sum_{j=0}^r(N_j-N)e_j,m[\varpi^a]^N\right),\]
    which forces that $N_j = N$ for all $j$ and $[\varpi^a] = m[\varpi^a]^N$. In particular, we get $N\leq 1$. If $N = 0$, we see that $N_j = 0$ and thus $n_{ij} = 0$ for any $0\leq i,j\leq r$; that is, $\beta(e_i)\in M_Z/\bA_{\inf}(Z)^{\times}$, which violates to that $\beta$ is an isomorphism. So we have $N = N_0 = \cdots = N_r = 1$ and $m=1$, which implies that $m_i=1$ for all $i$. Equivalently, we have
    \[\beta(e_i) = \sum_{j=0}^rn_{ij}e_j \text{ with }\sum_{j=0}^rn_{ij}=1.\]
    Since $\beta$ is an isomorphism, we see that there exists a permutation $\sigma$ of subscript such that $\beta(e_i) = e_{\sigma(i)}$ as desired.

    Put $t_i^{\prime} = \iota(t_i)$ and $v_i = \iota(z_i)$; then $z_i^{-1}v_i\in \big(\wtr_Z^{\prime}\widehat \otimes_{\bA_{\inf}(Z)}\bA_{\inf}(U)\big)^{\times}$ such that $|z_i^{-1}v_i|_I = 1$ for each $i$. One can conclude by noting that the endomorphism
    \[f:\big(\wtr_Z^{\prime}\widehat \otimes_{\bA_{\inf}(Z)}\bA_{\inf}(U)\big)\za X_0,\dots,X_r\ya\to\big(\wtr_Z^{\prime}\widehat \otimes_{\bA_{\inf}(Z)}\bA_{\inf}(U)\big)\za X_0,\dots,X_r\ya\]
    sending each $X_i$ to $z_i^{-1}v_iX_i$ is an isometric isomorphism and maps the ideal $(T_0-X_0t_0,\dots,T_r-X_rt_r,X_0\cdots X_r-1)$ onto the ideal $(T_0^{\prime}-X_0t_0^{\prime},\dots,T_r^{\prime}-X_rt_r^{\prime},X_0\cdots X_r-1)$.

    For item (2): Suppose that $\iota^{\prime}:\wtr_Z\to \bA_{\inf}(U)$ is another lifting of the natural map $\calR_Z\to S^+$. By deformation theory, for any $0\leq i\leq r$, we have 
    \[t_i^{\prime}=\iota^{\prime}(T_i) = t_iw_i\] 
    for some unit $w_i\in \bA_{\inf}(U)$ such that $w_0\cdots w_r=1$. As $|w_i|_I=1$ for any $0\leq i\leq r$, it remains only to observe that the endomorphism
    \[g:\big(\wtr_Z^{\prime}\widehat \otimes_{\bA_{\inf}(Z)}\bA_{\inf}(U)\big)\za X_0,\dots,X_r\ya\to\big(\wtr_Z^{\prime}\widehat \otimes_{\bA_{\inf}(Z)}\bA_{\inf}(U)\big)\za X_0,\dots,X_r\ya\]
    sending each $X_i$ to $w_iX_i$ is an isometric isomorphism and maps the ideal $(T_0-X_0t_0,\dots,T_r-X_rt_r,X_0\cdots X_r-1)$ onto the ideal $(T_0-X_0t_0^{\prime},\dots,T_r-X_rt_r^{\prime},X_0\cdots X_r-1)$.
\end{rmk}

\begin{exam}\label{exam:norm is one}
    For any $0\leq i\leq d$, consider $\frac{T_i}{t_i}\in (A_{\frakY_Z}^+)^{\times}$. We have
    \[|\frac{T_i}{t_i}|_I = |\frac{t_i}{T_i}|_I=1.\]
    Indeed, for $r+1\leq i\leq d$, this is obvious because $T_i$ and $t_i$ are units in $\wtr_Z^{\prime}$ and $t_i\in \bA_{\inf}(U)^{\times}$ respectively. For $0\leq i\leq r$, we have 
    \[|\frac{T_i}{t_i}|_I \leq |X_i|_G\leq 1 \text{ and } |\frac{t_i}{T_i}|_I\leq |X_1\cdots X_{i-1}X_{i+1}\cdots X_r|_G\leq 1.\]
    The claim follows from the chain of inequalities
    \[1\leq |\frac{T_i}{t_i}\cdot \frac{t_i}{T_i}|_I\leq |\frac{T_i}{t_i}|_I \cdot |\frac{t_i}{T_i}|_I\leq 1.\]
\end{exam}

Denote by $A_{\frakY_Z}^{\circ}\subset A_{\frakY_Z}$ the subset of power-bounded elements and let
\[\Spa(A_{\frakY_Z}^I,A_{\frakY_Z}^{I,+})\to\Spa(A_{\frakY_Z},A_{\frakY_Z}^{\circ})\]
be the rational localization with respect to the condition $[\varpi]\in|p|^I$. Again, we have $A_{\frakY_Z}^I = A_{\frakY_Z}^{I,+}[\frac{1}{p}]$ with $A_{\frakY_Z}^{I,+}$ being the $p$-adic completion of $A_{\frakY_Z}^{\circ}[\frac{[\varpi^{1/r}]}{p},\frac{p}{[\varpi^{1/s}]}]$; the norm $|\cdot|_I$ on $A_{\frakY_Z}$ extends to a norm which we still denote by $|\cdot|_I$ on $A_{\frakY_Z}^I$. Now, for any $t\in I$, one can define a norm $|\cdot|_t$ on $A_{\frakY_Z}^I[\frac{1}{\xi}]$ as follows: By Lemma \ref{lem:second description of C^I}(2) and Lemma \ref{lem:property on I-norm}(3), the norm $|\cdot|_t$ on $A_{\frakY_Z}^I$ extends uniquely to a norm on $A_{\frakY_Z}^I[\frac{1}{\xi}]$ by letting $|\frac{x}{\xi^n}|_t:=|x|_t\cdot |\xi|_t^{-n}$ for any $x\in A_{\frakY_Z}^I$ and any $n\in \bZ$. Using this, we can define a norm $|\cdot|_I$ on $A_{\frakY_Z}^I[\frac{1}{\xi}]$ by setting $|\cdot|_I:=\sup_{t\in I}|\cdot|_t$. By restriction, we get a norm $|\cdot|_I$ on the blow-up algebra $A_{\frakY_Z}^I[\frac{\Ker(\theta_{\frakY_Z}^{\prime})}{\xi}]$, and denote by $A_{\frakY_Z}^I[\frac{\Ker(\theta_{\frakY_Z}^{\prime})}{\xi}]^{\circ}$ the subset of power-bounded elements.

For any $1\leq \lambda<|\xi|_I^{-1}$, denote by
\[\Spa(A_{\frakY_Z}^{I,\lambda},A_{\frakY_Z}^{I,\lambda,+})\to \Spa(A_{\frakY_Z}^I[\frac{\Ker(\theta_{\frakY_Z}^{\prime})}{\xi}],A_{\frakY_Z}^I[\frac{\Ker(\theta_{\frakY_Z}^{\prime})}{\xi}]^{\circ})\]
the rational localization with respect to the condition $|\frac{\Ker(\theta_{\frakY_Z}^{\prime})}{\xi}|\leq \lambda$ (cf. Lemma \ref{lem:finite ideal}). Let 
\[\rd:\wtr_Z^{\prime}\to \Omega^{1,\log}_{\wtr_Z^{\prime}}\]
be the usual $(p,\xi)$-complete (log-)connection, which extends to a $\CCtilde^I(U)$-linear connection 
\[\rd:A_{\frakY_Z}^I\to A_{\frakY_Z}^{I}\otimes_{\wtr_Z^{\prime}}\Omega^{1,\log}_{\wtr_Z^{\prime}}\]
Observing that for any $n\geq 1$,
\[\rd((\frac{\Ker(\theta_{\frakY_Z}^{\prime})A_{\frakY_Z}^I}{\xi})^n)\subset (\frac{\Ker(\theta_{\frakY_Z}^{\prime})A_{\frakY_Z}^I}{\xi})^{n-1}\otimes_{\wtr_Z^{\prime}}\Omega^{1,\log}_{\wtr_Z^{\prime}}\cdot\xi^{-1}\]
we see that $\rd$ induces a connection
\[\nabla:A_{\frakY_Z}^{I,\lambda}\to A_{\frakY_Z}^{I,\lambda}\otimes_{\wtr_Z^{\prime}}\Omega^{1,\log}_{\wtr_Z^{\prime}}\{-1\}.\]
It follows directly that for any $1\leq \lambda_1\leq\lambda_2\leq |\xi|_I^{-1}$ and any inclusion of closed intervals $J\subset I$, we get the following commutative diagram which is compatible with the connections
\begin{equation*}
    \xymatrix@C=0.5cm{
      A_{\frakY_Z}^{I,\lambda_2} \ar[d]\ar[rr]&& A_{\frakY_Z}^{I,\lambda_1}\ar[d]\\
      A_{\frakY_Z}^{J,\lambda_2} \ar[rr]&& A_{\frakY_Z}^{J,\lambda_1}.
    }
\end{equation*}

\begin{dfn}\label{dfn:local OC^I}
    Suppose that $Z = \Spa(A,A^+)\in \Perfd$ and $\frakX_Z = \Spf(\calR_Z)$ is small semi-stable over $A^+$ with a fixed lifting $\wtx_Z = \Spf(\wtr_Z)$. Denote by $X_Z$ the generic fiber of $\frakX_Z$. Let $I = [r,s]$ be a closed sub-interval of $(0,+\infty)$ and let $\lambda\in [1,|\xi|_I^{-1})$. For each $\frakY_Z\in\Sigma_U$, the algebra $A_{\frakY_Z}^{I,\lambda}$ is the rational localization of $A_{\frakY_Z}^I[\frac{\Ker(\theta_{\frakY_Z}^{\prime})}{\xi}]$ on which $|\frac{\Ker(\theta_{\frakY_Z}^{\prime})}{\xi}|\leq \lambda$. Denote by $\calO\CCtilde^{I,\lambda}$ the sheaf on $X_{Z,v}$ corresponding to the sheafification of the presheaf 
    \[U\in X_{Z,v}\mapsto \colim_{\frakY_Z=\Spf(\wtr_Z^{\prime})\in\Sigma_U}A_{\frakY_Z}^{I,\lambda}.\]
    Then the usual differential $\rd$ on $\wtx_Z$ induces a $\CCtilde^I$-linear connection
    \[\nabla:\calO\CCtilde^{I,\lambda}\to\calO\CCtilde^{I,\lambda}\otimes_{\calO_{\wtx_Z}}\Omega^{1,\log}_{\wtx_Z}\{-1\}.\]
    It follows directly that for any $1\leq\lambda_1\leq\lambda_2<|\xi|_I^{-1}$, we have the following morphism which is compatible with the connection
    \[\calO\CCtilde^{I,\lambda_2}\to\calO\CCtilde^{I,\lambda_1}\to\calO\CCtilde^{I,1}_K.\]
    Denote by
    \[(\calO\CCtilde^{I,\dagger},\nabla:\calO\CCtilde^{I,\dagger}\to \calO\CCtilde^{I,\dagger}\otimes_{\calO_{\wtx_Z}}\Omega^{1,\log}_{\wtx_Z}\{-1\})\]
    the colimit 
    \[\colim_{\lambda\to 1^+}(\calO\CCtilde^{I,\lambda},\nabla:\calO\CCtilde^{I,\lambda}\to \calO\CCtilde^{I,\lambda}\otimes_{\calO_{\wtx_Z}}\Omega^{1,\log}_{\wtx_Z}\{-1\}).\]
    It follows directly that for any $1<\lambda_1\leq \lambda_2<|\xi|_I^{-1}$ and any $J\subset I$, the following diagram of morphisms of period sheaves with connection is commutative
    \begin{equation*}
        \xymatrix@C=0.5cm{
           (\calO\CCtilde^{I,\lambda_2},\nabla)\ar[r]\ar[d]&(\calO\CCtilde^{I,\lambda_1},\nabla)\ar[r]\ar[d]&(\calO\CCtilde^{I,\dagger},\nabla)\ar[r]\ar[d]&(\calO\CCtilde^{I,1}_K,\nabla)\ar[d]\\
          (\calO\CCtilde^{J,\lambda_2},\nabla)\ar[r]&(\calO\CCtilde^{J,\lambda_1},\nabla)\ar[r]&(\calO\CCtilde^{J,\dagger},\nabla)\ar[r]&(\calO\CCtilde^{J,1}_K,\nabla).
        }
    \end{equation*}
\end{dfn}

Note that the definition of $(\calO\CCtilde^{I,\lambda},\nabla)$ (and a fortiori $(\calO\CCtilde^{I,\dagger},\nabla)$) is independent of the chart on $\frakX_Z$.

\subsection{Local description of $(\calO\CCtilde^{I,\dagger},\nabla)$}\label{ssec:local description of OC^I}
 Keep notations as in the previous subsection. This subsection provides an explicit description of $(\calO\CCtilde^{I,\dagger},\nabla)$ when $\frakX_Z = \Spf(\calR_Z)$ is small semi-stable.

\begin{dfn}\label{dfn:overconvergent algebras}
     Let $(B,|\cdot|)$ be a $p$-adic Banach algebra. Fix a $b\in B$ such that $|b|<1$. Let $M$ be a Banach $B$-module with norm $|\cdot|$. For any $\lambda\in[1,|b|^{-1})$, denote by 
     \[M\{V_0/\lambda,\dots,V_d/\lambda\}\]
     the completion of $M[V_0,\dots,V_d]:=\bigoplus_{J=(j_0,\dots,j_d)\in\bN^{d+1}}M\cdot V_0^{j_0}\cdots V_d^{j_d}$ with respect to the norm
     \[|\sum_{J=(j_0,\dots,j_d)\in\bN^{d+1}}a_JV_0^{j_0}\cdots V_d^{j_d}|:=\sup_{J}\{|a_J|\cdot \lambda^{|J|}\}.\]
     Then $M\{V_0/\lambda,\dots,V_d/\lambda\}$ consists of formal power series
     \[\sum_{J\in\bN^{d+1}}a_J\underline V^J:=\sum_{J=(j_0,\dots,j_d)\in\bN^{d+1}}a_JV_0^{j_0}\cdots V_d^{j_d}\] 
     such that 
     \[\lim_{|J|=j_0+\cdots+j_d\to+\infty}|a_J|\cdot \lambda^{|J|} = 0\]
     and the norm
     \[|\sum_{J\in\bN^{d+1}}a_J\underline V^J|:=\sup_{J}|a_J|\cdot \lambda^{|J|}\]
     on $M\{V_0/\lambda,\dots,V_d/\lambda\}$ makes it a Banach $B$-module. In particular, setting $M = B$, $B\{V_0/\lambda,\dots,V_d/\lambda\}$ becomes a Banach $B$-algebra, and $M\{V_0/\lambda,\dots,V_d/\lambda\}$ is a module over it.
     
     Let $P\in\{\frac{\prod_{i=0}^r(1+bV_i)-1}{b}, V_0+\cdots+V_r\}$. Denote by 
     \[B\{\underline V/\lambda\}_{r,d,P}:=B\{V_0/\lambda,\dots,V_d/\lambda\}/(P)\]
     the quotient of $B\{V_0/\lambda,\dots,V_d/\lambda\}$ by the closed ideal generated by $P$, and endow it with the quotient norm $|\cdot|$. As $|b|<1$, it is easy to see that the map 
     \[B\{V_1/\lambda,\dots,V_d/\lambda\}\to B\{\underline V/\lambda\}_{r,d,P}\]
     is an isometric isomorphism of Banach $B$-algebras.
     Now, we put
     \[M\{\underline V/\lambda\}_{r,d,P}:=M\{V_0/\lambda,\dots,V_d/\lambda\}/(P)\]
     and equip it with the quotient norm. Again, the natural map
     \[M\{V_1/\lambda,\dots,V_d/\lambda\}\to M\{\underline V/\lambda\}_{r,d,P}\]
     is an isometric isomorphism of Banach $B$-modules.
     Define 
     \[M\{\underline V/\lambda\}:=M\{V_1/\lambda,\dots,V_d/\lambda\}.\]
     Denote respectively by $M\{\underline V\}_{r,d,P}^{\dagger}$ and $M\{\underline V\}^{\dagger}$ the colimits
     \[M\{\underline V\}_{r,d,P}^{\dagger}:= \colim_{\lambda\to 1^+}M\{\underline V/\lambda\}_{r,d,P} \text{ and }M\{\underline V\}^{\dagger}:=\colim_{\lambda\to 1^+}M\{\underline V/\lambda\}.\]
     When $M = B$, $B\{\underline V\}_{r,d,P}^{\dagger}$ and $B\{\underline V\}^{\dagger}$ are $B$-algebras.
     Let $M^+:=\{x\in M\mid |x|\leq 1\}$. We define
     \[\begin{split}
         M^+\{V_0/\lambda,\dots,V_d/\lambda\}&:=M\{V_0/\lambda,\dots,V_d/\lambda\}\cap M^+[[V_0,\dots,V_d]]\\
         &=\{\sum_{J\in\bN^{d+1}}a_J\underline V^J\in M\{V_0/\lambda,\dots,V_d/\lambda\}\mid a_J\in M^+,\forall J\}.
     \end{split}\]
     As $|b|<1$, we can still define $M^+\{\underline V/\lambda\}_{r,d,P}$ such that the natural map
     \[M^+\{V_1\lambda,\dots,V_d/\lambda\}\to M^+\{\underline V/\lambda\}_{r,d,P}\]
     is an isometric isomorphism. We can also define $M^+\{\underline V/\lambda\} = M^+\{V_1/\lambda,\dots,V_d/\lambda\}$ and define 
     \[M^+\{\underline V\}^{\dagger}_{r,d,P}:=\colim_{\lambda\to 1^+}M^+\{\underline V/\lambda\}_{r,d,P} \text{ and }M^+\{\underline V\}^{\dagger}:=\colim_{\lambda\to 1^+}M^+\{\underline V/\lambda\}\]
     in a similar way as above.
 \end{dfn}

 \begin{lem}\label{lem:induced derivation}
     Keep notations as in Definition \ref{dfn:overconvergent algebras}. Denote by 
     \[\Omega:=(\oplus_{i=0}^d\bZ\cdot e_i)/\bZ\cdot(\sum_{i=0}^re_i)\]
     the quotient of the free $\bZ$-module $\oplus_{i=0}^d\bZ\cdot e_i$ by $\bZ\cdot(\sum_{i=0}^re_i)$. For any $0\leq i\leq d$, denote by
     \[\partial_i:B\{V_0/\lambda,\dots,V_d/\lambda\}\to B\{V_0/\lambda,\dots,V_d/\lambda\}\]
     the derivation $\frac{\partial}{\partial V_i}$.
     \begin{enumerate}
         \item[(1)] For $P = \frac{\prod_{i=0}^r(1+bV_i)-1}{b}$, the $B$-linear connection
         \[\oplus_{i=0}^d(1+bV_i)\partial_i\otimes e_i:B\{V_0/\lambda,\dots,V_d/\lambda\}\to \oplus_{i=0}^dB\{V_0/\lambda,\dots,V_d/\lambda\}\cdot e_i\]
         induces a $B$-linear connection
         \[\nabla:B\{\underline V/\lambda\}_{r,d,P}\to B\{\underline V/\lambda\}_{r,d,P}\otimes\Omega.\]

         \item[(2)] For $P = V_0+\cdots+V_r$, the $B$-linear connection 
         \[\oplus_{i=0}^d\partial_i\otimes e_i:B\{V_0/\lambda,\dots,V_d/\lambda\}\to \oplus_{i=0}^dB\{V_0/\lambda,\dots,V_d/\lambda\}\cdot e_i\]
         induces a $B$-linear connection
         \[\nabla:B\{\underline V/\lambda\}_{r,d,P}\to B\{\underline V/\lambda\}_{r,d,P}\otimes\Omega.\]

         \item[(3)] Both the connections in item (1) and item (2) induce a $B$-linear connection 
         \[\nabla:B\{\underline V\}_{r,d,P}^{\dagger}\to B\{\underline V\}_{r,d,P}^{\dagger}\otimes\Omega \] 
         such that the following complex is exact:
         \[0\to B\to B\{\underline V\}_{r,d,P}^{\dagger}\to B\{\underline V\}_{r,d,P}^{\dagger}\otimes\Omega\to B\{\underline V\}_{r,d,P}^{\dagger}\otimes\wedge^2\Omega\to \cdots \to B\{\underline V\}_{r,d,P}^{\dagger}\otimes\wedge^d\Omega\to 0.\]
     \end{enumerate}
 \end{lem}
 \begin{proof}
     We treat the case where $P = \frac{\prod_{i=0}^r(1+bV_i)-1}{b}$; the other case follows similarly.
     Note that $(1+bV_i)\partial_i(P) = 1+bP$ for $0\leq i\leq r$ and $(1+bV_i)\partial_i(P) = 0$ for $r+1\leq i\leq d$.
     So for any $F\in B\{V_0/\lambda,\dots,V_d/\lambda\}$, we have
     \[\big(\sum_{i=0}^d(1+bV_i)\partial_i\otimes e_i\big)(PF) = P\sum_{i=0}^d(1+bV_i)\partial_i(F)\otimes e_i+(1+bP)F\otimes\sum_{i=0}^re_i.\]
     Thus, for any $f \in B\{\underline V/\lambda\}_{r,d,P}$ and any pre-image $F$ in $B\{V_0/\lambda,\dots,V_d/\lambda\}$, the image of 
     \[\big(\sum_{i=0}^d(1+bV_i)\partial_i\otimes e_i\big)(F)\in \oplus_{i=0}^dB\{V_0/\lambda,\dots,V_d/\lambda\}\cdot e_i\]
     via the natural surjection
     \[\oplus_{i=0}^dB\{V_0/\lambda,\dots,V_d/\lambda\}\cdot e_i\to B\{\underline V/\lambda\}_{r,d,P}\otimes\Omega,\]
     which is exactly $\nabla(f)$, is independent of the choice of $F$. Hence, $\nabla$ is well-defined. In particular, letting $\lambda\to 1^+$, we get a $B$-linear connection
     \[\nabla:B\{\underline V\}_{r,d,P}^{\dagger}\to B\{\underline V\}_{r,d,P}^{\dagger}\otimes\Omega\]
     as desired. It remains to show that the associated de Rham complex
     \[0\to B\to B\{\underline V\}_{r,d,P}^{\dagger}\to B\{\underline V\}_{r,d,P}^{\dagger}\otimes\Omega\to \cdots\to B\{\underline V\}_{r,d,P}^{\dagger}\otimes\wedge^d\Omega\to 0\]
     is exact.
     
     Now, consider the connection 
     \[\sum_{i=1}^d(1+bV_i)\partial_i\otimes e_i:B\{\underline V/\lambda\}\to \oplus_{i=1}^dB\{\underline V/\lambda\}\cdot e_i =  B\{\underline V/\lambda\}\otimes(\oplus_{i=1}^d\bZ\cdot e_i)\]
     which induces a connection
     \[\nabla:B\{\underline V\}^{\dagger}\to \oplus_{i=1}^dB\{\underline V\}^{\dagger}\cdot e_i\]
     by letting $\lambda\to 1^+$. Consider the map of $B$-algebras
     \[\iota:B\{\underline V/\lambda\}\to B\{\underline V/\lambda\}_{r,d,P}\]
     induced by the composite
     \[B\{\underline V/\lambda\}\xrightarrow{V_i\mapsto V_i,\forall 1\leq i\leq d} B\{V_0/\lambda,\dots,V_d/\lambda\}\to B\{\underline V/\lambda\}_{r,d,P},\]
     where the first map is given by the natural inclusion $B\{\underline V/\lambda\}\hookrightarrow B\{V_0/\lambda,\dots,V_d/\lambda\}$ sending each $V_i$ to $V_i$. Then $\iota$ is compatible with the connections in the sense that the following diagram is commutative
     \begin{equation*}
         \xymatrix@C=0.5cm{
         B\{\underline V/\lambda\}\ar[d]^{\iota}\ar[rrrrrr]^{\sum_{i=1}^d(1+bV_i)\partial_i\otimes e_i}&&&&&&  B\{\underline V/\lambda\}\otimes(\oplus_{i=1}^d\bZ\cdot e_i)\ar[d]^{\iota\otimes f}\\
         B\{\underline V/\lambda\}_{r,d,P}\ar[rrrrrr]^{\nabla}&&&&&& B\{\underline V/\lambda\}_{r,d,P}\otimes\Omega.}
     \end{equation*}
     where $f:\oplus_{i=1}^d\bZ\cdot e_i\to \Omega$ denotes the isomorphism of free $\bZ$-modules induced by the composite
     \[\oplus_{i=1}^d\bZ\cdot e_i\subset \oplus_{i=0}^d\bZ\cdot e_i\to \Omega.\]
     
     We claim $\iota$ is an isomorphism. Indeed, as $\lambda<|b|^{-1}$, for any $0\leq i\leq d$, we see that $1+bV_i$ is a unit in $B\{V_0/\lambda,\dots,V_d/\lambda\}$ and thus so is $1+bP =\prod_{i=0}^r(1+bV_i)$. So we see that the natural morphism
     \[i:B\{P/\lambda,V_1/\lambda,\dots,V_d/\lambda\}\to B\{V_0/\lambda,V_1/\lambda,\dots,V_d/\lambda\}\]
     sending $P$ to $\frac{\prod_{i=0}^r(1+bV_i)-1}{b}$ and each $V_i$ to $V_i$ for any $1\leq i\leq d$ is an isometric isomorphism of Banach $B$-algebras whose inverse is given by the morphism
     \[B\{V_0/\lambda,V_1/\lambda,\dots,V_d/\lambda\}\to B\{P/\lambda,V_1/\lambda,\dots,V_d/\lambda\}\]
     sending $V_0$ to $\frac{(1+bP)\cdot\prod_{i=1}^d(1+bV_i)^{-1}-1}{b}$ and each $V_i$ to $V_i$ for any $1\leq i\leq d$. Note that the isomorphism $i$ induces an isomorphism 
     \[\bar i: B\{P/\lambda,V_1/\lambda,\dots,V_d/\lambda\}/(P)\xrightarrow{\cong} B\{\underline V\}_{r,d,P}.\]
     The above claim follows as $\iota$ coincides with the composite
     \[B\{\underline V/\lambda\}\xrightarrow{\subset}B\{P/\lambda,V_1/\lambda,\dots,V_d/\lambda\}\to B\{P/\lambda,V_1/\lambda,\dots,V_d/\lambda\}/(P)\xrightarrow{\bar i} B\{\underline V\}_{r,d,P}.\]

     Consequently, $\iota$ induces an isomorphism
     \[\iota:B\{\underline V\}^{\dagger}\xrightarrow{\cong}B\{\underline V\}_{r,d,P}^{\dagger}\]
     which is compatible with the connections. To complete the proof, it suffices to show that the following complex associated to $(B\{\underline V\}^{\dagger},\nabla = \sum_{i=1}^d(1+bV_i)\partial_i\otimes e_i)$
     \[0\to B\to B\{\underline V\}^{\dagger} \to \oplus_{i=1}^dB\{\underline V\}^{\dagger}\cdot e_i\to \cdots \to \wedge^d(\oplus_{i=1}^dB\{\underline V\}^{\dagger}\cdot e_i)\to 0\]
     is exact. By the K\"unneth formula, we reduce to the case for $d = 1$. In other words, it is enough to show that the complex
     \[0\to B\to B\{V\}^{\dagger}\xrightarrow{(1+bV)\frac{\partial}{\partial V}}B\{V\}^{\dagger}\to 0\]
     is exact. As $1+bV$ is a unit in $B\{V\}^{\dagger}$, it suffices to show that the complex
     \[0\to B\to B\{V\}^{\dagger}\xrightarrow{\frac{\partial}{\partial V}}B\{V\}^{\dagger}\to 0\]
     is exact, which is well-known. This completes the proof.
\end{proof}

Now, we fix an affinoid perfectoid space $U = \Spa(S,S^+)\in X_{Z,v}$. Then $(\CCtilde^I(U),|\cdot|_I)$ is a $p$-adic Banach algebra and $\xi\in \CCtilde^I(U)$ such that $|\xi|_I<1$.

\begin{prop}\label{prop:local description}
    Fix a $\lambda\in[1,|\xi|_I^{-1})$. For any $\frakY_Z = \Spf(\wtr_Z^{\prime})\in \Sigma_U$, let $A_{\frakY_Z}^{I,\lambda}$ be as before. Then there is a well-defined isomorphism of $\CCtilde^I(U)$-algebras
    \[\iota:\CCtilde^I(U)\{\underline V/\lambda\}_{r,d,P}\xrightarrow{\cong}A_{\frakY_Z}^{I,\lambda}\]
    for $P = \frac{\prod_{i=0}^r(1+\xi V_i)-1}{\xi}$, defined by sending each $V_i$ to $\frac{\frac{T_i}{t_i}-1}{\xi}$ for any $0\leq i\leq d$. This isomorphism is compatible with the connection
    \[\nabla:A_{\frakY_Z}^{I,\lambda} \to A_{\frakY_Z}^{I,\lambda}\otimes\Omega^{1,\log}_{\wtr_Z^{\prime}}\{-1\}\]
    in the sense that it coincides with the connection 
    \[\nabla:\CCtilde^I(U)\{\underline V/\lambda\}_{r,d,P}\to \CCtilde^I(U)\{\underline V/\lambda\}_{r,d,P}\otimes\Omega\]
    via the isomorphism 
    \[\Omega^{1,\log}_{\wtr_Z^{\prime}}\{-1\}\cong \wtr_Z^{\prime}\otimes(\oplus_{i=0}^d\bZ\cdot e_i/\bZ\cdot\sum_{i=0}^re_i)\]
    where each $e_i$ is identified with the log-differential $\frac{\rd \log T_i}{\xi}$ for any $0\leq i\leq d$.
\end{prop}
\begin{proof}
    Without loss of generality, we may assume $\frakY_Z = \frakX_Z = \Spf(\calR_Z)$ and $\wtr_Z^{\prime} = \wtr_Z$. It suffices to show that the morphism $\iota$ is a well-defined isomorphism.
    Consider the morphism of $\CCtilde^I(U)$-algebras 
    \[\CCtilde^I(U)[V_0,\dots,V_d]\to A_{\frakX_Z}^{I,\lambda}\]
    sending each $V_i$ to $\frac{\frac{T_i}{t_i}-1}{\xi}$. Equip $\CCtilde^I(U)[V_0,\dots,V_d]$ with the norm $|\cdot|_{\lambda}$ as in Definition \ref{dfn:overconvergent algebras} (for $(B,|\cdot|) = (\CCtilde^I(U),|\cdot|_I)$). By definition of $A_{\frakX_Z}^{I,\lambda}$, we have $|\frac{\frac{T_i}{t_i}-1}{\xi}|\leq \lambda$ and thus the above morphism extends uniquely to a morphism
    \[\CCtilde^I(U)\{V_0/\lambda,\dots,V_d/\lambda\}\to A_{\frakX_Z}^{I,\lambda},\]
    which factors uniquely through the quotient morphism $\CCtilde^I(U)\{V_0/\lambda,\dots,V_d/\lambda\}\to \CCtilde^I(U)\{\underline V/\lambda\}_{r,d,P}$ because $A_{\frakY_Z}^{I,\lambda}$ is $\xi$-torsion free and $\prod_{i=0}^r\frac{T_i}{t_i} = 1$. This shows that the morphism $\iota$ is well-defined. 

    We now construct the inverse of $\iota$ explicitly. Keep notations as in Example \ref{exam:small affine is liftable}. Consider the morphism of $\bA_{\inf}(Z)$-algebras
    \[j_0:\bA_{\inf}(Z)\za T_0,\dots, T_r,T_{r+1}^{\pm 1},\dots,T_d^{\pm 1}\ya/(T_0\cdots T_r-[\varpi^a])\to \bA_{\inf}(U)\za \xi V_0,\dots, \xi V_d\ya/(\prod_{i=0}^r(1+\xi V_i)-1)\]
    sending each $T_i$ to $t_i(1+\xi V_i)$ for any $0\leq i\leq d$. Here, $\bA_{\inf}(U)\za \xi V_0,\dots, \xi V_d\ya$ denotes the $(p,\xi)$-adic completion of the polynomial ring $\bA_{\inf}(U)[\xi V_0,\dots, \xi V_d]$ over $\bA_{\inf}(U)$ with free variables $\xi V_0,\dots,\xi V_d$, and $(\prod_{i=0}^r(1+\xi V_i)-1)$ denotes the closed ideal generated by $\prod_{i=0}^r(1+\xi V_i)-1$. Note that the composite of $j_0$ with the morphism
    \[\bA_{\inf}(U)\za \xi V_0,\dots, \xi V_d\ya/(\prod_{i=0}^r(1+\xi V_i)-1)\to S^+\]
    sending each $\xi V_i$ to $0$ coincides with the composite 
    \[A^+\za T_0,\dots, T_r,T_{r+1}^{\pm 1},\dots,T_d^{\pm 1}\ya/(T_0\cdots T_r-p^a)\to \calR_Z\to S^+.\]
    By the \'etaleness of $\psi$, the morphism $j_0$ uniquely extends to a morphism
    \[\wtr_Z\to \bA_{\inf}(U)\za \xi V_0,\dots, \xi V_d\ya/(\prod_{i=0}^r(1+\xi V_i)-1)\]
    and thus, by linear extension, a morphism
    \[j_1:\wtr_Z\widehat \otimes_{\bA_{\inf}(Z)}\bA_{\inf}(U)\to \bA_{\inf}(U)\za \xi V_0,\dots, \xi V_d\ya/(\prod_{i=0}^r(1+\xi V_i)-1).\]
    As $j_1(T_i) = t_i(1+\xi V_i)$ is a unit-multiple of $t_i$, the morphism $j_1$ uniquely extends to a morphism
    \[j_2:A_{\frakY_Z}^+=\big(\wtr_Z\widehat \otimes_{\bA_{\inf}(Z)}\bA_{\inf}(U)\big)[(\frac{T_0}{t_0})^{\pm 1},\dots,(\frac{T_r}{t_r})^{\pm 1}]^{\wedge_{(p,\xi)}}\to \bA_{\inf}(U)\za \xi V_0,\dots, \xi V_d\ya/(\prod_{i=0}^r(1+\xi V_i)-1)\]
    sending each $\frac{T_i}{t_i}$ to $1+\xi V_i$ for any $0\leq i\leq d$. By construction, $j_2$ carries the ideal $\Ker(\theta_{\frakY_Z}^{\prime})$ into the ideal $(\xi V_0,\dots,\xi V_d)$. Note that the rational localization of
    \[\Spa(\big(\bA_{\inf}(U)\za \xi V_0,\dots, \xi V_d\ya/(\prod_{i=0}^r(1+\xi V_i)-1)\big)[\frac{1}{p}],\bA_{\inf}(U)\za \xi V_0,\dots, \xi V_d\ya/(\prod_{i=0}^r(1+\xi V_i)-1))\]
    with respect to the condition $|[\varpi]|\in |p|^I$ is exactly
    \[\Spa(\CCtilde^I(U)\za \xi V_0,\dots, \xi V_d\ya/(\prod_{i=0}^r(1+\xi V_i)-1),\CCtilde^{I,+}(U)\za \xi V_0,\dots, \xi V_d\ya/(\prod_{i=0}^r(1+\xi V_i)-1)),\]
    where 
    \[\CCtilde^{I,+}(U)\za \xi V_0,\dots, \xi V_d\ya/(\prod_{i=0}^r(1+\xi V_i)-1)) := \big(\bA_{\inf}(U)\za \xi V_0,\dots, \xi V_d\ya/(\prod_{i=0}^r(1+\xi V_i)-1))\big)\widehat \otimes_{\bA_{\inf}(U)}\CCtilde^{I,+}(U)\]
    coincides with the quotient of $\CCtilde^{I,+}(U)\za \xi V_0,\dots, \xi V_d\ya$ by the closed ideal generated by $\prod_{i=0}^r(1+\xi V_i)-1$, and 
    \[\CCtilde^I(U)\za \xi V_0,\dots, \xi V_d\ya/(\prod_{i=0}^r(1+\xi V_i)-1) = \big(\CCtilde^{I,+}(U)\za \xi V_0,\dots, \xi V_d\ya/(\prod_{i=0}^r(1+\xi V_i)-1))\big)[\frac{1}{p}].\]
    By the definition of $A_{\frakY_Z}$, $j_2$ uniquely extends to a morphism
    \[j_3:A_{\frakY_Z}^I\to \CCtilde^I(U)\za \xi V_0,\dots, \xi V_d\ya/(\prod_{i=0}^r(1+\xi V_i)-1).\]
    
    Now, equip $\CCtilde^I(U)\za \xi V_0,\dots, \xi V_d\ya$ with the Gauss norm $|\cdot|_G$ with respect to the norm $|\cdot|_I$ on $\CCtilde^I(U)$ such that $|\xi V_i|_G = 1$ for any $0\leq i\leq d$, and equip $\CCtilde^I(U)\za \xi V_0,\dots, \xi V_d\ya/(\prod_{i=0}^r(1+\xi V_i)-1)$ with the induced norm which we still denote by $|\cdot|_I$. It induces a norm $|\cdot|_I$ on 
    \[\big(\CCtilde^I(U)\za \xi V_0,\dots, \xi V_d\ya/(\prod_{i=0}^r(1+\xi V_i)-1)\big)[\frac{1}{\xi}]\]
    such that for any $x\in \CCtilde^I(U)\za \xi V_0,\dots, \xi V_d\ya/(\prod_{i=0}^r(1+\xi V_i)-1)$ and any $n\geq 0$, we have
    \[|\frac{x}{\xi^n}|_I:=\sup_{t\in I}|\frac{x}{\xi}|_t.\] 
    Equip the blow-up algebra 
    \[\big(\CCtilde^I(U)\za \xi V_0,\dots, \xi V_d\ya/(\prod_{i=0}^r(1+\xi V_i)-1)\big)[V_1,\dots,V_d]\subset\big(\CCtilde^I(U)\za \xi V_0,\dots, \xi V_d\ya/(\prod_{i=0}^r(1+\xi V_i)-1)\big)[\frac{1}{\xi}]\]
    with the norm $|\cdot|_I$ by restriction. Then for any $1\leq \lambda<|\xi|_I^{-1}$, the localization of 
    \[\Spa(\big(\CCtilde^I(U)\za \xi V_0,\dots, \xi V_d\ya/(\prod_{i=0}^r(1+\xi V_i)-1)\big)[V_1,\dots,V_d],\big(\CCtilde^I(U)\za \xi V_0,\dots, \xi V_d\ya/(\prod_{i=0}^r(1+\xi V_i)-1)\big)[V_1,\dots,V_d]^{\circ})\]
    with respect to the condition $|V_0|\leq \lambda,\dots,|V_d|\leq \lambda$ is exactly
    \[\Spa(\CCtilde^I(U)\{\underline V/\lambda\}_{r,d,P},\CCtilde^I(U)\{\underline V/\lambda\}_{r,d,P}^{\circ})\]
    where $\Box^{\circ}$ denotes the subring of power bounded elements in $\Box$. Therefore, by definition of $A_{\frakY_Z}^{I,\lambda}$, $j_3$ yields a morphism
    \[j:A_{\frakY_Z}^{I,\lambda}\to \CCtilde^I(U)\{\underline V/\lambda\}_{r,d,P}.\]
    Finally, it remains only to observe that $j$ is the inverse of $\iota$ as expected.
\end{proof}
\begin{cor}\label{cor:local description}
    We keep the notation introduced above. For any affinoid perfectoid $U\in X_{Z,v}$, the morphism of sheaves of $\CCtilde^I$-algebras on $X_{Z,v}/U$
        \[\CCtilde^I\{\underline V/\lambda\}_{r,d,P}\to \calO\CCtilde^{I,\lambda} \text{ (resp. $\CCtilde^I\{\underline V\}^{\dagger}_{r,d,P}\to \calO\CCtilde^{I,\dagger}$)}\]
        for $P=\frac{\prod_{i=0}^r(1+\xi V_i)-1}{\xi}$ which sends each $V_i$ to $\frac{\frac{T_i}{t_i}-1}{\xi}$ for any $0\leq i\leq d$ is an isomorphism such that the connection \[\nabla:\calO\CCtilde^{I,\lambda}\to \calO\CCtilde^{I,\lambda}\otimes_{\calO_{\wtx_Z}}\Omega^{1,\log}_{\wtx_Z}\{-1\} \text{ (resp. $\nabla:\calO\CCtilde^{I,\dagger}\to \calO\CCtilde^{I,\dagger}\otimes_{\calO_{\wtx_Z}}\Omega^{1,\log}_{\wtx_Z}\{-1\}$)}\] 
        is given by 
        \[\nabla = \sum_{i=0}^d(1+\xi V_i)\frac{\partial}{\partial V_i}\otimes \frac{\dlog T_i}{\xi}\]
        in the sense of Lemma \ref{lem:induced derivation} via the isomorphism 
        \[\Omega^{1,\log}_{\wtx_Z} \cong \oplus_{i=0}^d\calO_{\wtx_Z}\cdot e_i/\calO_{\wtx_Z}\cdot\sum_{i=0}^re_i,\] 
        where $e_i$ represents the log-differential $\dlog T_i $ for any $0\leq i\leq d$.
\end{cor}
\begin{proof}
    This follows from Proposition \ref{prop:local description} immediately.
\end{proof}
\begin{prop}\label{prop:local Poincare lemma for OC^I}
    We keep the notation introduced above. Then the following complex of $\CCtilde^I$-modules on $X_{Z,v}$ is exact:
    \[0\to \CCtilde^I\to \calO\CCtilde^{I,\dagger}\xrightarrow{\nabla}\calO\CCtilde^{I,\dagger}\otimes_{\calO_{\wtx_Z}}\Omega^{1,\log}_{\wtx_Z}\{-1\}\to\cdots\to\calO\CCtilde^{I,\dagger}\otimes_{\calO_{\wtx_Z}}\Omega^{d,\log}_{\wtx_Z}\{-d\}\to 0.\]
\end{prop}
\begin{proof}
    This follows directly from Corollary \ref{cor:local description} together with Lemma \ref{lem:induced derivation}(3).
\end{proof}

\subsection{Change of variables}\label{ssec:change of variable}
Let $\frakX_Z = \Spf(\calR_Z)$ be small and semi-stable and Keep the notation as above. We will give another description of $(\calO\CCtilde^{I,\dagger},\nabla)$. We start with some lemmas.

\begin{lem}\label{lem:change of variable}
    Let $(B,|\cdot|)$ be a $p$-adic Banach algebra and $b\in B$ be an element such that $|b|<p^{-\frac{1}{p-1}}$. Let $\lambda\in [1,|b|^{-1}p^{-\frac{1}{p-1}})$.
    Put $E(T):=\frac{\exp(T)-1}{T} = \sum_{n\geq 1}\frac{T^{n-1}}{n!}$ and $L(T):=\frac{\log(1+T)}{T} = \sum_{n\geq 1}\frac{(-1)^{n-1}T^{n-1}}{n}$.
    \begin{enumerate}
        \item The morphism of Banach $B$-algebras
        \[\iota:B\{V_0/\lambda,\dots,V_d/\lambda\}\to B\{W_0/\lambda,\dots,W_d/\lambda\}\]
        sending each $V_i$ to $W_i\cdot E(bW_i)$ for each $0\leq i\leq d$ is a well-defined isometric isomorphism (with respect to the norm $|\cdot|_{\lambda}$ on both sides) whose inverse is given by the morphism
        \[\iota^{-1}:B\{W_0/\lambda,\dots,W_d/\lambda\}\to B\{V_0/\lambda,\dots,V_d/\lambda\}\]
        sending each $W_i$ to $V_i\cdot L(bV_i)$ for any $0\leq i\leq d$. Moreover, the isomorphism $\iota$ fits into the following commutative diagram
        \begin{equation*}
            \xymatrix@C=0.5cm{
              B\{V_0/\lambda,\dots,V_d/\lambda\}\ar[rrrrrr]^{\sum_{i=0}^d(1+bV_i)\frac{\partial}{\partial V_i}\otimes e_i}\ar[d]^{\iota}&&&&&&\oplus_{i=0}^dB\{V_0/\lambda,\dots,V_d/\lambda\}\cdot e_i\ar[d]^{\sum_{i=0}^d\iota\otimes e_i}\\
              B\{W_0/\lambda,\dots,W_d/\lambda\}\ar[rrrrrr]^{\sum_{i=0}^d\frac{\partial}{\partial W_i}\otimes e_i}&&&&&&\oplus_{i=0}^dB\{W_0/\lambda,\dots,W_d/\lambda\}\cdot e_i.
            }
        \end{equation*}

        \item[(2)] Let $P:=\frac{\prod_{i=0}^r(1+bV_i)-1}{b}$ and $Q=\sum_{i=0}^rW_i$, and let $\nabla_P$ and $\nabla_Q$ respectively be the connections on $B\{\underline V/\lambda\}_{r,d,P}$ and $B\{\underline W/\lambda\}_{r,d,Q}$ in the sense of Lemma \ref{lem:induced derivation}. Then the isomorphism $\iota$ above induces an isomorphism 
        \[\iota:B\{\underline V/\lambda\}_{r,d,P}\to B\{\underline W/\lambda\}_{r,d,Q}\]
        which is compatible with the connections in the sense that the following diagram commutes
        \begin{equation*}
            \xymatrix@C=0.5cm{
             B\{\underline V/\lambda\}_{r,d,P}\ar[rrrrrr]^{\nabla_P}\ar[d]^{\iota}&&&&&&B\{\underline V/\lambda\}_{r,d,P}\otimes\Omega\ar[d]^{\iota\otimes\id}\\
             B\{\underline W/\lambda\}_{r,d,Q}\ar[rrrrrr]^{\nabla_Q}&&&&&&B\{\underline W/\lambda\}_{r,d,Q}\otimes\Omega.
            }
        \end{equation*}

        \item[(3)] The isomorphism $\iota$ above induces an isomorphism 
        \[\iota:B\{\underline V\}_{r,d,P}^{\dagger}\to B\{\underline W\}_{r,d,Q}^{\dagger}\]
        which is compatible with the connection as described in Lemma \ref{lem:induced derivation}.
    \end{enumerate}
\end{lem}
\begin{proof}
    For item (1): Note that for any $0\leq i\leq d$, we have 
    \[\frac{\partial}{\partial W_i}(W_iE(bW_i)) = \exp(bW_i)\]
    and thus 
    \[(1+bW_iE(W_i))\frac{\partial}{\partial(W_iE(W_i))} = \frac{\partial}{\partial W_i}.\]
    It suffices to show that that both $\iota$ and $\iota^{-1}$ are well-defined and bounded by $1$.
    To see this for $\iota$, it suffices to show that for any $0\leq i\leq d$, the series
    \[E(bW_i) = \sum_{n\geq 1}\frac{b^{n-1}}{n!}W_i^{n-1} = 1+\sum_{n\geq 1}\frac{b^n}{(n+1)!}W^n\]
    converges to a unit in $B\{\underline W/\lambda\}_{r,d,Q}$ satisfying $|E(bW_i)|_{\lambda}=1$.
    Recall that for any $n\geq 1$, if we denote by $e_n$ the maximal integer such that $p^{e_n}\mid n!$, then $e_n\leq \frac{n-1}{p-1}$.
    In particular, we have 
    \[\lambda^n\cdot|\frac{b^n}{(n+1)!}|\leq \lambda^n\cdot |b|^n\cdot p^{\frac{n}{p-1}} = (\lambda\cdot |b|\cdot p^{\frac{1}{p-1}})^n.\]
    As $\lambda<|b|^{-1}p^{-\frac{1}{p-1}}$, we have $\lambda\cdot |b|\cdot p^{\frac{1}{p-1}}<1$.
    As for any $n\geq 1$, 
    \[|\frac{b^n}{(n+1)!}W^n|_{\lambda}\leq (\lambda\cdot |b|\cdot p^{\frac{1}{p-1}})^n,\]
    we see that $E(bW_i)$ is a well-defined unit in $B\{\underline W/\lambda\}_{r,d,Q}$ satisfying $|E(bW_i)|_{\lambda}=1$ as desired.
    A similar argument shows that $\iota^{-1}$ is well-defined and bounded by $1$. This completes the proof of Item (1).

    For item (2): Note that $\iota(P) = \frac{\prod_{i=0}^r\exp(bW_i)-1}{b} = (\sum_{i=0}^rW_i)\cdot E(b\sum_{i=0}^rW_i)$.
    The assertion follows from Item (1) and the fact that $(\sum_{i=0}^rW_i)\cdot E(b\sum_{i=0}^rW_i)$ generates the same ideal of $B\{W_0/\lambda,\dots,W_d/\lambda\}$ as $\sum_{i=0}^rW_i$ (because $E(b\sum_{i=0}^rW_i)$ is a unit).

    For item (3): this follows from Item (2) by letting $\lambda\to 1^+$.
\end{proof}

Recall that for any closed sub-interval $I = [r,s]$ of $(\frac{1}{p-1},+\infty)$, we have $|\xi|_I<p^{-\frac{1}{p-1}}$ and $t=\log[\epsilon]\in \CCtilde^{I,+}(Z)$ is a unit multiple of $u\xi$ with $u = [\epsilon^{\frac{1}{p}}]-1$.

\begin{prop}\label{prop:change of variable}
    Let $I$ be a closed sub-interval of $(\frac{1}{p-1},+\infty)$ and Keep the notation as above. Let $\lambda\in[1,|\xi|_I^{-1}p^{-\frac{1}{p-1}})$ and $P = W_0+\cdots+W_r$. Then for any affinoid perfectoid $U\in X_{Z,v}$, the morphism of sheaves of $\CCtilde^I$-algebras on $X_{Z,v}/U$
    \[\iota:\CCtilde^{I,\lambda}\{\underline W/\lambda\}_{r,d,P}\to \calO\CCtilde^{I,\lambda},\]
    which sends each $W_i$ to $\frac{u}{t}\log\frac{T_i}{t_i}$ for any $0\leq i\leq d$, is an isomorphism such that via the isomorphism
    \[\calO\CCtilde^{I,\lambda}\otimes_{\calO_{\wtx_Z}}\Omega^{1,\log}_{\wtx_Z}\{-1\}\cong \calO\CCtilde^{I,\lambda}\otimes\big(\oplus_{i=0}^d\bZ\cdot e_i/\bZ\cdot\sum_{i=0}^re_i\big),\]
    where $e_i$ stands for $\frac{u}{t}\dlog T_i$ for any $0\leq i\leq d$, the connection 
    \[\nabla:\calO\CCtilde^{I,\lambda}\to\calO\CCtilde^{I,\lambda}\otimes_{\calO_{\wtx_Z}}\Omega^{1,\log}_{\wtx_Z}\{-1\}\]
    is given by
    \[\nabla = \sum_{i=0}^d\frac{\partial}{\partial W_i}\otimes e_i:\CCtilde^I\{\underline W/\lambda\}_{r,d,P}\to \CCtilde^I\{\underline W/\lambda\}_{r,d,P}\otimes\big(\oplus_{i=0}^d\bZ\cdot e_i/\bZ\cdot\sum_{i=0}^re_i\big)\]
    in the sense of Lemma \ref{lem:induced derivation}. The same result holds true for respectively replacing $\CCtilde^I\{\underline W/\lambda\}_{r,d,P}$ and $\calO\CCtilde^{I,\lambda}$ by $\CCtilde^I\{\underline W\}_{r,d,P}^{\dagger}$ and $\calO\CCtilde^{I,\dagger}$ instead.
\end{prop}
\begin{proof}
    As $\frac{t}{u}$ is a unit-multiple of $\xi$ in $\CCtilde^{I,+}(Z)$, the result follows from Corollary \ref{cor:local description} and Lemma \ref{lem:change of variable}(2) immediately.
\end{proof}

We end this subsection with the following example.

\begin{exam}\label{exam:Apply to Gamma torsor}
    Suppose that $\frakX_Z = \Spf(\calR_Z)$ is small semi-stable with generic fiber $X_Z = \Spa(R_Z,R_Z^+)$ and let $I=[r,s]$ be a closed sub-interval of $(\frac{1}{p-1},+\infty)$ as before. Let $U = X_{\infty,Z}\in X_{Z,v}$ be the covering of $X_{Z,v}$ as in Example \ref{exam:Gamma torsor} and keep notations there. For any $\lambda\in[1,|\xi|_I^{-1}p^{\frac{1}{p-1}})$, by Proposition \ref{prop:change of variable}, for $P = W_0+\cdots+W_r$, sending each $W_i$ to $\frac{u}{t}\log\frac{[T_i]^{\flat}}{T_i}$ gives rise to an isomorphism
    \[\iota:\CCtilde^I\{\underline W/\lambda\}_{r,d,P} = \CCtilde^I\{W_0/\lambda,\dots,W_d/\lambda\}/(W_0+\cdots+W_r)\to \calO\CCtilde^{I,\lambda}\]
    of sheaves of $\CCtilde^I$-modules on $X_{Z,v}/X_{\infty,Z}$ such that the connection $\nabla$ on $\calO\CCtilde^{I,\lambda}$ is given by
    \[\nabla = \sum_{i=0}^d-\frac{\partial}{\partial W_i}\otimes\frac{u}{t}\dlog T_i.\]
    Note that $\calO\CCtilde^{I,\lambda}$ admits a natural action of $\Gamma$ so that $\nabla$ is $\Gamma$-equivariant. Via the above isomorphism $\iota$, by \eqref{equ:gamma acts on T} and \eqref{equ:gamma acts on T^flat}, this action is determined by
    \begin{equation}\label{equ:gamma acts on W}
        \delta_i(W_j) = W_j+\delta_{ij}u
    \end{equation} 
    for any $0\leq i,j\leq d$.
    Note that the natural morphism
    \[\CCtilde^I\{\underline W/\lambda\} = \CCtilde^I\{W_1/\lambda,\dots,W_d/\lambda\}\to \CCtilde^I\{\underline W/\lambda\}_{r,d,P}\]
    sending each $W_i$ to $W_i$ is an isomorphism. So we get an isomorphism
    \[\CCtilde^I\{\underline W/\lambda\}\to \calO\CCtilde^{I,\lambda}\]
    of sheaves of $\CCtilde^I$-algebras on $X_{Z,v}/X_{\infty,Z}$,
    which sends each $W_i$ to $\frac{u}{t}\log\frac{[T_i^{\flat}]}{T_i}$, such that via the isomorphism
    \[\Omega^{1,\log}_{\wtx_Z}\cong \oplus_{i=1}^d\calO_{\wtx_Z}\cdot\dlog T_i,\]
    the connection $\nabla$ on $\calO\CCtilde^{I,\lambda}$ is given by
    \[\nabla = \sum_{i=1}^d-\frac{\partial}{\partial W_i}\otimes\frac{u}{t}\dlog T_i.\]
    Moreover, via the isomorphism \eqref{equ:gamma group-II}, the $\Gamma$-action on $\calO\CCtilde^{I,\lambda}$ is given by
    \begin{equation}\label{equ:gamma acts on W-II}
        \gamma_i(W_j) = W_j+\delta_{ij}u
    \end{equation}
    for any $1\leq i,j\leq d$.
\end{exam}

\subsection{The general case}\label{ssec:period sheaf-general case}
We now turn to the general case. Let $\frakX_Z$ be a semi-stable formal scheme over $A^+$ of relative dimension $d$ with a fixed lifting $\wtx_Z$ and the generic fiber $X_Z$ as before. Denote by $X_{Z,v}^{\Box}\subset X_{Z,v}$ the sub-site consisting of affinoid perfectoid $U = \Spa(S,S^+)$ such that the composite
\[U\to X_Z\to \frakX_Z\]
factors as
\[U\to Y_Z\to \frakY_Z\to \frakX_Z\]
for some \'etale localization $\frakY_Z\to \frakX_Z$ such that $\frakY_Z$ is small semi-stable with generic fiber $Y_Z$.
Since $\frakX_{Z,\et}$ admits a basis consisting of small affines, we see that $X_{Z,v}^{\Box}$ gives a basis for the topology of $X_{Z,v}$.

\begin{dfn}\label{dfn:OC^I}
    Let $I = [r,s]$ be a closed sub-interval of $(0,+\infty)$ and $\lambda\in [1,|\xi|_I^{-1})$ as before. Denote by $\calO\CCtilde^{I,\lambda}$ the period sheaf on $X_{Z,v}$ associated to the sheafification of the following presheaf on $X_{Z,v}^{\Box}$
    \[U\in X_{Z,v}^{\Box}\mapsto \calO\CCtilde^{I,\lambda}(U),\]
    where $\calO\CCtilde^{I,\lambda}(U)$ is defined in Definition \ref{dfn:local OC^I}, and let 
    \[\calO\CCtilde^{I,\dagger}:=\colim_{\lambda\to 1^+}\calO\CCtilde^{I,\lambda}.\]
    Then the usual differential $\rd:\calO_{\wtx_Z}\to \Omega^{1,\log}_{\wtx_Z}$ induces a $\CCtilde^I$-linear connection
    \[\nabla:\calO\CCtilde^{I,\lambda}\to \calO\CCtilde^{I,\lambda}\otimes_{\calO_{\wtx_Z}}\Omega^{1,\log}_{\wtx_Z}\{-1\}\]
    and thus a $\CCtilde^I$-linear connection
    \[\nabla:\calO\CCtilde^{I,\dagger}\to \calO\CCtilde^{I,\dagger}\otimes_{\calO_{\wtx_Z}}\Omega^{1,\log}_{\wtx_Z}\{-1\}.\]
\end{dfn}

\begin{thm}[Poincar\'e Lemma for $(\calO\CCtilde^{I,\dagger},\nabla)$]\label{thm:Poincare Lemma for OC^I}
    The following sequence of sheaves of $\CCtilde^I$-modules on $X_{Z,v}$ is exact:
    \[0\to \CCtilde^I\to \calO\CCtilde^{I,\dagger}\xrightarrow{\nabla} \calO\CCtilde^{I,\dagger}\otimes_{\calO_{\wtx_Z}}\Omega^{1,\log}_{\wtx_Z}\{-1\}\to \cdots \to \calO\CCtilde^{I,\dagger}\otimes_{\calO_{\wtx_Z}}\Omega^{d,\log}_{\wtx_Z}\{-d\}\to 0.\]
\end{thm}
\begin{proof}
    Since the problem is local on $\frakX_{Z,\et}$, we may assume that $\frakX_Z$ is small affine and conclude the result from Proposition \ref{prop:local Poincare lemma for OC^I}.
\end{proof}

\subsection{Relative Poincar\'e Lemma}\label{ssec:relative OC}
Let $f:\frakX^{\prime}_Z\to \frakX_Z$ be a smooth morphism of semi-stable formal schemes over $A^+$ of relative dimension $d^{\prime}-d$ with generic fiber $X^{\prime}_Z\to X_Z$. Here, $d$ and $d^{\prime}$ denotes the relative dimension of $\frakX_Z$ and $\frakX^{\prime}_{Z}$ over $A^+$ respectively. Suppose that both $\frakX_Z$ and $\frakX^{\prime}_Z$ admit liftings $\wtx_Z$ and $\wtx^{\prime}_Z$ over $\bA_{\inf}(Z)$ respectively such that $f$ extends to a smooth morphism $\widetilde f:\wtx^{\prime}_Z\to\wtx_Z$. Denote by 
\[\nabla_{\wtx^{\prime}_Z/\wtx_Z}:\calO_{\wtx^{\prime}_Z}\to \Omega^1_{\wtx^{\prime}_Z/\wtx_Z}\{-1\}\]
the composite
\[\calO_{\wtx^{\prime}_Z}\xrightarrow{\rd_{\wtx^{\prime}_Z/\wtx_Z}} \Omega^1_{\wty_Z/\wtx_Z}\subset \Omega^1_{\wtx^{\prime}_Z/\wtx_Z}\{-1\}\]
where $\rd_{\wtx^{\prime}_Z/\wtx_Z}:\calO_{\wtx^{\prime}_Z}\to \Omega^1_{\wtx^{\prime}_Z/\wtx_Z}$ denotes the usual $(p,\xi)$-complete relative derivation. For any $\calP\in \{\calO\CCtilde^{I,\lambda},\calO\CCtilde^{I,\dagger}\}$, denote by $(\calP_{\wtx_Z},\nabla_{\wtx_Z})$ (resp. $(\calP_{\wtx^{\prime}_Z},\nabla_{\wtx^{\prime}_Z})$) the corresponding period sheaf with connection associated to $\frakX_Z$ (resp. $\frakX^{\prime}_Z$) with respect to the given lifting $\wtx_Z$ (resp. $\wtx^{\prime}_Z$). The purpose of this subsection is to construct a morphism $\widetilde f^*\calO\CCtilde^{I,\dagger}_{\wtx_Z}\hookrightarrow\calO\CCtilde^{I,\dagger}_{\wtx^{\prime}_Z}$, by using the lifting $\widetilde f$, such that the relative Poincar\'e Lemma holds true; that is, the following sequence
\[0\to f^*\calO\CCtilde^{I,\dagger}_{\wtx_Z}\hookrightarrow\calO\CCtilde^{I,\dagger}_{\wtx^{\prime}_Z}\xrightarrow{\nabla_{\wtx^{\prime}_Z/\wtx_Z}}\calO\CCtilde^{I,\dagger}_{\wtx^{\prime}_Z}\otimes_{\calO_{\wtx^{\prime}_Z}}\Omega^1_{\wtx^{\prime}_Z/\wtx_Z}\{-1\}\to \cdots \]
is exact, where 
\[\nabla_{\wtx^{\prime}_Z/\wtx_Z}:\calO\CCtilde^{I,\dagger}_{\wtx^{\prime}_Z}\to\calO\CCtilde^{I,\dagger}_{\wtx^{\prime}_Z}\otimes_{\calO_{\wtx^{\prime}_Z}}\Omega^1_{\wtx^{\prime}_Z/\wtx_Z}\{-1\}\]
is induced by the relative differential $\rd_{\wtx^{\prime}_Z/\wtx_Z}$.

Since $f$ is smooth, \'etale locally on both $\frakX_{Z,\et}$ and $\frakX^{\prime}_{Z,\et}$, it is of the form
\begin{equation}\label{equ:relative chart-I}
    \xymatrix@C=0.45cm{
      A^+\za T_0,\dots,T_r,T_{r+1}^{\pm 1},\dots,T_d^{\pm 1}\ya/(T_0\cdots T_r-p^a)\ar[d]^{\psi}\ar[r]^{\subset}&A^+\za T_0,\dots,T_r,T_{r+1}^{\pm 1},\dots,T_{d^{\prime}}^{\pm 1}\ya/(T_0\cdots T_r-p^a)\ar[d]^{\psi^{\prime}}
      \\
      \calR_Z\ar[r]^{f}&\calR^{\prime}_Z}
\end{equation}
where $\psi$ and $\psi^{\prime}$ denotes the charts on $\frakY_Z=\Spf(\calR_Z)\in \frakX_{Z,\et}$ and $\frakY_Z^{\prime} = \Spf(\calR_Z^{\prime})\in \frakX^{\prime}_{Z,\et}$, respectively. So the lifting $\widetilde f:\wty_Z^{\prime}=\Spf(\wtr_Z^{\prime})\to\wty_Z=\Spf(\wtr_Z)$ fits into the following commutative diagram
\begin{equation}\label{equ:relative chart-II}
    \xymatrix@C=0.45cm{
      \bA_{\inf}(Z)\za T_0,\dots,T_r,T_{r+1}^{\pm 1},\dots,T_d^{\pm 1}\ya/(T_0\cdots T_r-[\varpi^a])\ar[d]^{\widetilde \psi}\ar[r]^{\subset}&\bA_{\inf}(Z)\za T_0,\dots,T_r,T_{r+1}^{\pm 1},\dots,T_{d^{\prime}}^{\pm 1}\ya/(T_0\cdots T_r-[\varpi^a])\ar[d]^{\widetilde \psi^{\prime}}
      \\
      \wtr_Z\ar[r]^{\widetilde f}&\wtr^{\prime}_Z},
\end{equation}
where $\widetilde \psi$ and $\widetilde \psi^{\prime}$ are the lifting of $\psi$ and $\psi^{\prime}$. Observe that for any $U = \Spa(S,S^+)\in X^{\prime}_{Z,v}$, we may regard it as an object in $X_{Z,v}$. Then we have the following commutative diagram
\[\xymatrix@C=0.5cm{\wtr_Z\widehat \otimes_{\bA_{\inf}(Z)}\bA_{\inf}(U)\ar[d]^{\theta_{\fraky_Z}}\ar[rr]^{\widetilde f\otimes\id}&&\wtr^{\prime}_Z\widehat \otimes_{\bA_{\inf}(Z)}\bA_{\inf}(U)\ar[d]^{\theta_{\fraky^{\prime}_Z}}\\
S^+\ar@{=}[rr]&& S^+},\]
yielding a commutative diagram by taking exactification \eqref{equ:exactification}
\[\xymatrix@C=0.5cm{A_{\frakY_Z}^+\ar[d]^{\theta^{\prime}_{\fraky_Z}}\ar[rrrr]^{\widetilde f^{\prime}}&&&&A_{\frakY_Z^{\prime}}^+\ar[d]^{\theta^{\prime}_{\fraky^{\prime}_Z}}\\
S^+\ar@{=}[rrrr]&&&& S^+}.\]
Thanks to Remark \ref{rmk:isometry for R^I}, for any $x\in A_{\frakY_Z}^+$, we have $|\widetilde f^{\prime}(x)|_I\leq |x|_I$. Therefore, the above morphism $\widetilde f^{\prime}:A_{\frakY_Z}^+\to A_{\frakY^{\prime}_Z}^+$ extends uniquely to a morphism $A_{\frakY_Z}^I\to A_{\frakY_Z^{\prime}}^I$, and thus a morphism
\[\widetilde f^{\prime\prime}: A_{\frakY_Z}^I[\frac{\ker(\theta_{\frakY_Z}^{\prime})}{\xi}]\to A_{\frakY_Z^{\prime}}^I[\frac{\ker(\theta_{\frakY^{\prime}_Z}^{\prime})}{\xi}]\]
such that for any $x\in \widetilde f^{\prime\prime}: A_{\frakY_Z}^I[\frac{\ker(\theta_{\frakY_Z}^{\prime})}{\xi}]$, $|\widetilde f^{\prime\prime}(x)|_I\leq |x|_I$. Thus, for any $1\leq \lambda<|\xi|_I^{-1}$, we get a morphism
\[\widetilde f^{I,\lambda}:A_{\frakY_Z}^{I,\lambda}\to A_{\frakY^{\prime}_Z}^{I,\lambda}.\]
To emphasize the dependence of $A_{\frakY_Z}^{I,\lambda}$ on the choice of the affinoid perfectoid $U$, in what follows, we denote by $A_{\frakY_Z}^{I,\lambda}(U)$ the corresponding $A_{\frakY_Z}^{I,\lambda}$ associated to $U$. Let $V\in Y_{Z,v}$ be affinoid perfectoid such that the morphism $U\to Y_Z$ factors as $U\to V\to Y_Z$. By Proposition \ref{prop:local description}, the scalar extension
\[A_{\frakY_Z}^{I,\lambda}(V)\widehat \otimes_{\CCtilde^I(V)}\CCtilde^I(U)\to A_{\frakY_Z}^{I,\lambda}(U)\]
is an isomorphism. Varying $U$ and $V$, we see that the presheaf
\[U\in Y^{\prime}_{Z,v}\to A_{\frakY_Z}^{\lambda,I}(U)\]
is indeed a sheaf which is exactly the pull-back $f^*\calO\CCtilde^{I,\lambda}_{\widetilde \frakY_Z}$. Therefore, by varying $U$, the morphism $A_{\frakY_Z}^{I,\lambda}(U)\to A_{\frakY^{\prime}_Z}^{I,\lambda}(U)$ induces a morphism of sheaves of $\CCtilde^I_{Y^{\prime}_Z}$-algebras
\[f^*\calO\CCtilde^{I,\lambda}_{\wty_Z}\to \calO\CCtilde^{I,\lambda}_{\wty_Z^{\prime}}.\]
Finally, varying $\frakY_Z$ and $\frakY_Z^{\prime}$ yields a morphism of sheaves of $\CCtilde^I_{X_Z^{\prime}}$-algebras
\begin{equation}\label{equ:relative morphism on OC^Ilambda}
    f^*\calO\CCtilde^{I,\lambda}_{\wtx_Z}\to \calO\CCtilde^{I,\lambda}_{\wtx_Z^{\prime}},
\end{equation}
and thus, by letting $\lambda\to 1^+$, a morphism of sheaves of $\CCtilde^I_{X_Z^{\prime}}$-algebras
\begin{equation}\label{equ:relative morphism on OC^Idagger}
    f^*\calO\CCtilde^{I,\dagger}_{\wtx_Z}\to \calO\CCtilde^{I,\dagger}_{\wtx_Z^{\prime}}.
\end{equation}
We claim this is the desired morphism.

\begin{prop}\label{prop:relative Poincare over OC^I}
    The following sequence of sheaves of $\CCtilde^I_{X_Z^{\prime}}$-algebras
    \[0\to f^*\calO\CCtilde^{I,\dagger}_{\wtx_Z}\hookrightarrow\calO\CCtilde^{I,\dagger}_{\wtx^{\prime}_Z}\xrightarrow{\nabla_{\wtx^{\prime}_Z/\wtx_Z}}\calO\CCtilde^{I,\dagger}_{\wtx^{\prime}_Z}\otimes_{\calO_{\wtx^{\prime}_Z}}\Omega^1_{\wtx^{\prime}_Z/\wtx_Z}\{-1\}\to \cdots \to \calO\CCtilde^{I,\dagger}_{\wtx^{\prime}_Z}\otimes_{\calO_{\wtx^{\prime}_Z}}\Omega^{d^{\prime}-d}_{\wtx^{\prime}_Z/\wtx_Z}\{-(d^{\prime}-d)\}\to 0\]
    is exact, where the connection 
    \[\nabla_{\wtx^{\prime}_Z/\wtx_Z}:\calO\CCtilde^{I,\dagger}_{\wtx^{\prime}_Z}\to\calO\CCtilde^{I,\dagger}_{\wtx^{\prime}_Z}\otimes_{\calO_{\wtx^{\prime}_Z}}\Omega^1_{\wtx^{\prime}_Z/\wtx_Z}\{-1\}\]
    is induced by the relative differential $\rd_{\wtx_Z^{\prime}/\wtx_Z}:\calO_{\wtx_Z^{\prime}}\to \Omega^1_{\wtx_Z^{\prime}/\wtx_Z}$.
\end{prop}
\begin{proof}
   The problem is local on $\frakX_{Z,\et}, \frakX^{\prime}_{Z,\et}$ and $X^{\prime}_{Z,v}$. We may assume that $\frakX_Z = \Spf(\calR_Z)$ and $\frakX_Z^{\prime} = \Spf(\calR_Z^{\prime})$ such that both \eqref{equ:relative chart-I} and \eqref{equ:relative chart-II} hold true. Then it suffices to show the above sequence remains exact when taking restriction to any affinoid perfectoid $U\in X_{Z,v}^{\prime}$. For a suitable choice of $t_0,\dots,t_d, t_{d+1},\dots,t_{d^{\prime}}$ as in Lemma \ref{lem:finite ideal}, it follows from Corollary \ref{cor:local description} that, if we put $V_i:=\frac{\frac{T_i}{t_i}-1}{\xi}$ for any $0\leq i\leq d^{\prime}$, then on $X_{Z,v}^{\prime}/U$, we have isomorphisms 
   \[f^*\calO\CCtilde^{I,\dagger}_{\wtx_Z}\cong \CCtilde^I\{\underline V\}_{r,d,P}^{\dagger} \text{ and }\calO\CCtilde^{I,\dagger}_{\wtx^{\prime}_Z}\cong \CCtilde^I\{\underline V\}_{r,d^{\prime},P}^{\dagger}\]
   such that the morphism \eqref{equ:relative morphism on OC^Idagger} becomes the natural inclusion
   \[\CCtilde^I\{\underline V\}_{r,d,P}^{\dagger}\to \CCtilde^I\{\underline V\}_{r,d^{\prime},P}^{\dagger}.\]
   In this case, via the isomorphism
   \[\Omega^1_{\wtx_Z^{\prime}/\wtx_Z}\{-1\}\cong \oplus_{i=d+1}^{d^{\prime}}\calO_{\wtx_Z^{\prime}}\cdot\frac{\rd T_i}{\xi},\]
   we deduce that 
   \[\nabla_{\wtx_Z^{\prime}/\wtx_Z} = \sum_{i=d+1}^{d^{\prime}}\frac{\partial}{\partial V_i}\otimes\frac{\rd T_i}{\xi}: \CCtilde^I\{\underline V\}_{r,d^{\prime},P}^{\dagger}\to \oplus_{i=d+1}^{d^{\prime}}\CCtilde^I\{\underline V\}_{r,d^{\prime},P}^{\dagger}\cdot \frac{\rd T_i}{\xi}.\]
   Now, the desired exactness follows easily (cf. the proof of Lemma \ref{lem:induced derivation}).
\end{proof}

\begin{exam}\label{exam:relative OC}
    Suppose that $f:\frakX_Z^{\prime} = \Spf(\calR_Z^{\prime})\to \frakX_Z = \Spf(\calR_Z)$ is of the form \eqref{equ:relative chart-I} such that its lifting $\widetilde f:\wtx_Z^{\prime} = \Spf(\wtr_Z^{\prime})\to\wtx_Z = \Spf(\wtr_Z)$ is of the form \eqref{equ:relative chart-II}. Denote by \[
    \text{$X_{\infty,Z}^{\prime} = \Spa(\widehat R_{\infty,Z}^{\prime},\widehat R_{\infty,Z}^{\prime+})\to X_Z^{\prime} = \Spa(R_Z^{\prime},R_Z^{\prime+})$ (resp. $X_{\infty,Z} = \Spa(\widehat R_{\infty,Z},\widehat R_{\infty,Z}^+)\to X_Z = \Spa(R_Z,R_Z^+)$)}\]
    the pro-finite-\'etale covering with Galois group $\Gamma^{\prime}$ (resp. $\Gamma$) considered in Example \ref{exam:Gamma torsor}. Then we have 
    \[\Gamma^{\prime} = \{\prod_{i=0}^{d^{\prime}}\delta_i^{n_i}\mid n_0,\dots,n_{d^{\prime}}\in \Zp\text{ such that }n_0+\cdots+n_r = 0\}\]
    and
    \[\Gamma = \{\prod_{i=0}^{d}\delta_i^{n_i}\mid n_0,\dots,n_d\in \Zp\text{ such that }n_0+\cdots+n_r = 0\}.\]
    Note that we have the following commutative diagram
    \[\xymatrix@C=0.5cm{
      X_{\infty,Z}^{\prime}\ar[d]\ar[rrr]&&& X_{\infty,Z}\ar[d]\\
      X_Z^{\prime}\ar[rrr]&&& X_Z
    }\]
    which induces a natural surjection 
    \[\Gamma^{\prime}\to \Gamma\]
    whose kernel $\Gamma_f$ is isomorphic to 
    \[\Gamma_f \cong \oplus_{i=d+1}^{d^{\prime}}\Zp\cdot\delta_i\subset \Gamma^{\prime}.\]

    Now, suppose that $I = [r,s]\subset (\frac{1}{p-1},+\infty)$ and $\lambda\in [1,|\xi|_I^{-1})$ as in Example \ref{exam:Apply to Gamma torsor}. Put $P = W_0+\cdots+W_r$. Then sending each $W_i$ to $\frac{u}{t}\log\frac{[T_i^{\flat}]}{T_i}$ for any $0\leq i\leq d^{\prime}$ induces isomorphisms of sheaves of $\CCtilde^I:=\CCtilde^I_{X_Z^{\prime}}$-algebras on $X_{Z,v}^{\prime}/X_{\infty,Z}^{\prime}$
    \[\CCtilde^I\{\underline W/\lambda\}_{r,d,P}\xrightarrow{\cong} f^*\calO\CCtilde^{I,\lambda}_{\wtx_Z}\text{ and }\CCtilde^I\{\underline W/\lambda\}_{r,d^{\prime},P}\xrightarrow{\cong} \calO\CCtilde^{I,\lambda}_{\wtx^{\prime}_Z}.\]
    Via these isomorphisms, the connection 
    \[\nabla_{\wtx_Z^{\prime}/\wtx_Z}:\calO\CCtilde^{I,\lambda}_{\wtx^{\prime}_Z}\to \calO\CCtilde^{I,\lambda}_{\wtx^{\prime}_Z}\otimes_{\calO_{\wtx_Z^{\prime}}}\Omega^1_{\wtx_Z^{\prime}/\wtx_Z}\{-1\}\]
    is given by 
    \[\nabla_{\wtx_Z^{\prime}/\wtx_Z} = \sum_{i={d+1}}^{d^{\prime}}-\frac{\partial}{\partial W_i}\otimes\frac{u}{t}\dlog T_i.\]
    Here, we trivialize $\Omega^1_{\wtx_Z^{\prime}/\wtx_Z}\{-1\}$ via the isomorphism
    \[\Omega^1_{\wtx_Z^{\prime}/\wtx_Z}\{-1\}\cong \oplus_{i=d+1}^{d^{\prime}}\calO_{\wtx_Z^{\prime}}\cdot\frac{u}{t}\dlog T_i.\]
\end{exam}

\section{Overconvergent de Rham period sheaf}\label{sec:overconvergent dR}
  Throughout this section, fix $Z = \Spa(A,A^+)\in \Perfd$, let $\frakX_Z$ be a semi-stable formal scheme over $A^+$ admitting a lifting $\wtx_Z$ over $\bA_{\inf}(Z)$, and let $X_Z$ be the generic fiber of $\frakX_Z$ with the associated lifting $\wtX_Z$ over $\BBdRp(Z)$ as before. Let $I=[r,s]\subset (0,+\infty)$ be a closed sub-interval. We always assume $1\in I$; then $|\xi|_I = p^{-r}$ (cf.\ Example \ref{exam:norm of u}). 
  
  In what follows, we will construction a period sheaf $\calO\bB_{\dR}^{\dagger,+}$ with connection satisfying the analogous properties as $\calO\CCtilde^{I,\dagger}$, and compare it with the known period sheaves $\calO\bB_{\dR}^{[0,\infty]}$, $\calO\bC^{\dagger}$ and $\calO\bB_{\dR,\pd}^+$ in \cite{DLLZ23}, \cite{Wan23} and \cite{LMNQW} respectively.

\subsection{Revisiting the de Rham period sheaves}\label{ssec:revisit to OBdR}
  For rigid spaces over $C$ with a smooth model over a finite extension of $K$, the period sheaves with connection $(\calO\BBdRp,\nabla)$ and $(\calO\BBdR,\nabla)$ were introduced in the literature, cf.\ \cite{Sch-Pi,DLLZ23}. We first reconstruct the corresponding sheaves for our $X_Z$ from the lifting $\wtx_Z$, following \cite[\S 2]{DLLZ23} and the notation of \S\ref{ssec:period sheaf-small case}.

  For any affinoid perfectoid $U = \Spa(S,S^+)\in X_{Z,v}^{\Box}$ and any $\frakY_Z = \Spf(\calR_Z^{\prime})\in \Sigma_U$ (which is small semi-stable and admits a lifting $\wty_Z = \Spf(\wtr_Z^{\prime})$ associated to $\wtx_Z$), there is a map
  \[\theta_{\frakY_Z}^{\prime}:A_{\frakY_Z}^+=\big(\wtr_Z^{\prime}\widehat \otimes_{\bA_{\inf}(Z)}\bA_{\inf}(U)\big)[\big(\frac{T_0}{t_0}\big)^{\pm 1},\dots,\big(\frac{T_r}{t_r}\big)^{\pm 1}]^{\wedge_{(p,\xi)}}\to S^+,\]
  which extends to a map
  \[\theta_{\frakY_Z}^{\prime}:A_{\frakY_Z}\to S\]
  after inverting $p$. Define
  \[A_{\frakY_Z,\dR}^+:=\varprojlim_{n}A_{\frakY_Z}/\Ker(\theta_{\frakY_Z}^{\prime})^n \text{ and }A_{\frakY_Z,\dR}^{[0,\infty]}:=\varprojlim_{n}A_{\frakY_Z}[\frac{\Ker(\theta_{\frakY_Z}^{\prime})}{\xi}]/(\xi^n).\]
  It follows directly that $A_{\frakY_Z,\dR}^+$ and $A_{\frakY_Z,\dR}^{[0,\infty]}$ are both $\BBdRp(U)$-algebras and there is a natural morphism of $\BBdRp(U)$-algebras
  \[A_{\frakY_Z,\dR}^+\to A_{\frakY_Z,\dR}^{[0,\infty]}.\]
  Moreover, the morphism $\theta_{\frakY_Z}^{\prime}$ and taking reduction modulo $\xi$ respectively induce surjections
  \[A_{\frakY_Z,\dR}^+\to S \text{ and }A_{\frakY_Z}^{[0,\infty]}\to S\]
  fitting into the commutative diagram
  \[\xymatrix@C=0.5cm{
    A_{\frakY_Z,\dR}^+\ar[rrrr]\ar[rrd]&&&& \ar[lld]A_{\frakY_Z,\dR}^{[0,\infty]}\\ &&S.
  }\]
  Note that the usual differential $\rd:\calO_{\wtx_Z}\to\Omega^{1,\log}_{\wtx_Z}$ induces $\BBdRp(U)$-linear differentials
  \[\nabla:A_{\frakY_Z,\dR}^+\to A_{\frakY_Z,\dR}^+\otimes_{\wtr_Z^{\prime}}\Omega^{1,\log}_{\wtr_Z^{\prime}} \text{ and }\nabla:A_{\frakY_Z,\dR}^{[0,\infty]}\to A_{\frakY_Z,\dR}^{[0,\infty]}\otimes_{\wtr_Z^{\prime}}\Omega^{1,\log}_{\wtr_Z^{\prime}}\{-1\}\]
  such that the morphism $A_{\frakY_Z,\dR}^+\to A_{\frakY_Z,\dR}^{[0,\infty]}$ above is compatible with the connections.
  \begin{dfn}\label{dfn:OBdR revisit}
      Keep the notation as above. Denote by $\calO\BBdRp$ and $\calO\bB_{\dR}^{[0,\infty]}$ the sheaves of $\BBdRp$-algebras on $X_{Z,v}$ associated to the presheaves on $X_{Z,v}^{\Box}$
      \[U\in X_{Z,v}^{\Box}\mapsto \colim_{\frakY_Z\in\Sigma_U}A_{\frakY_Z,\dR}^+ \text{ and }U\in X_{Z,v}^{\Box}\mapsto\colim_{\frakY_Z\in\Sigma_U}A_{\frakY_Z}^{[0,\infty]}\]
      respectively. Let $\widehat \calO_{X_Z}$ be the structure sheaf on $X_{Z,v}$, by construction, we get a commutative diagram of morphisms
      \[\xymatrix@C=0.5cm{
        \calO\BBdRp\ar[d]\ar[rrrr]&&&& \calO\bB_{\dR}^{[0,\infty]}\ar[d]\\
        \widehat \calO_{X_Z}\ar@{=}[rrrr]&&&&\widehat \calO_{X_Z}.
      }\]
      If we denote the left vertical morphism by $\theta^{\prime}:\calO\BBdRp\to \widehat \calO_{X_Z}$, then we have isomorphisms
      \[\calO\BBdRp/\Ker(\theta^{\prime})\cong\widehat \calO_{X_Z}\cong \calO\bB_{\dR}^{[0,\infty]}/\xi.\]
      The differential $\rd:\calO_{\wtx_Z}\to\Omega^{1,\log}_{\wtx_Z}$ induces $\BBdRp$-linear connections
      \[\nabla:\calO\BBdRp\to\calO\BBdRp\otimes_{\calO_{\wtx_Z}}\Omega^{1,\log}_{\wtx_Z} \text{ and }\nabla:\calO\bB_{\dR}^{[0,\infty]}\to\calO\bB_{\dR}^{[0,\infty]}\otimes_{\calO_{\wtx_Z}}\Omega^{1,\log}_{\wtx_Z}\{-1\}\]
      such that the morphism
      \[\calO\BBdRp\to \calO\bB_{\dR}^{[0,\infty]}\]
      is compatible with the connections.

      Equip $\calO\BBdRp\to \calO\bB_{\dR}^{[0,\infty]}$ with the $\Ker(\theta^{\prime})$-adic filtration and $\xi$-adic filtration respectively, define a filtration on $\calO\BBdRp[\frac{1}{\xi}]$ by letting 
      \[\Fil^n(\calO\BBdRp[\frac{1}{\xi}]) = \sum_{i+j = n,j\geq 0}\xi^i\Ker(\theta^{\prime})^j\]
      for each $n\in \bZ$. Then the natural morphism $\calO\BBdRp\to \calO\bB_{\dR}^{[0,\infty]}$ factors as
      \[\calO\BBdRp\to \Fil^0(\calO\BBdRp[\frac{1}{\xi}])\to\calO\bB_{\dR}^{[0,\infty]}\]
      Let $\calO\BBdR$ be the completion of $\calO\BBdRp[\frac{1}{\xi}]$ with respect to the above filtration and equip it with the induced filtration. As $\calO\bB_{\dR}^{[0,\infty]}$ is $\xi$-adic complete, we then get natural morphisms of sheaves of $\BBdRp$-algebras on $X_{Z,v}$
      \[\calO\BBdRp\to \Fil^0(\calO\BBdR)\to\calO\bB_{\dR}^{[0,\infty]}.\]
      Denote by $(\calO\bC,\Theta)$ the reduction of $(\calO\bB_{\dR}^{[0,\infty]},\nabla)$ modulo $\xi$. 
  \end{dfn}

  \begin{lem}\label{lem:local description of OBdRp}
      Suppose that $\frakX_Z = \Spf(\calR_Z)$ is small semi-stable and $U = \Spa(S,S^+)\in X_{Z,v}$ is  affinoid perfectoid. Keep notations as in Lemma \ref{lem:finite ideal}. 

      \begin{enumerate}
          \item[(1)] For any $\frakY_Z\in \Sigma_U$, the morphism of $\BBdRp(U)$-algebras
          \[\iota: \BBdRp(U)[[V_0,\dots,V_d]]/(\prod_{i=0}^r(1+V_i)-1)\to A_{\frakY_Z,\dR}^+,\]
          which sending each $V_i$ to $\frac{T_i}{t_i}-1$, is an isomorphism and it carries the ideal $(\xi,V_0,\dots,V_d)$ onto the kernel of $A_{\frakY_Z,\dR}^+\to S$.

          \item[(2)] For any $\frakY_Z\in \Sigma_U$, the morphism of $\BBdRp(U)$-algebras
          \[\iota: \BBdRp(U)[[U_0,\dots,U_d]]/(\sum_{i=0}^rU_i)\to A_{\frakY_Z,\dR}^+,\]
          which sending each $U_i$ to $\log\frac{T_i}{t_i}$, is an isomorphism and it carries the ideal $(\xi,U_0,\dots,U_d)$ onto the kernel of $A_{\frakY_Z,\dR}^+\to S$. Moreover, via the isomorphism
          \[\Omega^{1,\log}_{\wtr_Z}\cong \oplus_{i=0}^d\wtr_Z\cdot e_i/\wtr_Z\cdot \sum_{i=0}^re_i,\]
          the connection $\nabla$ on $A_{\frakY_Z,\dR}^+$ is given by 
          \[\sum_{i=0}^d\frac{\partial}{\partial U_i}\otimes e_i:\BBdRp(U)[[U_0,\dots,U_d]]/(\sum_{i=0}^rU_i)\to \BBdRp(U)[[U_0,\dots,U_d]]/(\sum_{i=0}^rU_i)\otimes\big(\oplus_{i=0}^d\bZ\cdot e_i/\bZ\cdot\sum_{i=0}^re_i\big)\]
          where $e_i$ stands for the log-differential $\dlog T_i\in \Omega^{1,\log}_{\wtr_Z}$.
      \end{enumerate}
  \end{lem}
  \begin{proof}
      For item (1): This follows from a similar argument in the proof of Proposition \ref{prop:local description}. 

      For item (2): The assertion follows from Item (1) together with a similar argument for the proof of Lemma \ref{lem:change of variable} 
  \end{proof}

  \begin{lem}\label{lem:local description of OBdRp-II}
      Suppose that $\frakX_Z = \Spf(\calR_Z)$ is small semi-stable and $U = \Spa(S,S^+)\in X_{Z,v}$ is  affinoid perfectoid. Keep notations as in Lemma \ref{lem:finite ideal}. 

      \begin{enumerate}
          \item[(1)] For any $\frakY_Z\in \Sigma_U$, the morphism of $\BBdRp(U)$-algebras
          \[\iota: \BBdRp(U)\za V_0,\dots,V_d\ya/(\xi^{-1}(\prod_{i=0}^r(1+\xi V_i)-1))\to A_{\frakY_Z,\dR}^{[0,\infty]},\]
          which sending each $V_i$ to $\frac{\frac{T_i}{t_i}-1}{\xi}$, is an isomorphism. Here, $\BBdRp(U)\za U_0,\dots,U_d\ya$ denotes the $\xi$-adic completion of the polynomial ring $\BBdRp(U)[U_0,\dots,U_d]$.

          \item[(2)] For any $\frakY_Z\in \Sigma_U$, the morphism of $\BBdRp(U)$-algebras
          \[\iota: \BBdRp(U)
          \za W_0,\dots,W_d\ya/(\sum_{i=0}^rW_i)\to A_{\frakY_Z,\dR}^{[0,\infty]},\]
          which sending each $W_i$ to $\frac{\log\frac{T_i}{t_i}}{\xi}$, is an isomorphism. Moreover, via the isomorphism
          \[\Omega^{1,\log}_{\wtr_Z}\{-1\}\cong \oplus_{i=0}^d\wtr_Z\cdot e_i/\wtr_Z\cdot \sum_{i=0}^re_i,\]
          the connection $\nabla$ on $A_{\frakY_Z,\dR}^{[0,\infty]}$ is given by 
          \[\sum_{i=0}^d\frac{\partial}{\partial W_i}\otimes e_i:\BBdRp(U)\za W_0,\dots,W_d\ya/(\sum_{i=0}^rW_i)\to \BBdRp(U)\za W_0,\dots,W_d\ya/(\sum_{i=0}^rW_i)\otimes\big(\oplus_{i=0}^d\bZ\cdot e_i/\bZ\cdot\sum_{i=0}^re_i\big)\]
          where $e_i$ stands for $\frac{\dlog T_i}{\xi}\in \Omega^{1,\log}_{\wtr_Z}\{-1\}$.
      \end{enumerate}
  \end{lem}
  \begin{proof}
      For item (1): This follows from a similar argument in the proof of Proposition \ref{prop:local description}. 

      For item (2): The assertion follows from Item (1) together with a similar argument for the proof of Lemma \ref{lem:change of variable} immediately.
  \end{proof}
  \begin{prop}\label{prop:local description of OBdRp}
      Suppose that $\frakX_Z = \Spf(\calR_Z)$ is small semi-stable and $U = \Spa(S,S^+)\in X_{Z,v}$ is  affinoid perfectoid. Keep notations in Lemma \ref{lem:finite ideal}. 

      \begin{enumerate}
          \item[(1)] The morphism of sheaves of $\BBdRp$-algebras on $X_{Z,v}/U$
          \[\BBdRp[[U_0,\dots,U_d]]/(\sum_{i=0}^rU_i)\to \calO\BBdRp\]
          sending each $U_i$ to $\log\frac{T_i}{t_i}$ for any $0\leq i\leq d$ is an isomorphism which carries the ideal $(\xi,U_0,\dots,U_d)$ onto the ideal $\Ker(\theta^{\prime}:\calO\BBdRp\to \widehat \calO_{X_Z})$. Moreover, via the isomorphism 
          \[\Omega^{1,\log}_{\wtx_Z}\cong \oplus_{i=0}^d\calO_{\wtx_Z}\cdot e_i/\calO_{\wtx_Z}\cdot\sum_{i=0}^re_i,\]
          the connection $\nabla$ on $\calO\BBdRp$ is given by $\sum_{i=0}^d\frac{\partial}{\partial U_i}\otimes e_i$, where $e_i$ stands for $\dlog T_i$ for any $0\leq i\leq d$.

          \item[(2)] The morphism of sheaves of $\BBdRp$-algebras on $X_{Z,v}/U$
          \[\BBdRp\za W_0,\dots,W_d\ya/(\sum_{i=0}^rW_i)\to \calO\bB_{\dR}^{[0,\infty]}\]
          sending each $W_i$ to $\frac{\log\frac{T_i}{t_i}}{\xi}$ for any $0\leq i\leq d$ is an isomorphism. Moreover, via the isomorphism 
          \[\Omega^{1,\log}_{\wtx_Z}\{-1\}\cong \oplus_{i=0}^d\calO_{\wtx_Z}\cdot e_i/\calO_{\wtx_Z}\cdot\sum_{i=0}^re_i,\]
          the connection $\nabla$ on $\calO\bB_{\dR}^{[0,\infty]}$ is given by $\sum_{i=0}^d\frac{\partial}{\partial W_i}\otimes e_i$, where $e_i$ stands for $\frac{\dlog T_i}{\xi}$ for any $0\leq i\leq d$.
      \end{enumerate}
  \end{prop}
  \begin{proof}
      The assertion follows from Lemmas \ref{lem:local description of OBdRp} and \ref{lem:local description of OBdRp-II} directly.
  \end{proof}
  \begin{prop}\label{prop:Poincare's Lemma for OBdR}
      Let $\frakX_Z$ be a semi-stable formal scheme over $A^+$ of relative dimension $d$ with a fixed lifting $\wtx_Z$ over $\bA_{\inf}(Z)$. Let $X_Z$ be the generic fiber of $\frakX_Z$.
      \begin{enumerate}
          \item[(1)] The natural morphism of sheaves of $\BBdRp$-algebras
          \[\Fil^0\calO\BBdR\to \calO\bB_{\dR}^{[0,\infty]}\]
          is an isomorphism.

          \item[(2)] The following complex of sheaves of $\BBdRp$-modules on $X_{Z,v}$ is exact
          \[0\to \BBdRp\to \calO\bB_{\dR}^{[0,\infty]}\to \calO\bB_{\dR}^{[0,\infty]}\otimes_{\calO_{\wtx_Z}}\Omega^{1,\log}_{\wtx_Z}\{-1\}\to\cdots\to \calO\bB_{\dR}^{[0,\infty]}\otimes_{\calO_{\wtx_Z}}\Omega^{d,\log}_{\wtx_Z}\{-d\}.\]
      \end{enumerate}
  \end{prop}
  \begin{proof}
      Since the problems in items (1) and (2) are both local on $\frakX_{Z,\et}$ and $X_{Z,v}$, we may assume $\frakX_Z = \Spf(\calR_Z)$ is small semi-stable and check everything on $X_{Z,v}/U$ for $U$ affinoid perfectoid.

      For item (1): By Proposition \ref{prop:local description of OBdRp} (working on the chart with coordinates $U_1,\dots,U_d$ and $W_1,\dots,W_d$ after eliminating $U_0$ and $W_0$), we have the following isomorphisms of sheaves of $\BBdRp$-algebras on $X_{Z,v}/U$
      \[i^+:\BBdRp[[U_1,\dots,U_d]]\to \calO\BBdRp\]
      and 
      \[i^{[0,\infty]}:\BBdRp\za W_1,\dots,W_d\ya\to \calO\bB_{\dR}^{[0,\infty]},\]
      which send $U_i$ to $\log\frac{T_i}{t_i}$ and $W_i$ to $\frac{\log\frac{T_i}{t_i}}{\xi}$, respectively, for any $1\leq i\leq d$. Via the above two isomorphisms, the natural map $\calO\BBdRp\to \calO\bB_{\dR}^{[0,\infty]}$ is the morphism
      \[\BBdRp[[U_1,\dots,U_d]]\to \BBdRp\za W_1,\dots,W_d\ya\]
      sending each $U_i$ to $\xi W_i$ for any $1\leq i\leq d$.

      Note that the isomorphism $i^+$ carries the ideal of $(\xi, U_1,\dots,U_d)$ to $\Ker(\theta^{\prime})$. By definition of $\Fil^0\calO\BBdR$, it is easy to see that $i^+$ induces an isomorphism of sheaves of $\BBdRp$-algebras
      \[\BBdRp[[U_1,\dots,U_d]][\frac{U_1}{\xi},\dots,\frac{U_d}{\xi}]^{\wedge_{\xi}}\to \Fil^0\calO\BBdR,\]
      where the left hand side of the above isomorphism denotes the $\xi$-adic completion of the blow-up algebra 
      \[\BBdRp[[U_1,\dots,U_d]][\frac{U_1}{\xi},\dots,\frac{U_d}{\xi}].\]
      As a consequence, the morphism $\Fil^0\calO\BBdR\to \calO\bB_{\dR}^{[0,\infty]}$ is given by the morphism
      \[\BBdRp[[U_1,\dots,U_d]][\frac{U_1}{\xi},\dots,\frac{U_d}{\xi}]^{\wedge_{\xi}}\to \BBdRp\za W_1,\dots,W_d\ya\]
      which sends each $\frac{U_i}{\xi}$ to $W_i$ for any $1\leq i\leq d$. One can conclude by noting that this map is an isomorphism.

      For item (2): Note that via the isomorphism
      \[\Omega^{1,\log}_{\wtx_Z}\{-1\}\cong \oplus_{i=1}^d\calO_{\wtx_Z}\frac{\dlog T_i}{\xi}\]
      and the isomorphism $i^{[0,\infty]}$, the connection $\nabla$ on $\calO\bB_{\dR}^{[0,\infty]}$ is given by
      \[\sum_{i=1}^d\frac{\partial}{\partial W_i}\otimes\frac{\dlog T_i}{\xi}:\BBdRp\za W_1,\dots,W_d\ya\to \oplus_{i=1}^d\BBdRp\za W_1,\dots,W_d\ya\cdot\frac{\dlog T_i}{\xi}.\]
      So one can conclude by \cite[Cor. 2.3.17]{DLLZ23}.
  \end{proof}

  \begin{cor}\label{cor:description of OC}
      Keep notations as in Proposition \ref{prop:Poincare's Lemma for OBdR}.
      \begin{enumerate}
          \item[(1)] The following complex of sheaves of $\widehat \calO_{X_Z}$-modules on $X_{Z,v}$ is exact
          \[0\to \widehat \calO_{X_Z}\to \calO\bC\xrightarrow{\Theta}\calO\bC\otimes_{\calO_{\frakX_Z}}\Omega^{1,\log}_{\frakX_Z}\{-1\}\to \cdots\to \calO\bC\otimes_{\calO_{\frakX_Z}}\Omega^{d,\log}_{\frakX_Z}\{-d\}\to 0.\]

          \item[(2)] Suppose that $\frakX_Z = \Spf(\calR_Z)$ is small semi-stable and keep notations as in Proposition \ref{prop:local description of OBdRp}. For any affinoid perfectoid $U\in X_{Z,v}$, the morphism of sheaves of $\widehat \calO_{X_Z}$-algebras on $X_{Z,v}/U$
          \[\widehat \calO_{X_Z}[W_0,\dots,W_d]/(W_0+\cdots+W_d)\to \calO\bC\]
          sending each $W_i$ to $\frac{\log \frac{T_i}{t_i}}{\xi}$ for any $0\leq i\leq d$ is an isomorphism such that via the isomorphism, 
          \[\Omega^{1,\log}_{\wtx_Z}\{-1\}\cong \oplus_{i=0}^d\calO_{\wtx_Z}\cdot e_i/\calO_{\wtx_Z}\cdot\sum_{i=0}^re_i,\]
          the Higgs field $\Theta$ on $\calO\bC$ is given by $\sum_{i=0}^d\frac{\partial}{\partial W_i}\otimes e_i$, where $e_i$ stands for $\frac{\dlog T_i}{\xi}$ for any $0\leq i\leq d$.
      \end{enumerate}
  \end{cor}
  \begin{proof}
      The assertion follows from Propositions \ref{prop:local description of OBdRp} and \ref{prop:Poincare's Lemma for OBdR} by taking reduction of $(\calO\bB_{\dR}^{[0,\infty]},\nabla)$ modulo $\xi$.
  \end{proof}

  At the end of this subsection, we describe the integral Faltings' extension appeared in \cite{Wan23,SW24}.

  \begin{dfn}\label{dfn:integral Faltings extension}
      Fix a $Z = \Spa(A,A^+)\in \Perfd$, and let $\frakX_Z$ be a semi-stable formal scheme over $A^+$ that admits a lifting $\wtx_Z$ over $\bA_{\inf}(Z)$ and $X_Z$ be the generic fiber of $\frakX_Z$ with the associated lifting $\wtX_Z$ over $\BBdRp(Z)$ as before. Denote by $\calO\bA_{\inf}$ the $(p,\xi)$-adic completion of the exactification with respect to the natural surjection
      \[\calO_{\wtx_Z}\widehat \otimes_{\bA_{\inf}(Z)}\bA_{\inf}\to \widehat \calO^+_{X_Z}\]
      and let 
      \[\theta^{\prime}:\calO\bA_{\inf}\to \widehat \calO^+_{X_Z}\]
      be the induced surjection. Note that $\calO\bA_{\inf}$ is exactly the sheaf associated to the presheaf on $X_{Z,v}^{\Box}$ sending each affinoid perfectoid $U = \Spa(S,S^+)$ to $\colim_{\frakY_Z\in \Sigma_U}A_{\frakY_Z}^+$ and the morphism $\theta^{\prime}$ is induced by the surjections $\theta^{\prime}_{\frakY_Z}:A_{\frakY_Z}^+\to S^+$. So we have natural morphisms
      \[\calO\bA_{\inf}[\frac{\Ker(\theta^{\prime})}{\xi}]\to \calO\bB_{\dR}^{[0,\infty]}\xrightarrow{\mod \xi}\calO\bC.\]
      Denote by $\calE^+\subset \calO\bC$ the sub-sheaf of $\widehat \calO_{X_Z}^+$-modules generated by the image of $\frac{\Ker(\theta^{\prime})}{\xi}$ via the composite of the above morphisms. As $\xi\in \Ker(\theta^{\prime})$, $\calE^+$ contains $\widehat \calO_{X_Z}^+$ as a sub-module.
  \end{dfn}
  \begin{lem}\label{lem:local description of integral Faltings extension}
      Suppose that $\frakX_Z$ is small semi-stable and $U = \Spa(S,S^+)\in X_{Z,v}$ is affinoid perfectoid, and Keep notations in Corollary \ref{cor:description of OC}. Then for any $\frakY_Z = \Spf(\calR_Z^{\prime})\in \Sigma_U$, the $S^+$-module generated by the image of $\frac{\Ker(\theta_{\frakY_Z}^{\prime})}{\xi}$ via the composite
      \[A_{\frakY_Z}^+[\frac{\Ker(\theta_{\frakY_Z}^{\prime})}{\xi}]\to A_{\frakY_Z,\dR}^{[0,\infty]} \to S[W_0,\dots,W_d]/(W_0+\cdots+W_d) \]
      coincides with the sub-$S^+$-module of $S[W_0,\dots,W_d]/(W_0+\cdots+W_d)$ generated by $1,W_0,\dots,W_d$ and thus is isomorphic to the free $S^+$-module of rank $d+1$
      \[S^+\cdot 1\oplus(\oplus_{i=0}^dS^+\cdot W_i/S^+\cdot\sum_{i=0}^rW_i).\]
  \end{lem}
  \begin{proof}
      Denote by $E^+$ the sub-$S^+$-module of $S[W_0,\dots,W_d]/(W_0+\cdots+W_d)$ generated by the image of $\frac{\Ker(\theta_{\frakY_Z}^{\prime})}{\xi}$.
      Note that \ref{lem:local description of OBdRp-II}, for any $0\leq i\leq d$, we have
      \[\frac{\frac{T_i}{t_i}-1}{\xi}\equiv \frac{\log\frac{T_i}{t_i}}{\xi}\mod \xi A_{\frakY_Z,\dR}^{[0,\infty]}.\]
      So it follows from Lemma \ref{lem:finite ideal} that $E^+$ is the $S^+$-module generated by the image of $1,W_0,\dots,W_d$ as desired.
  \end{proof}
  \begin{prop}\label{prop:Integral Faltings extension}
      The restriction of $\Theta$ on $\calO\bC$ to $\calE^+$ takes values in $\widehat \calO_{X_Z}^+\otimes_{\calO_{\frakX_Z}}\Omega^{1,\log}_{\frakX_Z}\{-1\}$ and induces a short exact sequence of sheaves of $\widehat \calO_{X_Z}^+$-modules on $X_{Z,v}$
      \[0\to \widehat \calO_{X_Z}^+\to \calE^+\to \widehat \calO^+_{X_Z}\otimes_{\calO_{\frakX_Z}}\Omega^{1,\log}_{\frakX_Z}\{-1\}\to 0.\]
  \end{prop}
  \begin{proof}
      Since the problem is local on $\frakX_{Z,\et}$ and $X_{Z,v}$, we may assume $\frakX_Z$ is small semi-stable and work on $X_{Z,v}/U$ for any fixed affinoid perfectoid $U = \Spa(S,S^+)\in X_{Z,v}$. By Corollary \ref{cor:description of OC}, the Higgs field $\Theta$ is given by $\sum_{i=0}^d\frac{\partial}{\partial W_i}\otimes\frac{\dlog T_i}{\xi}$. So the result follows from Lemma \ref{lem:local description of integral Faltings extension} immediately.
  \end{proof}

\subsection{The overconvergent de Rham period sheaf $\calO\bB_{\dR}^{\dagger,+}$}
  Recall that $I = [r,s]\subset (0,+\infty)$ contains $1$. One may ask whether there exists a morphism from $\calO\CCtilde^{I,\dagger}$ to $\calO\bB_{\dR}^{[0,\infty]}$ compatible with the connections and with $\CCtilde^I\to\BBdRp$. This fails in general; instead we introduce a new period sheaf with connection $(\calO\bB_{\dR}^{\dagger,+},\nabla)$---the natural target for the Riemann--Hilbert correspondence in this paper---for which there is a morphism $\calO\CCtilde^{I,\dagger}\to\calO\bB_{\dR}^{\dagger,+}$ compatible with the connections and with $\CCtilde^I\to\BBdRp$.

  Recall that for any affinoid perfectoid $U = \Spa(S,S^+)\in X_{Z,v}$ and any $\frakY_Z = \Spf(\wtr_Z^{\prime})$, we have introduced the $\CCtilde^I(U)$-algebra $A_{\frakY_Z}^{I,\lambda}$ (cf. the paragraph above Definition \ref{dfn:local OC^I}) for any $\lambda\in[1,|\xi|_I^{-1})$. By its construction, there exist a natural map
  \[A_{\frakY_Z}[\frac{\Ker(\theta_{\frakY_Z})}{\xi}]\to A_{\frakY_Z}^{I,\lambda}.\]
  Denote by $A_{\frakY_Z,\dR}^{I,\lambda}$ the $\xi$-adic completion of $A_{\frakY_Z}^{I,\lambda}$ so we get a morphism of $\BBdRp(U)$-algebras
  \[A_{\frakY_Z,\dR}^{[0,\infty]}\to A_{\frakY_Z,\dR}^{I,\lambda}.\]

  \begin{dfn}\label{dfn:Overconvergent de Rham period sheaf}
      Keep the notation as above. Denote by $\calO\bB_{\dR}^{I,\lambda,+}$ the sheaf on $X_{Z,v}$ associated to the sheafification of the presheaf 
      \[U = \Spa(S,S^+)\in X_{Z,v}^{\Box}\mapsto \colim_{\frakY_Z\in \Sigma_U}A_{\frakY_Z,\dR}^{I,\lambda}.\]
      Define
      \[\calO\bB_{\dR}^{I,\dagger,+}:=(\colim_{\lambda\to 1^+}\calO\bB_{\dR}^{I,\lambda,+})^{\wedge_{\xi}}.\]
      Then the usual $(p,\xi)$-adic continuous differential $\rd:\calO_{\wtx_Z}\to \Omega^{1,\log}_{\wtx_Z}$ induces connections
      \[\nabla:\calO\bB_{\dR}^{I,\lambda,+}\to\calO\bB_{\dR}^{I,\lambda,+}\otimes_{\calO_{\wtx_Z}}\Omega^{1,\log}_{\wtx_Z}\{-1\}\]
      and 
      \[\nabla:\calO\bB_{\dR}^{I,\dagger,+}\to\calO\bB_{\dR}^{I,\dagger,+}\otimes_{\calO_{\wtx_Z}}\Omega^{1,\log}_{\wtx_Z}\{-1\}.\]
      Moreover, we get a morphism of period sheaves with connections
      \[(\calO\bB_{\dR}^{[0,\infty]},\nabla)\to(\calO\bB_{\dR}^{I,\lambda,+},\nabla).\]
      For any closed sub-interval $J\subset I$ with $1\in J$ and any $\lambda\in [1,|\xi|_I^{-1})\subset[1,|\xi|_J^{-1})$, there is a natural morphism of sheaves with connection
      \[(\calO\bB_{\dR}^{I,\lambda,+},\nabla)\to (\calO\bB_{\dR}^{J,\lambda,+},\nabla)\]
      and thus a natural morphism
      \[(\calO\bB_{\dR}^{I,\dagger,+},\nabla)\to (\calO\bB_{\dR}^{J,\dagger,+},\nabla).\]
      We will see both of them are isomorphisms and thus just denote $(\calO\bB_{\dR}^{I,\lambda,+},\nabla)$ and $(\calO\bB_{\dR}^{I,\dagger,+},\nabla)$ by $(\calO\bB_{\dR}^{\lambda,+},\nabla)$ and $(\calO\bB_{\dR}^{\dagger,+},\nabla)$, respectively, when context is clear. Denote by $(\calO\widehat \bC_{\lambda},\Theta)$ and $(\calO\bC^{\dagger},\Theta)$ the reduction of $(\calO\bB_{\dR}^{\lambda,+},\nabla)$ and $(\calO\bB_{\dR}^{\dagger,+},\nabla)$ modulo $\xi$ respectively.
  \end{dfn}
  \begin{lem}\label{lem:independent of I}
      Suppose that $\frakX_Z$ is small semi-stable. Let $J\subset I$ be a closed sub-interval $J$ containing $1$ and $\lambda\in [1,|\xi|_I^{-1})$. For any affinoid perfectoid $U = \Spa(S,S^+)\in X_{Z,v}$ and any $\frakY_Z\in \Sigma_U$, the natural morphism
      \[A_{\frakY_Z,\dR}^{I,\lambda}\to A_{\frakY_Z,\dR}^{J,\lambda}\]
      is an isomorphism.
  \end{lem}
  \begin{proof}
      Fix an $n\geq 1$.
      Denote by $|\cdot|_I$ and $|\cdot|_J$ the induced norm on 
      \[\BBdRp(U)/\xi^n\cong \CCtilde^I(U)/\xi^n\cong \CCtilde^J(U)/\xi^n\]
      and then by the proof of Lemma \ref{lem:compare with BdR}, we have
      \[|\cdot|_J\leq|\cdot|_I\leq p^{n-1}|\cdot|_J.\]
      This forces that
      \[\CCtilde^I(U)\{V_0/\lambda,\dots,V_d/\lambda\}/\xi^n\cong \CCtilde^J(U)\{V_0/\lambda,\dots,V_d/\lambda\}/\xi^n\]
      as $p$-adic Banach algebras. Thus, by Proposition \ref{prop:local description}, we have
      \[A_{\frakY_Z}^{I,\lambda}/\xi^n\cong A_{\frakY_Z}^{J,\lambda}/\xi^n\]
      yielding the desired isomorphism
      \[A_{\frakY_Z,\dR}^{I,\lambda}\to A_{\frakY_Z,\dR}^{J,\lambda}\]
      by letting $n\to +\infty$.
  \end{proof}
  \begin{prop}\label{prop:independent of I}
      Let $J\subset I$ be a closed sub-interval $J$ containing $1$ and $\lambda\in [1,|\xi|_I^{-1})$. The morphisms of period sheaves with connection on $X_{Z,v}$
      \[(\calO\bB_{\dR}^{I,\lambda,+},\nabla)\to (\calO\bB_{\dR}^{J,\lambda,+},\nabla)\]
      and
      \[(\calO\bB_{\dR}^{I,\dagger,+},\nabla)\to (\calO\bB_{\dR}^{J,\dagger,+},\nabla)\]
      are both isomorphisms. Moreover, we have natural morphisms of period sheaves with connection
      \[(\calO\CCtilde^{I,\lambda},\nabla)\to(\calO\bB_{\dR}^{I,\lambda,+},\nabla)\]
      and
      \[(\calO\CCtilde^{I,\dagger},\nabla)\to(\calO\bB_{\dR}^{I,\dagger,+},\nabla)\]
      such that $\calO\bB_{\dR}^{I,\lambda,+}$ and $\calO\bB_{\dR}^{I,\dagger,+}$ are the $\xi$-adic completions of $\calO\CCtilde^{I,\lambda}$ and $\calO\CCtilde^{I,\dagger}$ respectively.
  \end{prop}
  \begin{proof}
      By Lemma \ref{lem:independent of I} (and Proposition \ref{prop:local description} for the ``moreover'' part), the result follows by unwinding definitions.
  \end{proof}
  Thanks to Proposition \ref{prop:independent of I}, from now on we always assume $\lambda\in [1,p^{\frac{1}{p-1}})\cap p^{\bQ}$ and regard $\calO\bB_{\dR}^{I,\lambda,+}$ and $\calO\bB_{\dR}^{I,\dagger,+}$ as the $\xi$-adic completions of $\calO\CCtilde^{1,\lambda}$ and $\calO\CCtilde^{1,\dagger}$ respectively. 

  \begin{lem}\label{lem:local description of OBdRp^lambda}
      Suppose that $\frakX_Z = \Spf(\calR_Z)$ is small semi-stable and let $\lambda = p^{\rho}$ for some $\rho\in\bQ\cap[0,\frac{1}{p-1})$. Let $U = \Spa(S,S^+)\in X_{Z,v}$ be affinoid perfectoid and keep notations in Proposition \ref{prop:change of variable}. Then the map of sheaves of $\BBdRp$-algebras on $X_{Z,v}$
      \[\iota:\BBdRp\za[\varpi^{\rho}]W_0,\dots,[\varpi^{\rho}]W_d\ya/(W_0+\cdots+W_d)\to \calO\bB_{\dR}^{\lambda,+},\]
      which sends each $W_i$ to $\frac{\log\frac{T_i}{t_i}}{\xi}$, is an isomorphism, where $\BBdRp\za[\varpi^{\rho}]W_0,\dots,[\varpi^{\rho}]W_d\ya$ denotes the $\xi$-adic completion of 
      \[\CCtilde^{1,+}\za [\varpi^{\rho}]W_0,\dots,[\varpi^{\rho}]W_d\ya[\frac{1}{p}]\]
      and $\CCtilde^{1,+}\za [\varpi^{\rho}]W_0,\dots,[\varpi^{\rho}]W_d\ya$ denotes the $p$-adic completion of $\CCtilde^{1,+}[ [\varpi^{\rho}]W_0,\dots,[\varpi^{\rho}]W_d]$, and $(W_0+\dots+W_r)$ denotes the closed ideal generated by $W_0+\dots+W_r$. Moreover, via the isomorphism
      \[\Omega^{1,\log}_{\wtx_Z}\{-1\}\cong (\oplus_{i=0}^d\calO_{\wtx_Z}\cdot\frac{\dlog T_i}{\xi})/\calO_{\wtx_Z}\cdot\sum_{i=0}^r\frac{\dlog T_i}{\xi}\]
      the connection $\nabla$ on $\calO\bB_{\dR}^{\lambda,+}$ is given by 
      \[\nabla = \sum_{i=0}^d\frac{\partial}{\partial W_i}\otimes\frac{\dlog T_i}{\xi}.\]
  \end{lem}
  \begin{proof}
      As $|[\varpi^{\rho}]|_1 = p^{-\rho} = \lambda^{-1}$, we see that
      \[\CCtilde^1\{W_0/\lambda,\dots,W_d/\lambda\} = \CCtilde^{1,+}\za [\varpi^{\rho}]W_0,\dots,[\varpi^{\rho}]W_d\ya[\frac{1}{p}].\]
      Thus, the result follows from Propositions \ref{prop:change of variable} and \ref{prop:independent of I} immediately.
  \end{proof}
  \begin{cor}\label{cor:local description of OBdRp^lambda}
      Suppose that $\frakX_Z = \Spf(\calR_Z)$ is small semi-stable and let $\lambda = p^r$ for some $\rho\in\bQ\cap[0,\frac{1}{p-1})$ and $U = \Spa(S,S^+)\in X_{Z,v}$ is affinoid perfectoid as in Lemma \ref{lem:local description of OBdRp^lambda}. Put $P=W_0+\cdots+W_r$ and keep notations in Definition \ref{dfn:overconvergent algebras}.
      \begin{enumerate}
          \item[(1)] The map of sheaves of $\widehat \calO_{X_Z}$-algebras on $X_{Z,v}/U$
          \[\widehat \calO_{X_Z}\{\underline W/\lambda\}_{r,d,P}\to\calO\widehat \bC_{\lambda}\]
          which sends each $W_i$ to $\frac{\log\frac{T_i}{t_i}}{\xi}$ for any $0\leq i\leq d$, is an isomorphism. 

          \item[(2)] The map of sheaves of $\widehat \calO_{X_Z}$-algebras on $X_{Z,v}/U$ 
          \[\widehat \calO_{X_Z}\{\underline W\}_{r,d,P}^{\dagger}\to\calO\bC^{\dagger}\]
          which sends each $W_i$ to $\frac{\log\frac{T_i}{t_i}}{\xi}$ for any $0\leq i\leq d$, is an isomorphism. 
          
          \item[(3)] Via the isomorphism
          \[\Omega^{1,\log}_{\wtx_Z}\{-1\}\cong (\oplus_{i=0}^d\calO_{\wtx_Z}\cdot\frac{\dlog T_i}{\xi})/\calO_{\wtx_Z}\cdot\sum_{i=0}^r\frac{\dlog T_i}{\xi}\]
          the Higgs field $\Theta$ on $\calO\widehat \bC_{\lambda}$ (resp. $\calO\bC^{\dagger}$) is given by $\Theta = \sum_{i=0}^d\frac{\partial}{\partial W_i}\otimes\frac{\dlog T_i}{\xi}$.
      \end{enumerate}
  \end{cor}
  \begin{proof}
      Keep notations as in Lemma \ref{lem:local description of OBdRp^lambda}. Note that 
      \[\BBdRp\za[\varpi^{\rho}]W_0,\dots,[\varpi^{\rho}]W_d\ya/\xi = \widehat \calO_{X_Z}\za p^{\rho}W_0,\dots,p^{\rho} W_d\ya:=\widehat \calO_{X_Z}^+\za p^{\rho}W_0,\dots,p^{\rho} W_d\ya[\frac{1}{p}]\]
      where $\widehat \calO_{X_Z}^+\za p^{\rho}W_0,\dots,p^{\rho} W_d\ya$ denotes the $p$-adic completion of $\widehat \calO_{X_Z}^+[p^{\rho}W_0,\dots,p^{\rho} W_d]$. So we see that
      \[\widehat \calO_{X_Z}\za p^{\rho}W_0,\dots,p^{\rho} W_d\ya = \widehat \calO_{X_Z}\{W_0/\lambda,\dots,W_d/\lambda\}.\]
      So everything follows from Lemma \ref{lem:local description of OBdRp^lambda} directly.
  \end{proof}

  \begin{prop}\label{prop:Poincare's Lemma for overconvergent dR}
    Suppose that $\frakX_Z$ is semi-stable of relative dimension $d$ with a fixed lifting $\wtx_Z$ over $A^+$ as before.
    \begin{enumerate}
        \item[(1)] The Higgs complex associated to $(\calO\bC^{\dagger},\Theta)$
        \[0\to \widehat \calO_{X_Z}\to\calO\bC^{\dagger}\xrightarrow{\Theta}\calO\bC^{\dagger}\otimes_{\calO_{\frakX_Z}}\Omega^{1,\log}_{\frakX_Z}\{-1\}\to\cdots\to\calO\bC^{\dagger}\otimes_{\calO_{\frakX_Z}}\Omega^{d,\log}_{\frakX_Z}\{-d\}\to 0\]
        is exact.

        \item[(2)] The de Rham complex associated to $(\calO\bB_{\dR}^{\dagger,+},\nabla)$
        \[0\to \BBdRp\to\calO\bB_{\dR}^{\dagger,+}\xrightarrow{\nabla}\calO\bB_{\dR}^{\dagger,+}\otimes_{\calO_{\wtx_Z}}\Omega^{1,\log}_{\wtx_Z}\{-1\}\to\cdots\to\calO\bB_{\dR}^{\dagger,+}\otimes_{\calO_{\wtx_Z}}\Omega^{d,\log}_{\wtx_Z}\{-d\}\to 0\]
        is exact.
    \end{enumerate}
  \end{prop}
  \begin{proof}
      It suffices to prove Item (1) as Item (2) follows by the derived Nakayama's Lemma. For Item (1), as the problem is local on $\frakX_{Z,\et}$ and $X_{Z,v}$, we may assume that $\frakX_Z$ is small semi-stable and check the desired exactness on $X_{Z,v}/U$ for any affinoid perfectoid $U\in X_{Z,v}$. By Corollary \ref{cor:local description of OBdRp^lambda}, it is easy to see that the morphism of sheaves of $\widehat \calO_{X_Z}$-algebras on $X_{Z,v}/U$
      \[\widehat \calO_{X_Z}\{\underline W\}^{\dagger}\to\calO\bC^{\dagger},\]
      which sends each $W_i$ to $\frac{\log\frac{T_i}{t_i}}{\xi}$ for any $1\leq i\leq d$, is an isomorphism. Moreover, via the isomorphism
      \[\Omega^{1,\log}_{\wtx_Z}\{-1\}\cong \oplus_{i=1}^d\calO_{\wtx_Z}\cdot\frac{\dlog T_i}{\xi}\]
      the Higgs field $\Theta$ on $\calO \bC^{\dagger}$ is given by $\Theta = \sum_{i=1}^d\frac{\partial}{\partial W_i}\otimes\frac{\dlog T_i}{\xi}$. Then one can conclude by the same proof of \cite[Cor. 2.23]{Wan23}.
  \end{proof}

  At the end of this subsection, similar to \S\ref{ssec:relative OC}, we will give a relative version of the above result. So in what follows, suppose we are in the setting of \S\ref{ssec:relative OC} and keep notations therein. Suppose that $\frakY_Z = \Spf(\calR_Z)\in\frakX_{Z,\et}$ and $\frakY_Z^{\prime} = \Spf(\calR_Z^{\prime})\in\frakX_{Z,\et}$ with corresponding lifting $\wty_Z = \Spf(\wtr_Z)$ and $\wty_Z^{\prime} = \Spf(\wtr_Z^{\prime})$ fit into the commutative diagrams \eqref{equ:relative chart-I} and \eqref{equ:relative chart-II}, respectively. For any affinoid perfectoid $U\in Y^{\prime}_{Z,v}$, consider the morphism
  \[A_{\frakY_Z}^{I,\lambda}(U)\to A_{\frakY_Z^{\prime}}^{I,\lambda}(U)\]
  for $I = \{1\}$ and $\lambda = p^{\rho}$ with $\rho\in \bQ\cap[0,\frac{1}{p-1})$.

  Taking $\xi$-adic completion, we get a morphism of $\BBdRp(U)$-algebras
  \[A_{\frakY_Z,\dR}^{1,\lambda}(U)\to A_{\frakY_Z^{\prime},\dR}^{1,\lambda}(U).\]
  By letting $U$ and then $\frakY_Z$ and $\frakY_Z^{\prime}$ vary, we get a morphism of sheaves of $\BBdRp$-algebras on $X_{Z,v}^{\prime}$
  \[f^*\calO\bB_{\dR,\wtx_Z}^{1,\lambda,+}\to\calO\bB_{\dR,\wtx^{\prime}_Z}^{1,\lambda,+},\]
  where $\calO\bB_{\dR,\wtx_Z}^{1,\lambda,+}$ (resp. $\calO\bB_{\dR,\wtx^{\prime}_Z}^{1,\lambda,+}$) denotes the overconvergent de Rham period sheaf $\calO\bB_{\dR}^{1,\lambda,+}$ associated to $\wtx_Z$ (resp. $\wtx_Z^{\prime}$). By letting $\lambda\to 1^+$ and taking $\xi$-adic completion, we then get a morphism  of sheaves of $\BBdRp$-algebras on $X_{Z,v}^{\prime}$
  \begin{equation}\label{equ:relative morphism on OBdRdagger}
      f^*\calO\bB_{\dR,\wtx_Z}^{\dagger,+}\to \calO\bB_{\dR,\wtx^{\prime}_Z}^{\dagger,+}.
  \end{equation}
 \begin{prop}\label{prop:relative Poincare over OBdR}
     The following sequence of sheaves of $\bB_{\dR}^+$-algebras on $X_{Z,v}^{\prime}$
     \[0\to f^*\calO\bB_{\dR,\wtx_Z}^{\dagger,+}\to \calO\bB_{\dR,\wtx^{\prime}_Z}^{\dagger,+}\xrightarrow{\nabla_{\wtx^{\prime}_Z/\wtx_Z}}\calO\bB_{\dR,\wtx^{\prime}_Z}^{\dagger,+}\otimes_{\calO_{\wtx^{\prime}_Z}}\Omega^1_{\wtx^{\prime}_Z/\wtx_Z}\{-1\}\to \cdots\]
     is exact, where the connection 
     \[\nabla_{\wtx^{\prime}_Z/\wtx_Z}:\calO\bB_{\dR,\wtx^{\prime}_Z}^{\dagger,+}\to\calO\bB_{\dR,\wtx^{\prime}_Z}^{\dagger,+}\otimes_{\calO_{\wtx^{\prime}_Z}}\Omega^1_{\wtx^{\prime}_Z/\wtx_Z}\{-1\}\]
     is induced by the relative differential $\rd_{\wtx_Z^{\prime}/\wtx_Z}:\calO_{\wtx_Z^{\prime}}\to \Omega^1_{\wtx_Z^{\prime}/\wtx_Z}$.
 \end{prop}
 \begin{proof}
     We may assume $f:\frakX_Z^{\prime}\to \frakX_Z$ fits into the diagram \eqref{equ:relative chart-I} such that its lifting $\widetilde f$ satisfies \eqref{equ:relative chart-II}.
     By Nakayama's Lemma, we can check the desired exactness modulo $\xi$ and then conclude by (the same proof of) \cite[Prop. 2.29]{Wan23} (together with the compatibility of period sheaves, cf. Proposition \ref{prop:compare with Wang}).
 \end{proof}

\subsection{Compared with period sheaves in \cite{Wan23} and \cite{LMNQW}}\label{ssec:compared with known period sheaves}
  In this subsection, we still assume that $\lambda = p^{\rho}$ for some $\rho\in [0,\frac{1}{p-1})$ and regard $\calO\bB_{\dR}^{I,\lambda,+}$ and $\calO\bB_{\dR}^{I,\dagger,+}$ as the $\xi$-adic completions of $\calO\CCtilde^{1,\lambda}$ and $\calO\CCtilde^{1,\dagger}$ respectively. Our purpose is to compare $(\calO\bC^{\dagger},\Theta)$ and $(\calO\bB_{\dR}^{\dagger,+},\nabla)$ constructed above with the period sheaves in \cite{Wan23} and \cite{LMNQW} respectively.

  We first compare our $(\calO\bC^{\dagger},\Theta)$ with that in \cite{Wan23}. Let us recall the construction in \cite[\S2.3]{Wan23}. Noting that by letting $k$ be absolutely unramified and that $\frac{t}{\xi}\equiv \zeta_p-1\mod t\BBdRp(Z)$, if we choose $\rho_k = \zeta_p-1$ in \cite{Wan23}, then $\rho_k\widehat \calO_{X_Z}^+\otimes_{\calO_{\wtx_Z}}\Omega^{1,\log}_{\frakX_Z}(-1)$ (for $Z = \Spa(C,\calO_C)$ and $\frakX_Z$ smooth) in \cite{Wan23} is exactly our $\widehat \calO_{X_Z}^+\otimes_{\calO_{\frakX_Z}}\Omega^{1,\log}_{\frakX_Z}\{-1\}$.
  \begin{construction}\label{construction:Wang}
      Let $\calE^+$ be the integral Faltings' extension introduced in Definition \ref{dfn:integral Faltings extension}. For any $\rho\in [0,\frac{1}{p-1})$ with $\lambda = p^{\rho}$, denote by $\calE^+_{\lambda}$ the pull-back the $\calE^+$ along the natural inclusion 
      \[p^{\rho}\widehat \calO_{X_Z}^+\otimes_{\calO_{\frakX_Z}}\Omega^{1,\log}_{\frakX_Z}\{-1\}\subset \widehat \calO_{X_Z}^+\otimes_{\calO_{\frakX_Z}}\Omega^{1,\log}_{\frakX_Z}\{-1\}\]
      with respect to the short exact sequence in Proposition \ref{prop:Integral Faltings extension}. In other words, $\calE_{\lambda}^+$ is a sub-$\widehat \calO_{X_Z}^+$-module of $\calE^+$ such that the following diagram is commutative
      \begin{equation*}
          \xymatrix@C=0.5cm{
            0\ar[r] &\widehat \calO_{X_Z}^+\ar@{=}[d]\ar[rr]&&\calE^+_{\lambda}\ar[rr]\ar[d]^{\cap}&& p^{\rho}\widehat \calO_{X_Z}^+\otimes_{\calO_{\frakX_Z}}\Omega^{1,\log}_{\frakX_Z}\{-1\}\ar[d]^{\cap}\ar[r]&0\\
            0\ar[r] &\widehat \calO_{X_Z}^+\ar[rr]&&\calE^+\ar[rr]&& \widehat \calO_{X_Z}^+\otimes_{\calO_{\frakX_Z}}\Omega^{1,\log}_{\frakX_Z}\{-1\}\ar[r]&0.
          }
      \end{equation*}
      In particular, we have $\calE_0^+ = \calE^+$. Denote by $e$ the global section of $\calE^+_{\lambda}$ corresponding to the image of $1\in \widehat \calO_{X_Z}^+$ via the inclusion $\widehat \calO_{X_Z}^+\subset \calE_{\lambda}^+$.
      Denote by $\calO\bC_{\lambda}^+$ the following colimit
      \[\calO\bC_{\lambda}^+:=\colim_{n}\Sym^n_{\widehat \calO_{X_Z}^+}\calE_{\lambda}^+\]
      where the translation map $\Sym^n_{\widehat \calO_{X_Z}^+}\calE_{\lambda}^+\to \Sym^{n+1}_{\widehat \calO_{X_Z}^+}\calE_{\lambda}^+$
      is given by sending each local section 
      \[x_1\otimes\cdots\otimes x_n\]
      to 
      \[e\otimes x_1\otimes\cdots\otimes x_n\]
      for $n\geq 1$ and is given by the natural inclusion of $\widehat \calO_{X_Z}^+\subset \calE_{\lambda}^+$ for $n=0$. This is a sheaf of $\widehat \calO_{X_Z}^+$-algebras on $X_{Z,v}$ and is equiped with a Higgs field
      \[\Theta:\calO\bC_{\lambda}^+\to\calO\bC_{\lambda}^+\otimes_{\calO_{\frakX_Z}}\Omega^{1,\log}_{\frakX_Z}\{-1\}\]
      by taking the colimit of the map
      \[\Sym_{\widehat \calO_{X_Z}^+}^n\calE_{\lambda}^+\to\Sym_{\widehat \calO_{X_Z}^+}^{n-1}\calE_{\lambda}^+\otimes_{\calO_{\frakX_Z}}\Omega^{1,\log}_{\frakX_Z}\{-1\}\]
      which sends each local section 
      \[x_1\otimes\cdots\otimes x_n\]
      to
      \[\sum_{i=1}^nx_1\otimes\cdots\otimes x_{i-1}\otimes x_{i+1}\otimes\cdots\otimes x_n\otimes\overline{x_i},\]
      where $\overline{x_i}$ denotes the image of $x_i$ via the projection $\calE^+_{\lambda}\to p^{\rho}\widehat \calO_{X_Z}^+\otimes_{\calO_{\frakX_Z}}\Omega^{1,\log}_{\frakX_Z}\{-1\}$. 

      Let $(\calO\widehat \bC_{\lambda}^+,\Theta)$ be the $p$-adic completion of $(\calO\bC_{\lambda}^+,\Theta)$. Then the colimit 
      \[(\colim_{\lambda\to 1^+}\calO\widehat \bC_{\lambda}^+[\frac{1}{p}],\Theta)\]
       is exactly the period sheaf with Higgs field in \cite[Def. 2.27]{Wan23}.
  \end{construction}
  \begin{lem}\label{lem:OC_lambda}
      Let $\Sym_{\widehat \calO_{X_Z}^+}\calE_{\lambda}^+$ be the tensor algebra over $\widehat \calO_{X_Z}^+$ generated by $\calE_{\lambda}^+$. Then the morphism 
      \[\Sym_{\widehat \calO_{X_Z}^+}\calE_{\lambda}^+\to\calO\bC_{\lambda}^+\]
      of sheaves of $\widehat \calO_{X_Z}^+$-algebras on $X_{Z,v}$ associated to the natural inclusion $\calE_{\lambda}^+\to \calO\bC_{\lambda}^+$ uniquely factors through the quotient map
      \[\Sym_{\widehat \calO_{X_Z}^+}\calE_{\lambda}^+\to\Sym_{\widehat \calO_{X_Z}^+}\calE_{\lambda}^+/(e-1)\]
      and induces an isomorphism of sheaves of $\widehat \calO_{X_Z}^+$-algebras
      \[\Sym_{\widehat \calO_{X_Z}^+}\calE_{\lambda}^+/(e-1)\xrightarrow{\cong}\calO\bC_{\lambda}^+.\]
  \end{lem}
  \begin{proof}
      Since the problem is local on $\frakX_{Z,\et}$ and $X_{Z,v}$, we may assume that $\frakX_Z$ is small semi-stable and prove the result on $X_{Z,v}$ for any affinoid perfectoid $U\in X_{Z,v}$. Keep notations as in Lemma \ref{lem:local description of integral Faltings extension}. By construction of $\calE_{\lambda}^+$, we see that 
      \[\calE_{\lambda}^+\cong \widehat \calO_{X_Z}^+\cdot e\oplus(\oplus_{i=0}^d\widehat \calO_{X_Z}^+\cdot p^{\rho}W_i/\widehat \calO_{X_Z}^+\cdot \sum_{i=0}^rp^{\rho}W_i)\]
      such that the natural projection $\calE_{\lambda}^+\to p^{\rho}\widehat \calO_{X_Z}^+\otimes_{\calO_{\frakX_Z}}\Omega^{1,\log}_{\frakX_Z}\{-1\}$ is given by sending $e$ to $0$ and each $p^{\rho}W_i$ to $p^{\rho}\frac{\dlog T_i}{\xi}$ for any $0\leq i\leq d$. In particular, via the isomorphisms
      \[\calE_{\lambda}^+\cong \widehat \calO_{X_Z}^+\cdot e\oplus (\oplus_{i=1}^d\widehat \calO_{X_Z}^+\cdot p^{\rho}W_i)\]
      and
      \[\widehat \calO_{X_Z}^+\otimes_{\calO_{\frakX_Z}}\Omega_{\frakX_Z}^{1,\log}\{-1\}\cong\oplus_{i=1}^d\widehat \calO_{X_Z}^+\cdot\frac{\dlog T_i}{\xi},\]
      as mentioned in the paragraph above \cite[Cor. 2.22]{Wan23}, we see that 
      \[\calO\bC_{\lambda}^+\cong \widehat \calO_{X_Z}^+[p^{\rho}W_1,\dots,p^{\rho}W_d]\]
      such that the Higgs field $\Theta$ on $\calO\bC_{\lambda}^+$ is given by
      \[\Theta = \sum_{i=1}^d\frac{\partial}{\partial (p^{\rho}W_i)}\otimes p^{\rho}\frac{\dlog T_i}{\xi}: \widehat \calO_{X_Z}^+[p^{\rho}W_1,\dots,p^{\rho}W_d]\to \oplus_{i=1}^d \widehat \calO_{X_Z}^+[p^{\rho}W_1,\dots,p^{\rho}W_d]\cdot p^{\rho}\frac{\dlog T_i}{\xi}.\]
      On the other hand, it is easy to see that
      \[\Sym_{\widehat \calO_{X_Z}^+}\calE_{\lambda}^+\cong \widehat \calO_{X_Z}^+[e,p^{\rho}W_1,\dots,p^{\rho}W_d]\]
      and that the natural map $\Sym_{\widehat \calO_{X_Z}^+}\calE_{\lambda}^+\to \calO\bC_{\lambda}^+$ sends $e$ to $1$, by the construction of $\calO\bC_{\lambda}^+$. So the result follows trivially.
  \end{proof}

  Now, we are able to prove that our period sheaf with connection $(\calO\bC^{\dagger},\Theta)$ is compatible with that in \cite{Wan23} (which is $(\colim_{\lambda\to 1^+}\calO\widehat \bC_{\lambda}^+[\frac{1}{p}],\Theta)$ in Construction \ref{construction:Wang}).

  \begin{prop}\label{prop:compare with Wang}
      Consider the morphism of sheaves of $\widehat \calO_{X_Z}^+$-algebras on $X_{Z,v}$
      \[\Sym_{\widehat \calO_{X_Z}^+}\calE_{\lambda}^+\to \calO\widehat \bC_{\lambda}\]
      associated to the composite of morphisms of $\widehat \calO_{X_Z}^+$-modules
      \[\calE_{\lambda}^+\hookrightarrow \calE^+\hookrightarrow\calO\bC\hookrightarrow\calO\widehat \bC_{\lambda}.\]
      \begin{enumerate}
          \item[(1)] The morphism $\Sym_{\widehat \calO_{X_Z}^+}\calE_{\lambda}^+\to \calO\widehat \bC_{\lambda}$ uniquely factors through the quotient
          \[\Sym_{\widehat \calO_{X_Z}^+}\calE_{\lambda}^+\to \Sym_{\widehat \calO_{X_Z}^+}\calE_{\lambda}^+/(e-1)\]
          and thus induces a morphism of sheaves of $\widehat \calO_{X_Z}^+$-algebras
          \[\calO\bC_{\lambda}^+\to \calO\widehat \bC_{\lambda}\]
          which is compatible with Higgs fields.

          \item[(2)] The morphism $\calO\bC_{\lambda}^+\to \calO\widehat \bC_{\lambda}$ canonically factors as
          \[\calO\bC_{\lambda}^+\to\calO\widehat \bC_{\lambda}^+\to\calO\widehat \bC_{\lambda}^+[\frac{1}{p}]\to\calO\widehat \bC_{\lambda}\]
          and induces an isomorphism
          \[\calO\widehat \bC_{\lambda}^+[\frac{1}{p}]\xrightarrow{\cong}\calO\widehat \bC_{\lambda}\]
          which is compatible with Higgs fields.

          \item[(3)] There is a natural isomorphism of period sheaves with Higgs field
          \[(\colim_{\lambda\to 1^+}\calO\widehat \bC_{\lambda}^+[\frac{1}{p}],\Theta)\xrightarrow{\cong}(\calO\bC^{\dagger},\Theta).\]
      \end{enumerate}
  \end{prop}
  \begin{proof}
      For item (1): By the definition of $\calE^+$, we see that the morphism 
      \[\calE_{\lambda}^+\hookrightarrow \calE^+\hookrightarrow\calO\bC\]
      maps $e$ to $1$. This yields the desired factorization and thus gives rise to a morphism
      \[\calO\bC_{\lambda}^+\to \calO\widehat \bC_{\lambda}.\]
      It suffices to show that this morphism is compatible with Higgs fields. 
      
      Since the problem is local on $\frakX_{Z,\et}$ and $X_{Z,v}$, we may assume $\frakX_Z$ is small semi-stable and work on $X_{Z,v}/U$ for any affinoid perfectoid $U\in X_{Z,V}$. Keep notations as in the proof of Lemma \ref{lem:OC_lambda}. By Corollary \ref{cor:local description of OBdRp^lambda}, the morphism of sheaves of $\widehat \calO_{X_Z}$-algebras on $X_{Z,v}/U$
      \[\widehat \calO_{X_Z}\{W_1/\lambda,\dots,W_d/\lambda\}\to\calO\widehat \bC_{\lambda},\]
      which sends each $W_i$ to $\frac{\log\frac{T_i}{t_i}}{\xi}$ for any $1\leq i\leq d$, is an isomorphism such that Higgs field on $\calO\widehat \bC_{\lambda}$ is given by 
      \[\Theta = \sum_{i=1}^d\frac{\partial}{\partial W_i}\otimes\frac{\dlog T_i}{\xi}.\]

      On the other hand, the morphism of sheaves of $\widehat \calO_{X_Z}$-algebras on $X_{Z,v}/U$
      \[\widehat \calO_{X_Z}[p^{\rho}W_1,\dots,p^{\rho}W_d]\to\calO\bC_{\lambda}^+,\]
      which sends each $W_i$ to $\frac{\log\frac{T_i}{t_i}}{\xi}$ for any $1\leq i\leq d$, is again an isomorphism such that Higgs field on $\calO\bC_{\lambda}^+$ is given by
      \[\Theta = \sum_{i=1}^d\frac{\partial}{\partial (p^{\rho}W_i)}\otimes p^{\rho}\frac{\dlog T_i}{\xi} = \sum_{i=1}^d\frac{\partial}{\partial W_i}\otimes\frac{\dlog T_i}{\xi}.\]

      Via the above two isomorphisms, the map $\calO\bC^+_{\lambda}\to \calO\widehat \bC_{\lambda}$ is given by the morphism of sheaves of $\widehat \calO_{X_Z}^+$-algebras
      \[\widehat \calO_{X_Z}^+[p^{\rho}W_1,\dots,p^{\rho}W_d]\to \widehat \calO_{X_Z}\za W_1/\lambda,\dots,W_d/\lambda\ya\]
      mapping each $W_i$ to $W_i$ for any $1\leq i\leq d$. So we get the compatibility of Higgs fields.

      For item (2): As $\lambda = p^{\rho}$, we see that  
      \[\widehat \calO_{X_Z}\za W_1/\lambda,\dots,W_d/\lambda\ya = \widehat \calO_{X_Z}^+\za p^{\rho }W_1,\dots,p^{\rho }W_d\ya[\frac{1}{p}],\]
      where $\widehat \calO_{X_Z}^+\za p^{\rho }W_1,\dots,p^{\rho }W_d\ya$ is the $p$-adic completion of $\widehat \calO_{X_Z}^+[p^{\rho}W_1,\dots,p^{\rho}W_d]$. By the proof of Item (1), we see that the morphism $\calO\bC^+_{\lambda}\to \calO\widehat \bC_{\lambda}$ factors through $\calO\widehat \bC^+_{\lambda}$, and then factors through $\calO\widehat \bC^+_{\lambda}[\frac{1}{p}]$ because $p$ is invertible in $\calO\widehat \bC_{\lambda}$. The argument in Item (1) also implies that the morphism $\calO\widehat \bC^+_{\lambda}[\frac{1}{p}]\to \calO\widehat \bC_{\lambda}$ is an isomorphism and compatible with Higgs fields.

      For item (3): It follows from Item (2) by letting $\lambda\to 1^+$ directly.
  \end{proof}

  The rest of this subsection is devoted to comparing $(\calO\bB_{\dR}^{\dagger,+},\nabla)$ with $(\calO\bB_{\dR,\pd}^+,\nabla)$ in \cite{LMNQW}. Let us recall the construction of $(\calO\bB_{\dR,\pd}^+,\nabla)$ for the convenience of the reader.
  \begin{construction}\label{construction:LMNQW}
      Fix an affinoid perfectoid $U = \Spa(S,S^+)\in X_{Z,v}^{\Box}$ and a $\frakY_Z = \Spf(\calR^{\prime}_Z)$ and let $\theta^{\prime}_{\frakY_Z}:A_{\frakY_Z}^+\to S^+$ be the surjection as before. Consider the blow-up algebra
      \[(A_{\frakY_Z}^+)^{\wedge}[\frac{u\cdot\Ker(\theta_{\frakY_Z}^{\prime})}{\xi}]\subset (A_{\frakY_Z}^+)^{\wedge}[\frac{1}{\xi}],\] 
      where $(A_{\frakY_Z}^+)^{\wedge}$ denotes the $(p,\Ker(\theta_{\frakY_Z}^{\prime}))$-adic completion of $A_{\frakY_Z}^+$. Define
      \[B_{\frakY_Z,\dR,\pd}^+:=\varprojlim_n\big((A_{\frakY_Z}^+)^{\wedge}[\frac{u\cdot\Ker(\theta_{\frakY_Z}^{\prime})}{\xi}]^{\wedge}_{\pd}\widehat \otimes_{\bA_{\inf}(U)}\CCtilde^1(U)\big)/(\xi^n),\]
      where $(A_{\frakY_Z}^+)^{\wedge}[\frac{u\cdot\Ker(\theta_{\frakY_Z}^{\prime})}{\xi}]^{\wedge}_{\pd}$ is the $(p,\Ker(\theta_{\frakY_Z}^{\prime}))$-adic completion of $(A_{\frakY_Z}^+)^{\wedge}[\frac{u\cdot\Ker(\theta_{\frakY_Z}^{\prime})}{\xi}]_{\pd}$, the pd-envelope of $(A_{\frakY_Z}^+)^{\wedge}[\frac{u\cdot\Ker(\theta_{\frakY_Z}^{\prime})}{\xi}]$ with respect to the ideal generated by $\frac{u\cdot\Ker(\theta_{\frakY_Z}^{\prime})}{\xi}$, and the complete tensor product is with respect to the $(p,\xi)$-adic topology on $\bA_{\inf}(U)$. Then $\calO\bB_{\dR,\pd}^+$ is the sheaf of $\BBdRp$-algebras on $X_{Z,v}$ associated to the presheaf
      \[U\in X_{Z,v}^{\Box}\mapsto \colim_{\frakY_Z}B_{\frakY_Z,\dR,\pd}^+\]
      and the connection 
      \[\nabla:\calO\bB_{\dR,\pd}^+\to \calO\bB_{\dR,\pd}^+\otimes_{\calO_{\wtx_Z}}\Omega^{1,\log}_{\frakX_Z}\{-1\}\]
      is induced by the usual $(p,\xi)$-adically continuous differential on $\wtx_Z$.
  \end{construction}
  \begin{lem}\label{lem:local description of OBdRpd}
      Suppose that $\frakX_Z$ is small affine and $U\in X_{Z,v}$ is affinoid perfectoid. keep notations in Lemma \ref{lem:finite ideal}. Then the morphism of sheaves of $\BBdRp$-algebras on $X_{Z,v}/U$
      \[\BBdRp\za uW_1,\dots,uW_d\ya_{\pd}\to \calO\bB_{\dR,\pd}^+,\]
      which sends each $W_i$ to $\frac{\log\frac{T_i}{t_i}}{\xi}$ for any $1\leq i\leq d$, is an isomorphism, where \[\BBdRp\za uW_1,\dots,uW_d\ya_{\pd} = \big(\CCtilde^{1,+}[uW_1,\dots,uW_d]^{\wedge_p}_{\pd}[\frac{1}{p}]\big)^{\wedge_{\xi}}\]
      is the $\xi$-adic completion of $\CCtilde^{1,+}[uW_1,\dots,uW_d]^{\wedge_p}_{\pd}[\frac{1}{p}]$ while $\CCtilde^{1,+}[uW_1,\dots,uW_d]^{\wedge_p}_{\pd}$ is the $p$-adic completion of the free pd-algebras over $\CCtilde^{1,+}$ generated by $W_1,\dots,W_d$. Moreover, via the isomorphism
      \[\Omega^{1,\log}_{\wtx_Z}\{-1\}\cong \oplus_{i=1}^d\calO_{\wtx_Z}\cdot\frac{\dlog T_i}{\xi},\]
      the connection $\nabla$ on $\calO\bB_{\dR,\pd}^+$ is given by
      \[\nabla = \sum_{i=1}^d\frac{\partial}{\partial W_i}\otimes\frac{\dlog T_i}{\xi}.\]
  \end{lem}
  \begin{proof}
      This follows from \cite[Prop. 3.11]{LMNQW} directly.
  \end{proof}
  In order to compare $(\calO\bB_{\dR}^{\dagger,+},\nabla)$ with $(\calO\bB_{\dR,\pd}^+,\nabla)$, we have to compare $(\calO\bB_{\dR,\pd}^+,\nabla)$ with $(\calO\bB_{\dR}^{1,+},\nabla)$ at first. (Recall that $\calO\bB_{\dR}^{1,+} = \calO\bB_{\dR}^{I,\lambda,+}$ for $I = \{1\}$ and $\lambda=1$.) Note that $u=[\epsilon^{\frac{1}{p}}]-1$ admits arbitrary pd-powers in $\CCtilde^{1,+}(Z)$ as $|u|_1 = p^{-\frac{1}{p-1}}$, cf. Example \ref{exam:norm of u}.
  \begin{construction}\label{construction:compare with LMBQW}
      Fix an affinoid perfectoid $U\in X_{Z,v}$ and a $\frakY = \Spf(\calR_Z^{\prime})\in \Sigma_U$. Then the natural map $A_{\frakY_Z}^+\to B_{\frakY_Z,\dR}^{1,1,+}$ factors as
      \[A_{\frakY_Z}^+\to A_{\frakY_Z}^{1,1,+}\to A_{\frakY_Z}^{1,1}\to B_{\frakY_Z,\dR}^{1,1,+}.\]
      By definition of $A_{\frakY_Z}^{1,1}$, for any $x\in \Ker(\theta_{\frakY_Z}^{\prime})$, we have $|\frac{x}{\xi}|_1 \leq 1$, so the natural morphism $A_{\frakY_Z}^+\to A_{\frakY_Z}^{1,1,+}$ uniquely factors as
      \[A_{\frakY_Z}^+\to (A_{\frakY_Z}^+)^{\wedge}\to(A_{\frakY_Z}^+)^{\wedge}[\frac{u\cdot\Ker(\theta_{\frakY_Z}^{\prime})}{\xi}]^{\wedge}_{\pd}\to \big((A_{\frakY_Z}^+)^{\wedge}[\frac{u\cdot\Ker(\theta_{\frakY_Z}^{\prime})}{\xi}]^{\wedge}_{\pd}\big)\widehat \otimes_{\bA_{\inf}(U)}\CCtilde^{1,+}(U)\to A_{\frakY_Z}^{1,1,+}.\]
      Inverting $p$ and taking $\xi$-adic completion, we get a natural morphism of $\BBdRp(U)$-algebras
      \[B_{\frakY_Z,\dR,\pd}^+\to B_{\frakY_Z,\dR}^{1,+}.\]
      Letting $U\in X_{Z,v}^{\Box}$ and $\frakY_Z\in \Sigma_U$ vary, we obtain a natural morphism of sheaves of $\BBdRp$-algebras on $X_{Z,v}$ which is compatible with the connections
      \[(\calO\bB_{\dR,\pd}^+,\nabla)\to (\calO\bB_{\dR}^{1,+},\nabla).\]
  \end{construction}

  \begin{lem}\label{lem:compare with LMNQW}
      Suppose that $\frakX_{Z,v}$ is small semi-stable. Then for any affinoid perfectoid $U\in X_{Z,v}$, via the isomorphisms
      \[\calO\bB_{\dR,\pd}^+\cong \BBdRp\za uW_1,\dots,uW_d\ya_{\pd}\]
      and 
      \[\calO\bB_{\dR}^{1,+}\cong\BBdRp\za W_1,\dots,W_d\ya\]
      respectively described in Lemmas \ref{lem:local description of OBdRp^lambda} and \ref{lem:local description of OBdRpd}, the morphism $\calO\bB_{\dR,\pd}^+\to \calO\bB_{\dR}^{1,+}$ constructed above is given by
      \[\BBdRp\za uW_1,\dots,uW_d\ya_{\pd}\to\BBdRp\za W_1,\dots,W_d\ya \]
      sending each $uW_i$ to $uW_i$ for any $1\leq i\leq d$. In particular, the morphism $\calO\bB_{\dR,\pd}^+\to \calO\bB_{\dR}^{1,+}$ is an injection.
  \end{lem}
  \begin{proof}
      This follows from Lemma \ref{lem:local description of OBdRp^lambda} and Lemma \ref{lem:local description of OBdRpd} directly.
  \end{proof}
  Now, we are able to compare $(\calO\bB_{\dR}^{\dagger,+},\nabla)$ with $(\calO\bB_{\dR,\pd}^+,\nabla)$.
  \begin{prop}\label{prop:compare with LMNQW}
      For any $\lambda\in (1,p^{\frac{1}{p-1}})$, the natural morphism $\calO\bB_{\dR}^{\lambda,+}\to \calO\bB_{\dR}^{1,+}$ factors through the image of the map $\calO\bB_{\dR,\pd}^+\to \calO\bB_{\dR}^{1,+}$ and thus induces a unique morphism of sheaves of $\BBdRp$-algebras on $X_{Z,v}$
      \[\calO\bB_{\dR}^{\lambda,+}\to \calO\bB_{\dR,\pd}^+\]
      which is compatible with the connections. As a consequence, we get a morphism of period sheaves with connections
      \[(\calO\bB_{\dR}^{\dagger,+},\nabla)\to (\calO\bB_{\dR,\pd}^+,\nabla).\]
  \end{prop}
  \begin{proof}
      It suffices to show that the morphism $\calO\bB_{\dR}^{\lambda,+}\to \calO\bB_{\dR}^{1,+}$ factors through the image of the map $\calO\bB_{\dR,\pd}^+\to \calO\bB_{\dR}^{1,+}$ for any $\lambda = p^{\rho}$ for $\rho\in (0,\frac{1}{p-1})\cap \bQ$. Since this is a local problem on $\frakX_{Z,\et}$ and $X_{Z,v}$, we may assume that $\frakX_Z$ is small semi-stable and work on $X_{Z,v}/U$ for any affinoid perfectoid $U\in X_{Z,v}$. Note that $\frac{u}{[\varpi^{\frac{1}{p-1}}]}$ is a unit in $\CCtilde^{1,+}(Z)$, cf. Example \ref{exam:t in C^I}(3). Using Lemmas \ref{lem:local description of OBdRp^lambda} and \ref{lem:compare with LMNQW}, it suffices to show that the morphism
      \[\BBdRp\za [\varpi^{\rho}]W_1,\dots,[\varpi^{\rho}]W_d\ya\to \BBdRp\za W_1,\dots,W_d\ya\]
      uniquely factors through the inclusion
      \[\BBdRp\za [\varpi^{\frac{1}{p-1}}]W_1,\dots,[\varpi^{\frac{1}{p-1}}]W_d\ya_{\pd}\to \BBdRp\za W_1,\dots,W_d\ya.\]
      
      For any $n\geq 0$, define $e_n\in\bN$ such that it is $0$ when $n = 0$ and is the maximal integer such that $p^{e_n}\mid n$ when $n\geq 1$. It is well-known that 
      \[e_n = \frac{n-\sum_{i=0}^ra_i}{p-1}\geq\frac{n}{p-1}-\log_p(n)-1,\]
      where $a_i\in \{0,1,\dots,p-1\}$ for each $i$ such that $n = \sum_{i=0}^ra_ip^i$ with $a_r\neq 0$. So for any $(n_1,\dots,n_d)\in\bN^d$, we have
      \[\sum_{i=1}^de_{n_i}-(\frac{1}{p-1}-\rho)\sum_{i=1}^dn_i\geq \rho\cdot\sum_{i=1}^dn_i-\sum_{i=1}^d\log_pn_i-d\to +\infty,\]
      yielding that 
      \[|[\varpi^{(\rho-\frac{1}{p-1})\cdot\sum_{i=1}^dn_i}]\cdot\prod_{i=1}^dn_i!|_1 \to 0\]
      as $n_1+\dots+n_d\to +\infty$.

      So we obtain the desired inclusion
      \[\BBdRp\za [\varpi^{\rho}]W_1,\dots,[\varpi^{\rho}]W_d\ya\subset\BBdRp\za [\varpi^{\frac{1}{p-1}}]W_1,\dots,[\varpi^{\frac{1}{p-1}}]W_d\ya_{\pd}\]
      by
      noting that
      \[\prod_{i=1}^d([\varpi^{\rho}]W_i)^{n_i} = [\varpi^{(\rho-\frac{1}{p-1})\cdot\sum_{i=1}^dn_i}]\cdot\prod_{i=1}^dn_i!\cdot \prod_{i=1}^d\frac{([\varpi^{\frac{1}{p-1}}]W_i)^{n_i}}{n_i!}\]
      for any $(n_1,\dots,n_d)\in\bN^d$.
  \end{proof}
  \begin{rmk}
      Denote by $(\calO\widehat \bC_{\pd},\Theta)$ the reduction of $(\calO\bB_{\dR,\pd}^+,\nabla)$ modulo $\xi$. Then Proposition \ref{prop:compare with LMNQW} gives rise to a morphism of period sheaves with Higgs fields
      \[(\calO\bC^{\dagger},\Theta)\to(\calO\widehat \bC_{\pd},\Theta).\]
      This map essentially occurs in \cite[\S2.3]{MW-AIM}.
  \end{rmk}

\section{A decompletion theorem}\label{sec:decompletion}
  This section is devoted to studying the decompletion theory mentioned in \S\ref{sec:Intro}. Before we move on, let us fix some notation.
  Fix an affinoid perfectoid $Z = \Spa(A,A^+)\in \Perfd$ and a closed interval $I=[r,s]\subset (0,+\infty)$. Let $\frakX_Z = \Spf(\calR_Z)$ be a small semi-stable formal scheme of relative dimension $d$ over $A^+$ with generic fiber $X_Z = \Spa(R_Z,R_Z^+)$ and fixed prismatic lifting $\wtx_Z = \Spf(\wtr_Z)$ over $\bA_{\inf}(Z)$. Let $X_{\infty,Z} = \Spa(\widehat R_{\infty,Z},\widehat R_{\infty,Z}^+)\to X_Z$ be the pro-finite-\'etale covering with Galois group $\Gamma$ from Example \ref{exam:Gamma torsor}. The chart $\psi$ on $\frakX_Z$ induces a morphism 
  \[\widetilde \psi:\wtr_Z\to \bA_{\inf}(X_{\infty,Z})\]
  uniquely determined by sending each $T_i$ to $[T_i^{\flat}]$ such that the decomposition \eqref{equ:decomposition over Ainf} holds:
  \[\bA_{\inf}(X_{\infty,Z}) = \widehat \bigoplus_{\underline \alpha\in J_{r,d}}\wtr_Z\cdot[(\underline T^{\flat})^{\underline \alpha}],\]
  where for any $\underline \alpha = (\alpha_0,\dots,\alpha_d)\in J_{r,d}$, $(\underline T^{\flat})^{\underline \alpha}=\prod_{i=0}^d(T_i^{\flat})^{\alpha_i}$. The ring $\bA_{\inf}(X_{\infty,Z})$ carries a $\Gamma$-action such that for any $\delta = \prod_{i=0}^d\delta_i^{n_i}\in \Gamma$ and any $\underline \alpha = (\alpha_0,\dots,\alpha_d)\in J_{r,d}$, 
  \[\delta([(\underline T^{\flat})^{\underline \alpha}]) = [\epsilon^{\sum_{i=0}^dn_i\alpha_i}]\cdot[(\underline T^{\flat})^{\underline \alpha}],\]
  and this decomposition is $\Gamma$-equivariant. The morphism $\widetilde \psi$ induces
  \[\widetilde \psi^I:\wtr_Z^I\to \CCtilde^I(X_{\infty,Z}),\]
  and one obtains the $\Gamma$-equivariant decompositions \eqref{equ:decomposition over C^Iplus} and \eqref{equ:decomposition over C^I} of \S\ref{sec:period sheaf without connection}:
  \[\CCtilde^{I,+}(X_{\infty,Z}) = \widehat \bigoplus_{\underline \alpha\in J_{r,d}}\wtr_Z^{I,+}\cdot[(\underline T^{\flat})^{\underline \alpha}],\quad \CCtilde^I(X_{\infty,Z}) = \widehat \bigoplus_{\underline \alpha\in J_{r,d}}\wtr_Z^I\cdot[(\underline T^{\flat})^{\underline \alpha}].\]
  We write $\CCpsi:=\CCtilde^I(X_{\infty,Z})$ and $\Cpsi:=\widetilde \psi^I(\wtr_Z^I)$ for the completed and decompleted period rings on $X_{\infty,Z}$, and $\CCppsi$, $\Cppsi$ for their integral models.
  \begin{dfn}\label{dfn:Cpsi}
      Put $(\CCpsi,\CCppsi):=(\CCtilde^I(X_{\infty,Z}),\CCtilde^{I,+}(X_{\infty,Z}))$ and $(\Cpsi,\Cppsi):=(\widetilde \psi^I(\wtr_Z^I),\widetilde \psi^I(\wtr_Z^{I,+}))$, equipped with the $\Gamma$-action above, so that
      \[\CCppsi = \widehat \bigoplus_{\underline \alpha\in J_{r,d}}\Cppsi\cdot[(\underline T^{\flat})^{\underline \alpha}]\text{ and }\CCpsi = \widehat \bigoplus_{\underline \alpha\in J_{r,d}}\Cpsi\cdot[(\underline T^{\flat})^{\underline \alpha}].\]
  \end{dfn}
\subsection{Smallness and iterative refinement}\label{ssec:basic facts in decompletion}
  \begin{lem}\label{lem:isometric decomposition}
      The morphism $\widetilde \psi^I:\wtr_Z^{I,+}\to \CCppsi$ is isometric and thus the norm on $\Cppsi$ induced by the norm $|\cdot|_I$ on $\wtr_Z^{I,+}$ via the morphism $\widetilde \psi^I$ coincides with the restriction of the norm $|\cdot|_I$ on $\CCppsi$. Moreover, if we equip $\CCppsi/\Cppsi$ with the norm induced by $|\cdot|_I$ on $\CCppsi$, then the natural surjection of $\Cppsi$-modules $\CCppsi\to\CCppsi/\Cppsi$ admits an isometric $\Gamma$-equivariant section 
      \[s:\CCppsi/\Cppsi\to \CCppsi\]
      given by the inclusion
      \[\CCppsi/\Cppsi = \widehat \bigoplus_{0\neq \underline \alpha\in J_{r,d}}\Cppsi\cdot[(\underline T^{\flat})^{\underline \alpha}]\subset \widehat \bigoplus_{\underline \alpha\in J_{r,d}}\Cppsi\cdot[(\underline T^{\flat})^{\underline \alpha}] = \CCppsi.\]
  \end{lem}
  \begin{proof}
      This follows trivially from Corollary \ref{cor:Independence of embedding}.
  \end{proof}
  \begin{lem}\label{lem:smallness for Cpsi}
      Put $u = [\epsilon^{\frac{1}{p}}]-1$.
      For any $\gamma\in \Gamma$ and any $x\in \Cppsi$, we have $(\gamma-1)(x)\in \phi(u)\Cppsi$ such that
      \[|\frac{(\gamma-1)(x)}{u}|_I \leq |\xi|_I|x|_I.\]
  \end{lem}
  \begin{proof}
      As shown in \cite[\S9.2, pp. 340]{BMS18}, we know that for any $\delta\in \Gamma$ and any $x\in \wtr_Z$, we have 
      \[(\gamma-1)(x)\in \phi(u)\wtr_Z.\]
      As $\Cppsi = \wtr_Z\widehat \otimes_{\bA_{\inf}(Z)}\CCtilde^{I,+}(Z)$, we see that the same holds true for any $x\in \Cppsi$ as expected. It remains to show that $|\frac{(\gamma-1)(x)}{u}|_I \leq |\xi|_I|x|_I$. It suffices to show that for any $t\in I$, we have
      \[|\frac{(\gamma-1)(x)}{u}|_t \leq |\xi|_t|x|_t.\]
      Granting this, we see that 
      \[|\frac{(\gamma-1)(x)}{u}|_t \leq |\xi|_t|x|_t\leq |\xi|_I|x|_t\]
      which implies $|\frac{(\gamma-1)(x)}{u}|_I \leq |\xi|_I|x|_I$ by furthermore taking $\sup_{t\in I}$ on both sides.

      Now, we are going to prove $|\frac{(\gamma-1)(x)}{u}|_t \leq |\xi|_t|x|_t$. For any $0<\varepsilon<1$, choose $a\in C^{\flat\times}$ such that $\varepsilon<|[a]x|_t\leq 1$ and write $(\gamma-1)([a]x) = \phi(u)y = u\xi y$ for some $y\in \Cppsi$. Then we have 
      \[|[a]|_t|\frac{(\gamma-1)(x)}{u}|_t = |[a]\frac{(\gamma-1)(x)}{u}|_t=|\frac{(\gamma-1)([a]x)}{u}|_t = |\xi y|_t\leq |\xi|_I|y|_t\leq \varepsilon^{-1}|\xi|_t|[a]x|_t = |[a]|_t\varepsilon^{-1}|\xi|_t|x|_t,\]
      yielding that $|\frac{(\gamma-1)(x)}{u}|_t\leq \varepsilon^{-1}|\xi|_t|x|_t$ because $|[a]|_t\neq 0$.       Now, one can conclude by letting $\varepsilon\to 1$.
  \end{proof}
  \begin{rmk}
      Lemma \ref{lem:smallness for Cpsi} is the $\Gamma$-equivariant estimate that makes the $u$-adic error terms in Proposition \ref{prop:u-error} controllable: the action on $\Cppsi$ is ``$p$-adically small'' relative to $|\xi|_I$, which is exactly what later enters the $\lambda$-triviality inequalities for decompleted representations.
  \end{rmk}

\subsection{Decompletion of small representations}\label{ssec:decompletion theory}

 \begin{dfn}
     Let $G$ be a topological ring and let $B$ be a topological ring equipped with a continuous $G$-action. By a representation of $G$ over $B$ of rank $r$, we mean a finite free $B$-module $M$ of rank $r$ equipped with a continuous semi-linear action of $G$. Denote by $\Rep_G(B)$ the category of representations of $G$ over $B$.
 \end{dfn}
 In this subsection, we will consider the case for $G = \Gamma$ and $B\in \{\Cppsi,\CCppsi,\Cpsi,\CCpsi\}$.
 
 \begin{dfn}\label{dfn:small representation}
     Let $(B,B^+)$ be either $(\CCpsi,\CCppsi)$ or $(\Cpsi,\Cppsi)$ and equip $B$ with the norm $|\cdot|_I$.
     \begin{enumerate}
         \item Let $\lambda>1$ be a real number and $M$ be a representation of $\Gamma$ over $B^+$ of rank $r$. We say $M$ is \emph{$\lambda$-trivial} if there exists a $B^+$-basis $e_1,\dots,e_r$ of $M$ such that if we equip $M$ with the supremum norm, denoted by $|\cdot|_I$, via the isomorphism $M\cong \oplus_{i=1}^rB^+\cdot e_i$ corresponding to the given basis, then for any $\gamma\in\Gamma$ and any $1\leq i\leq r$, we have 
         \begin{equation}\label{equ:small representation}
             (\gamma-1)(e_i)\subset uM \text{ such that } |\frac{(\gamma-1)}{u}(e_i)|\leq \lambda^{-1}.
         \end{equation}
         A representation $M$ of $\Gamma$ over $B$ is called \emph{$\lambda$-trivial} if it admits a $\Gamma$-stable $B^+$-lattice which is $\lambda$-small. Denote by $\Rep_{\Gamma}^{\lambda}(B^+)$ (resp. $\Rep_{\Gamma}^{\lambda}(B)$) the category of $\lambda$-trivial representations of $\Gamma$ over $B^+$ (resp. $B$).

         \item Let $\lambda_0\geq 1$. A representation $M$ of $\Gamma$ over $B^+$ (resp. $B$) is called \emph{$\lambda_0$-small} if it is $(\lambda p^{\frac{1}{p-1}})$-trivial for some $\lambda>\lambda_0$.
         Denote by $\Rep_{\Gamma}^{\lambda_0,\dagger}(B^+)$ (resp. $\Rep_{\Gamma}^{\lambda_0,\dagger}(B)$) the category of $\lambda_0$-small representations of $\Gamma$ over $B^+$ (resp. $B$). We also say a representation $M$ of $\Gamma$ over $B^+$ (resp. $B$) is \emph{small} if it is $1$-small and denote by $\Rep_{\Gamma}^{\dagger}(B^+)$ (resp. $\Rep_{\Gamma}^{\dagger}(B)$) the corresponding category $\Rep_{\Gamma}^{1,\dagger}(B^+)$ (resp. $\Rep_{\Gamma}^{1,\dagger}(B)$) for short.
     \end{enumerate}
 \end{dfn}
 \begin{rmk}\label{rmk:small cocycle}
     Keep notations in Definition \ref{dfn:small representation}(1) and equip $\Mat_r(B)$ with the supremum norm $|\cdot|_I$. Let $f:\Gamma \to \GL_r(B^+)$ be the $1$-cocycle determined by $M$ with respect to the given $B^+$-basis $e_i$'s. Then the condition (\ref{equ:small representation}) amounts to that the function $f-1:\Gamma\to \Mat_r(B^+)$ takes values in $u\cdot\Mat_r(B^+)$ with the norm $|\frac{f-1}{u}|_I\leq \lambda^{-1}$. In particular, $\Gamma$ acts isometrically on $M$.
 \end{rmk}

 \begin{lem}\label{lem:lambda(M)}
     Let $\lambda>1$. For any $M\in \Rep^{\lambda}_{\Gamma}(\Cppsi)$, put $\lambda(M):=\inf(\lambda,|\xi|_I^{-1})$. Then for any $x\in M$ and any $\gamma\in\Gamma$, we have
     \[(\gamma-1)x \in uM \text{ such that }|\frac{(\gamma-1)}{u}(x)|_I\leq \lambda(M)^{-1}|x|_I.\]
 \end{lem}
 \begin{proof} 
     Fix an $x\in M$, and write $x = \sum_{i=1}^rc_ie_i$. For any $\delta\in \Gamma$, we have
     \[(\gamma-1)(x) = \sum_{i=1}^d(\gamma-1)(c_i)\gamma(e_i)+c_i(\gamma-1)(e_i).\]
     As $(\gamma-1)(c_i)\in \phi(u)\Cppsi$ by Lemma \ref{lem:smallness for Cpsi}, we see that $(\gamma-1)(x)\in uM$ so that
     \[\frac{(\gamma-1)(x)}{u} = \sum_{i=1}^d\big(\frac{(\gamma-1)(c_i)}{u}\gamma(e_i)+c_i\frac{(\gamma-1)(e_i)}{u}\big).\]
     Applying Lemma \ref{lem:smallness for Cpsi} together with \eqref{equ:small representation} again, we see that 
     \[|\frac{(\gamma-1)(c_i)}{u}\gamma(e_i)+c_i\frac{(\gamma-1)(e_i)}{u}|_I\leq \sup\{|\xi|_I|c_i|_I,\lambda^{-1}|c_i|_I\}= \lambda(M)^{-1}|c_i|_I,\]
     yielding that $|\frac{(\gamma-1)(x)}{u}|_I\leq\lambda(M)^{-1}\sup_{1\leq i\leq r}\{|c_i|_I\} = \lambda(M)^{-1}|x|_I$ as desired.
 \end{proof}
 Note that for any $\lambda>1$, the base-change $M\mapsto M_{\infty}:=M\otimes_{\Cpsi}\CCpsi$ gives rise to a functor 
 \[\Rep^{\lambda}_{\Gamma}(\Cpsi)\to\Rep^{\lambda}_{\Gamma}(\CCpsi),\]
 such that, via the decomposition in Definition \ref{dfn:Cpsi}, we have a $\Gamma$-equivariant decomposition
 \begin{equation}\label{equ:decomposition for M}
     M_{\infty} = \widehat \bigoplus_{\underline \alpha\in J_{r,d}}M\cdot[(\underline T^{\flat})^{\underline \alpha}].
 \end{equation}
 The following proposition tells us that this functor is fully faithful.

 \begin{prop}\label{prop:u-error}
    Fix a $\lambda> 1$ and put $u = [\epsilon^{\frac{1}{p}}]-1$.
    For any $\lambda$-trivial representation $M$ of $\Gamma$ over $\Cppsi$, we have that for any $n\geq 0$, the cohomological group
    \[\rH^n(\Gamma,M\otimes_{\Cppsi}(\CCppsi/\Cppsi))\]
    is killed by $u$ and vanishes when $n = 0$.
 \end{prop}
 \begin{proof}
     It suffices to show that for any $\underline \alpha = (\alpha_0,\dots,\alpha_d)\in J_{r,d}$, the cohomological group 
     \[\rH^n(\Gamma,M\cdot[(\underline T)^{\flat})^{\underline \alpha}])\]
     is killed by $u$ and vanishes when $n=0$.
     Without loss of generality, we may assume that $\alpha_0 = 0$ and furthermore assume $\alpha_1\neq 0$. Note that we have the isomorphism \eqref{equ:gamma group-II}
     \[\Gamma\cong \oplus_{i=1}^d\Zp\cdot \gamma_i\]
     such that for any $1\leq i\leq d$, 
     \[\gamma_i(T_j^{\flat}) = \epsilon^{\delta_{ij}}T_j^{\flat}.\]
     So for $b = \frac{u}{[\epsilon^{\alpha_1}]-1}\in \Ainf$, we have
     \[(\gamma_1-1)(b[(\underline T^{\flat})^{\underline \alpha}]) = u[(\underline T^{\flat})^{\underline \alpha}].\]

     Fix a basis $e_1,\dots,e_r$ of $M$ and let $\lambda(M)$ be as in Lemma \ref{lem:lambda(M)}. Then for any $x=x_0\in M$, $y_0:=\frac{(\gamma_1-1)}{u}(x_0)$ is a well-defined element in $M$ satisfying $|y_0|_I\leq \lambda(M)^{-1}|x_0|_I$. Since
     \[ux_0[(\underline T^{\flat})^{\underline \alpha}] - (\gamma_1-1)(bx_0[(\underline T^{\flat})^{\underline \alpha}]) = -u[\epsilon^{\alpha_1}]by_0[(\underline T^{\flat})^{\underline \alpha}],\]
     put $x_1= -[\epsilon^{j_1}]by_0\in M$ and then we have that $|x_1|_I\leq \lambda(M)^{-1}|x_0|_I$ and
     \[ux_0[(\underline T^{\flat})^{\underline \alpha}] - (\gamma_1-1)(bx_0[(\underline T^{\flat})^{\underline \alpha}]) = ux_1[(\underline T^{\flat})^{\underline \alpha}].\]
     Repeating the above process, we then have $x_n\in  M$ satisfying 
     \[|x_n|_I\leq \lambda(M)^{-1}|x_{n-1}|_I\leq \cdots\leq \lambda(M)^{-n}|x_0|_I\]
     for any $n\geq 1$ such that
     \[ux_{n-1}[(\underline T^{\flat})^{\underline \alpha}] - (\gamma_1-1)(bx_{n-1}[(\underline T^{\flat})^{\underline \alpha}]) = ux_n[(\underline T^{\flat})^{\underline \alpha}].\]
     Now $z:=\sum_{n\geq 0}x_n$ is a well-defined element in $M$ with norm $|z|_I = |x|_I$ such that
     \[ux[(\underline T^{\flat})^{\underline \alpha}] = (\gamma_1-1)(bz[(\underline T^{\flat})^{\underline \alpha}]).\]
     This implies that $\rH^1(\Zp\gamma_1,M\cdot [(\underline T^{\flat})^{\underline \alpha}])$ is killed by $u$, yielding that $\rH^n(\Gamma,M\cdot [(\underline T^{\flat})^{\underline \alpha}])$ is also killed by $u$, by Hochschild--Serre spectral sequence, for any $n\geq 1$. 

     It remains to show that $\rH^0(\Gamma,M\cdot [(\underline T^{\flat})^{\underline \alpha}]) = 0$.
     Let $x\in M$ be an element such that $x[(\underline T^{\flat})^{\underline \alpha}]\in \rH^0(\Gamma,M\cdot [(\underline T^{\flat})^{\underline \alpha}])$. We have to show that $x = 0$. Applying Lemma \ref{lem:lambda(M)} again, we see that $y = \frac{(\gamma_1-1)}{u}(x)$ is a well-defined element in $M$ such that $|y|_I\leq \lambda(M)^{-1}|x|_I$. As $\gamma_1(x[(\underline T^{\flat})^{\underline \alpha}])=x[(\underline T^{\flat})^{\underline \alpha}]$, we have
     \[x = -\frac{u}{[\epsilon^{j_1}]-1}[\epsilon^{\alpha_1}]y.\]
     Note that $-\frac{u}{[\epsilon^{\alpha_1}]-1}[\epsilon^{j_1}]\in\bfA_{\inf}$ because $\alpha\in \bN[\frac{1}{p}]\cap(0,1)$. So we have
     \[|x|_I\leq |y|_I\leq \lambda(M)^{-1}|x|_I,\]
     which forces $|x|_I = 0$ (as $\lambda(M)>1$), and a fortiori $x = 0$ as desired.
 \end{proof}
 \begin{cor}\label{cor:u-error}
     The cohomological group $\rH^n(\Gamma,\CCppsi/\Cppsi)$ is killed by $u$ for any $n\geq 0$, and vanishes when $n = 0$.
 \end{cor}
 \begin{proof}
     This is a special case of Proposition \ref{prop:u-error} for $M = \Cpsi$.
 \end{proof}
 \begin{cor}\label{cor:u-error implies fully faithful}
     Let $\lambda>1$. 
     \begin{enumerate}
         \item[(1)] The base-change functor
         \[\Rep_{\Gamma}^{\lambda}(\Cppsi)\to\Rep_{\Gamma}^{\lambda}(\CCppsi)\]
         sending each $M$ to $M_{\infty}:=M\otimes_{\Cppsi}\CCppsi$ is fully faithful.

         \item[(2)] The base-change functor
         \[\Rep_{\Gamma}^{\lambda}(\Cpsi)\to\Rep_{\Gamma}^{\lambda}(\CCpsi)\]
         sending each $M$ to $M_{\infty}:=M\otimes_{\Cppsi}\CCppsi$ is fully faithful. If moreover $u$ is invertible in $\CCtilde^I(Z)$ (cf. Example \ref{exam:t in C^I}(2)), then the inclusion $M\subset M_{\infty}$ induces a quasi-isomorphism
         \[\rR\Gamma(\Gamma,M)\simeq \rR\Gamma(\Gamma,M_{\infty}).\]
     \end{enumerate}
 \end{cor}
 \begin{proof}
     As the base-change functor $M\mapsto M_{\infty} = M\otimes_{\Cpsi}\CCpsi$ preserves duals and tensor products, the result follows directly from Proposition \ref{prop:u-error}.
 \end{proof}
 Now, let $S$ be a topological space which is endowed with a continuous $\Gamma$-action. For any $n\geq 0$, let $\rC(\Gamma^n,S)$ be the set of continuous functions from $\Gamma^n$ to $S$. When $S$ is furthermore a Banach space over $\Qp$ with norm $|\cdot|$, we always equip $\rC(\Gamma^n,S)$ with the supremum norm $|f|:=\sup_{g\in\Gamma^n}|f(g)|$ for any $f\in \rC(\Gamma^n,S)$. In this case, we know that the complex $\rC(\Gamma^{\bullet},S)$ computes the (continuous) $\Gamma$-cohomology $\rH^n(\Gamma,S)$ of $S$.
 \begin{prop}\label{prop:essentially surjective}
     Let $I\subset(0,+\infty)$ be a closed interval such that $u$ is invertible in $\CCtilde^I(Z)$ (cf. Example \ref{exam:t in C^I}(2)). Fix a $\lambda> p^{\frac{1}{p-1}}$. Then for any $M_{\infty}\in \Rep^{\lambda}_{\Gamma}(\CCppsi)$, there exists some $M\in \Rep^{\lambda}_{\Gamma}(\Cppsi)$ such that 
     \[M_{\infty}\cong M\otimes_{\Cppsi}\CCppsi.\]
 \end{prop}
 \begin{proof}
     We will proceed as in the proof \cite[Lem. A.1.15]{DLLZ23}. 
     Let $e_1,\dots,e_r$ be a $\CCppsi$-basis of $M_{\infty}$ as in Definition \ref{dfn:small representation}. For any $\gamma\in \Gamma$, let $f(\gamma)\in \Mat_r(\CCppsi)$ be the matrix of $\gamma$ corresponding to the chosen basis $\underline e:=(e_1,\dots,e_r)$; that is, we have
     \[\gamma(\underline e) = \underline e\cdot f(\gamma).\]
     Then $f\in \rC(\Gamma,\GL_r(\CCppsi))$ is an $1$-cocycle such that 
     \begin{equation}\label{equ:another description of lambda-small}
        (f-1)(\Gamma)\subset u\Mat_r(\CCppsi)\text{ with }|\frac{f-1}{u}|_I\leq\lambda^{-1}.
     \end{equation}
     as mentioned in Remark \ref{rmk:small cocycle}. 
     One can also regard $\rC(\Gamma^n,\GL_r(\CCppsi))$ as a subset of $\rC(\Gamma^n,\Mat_r(\CCppsi))$ (so that $f$ is just a cochain in the latter). For any $g\in \rC(\Gamma^n,\Mat_r(\Cppsi))$, denote by $\overline g$ its image in $\rC(\Gamma^n,\Mat_r(\CCppsi)/\Mat_r(\Cppsi))$.
     
     We claim for any $f\in \rC(\Gamma,\GL_r(\CCppsi))$ satisfying (\ref{equ:another description of lambda-small}), there is an 
     \[h\in \rC(\Gamma,\Mat_r(\CCppsi)/\Mat_r(\Cppsi)),\]
     viewed also as an element in $\rC(\Gamma,\Mat_r(\CCppsi))$, via the section $s$ given in Lemma \ref{lem:isometric decomposition}, such that $|h|_I\leq |\frac{\overline{f}}{u}|_I\leq |\frac{f-1}{u}|_I\leq \lambda^{-1}<1$ (which implies that $1+h\in \GL_r(\CCppsi)$) and that the $1$-cocycle
     \[f^{\prime}:\Gamma\to \GL_r(\CCppsi),\quad \gamma\mapsto (1+h)^{-1}f(\gamma)\gamma(1+h)\]
     satisfies the condition (\ref{equ:another description of lambda-small}) and $|\frac{\overline{f^{\prime}}}{u}|_I\leq\lambda^{-1}|\frac{\overline{f}}{u}|_I$. 
     
     Granting the claim, we can iteratively find a sequence $h_1, h_2,\ldots$ in $\Mat_r(\CCppsi/\Cppsi)$ with $|h_n|_I\leq \lambda^{1-n}|\frac{\overline{f}}{u}|_I\leq \frac 1{\lambda^n}$ such that
        \[
            |\overline{\left(\prod_{i = 1}^n (1 + h_i)^{-1}\right)f(\gamma)\gamma\left(\prod_{i = 1}^n (1 + h_i)\right)}/u|_I\leq \lambda^{-n}|\frac{\overline{f}}{u}|_I.
        \]
     Then the infinite product $h_{\infty}:=\prod_{i = 1}^\infty (1 + h_i)$ is a well-defined element in $\GL_r(\CCppsi)$. Consider the $1$-cocycle
     
     \[\tilde{f} : \Gamma\to \GL_r(\CCppsi),\quad \gamma \mapsto h_\infty^{-1}f(\gamma)\gamma(h_\infty).\]
     It satisfies the condition \eqref{equ:another description of lambda-small} by construction such that 
     \[|\frac{\overline{\tilde{f}}}{u}|_I\leq \lambda^{-n}|\frac{\overline f}{u}|_I\]
     for any $n\geq 0$. Letting $n\to+\infty$, we get 
     \[|\frac{\overline{\tilde{f}}}{u}|_I = 0\]
     and thus 
     \[\overline{\tilde{f}} = 0.\]
     In other words, $\tilde{f}$ takes values in $\Mat_r(\Cppsi)\cap \GL_r(\CCppsi)$. So $\tilde{f}$ is an $1$-cocycle in $\rC(\Gamma, \GL_r(\Cppsi))$ satisfying the condition \eqref{equ:another description of lambda-small}. By Remark \ref{rmk:small cocycle}, $\tilde{f}$ defines a $\lambda$-small representation $M$ of $\Gamma$ over $\Cppsi$ and by construction we have
     \[M\otimes_{\Cppsi}\CCppsi \cong M_{\infty}\]
     as desired.
     
     Now, we are going to prove the claim. Regarding $f$ as a cochain in $\rC(\Gamma,\Mat_r(\CCppsi))$, for any $\gamma,\gamma^{\prime}\in \Gamma$, we have 
        \[\begin{split}
           (\rd f)(\gamma,\gamma^{\prime}) = & \gamma(f(\gamma^{\prime}))-f(\gamma\gamma^{\prime})+f(\gamma)\\
           = &\gamma(f(\gamma^{\prime}))-f(\gamma)\gamma(f(\gamma^{\prime}))+f(\gamma)\\
           =& (1-f(\gamma))(\gamma(f(\gamma^{\prime}))-1)-1.
        \end{split}\]
        This implies that 
        \[\overline{\rd f}(\gamma,\gamma^{\prime}) = \overline{(1-f(\gamma))(\gamma(f(\gamma^{\prime}))-1)-1} = \overline{(1-f(\gamma))(\gamma(f(\gamma^{\prime}))-1)}.\]
        By \eqref{equ:another description of lambda-small}, we have $\rd\overline{f}\in u^2\rC(\Gamma^2,\Mat_r(\CCppsi)/\Mat_r(\Cppsi))$. Using \cite[Lem. A.1.14]{DLLZ23} (and the section $s$ given in Lemma \ref{lem:isometric decomposition}), we obtain that
        \[|\frac{\rd\overline{f}}{u^2}|_I\leq |\frac{f-1}{u}|_I|\frac{\overline f}{u}|_I\leq \lambda^{-1}|\frac{\overline f}{u}|_I,\]
        where the last inequality follows as $|\frac{f-1}{u}|_I\leq\lambda^{-1}$ by  \eqref{equ:another description of lambda-small}. Since $\frac{\rd \overline f}{u^2}\in\rC^2(\Gamma,\Mat_r(\CCppsi)/\Mat_r(\Cppsi))$ is a cocycle, by Corollary \ref{cor:u-error}, there exists some $\overline g\in \rC^1(\Gamma,\Mat_r(\CCppsi)/\Mat_r(\Cppsi))$ with norm
        \[|\overline g|_I\leq |\frac{\rd\overline{f}}{u^2}|_I\leq \lambda^{-1}|\frac{\overline f}{u}|_I\]
        such that $\rd\overline g = \frac{\rd \overline f}{u}$. This implies that $\overline{g}-\frac{\overline f}{u}$ is a cocycle in $\rC^1(\Gamma,\Mat_r(\CCppsi)/\Mat_r(\Cppsi))$ and we have 
        \[|\overline{g}-\frac{\overline f}{u}|_I = |\frac{\overline f}{u}|_I.\] Again by Corollary \ref{cor:u-error}, there exists some $\overline h\in \Mat_r(\CCppsi)/\Mat_r(\Cppsi)$ with norm 
        \[|\overline h|_I\leq |\frac{\overline f}{u}|_I\leq \lambda^{-1}\] 
        such that for any $\gamma\in\Gamma$,
        \[\gamma(\overline h)-\overline h = \rd\overline h(\gamma) = u\overline g(\gamma)-\overline f(\gamma).\]
        Regard $\overline h$ as an element in $\Mat_r(\CCppsi)$ via the section $s$ in Lemma \ref{lem:isometric decomposition} and put $h:=s(\overline h)$. Then for any $\gamma\in\Gamma$, we have
        \[\begin{split}
          f^{\prime}(\gamma)-1 = & (1+h)^{-1}f(\gamma)\gamma(1+h)-1\\
          = & (1+h)^{-1}(f(\gamma)-1)(1+\gamma(h))+(1+h)^{-1}(1+\gamma(h))-1\\
          = & (1+h)^{-1}(f(\gamma)-1)(1+\gamma(h))+(1+h)^{-1}\rd h(\gamma).
        \end{split}\]
        Note that 
        \[|\frac{\rd h}{u}|_I = |\frac{\rd \overline h}{u}|_I\leq |\overline g - \frac{\overline f}{u}|_I = |\frac{\overline f}{u}|_I\leq \lambda^{-1}.\] 
        As $f$ satisfies the condition \eqref{equ:another description of lambda-small}, so is $f^{\prime}$. Note that for any $\gamma\in\Gamma$,
        \[\begin{split}
            (\frac{\overline{f'}}{u} - \overline g)(\gamma) = & \overline{(1+h)^{-1}\frac{f(\gamma)-1}{u}(1+\gamma(h))+(1+h)^{-1}\frac{\rd h(\gamma)}{u}} - \overline{\frac{\rd h(\gamma)+f(\gamma)}{u}}\\
            = & \overline{(1+h)^{-1}\frac{f(\gamma)-1}{u}(1+\gamma(h))-\frac{f(\gamma)-1}{u}}+\overline{((1+h)^{-1}-1)\frac{\rd h(\gamma)}{u}}.
        \end{split}\]
        By \cite[Lem. A.1.14]{DLLZ23}, we then have
        \[|\frac{\overline{f'}}{u} - \overline g|_I\leq \sup(|h|_I\cdot |\frac{\overline f}{u}|_I,|\frac{f-1}{u}|_I\cdot |\overline h|_I,|h|_I\cdot|\frac{\overline{\rd h}}{u}|_I,|\overline h|_I\cdot|\frac{\rd h}{u}|_I).\]
        As $|h|_I = |\overline h|_I\leq |\frac{\overline f}{u}|_I\leq \lambda^{-1}$ and $|\overline g|_I\leq \lambda^{-1}|\frac{\overline f}{u}|_I$, we finally obtain that $|\frac{\overline f^{\prime}}{u}|_I\leq \lambda^{-1}|\frac{\overline f}{u}|_I$ as desired. This proves the claim and thus completes the proof.
 \end{proof}

 Combining Corollary \ref{cor:u-error implies fully faithful} with Proposition \ref{prop:essentially surjective}, we obtain the following decompletion theorem.
 \begin{thm}\label{thm:decompletion}
     Let $I = [r,s]\subset(0,p)$ be a closed interval and $\lambda> p^{\frac{1}{p-1}}$ be a real number.
     \begin{enumerate}
         \item[(1)] The base-change $M\mapsto M_{\infty}:=M\otimes_{\Cppsi}\CCppsi$ induces an equivalence of categories
         \[\Rep^{\lambda}_{\Gamma}(\Cppsi)\simeq \Rep^{\lambda}_{\Gamma}(\CCppsi)\]
         such that the morphism
         \[\rR\Gamma(\Gamma,M)\to\rR\Gamma(\Gamma,M_{\infty})\]
         induced by the inclusion $M\subset M_{\infty}$ identifies $\rR\Gamma(\Gamma,M)$ as a direct summand of $\rR\Gamma(\Gamma,M_{\infty})$ whose complement is concentrated in degree $\geq 1$ and is killed by $u$.

         \item[(2)] The base-change $M\mapsto M_{\infty}:=M\otimes_{\Cppsi}\CCppsi$ induces an equivalence of categories
         \[\Rep^{\lambda}_{\Gamma}(\Cpsi)\simeq \Rep^{\lambda}_{\Gamma}(\CCpsi)\]
         such that the morphism
         \[\rR\Gamma(\Gamma,M)\to\rR\Gamma(\Gamma,M_{\infty})\]
         induced by the inclusion $M\subset M_{\infty}$ is a quasi-isomorphism.
     \end{enumerate}
 \end{thm}
 \begin{proof}
     The assertion follows from Proposition \ref{prop:u-error}, Corollary \ref{cor:u-error implies fully faithful} and Proposition \ref{prop:essentially surjective}.
 \end{proof}
 \begin{cor}\label{cor:decompletion}
     Let $I = [r,s]\subset(0,p)$ be a closed interval and $\lambda_0\geq 1$ be a real number.
     \begin{enumerate}
         \item[(1)] The base-change $M\mapsto M_{\infty}:=M\otimes_{\Cppsi}\CCppsi$ induces an equivalence of categories
         \[\Rep^{\lambda_0,\dagger}_{\Gamma}(\Cppsi)\simeq \Rep^{\lambda_0,\dagger}_{\Gamma}(\CCppsi)\]
         such that the morphism
         \[\rR\Gamma(\Gamma,M)\to\rR\Gamma(\Gamma,M_{\infty})\]
         induced by the inclusion $M\subset M_{\infty}$ identifies $\rR\Gamma(\Gamma,M)$ as a direct summand of $\rR\Gamma(\Gamma,M_{\infty})$ whose complement is concentrated in degree $\geq 1$ and is killed by $u$.

         \item[(2)] The base-change $M\mapsto M_{\infty}:=M\otimes_{\Cppsi}\CCppsi$ induces an equivalence of categories
         \[\Rep^{\lambda_0,\dagger}_{\Gamma}(\Cpsi)\simeq \Rep^{\lambda_0,\dagger}_{\Gamma}(\CCpsi)\]
         such that the morphism
         \[\rR\Gamma(\Gamma,M)\to\rR\Gamma(\Gamma,M_{\infty})\]
         induced by the inclusion $M\subset M_{\infty}$ is a quasi-isomorphism.
     \end{enumerate}
 \end{cor}
 \begin{proof}
     The assertion follows from Theorem \ref{thm:decompletion} by noting that
     \[\Rep_{\Gamma}^{\lambda_0,\dagger}(B) = \colim_{\lambda>\lambda_0 p^{\frac{1}{p-1}}}\Rep_{\Gamma}^{\lambda}(B)\]
     for any $B\in \{\Cppsi,\CCppsi,\Cpsi,\CCpsi\}$.
 \end{proof}

 Then above equivalence behaves well when shrinking $I$.
 \begin{lem}\label{lem:shrink I-local representation}
     Let $J\subset I$ be a closed sub-interval. Then for any $\lambda>1$, the base-change \[ \text{$M\mapsto M\otimes_{\Cppsi}\widetilde \bfC^{J,+}_{\widetilde \psi} $(resp. $M_{\infty}\mapsto M_{\infty}\otimes_{\CCppsi}\widetilde \bfC^{J,+}_{\widetilde \psi,\infty}$)}\]
     induces a functor
     \[\text{$\Rep_{\Gamma}^{\lambda}(\Cppsi)\to \Rep_{\Gamma}^{\lambda}(\widetilde \bfC^{J,+}_{\widetilde \psi})$ (resp. $\Rep_{\Gamma}^{\lambda}(\CCppsi)\to \Rep_{\Gamma}^{\lambda}(\widetilde \bfC^{J,+}_{\widetilde \psi,\infty})$)}\]
     such that the following diagram commutes
     \[\xymatrix@C=0.5cm{
       \Rep_{\Gamma}^{\lambda}(\Cppsi)\ar[d]\ar[rr]&& \Rep_{\Gamma}^{\lambda}(\widetilde \bfC^{J,+}_{\widetilde \psi})\ar[d]\\
       \Rep_{\Gamma}^{\lambda}(\CCppsi)\ar[rr]&& \Rep_{\Gamma}^{\lambda}(\widetilde \bfC^{J,+}_{\widetilde \psi,\infty})
     }\]
     where the vertical arrows denotes the functors induced by the base-changes along $\Cppsi\to\CCppsi$ and $\widetilde \bfC^{J,+}_{\widetilde \psi}\to\widetilde \bfC^{J,+}_{\widetilde \psi,\infty}$ respectively.
 \end{lem}
 \begin{proof}
     The base-change functor is well-defined following from the definition of $\lambda$-trivial representations \eqref{equ:small representation} because for any $x\in \CCppsi$, it always holds true that $|x|_J\leq |x|_I$. The above diagram is clearly commutative because every functor involved is induced by base-change. 
 \end{proof}

\section{Local Riemann--Hilbert Correspondence: The Robba Case}\label{sec:local RH}
   Throughout this section, let $Z = \Spa(A,A^+)\in\Perfd$ and let $\frakX_Z = \Spf(\calR_Z)$ be a small semi-stable formal scheme of relative dimension $d$ over $A^+$, equipped with a chart $\psi$ and a fixed lifting $\wtx_Z = \Spf(\wtr_Z)$ over $\bA_{\inf}(Z)$. Let $I = [r,s]\subset (\frac{1}{p-1},p)$ be a closed sub-interval. We trivialize the Breuil--Kisin twist $\CCtilde^{I,+}\{1\}$ using $\frac{t}{u}$, where $t = \log[\epsilon]$ is a unit multiple of $u\xi$ with $u = [\epsilon^{1/p}]-1$ (cf.\ Example \ref{exam:t in C^I}(1)). Then we have $\CCtilde^{I,+}(Z)$-algebras $\wtr_Z^{I,+}$ and $\wtr_Z^I$, and as in \S\ref{ssec:ringed space X^I} there is a $\CCtilde^{I,+}(Z)$-linear connection
   \[\nabla:\wtr_Z^{I,+}\to\Omega^{1,\log}_{\wtr_Z^{I,+}}\{-1\}.\]
   The chart $\psi$ (together with $\widetilde \psi$, cf.\ Example \ref{exam:small affine is liftable}) induces an isomorphism 
   \[\Omega^{1,\log}_{\wtr_Z^{I,+}}\{-1\}\cong (\oplus_{i=0}^d\wtr_Z^{I,+}\cdot e_i)/\wtr_Z^{I,+}\cdot\sum_{i=0}^re_i,\]
   where $e_i = \frac{u}{t}\rd\log T_i$ for $0\leq i\leq d$,
   so one can equip $\Omega^{1,\log}_{\wtr_Z^{I,+}}\{-1\}$ with the norm $|\cdot|_I$ given by the quotient norm from the supremum norm on $\oplus_{i=0}^d\wtr_Z^{I,+}\cdot e_i$. Via
   \[(\oplus_{i=0}^d\wtr_Z^{I,+}\cdot e_i)/\wtr_Z^{I,+}\cdot\sum_{i=0}^re_i\cong \oplus_{i=1}^d\wtr_Z^{I,+}\cdot e_i\]
   this norm agrees with the supremum norm on $\oplus_{i=1}^d\wtr_Z^{I,+}\cdot e_i$. 
   In this section, we will prove a local version of $p$-adic Riemann--Hilbert over Robba rings (cf. Theorem \ref{thm:intro-C}).

  \subsection{Small local connections}\label{ssec:small local connections Robba}
    \begin{dfn}\label{dfn:small local connection}
        We keep the notation introduced above.
        \begin{enumerate}
            \item By a \emph{flat (or integrable) connection} of rank $r$ over $\wtr_Z^{I,+}$, we mean a finite free module $D$ of rank $r$ over $\wtr_Z^{I,+}$ together with a $\CCtilde^{I,+}(Z)$-linear map
            \[\nabla_D:D\to D\otimes_{\wtr_Z^{I,+}}\Omega^{1,\log}_{\wtr_Z^{I,+}}\{-1\}\]
            satisfying the Leibniz rule with respect to $\nabla:\wtr_Z^{I,+}\to\Omega^{1,\log}_{\wtr_Z^{I,+}}\{-1\}$ such that $\nabla_D\wedge\nabla_D = 0$. Denote by $\MIC(\wtr_Z^{I,+})$ the category of flat connections over $\wtr_Z^{I,+}$. 

            \item Let $\lambda>1$. A flat connection $(D,\nabla_D)$ over $\wtr_Z^{I,+}$ is called \emph{$\lambda$-trivial}, if there exists a basis $e_1,\dots,e_r$ of $D$ such that if we equip $D$ with the supremum norm $|\cdot|_I$ associated to the given basis, then for any $1\leq i\leq d$, we have
            \[|\nabla(e_i)|_I\leq \lambda^{-1}.\]
            Denote by $\MIC^{\lambda}(\wtr_Z^{I,+})$ the full sub-category of $\MIC(\wtr_Z^{I,+})$ consisting of $\lambda$-trivial objects.

            \item Let $\lambda_0\geq 1$. A flat connection $(D,\nabla_D)$ over $\wtr_Z^{I,+}$ is called \emph{$\lambda_0$-small} if there exists some $\lambda>\lambda_0$ such that $(D,\nabla_D)$ is $(\lambda p^{\frac{1}{p-1}})$-trivial. Denote by $\MIC^{\lambda_0,\dagger}(\wtr_Z^{I,+})$ the full sub-category of $\MIC(\wtr_Z^{I,+})$ consisting of $\lambda_0$-small objects. We also say $(D,\nabla_D)$ is \emph{small} if it is $1$-small, and denote the corresponding category by $\MIC^{\dagger}(\wtr_Z^{I,+})$.

            \item One can similarly introduce the category $\MIC(\wtr_Z^I)$ of flat connections over $\wtr_Z^I$ as in Item (1). For such a flat connection $(D,\nabla_D)\in \MIC(\wtr_Z^I)$, we call it $\lambda$-trivial (resp. $\lambda_0$-small, small) if it admits a $\lambda$-trivial (resp. $\lambda_0$-small, small) lattice over $\wtr_Z^{I,+}$, and denote by $\MIC^{\lambda}(\wtr_Z^I)$ (resp. $\MIC^{\lambda_0,\dagger}(\wtr_Z^I)$, $\MIC^{\dagger}(\wtr_Z^I)$) the corresponding category. 

            \item For any flat connection $(D,\nabla_D)$ over $\wtr_Z^{I,+}$ or $\wtr_Z^I$, denote by $\rD\rR(D,\nabla_D)$ its associated de Rham complex:
            \[0\to D\xrightarrow{\nabla_D}D\otimes_{\wtr_Z^{I,+}}\Omega^{1,\log}_{\wtr_Z^{I,+}}\{-1\}\xrightarrow{\nabla_D}D\otimes_{\wtr^{I,+}}\Omega^{2,\log}_{\wtr_Z^{I,+}}\{-2\}\to\cdots\]
            and for any $n\geq0$, define
            \[\rH^n_{\dR}(D,\nabla_D):=\rH^n(\rD\rR(D,\nabla_D)).\]
        \end{enumerate}
    \end{dfn}
    It is worth pointing out that, by Remark \ref{rmk:xi-connection}, a flat connection $(D,\nabla_D)$ in Definition \ref{dfn:small local connection} is actually a $\xi$-connection. In what follows, we list some obvious lemmas.

    \begin{lem}\label{lem:nabla is small}
        For any $f\in \wtr_Z^I$, we have $|\nabla(f)|_I\leq |\xi|_I\cdot|f|_I$.
    \end{lem}
    \begin{proof}
        Without loss of generality, we may always assume $f\neq 0$.
        We first deal with the case for $I=\{r\}$. In this case, for any $0<\varepsilon<1$, we can choose some $a\in C^{\flat\times}$ such that $\varepsilon<|[a]f|_r\leq 1$. So we have that
        \[|[a]|_r\cdot|\nabla(f)|_r = |\nabla([a]f)|_r \leq |\xi|_r\]
        where the last inequality holds true because $\nabla(\wtr_Z^{I,+})\subset \xi\cdot \Omega^{1,\log}_{\wtr_Z^{I,+}}\{-1\}$. This implies that
        \[|\nabla(f)|_r\leq |[a]|_r^{-1}|\xi|_r\leq \varepsilon^{-1}|f|_r\]
        by the choice of $a$. So we get $|\nabla(f)|_r\leq |\xi|_r|f|_r$ by letting $\varepsilon\to 1$. For general case, it remains only to observe that
        \[|\nabla(f)|_I = \sup_{t\in I}|\nabla(f)|_t \leq \sup_{t\in I}\{|\xi|_t|f|_t\}\leq |\xi|_I|f|_I.\]
    \end{proof}
    \begin{rmk}
        Lemma \ref{lem:nabla is small} implies that the universal connection $\nabla$ on $\wtr_Z^I$ is Lipschitz (with constant $|\xi|_I$) for the Robba-type norm; in particular, it respects the overconvergent structure used to define small connections.
    \end{rmk}
    \begin{lem}\label{lem:lambda(D)}
        Let $\lambda>1$. For any $(D,\nabla_D)\in \MIC^{\lambda}(\wtr_Z^{I,+})$, put $\lambda(D) := \inf(\lambda,|\xi|_I^{-1})$. Then for any $x\in D$, we have 
        \[|\nabla_D(x)|_I\leq \lambda(D)^{-1}|x|_I.\]
    \end{lem}
    \begin{proof}
        Let $e_1,\dots,e_r$ be an $\wtr_Z^{I,+}$-basis of $D$ as in Definition \ref{dfn:small local connection}. For any $x\in D$, write $x = \sum_{i=1}^rf_ie_i$ and then we have
        \[\nabla_D(x) = \sum_{i=1}^rf_i\nabla_D(e_i)+\sum_{i=1}^r\nabla(f)_i\otimes e_i.\]
        By the choice of $e_i$'s, we have
        \[|\sum_{i=1}^rf_i\nabla_D(e_i)|_I\leq \lambda^{-1}\sup_{1\leq i\leq r}\{|f_i|_I|e_i|_I\} = \lambda^{-1}|x|_I\]
        and by Lemma \ref{lem:nabla is small}, we have
        \[|\sum_{i=1}^r\nabla(f)_i\otimes e_i|_I\leq |\xi|_I\sup_{1\leq i\leq r}\{|f_i|_I|e_i|_I\} = |\xi|_I|x|_I.\]
        So we finally obtain that
        \[|\nabla_D(x)|_I\leq \sup\{|\sum_{i=1}^rf_i\nabla_D(e_i)|_I,|\sum_{i=1}^r\nabla(f)_i\otimes e_i|_I\}\leq \sup\{\lambda^{-1},|\xi|_I\}|x|_I = \lambda(D)^{-1}|x|_I\]
        as desired.
    \end{proof}
    \begin{lem}\label{lem:shrink I-local connection}
        Let $J\subset I$ be a closed sub-interval. Then for any $\lambda>1$, the base-change $M\mapsto M\otimes_{\wtr_Z^{I,+}}\wtr_Z^{J,+}$ gives rise to a functor
        \[\MIC^{\lambda}(\wtr_Z^{I,+})\to \MIC^{\lambda}(\wtr_Z^{J,+}),\]
        which is compatible with cohomology in the sense that for any $(D,\nabla_D)\in \MIC^{\lambda}(\wtr_Z^{I,+})$, there is a quasi-isomorphism
        \[\rD\rR(D,\nabla_D)\widehat \otimes_{\CCtilde^{I,+}(Z)}\CCtilde^{J,+}(Z)\simeq \rD\rR(D\otimes_{\wtr_Z^{I,+}}\wtr_Z^{J,+},\nabla_D).\]
    \end{lem}
    \begin{proof}
        The well-definedness of the base-change functor follows from the definition of $\lambda$-trivial connections together with the fact that for any $f\in \wtr_Z^{I,+}$, $|x|_J\leq |x|_I$. The desired cohomological comparison follows from the fact that the base-change map
        \[\wtr_Z^{I,+}\widehat \otimes_{\CCtilde^{I,+}(Z)}\CCtilde^{J,+}(Z)\cong \wtr_Z^{J,+}\]
        is an isomorphism, which forces that 
        \[D\widehat \otimes_{\CCtilde^{I,+}(Z)}\CCtilde^{J,+}(Z)\cong D\widehat \otimes_{\wtr_Z^{I,+}}\wtr_Z^{J,+}\cong  D\otimes_{\wtr_Z^{I,+}}\wtr_Z^{J,+},\]
        where the last isomorphism is true because $D$ is finite free over $\wtr_Z^{I,+}$.
    \end{proof}

  \subsection{Logarithmic coordinates and the local theorem}\label{ssec:local RH W and theorem}
    Put $P=W_0+\cdots+W_r$ as in Example \ref{exam:Apply to Gamma torsor}. Via the isomorphism
    \[\Omega^{1,\log}_{\wtr_Z^{I,+}}\{-1\}\cong \oplus_{i=0}^d\wtr_Z^{I,+}\cdot \frac{u}{t}\dlog T_i/\wtr_Z^{I,+}\sum_{i=0}^r\frac{u}{t}\dlog T_i,\]
    we have a $\Gamma$-equivariant isomorphism
    \[(\CCpsi\{\underline W/\lambda\}_{r,d,P} = \CCpsi\{W_0/\lambda,\dots,W_d/\lambda\}^{\dagger}/(P),\nabla = \sum_{i=0}^d-\frac{\partial}{\partial W_i}\otimes\frac{u}{t}\dlog T_i)\cong (\calO\bC^{I,\lambda},\nabla)(X_{\infty,Z})\]
    yielding a $\Gamma$-equivariant isomorphism
    \[(\CCpsi\{\underline W\}_{r,d,P}^{\dagger} = \CCpsi\{W_0,\dots,W_d\}^{\dagger}/(P),\nabla = \sum_{i=0}^d-\frac{\partial}{\partial W_i}\otimes\frac{u}{t}\dlog T_i)\cong (\calO\bC^{I,\dagger},\nabla)(X_{\infty,Z})\]
    by letting $\lambda\to 1^+$, where the $\Gamma$-action on $\CCpsi\{\underline W/\lambda\}_{r,d,P}$ (and thus on $\CCpsi\{\underline W\}_{r,d,P}^{\dagger}$) is given by \eqref{equ:gamma acts on W}. Note that for any $\bfB\in \{\CCppsi,\Cpsi,\Cppsi\}$, 
    \[\bfB\{\underline W\}_{r,d,P}^{\dagger}:=\bfB\{W_0,\dots,W_d\}^{\dagger}/(P)\subset \CCpsi\{\underline W\}_{r,d,P}^{\dagger}\]
    is $\nabla$-preserving and stable under $\Gamma$-action. Similarly, for any $\lambda>1$,
    \[\bfB\{\underline W/\lambda\}_{r,d,P}:=\bfB\{W_0/\lambda,\dots,W_d/\lambda\}/(P)\subset\CCpsi\{\underline W/\lambda\}_{r,d,P}\]
    is also $\nabla$-preserving and stable under $\Gamma$-action.

    Now, we are able to state our main theorem in this section.
    \begin{thm}[Local Riemann--Hilbert correspondence over Robba rings]\label{thm:local correspondence over Robba ring}
        Keep the notation as above.
        \begin{enumerate}
            \item[(1)] For any $M_{\infty}\in \Rep^{\dagger}_{\Gamma}(\CCpsi)$ of rank $r$, we have
            \[
              \rH^n(\Gamma,M_{\infty}\otimes_{\CCpsi}\Ctilde_{\psi,\infty}^I\{\underline W\}^{\dagger}_{r,d,P}) = 
                 \left\{
                    \begin{array}{rcl}
                     D_{\psi}^I(M_{\infty}), & n = 0 \\
                      0, & n\geq 1,
                    \end{array}
                 \right.
            \]
        where $D_{\psi}^I(M_{\infty})$ is a finite free $\wtr^I$-module of rank $r$. Moreover, the connection $\nabla$ on $\calO\CCtilde^{I,\dagger}$ induces a flat connection $\nabla_{D_{\psi}^I(M_{\infty})}$ on $D_{\psi}^I(M_{\infty})$ making $(D_{\psi}^I(M_{\infty}),\nabla_{D_{\psi}^I(M_{\infty})})$ an object in $\MIC^{\dagger}(\wtr_Z^I)$.

        \item[(2)] For any $(D,\nabla_D)\in \MIC^{\dagger}(\wtr_Z^I)$ of rank $r$, define
        \[\rd_D = \nabla_D\otimes\id+\id\otimes\nabla: D\otimes_{\wtr_Z^I}\Ctilde_{\psi,\infty}^I\{\underline W\}^{\dagger}_{r,d,P}\to D\otimes_{\wtr^I}\Ctilde_{\psi,\infty}^I\{\underline W\}^{\dagger}_{r,d,P}\otimes_{\wtr_Z^I}\Omega^1_{\wtr_Z^I}\{-1\}\]
        which satisfies $\rd_D\wedge\rd_D = 0$. Then we have
        \[
          \rH^n_{\dR}(D\otimes_{\wtr_Z^I}\Ctilde_{\psi,\infty}^I\{\underline W\}^{\dagger}_{r,d,P},\rd_D) = 
            \left\{
              \begin{array}{rcl}
                  M_{\infty}(D,\nabla_D), & n=0 \\
                  0, & n\geq 1,
              \end{array}
            \right.
        \]
        where $M_{\infty}(D,\nabla_D)$ is a finite free $\CCpsi$-module of rank $r$. Moreover, the $\Gamma$-action on $\CCpsi\{\underline W\}^{\dagger}$ induces a $\Gamma$-action on $M_{\infty}(D,\nabla_D)$ making it an object in $\Rep_{\Gamma}^{\dagger}(\CCpsi)$.

        \item[(3)] The functors $M_{\infty}\mapsto D^I_{\psi}(M_{\infty})$ and $(D,\nabla_D)\mapsto M_{\infty}(D,\nabla_D)$ induce an equivalence of categories
        \[\Rep_{\Gamma}^{\dagger}(\CCpsi)\simeq \MIC^{\dagger}(\wtr_Z^I)\]
        preserving ranks, tensor products and duals, such that for any $M_{\infty}\in \Rep_{\Gamma}^{\dagger}(\CCpsi)$ corresponding to $(D,\nabla_D)\in \MIC^{\dagger}(\wtr_Z^I)$, we have a quasi-isomorphism
        \[\RGamma(\Gamma,M_{\infty})\simeq \rD\rR(D,\nabla_D).\]
        \end{enumerate}
    \end{thm}

\subsection{$\Gamma$-cohomology of small representations}\label{ssec:Gamma cohomology of small representation}
    We first give another description of $\CCpsi\{\underline W\}^{\dagger}_{r,d,P}$, which will simplify some local calculations. 
    \begin{convention}\label{convention:good form for calculation}
        Let $B$ be a Banach $\Qp$-algebra with a unit $u\in B^{\times}$ such that $|u|<1$.
        For any $j\geq 0$ and any formal variable $V$, we put
        \[\binom{V}{j}_u:=\binom{\frac{V}{u}}{j} = \frac{\prod_{i=0}^{j-1}(\frac{V}{u}-i)}{j!} = \frac{\prod_{i=0}^{j-1}(V-iu)}{u^jj!}.\]
        For any $J=(j_1,\dots,j_d)\in\bN^d$ and any formal variables $V_1,\dots,V_d$, we put
        \[\binom{\underline V}{J}_u:=\prod_{i=1}^d\binom{V_i}{j_i}_u = \frac{\prod_{i=1}^d\prod_{l_i=0}^{j_i}(V_i-l_iu)}{u^{|J|}J!},\]
        where $|J|:=j_1+\cdots+j_d$ and $J!:=j_1!\cdots j_d!$.
    \end{convention}
    \begin{lem}\label{lem:good form for calculation}
        Let $B$ be a Banach $\Qp$-algebra with a unit $u\in B^{\times}$ such that $|u|<1$, and let $\lambda>1$. For any Banach $B$-module $M$ with $M^+:=\{x\in M\mid |x|\leq 1\}$, let $M^+\{\underline V/\lambda\} = M^+\{V_1/\lambda,\dots,V_d/\lambda\}$ be as in Definition \ref{dfn:overconvergent algebras}.
        Then we have
        \[M^+\{\underline V/\lambda\} = \{\sum_{J\in \bN^d}a_Ju^{|J|}J!\binom{\underline V}{J}_u\mid |a_J|\leq 1,~\lim_{|J|\to +\infty\}}\lambda^{|J|}|a_J| = 0\}\]
        As a consequence, we have
        \[M\{\underline V/\lambda\} = \{\sum_{J\in \bN^d}a_J\binom{\underline V}{J}_u\mid \lim_{|J|\to +\infty\}}\lambda^{|J|}|\frac{a_J}{u^{|J|}J!}| = 0\}.\]
    \end{lem}
    \begin{proof}
        Recall that by definition, $M^+\{\underline V/\lambda\}$ is the completion of $M^+[V_1,\dots,V_d] = \bigoplus_{J\in\bN^d}M^+\cdot \underline V^J$ with respect to the norm
        \[|\sum_{J\in\bN^d}a_J \underline V^J|:=\sup_{J\in\bN^d}\lambda^{|J|}|a_J|,\]
        where for any $J=(j_1,\dots,j_d)\in\bN^d$, $\underline V^J= \prod_{i=1}^dV_i^{j_i}$.
        
        Note that for any $j\geq 0$, there are rational integers $a_{j,1},\dots,a_{j,j}$ and $b_{j,1},\dots,b_{j,j}$ such that
        \[V^j = u^jj!\binom{V}{j}_u+\sum_{i=1}^ju^ia_{j,i}\cdot u^{j-i}(j-i)!\binom{V}{j-i}_u\]
        and 
        \[u^jj!\binom{V}{j}_u = V^j+\sum_{i=1}^ju^ib_{j,i}\cdot V^{j-i}.\]
        Therefore, we have rational integers $\{z_{J,I}\}_{I,J\in\bN^d,|I|<|J|}$ and $\{w_{J,I}\}_{I,J\in\bN^d,|I|<|J|}$ such that
        \[\underline V^J = u^{|J|}J!\binom{\underline V}{J}_u+\sum_{i=0}^{|J|-1}u^{|J|-i}\sum_{|I|=i}z_{J,I}\cdot u^{|I|}I!\binom{\underline V}{I}_u\]
        and
        \[u^{|J|}J!\binom{\underline V}{J}_u = \underline V^J+\sum_{i=0}^{|J|-1}u^{|J|-i}\sum_{|I|=i}w_{J,I}\cdot\underline V^I.\]
        In particular, we have $M^+[V_1,\dots,V_d] = \bigoplus_{J\in\bN^d}M^+\cdot u^{|J|}J!\binom{\underline V}{J}_u$.
        Now, equip $M^+[V_1,\dots,V_d]$ with another norm 
        \[|\sum_{J\in\bN^d}a_Ju^{|J|}J!\binom{\underline V}{J}_u|^{\prime}:=\sup_{J\in\bN^d}\lambda^{|J|}|a_J|.\]
        It suffices to show that the norm $|\cdot|^{\prime}$ conicides with $|\cdot|$ on $M^+[V_1,\dots,V_d]$. 
        Using that
        \[\begin{split}
            \sum_{J}a_J\underline V^J = & \sum_Ja_J\big(u^{|J|}J!\binom{\underline V}{J}_u+\sum_{|I|<|J|}u^{|J|-|I|}z_{J,I}\cdot u^{|I|}I!\binom{\underline V}{I}_u\big)\\
            = & \sum_{J}\big(a_J+\sum_{|I|>|J|}a_Iu^{|I|-|J|}z_{I,J}\big)u^{|J|}J!\binom{\underline V}{J}_u,
        \end{split}\]
        we obtain that 
        \[|\sum_{J}a_J\underline V^J|^{\prime} = \sup_{J}\lambda^{|J|}|a_J+\sum_{|I|>|J|}a_Iu^{|I|-|J|}z_{I,J}|\leq \sup_J\lambda^J|a_J| = |\sum_{J}a_J\underline V^J|\]
        because $|u|< 1$, $\lambda>1$ and $z_{I,J}\in\bZ$. In other words, for any $f\in M^+[V_1,\dots,V_d]$, we have 
        \[|f|^{\prime}\leq |f|.\]
        By a similar argument, we obtain that  for any $f\in M^+[V_1,\dots,V_d]$, we have 
        \[|f|\leq |f|^{\prime},\]
        which implies that $|\cdot| = |\cdot|^{\prime}$ as expected.
    \end{proof}
    The following two lemmas play the central role in this subsection.
    \begin{lem}\label{lem:Zp cohomology-principal}
        Let $B$ be a Banach $\Qp$-algebra and $B^+:=\{x\in B\mid |x|\leq 1\}$. Let $u\in B^{\times}$ with $|u|<1$. Put $G = \Zp\gamma$. Let $M$ be a Banach $B$-module which is endowed with a continuous $B$-linear $G$-action such that $M^+:=\{x\in M\mid |x|\leq 1\}$ is stable under this action. Let $\lambda_0>1$ be a real number. Suppose that for any $x\in M^+$ and any $g\in G$, we have $(g-1)(x)\in uM^+$ such that
        \begin{equation}\label{equ:Zp cohomology-principal}
            |\frac{(g-1)(x)}{u}|\leq \lambda_0^{-1}p^{-\frac{1}{p-1}}|x|.
        \end{equation}
        For any $1\leq \lambda<\lambda_0$, equip $B\{W/\lambda\}$ with the $G$-action such that $\gamma(W) = W+u$.
        \begin{enumerate}
            \item[(1)] For any $1\leq\lambda\leq\lambda_0$, the morphism $B^+\{W/\lambda\}\to B^+$ sending $W$ to $0$ induces an isometric isomorphism
            \[\rH^0(G,M^+\widehat \otimes_{B^+}B^+\{W/\lambda\})\xrightarrow{\cong}M^+.\]
            More precisely, we have
            \[\rH^0(G,M^+\widehat \otimes_{B^+}B^+\{W/\lambda\}) = \gamma^{-\frac{W}{u}}(M^+):= \{\gamma^{-\frac{W}{u}}(x):=\sum_{n\geq 0}(\gamma^{-1}-1)^n(x)\cdot\binom{W}{n}_u\mid x\in M^+\}.\]

            \item[(2)] For any $n\geq 1$, we have $\rH^n(G,\colim_{\lambda\to 1^+}M\widehat \otimes_BB\{W/\lambda\}) = 0$.
        \end{enumerate}
    \end{lem}
    Before proving the result, let us recall in $\bQ[[X,Y]]$, we always have 
    \[(1+X)^Y = \sum_{n\geq 0}X^n\binom{Y}{n}.\] 
    Letting $X = \gamma^{-1}-1$ and $Y = \frac{W}{u}$, informally, we get
    \[\gamma^{-\frac{W}{u}} = \sum_{n\geq 0}(\gamma^{-1}-1)^n\binom{\frac{W}{u}}{n} = \sum_{n\geq 0}(\gamma^{-1}-1)^n\binom{W}{n}_u.\]
    \begin{proof}
       For item (1): Note that $G$ acts on $M^+$ isometrically by \eqref{equ:Zp cohomology-principal}. By Lemma \ref{lem:good form for calculation}, we see that $M^+\widehat \otimes_{B^+} B^+\{W/\lambda\}$ consists of elements
       \[f(W) = \sum_{n\geq 0}x_n\binom{W}{n}_u\]
       satisfying the condition that
       $x_n\in u^nn!M^+$ for any $n\geq 0$ such that $\lim_{n\to \infty}|\frac{x_n}{u^nn!}|\lambda^n = 0$.
       Note that for any $n\geq 1$, 
       \[\gamma(\binom{W}{n}_u) = \binom{W}{n}_u+\binom{W}{n-1}_u.\]
       Thus, we have
       \begin{equation}\label{equ:principal term gamma-1}
           (\gamma-1)(f(W)) = \sum_{n\geq 0}\left(\gamma(x_n)-x_n+\gamma(x_{n+1})\right)\binom{W}{n}_u.
       \end{equation}
       In particular, $f$ is fixed by $G$ if and only if for any $n\geq 0$,
       \[x_{n+1} = (\gamma^{-1}-1)(x_n), \quad \forall n\geq 0.\]
       By iteration, this is equivalent to that
       \[x_n = (\gamma^{-1}-1)^n(x_0), \quad \forall n\geq 0.\]
       So the map $f(W)\mapsto f(0)$ gives rise to an injection 
       \[\rH^0(G,M^+\{W/\lambda\})\to M^+.\] 
       We claim this map is also surjective; that is, for any $x\in M^+$, 
       \[f_x(W):=\sum_{n\geq 0}(\gamma^{-1}-1)^n(x)\binom{W}{n}_u\]
       is a well-defined element in $M^+\{W/\lambda\}$. To do so, it suffices to show that for any $n\geq 0$, we have
       \[(\gamma^{-1}-1)^n(x)\in u^nn!M^+ \text{ such that }\lim_{n\to +\infty}|\frac{(\gamma^{-1}-1)^n(x)}{u^nn!}|\lambda^n = 0.\]
       Using \eqref{equ:Zp cohomology-principal}, we see that 
       \[(\gamma^{-1}-1)^n(x)\in u^nM^+ \text{ such that }|\frac{(\gamma^{-1}-1)^n(x)}{u^n}|\leq (\lambda_0^{-1}p^{-\frac{1}{p-1}})^n|x|.\]
       Combining this with the fact that $|n!|\leq p^{-\frac{n-1}{p-1}}$, we finally obtain that 
       \[|\frac{(\gamma^{-1}-1)^n(x)}{u^nn!}|\leq \lambda_0^{-n}p^{-\frac{1}{p-1}}|x|\leq 1,\]
       which amounts to that $(\gamma^{-1}-1)^n(x)\in u^nn!M^+$. Moreover, as $\lambda<\lambda_0$, we also have 
       \[|\frac{(\gamma^{-1}-1)^n(x)}{u^nn!}|\lambda^n\leq (\frac{\lambda}{\lambda_0})^np^{-\frac{1}{p-1}}|x|\to 0\text{ as }n\to +\infty,\]
       yielding that $f_x(W)\in M^+\{W/\lambda\}$ as desired.
       All in all, we have 
       \[\rH^0(G,M^+\widehat \otimes_{B^+}B^+\{W/\lambda\}) = \gamma^{-\frac{W}{u}}(M^+),\]
       which completes the proof of Item (1).

       For item (2): 
       Put $M\{W\}^{\dagger}:= \colim_{\lambda\to 1^+}M\widehat \otimes_BB\{W/\lambda\}$ for short.
       The only thing we have to prove is that
       \[\rH^1(G,M\{W\}^{\dagger}) = 0.\]
       For this purpose, it is enough to show that for any 
       \[g(W) = \sum_{n\geq 0}y_n\binom{W}{n}_u\in M\{W\}^{\dagger},\]
       there exists some 
       \[f(W) = \sum_{n\geq 0}x_n\binom{W}{n}_u\in M\{W\}^{\dagger}\]
       such that 
       \[(\gamma-1)(f(W)) = g(W).\]
       Note that there exists some $1<\lambda_1<\lambda_0$ such that $\lim_{n\to+\infty}|\frac{y_n}{u^nn!}|\lambda_1^n = 0$. 
       
       Suppose that 
       $(\gamma-1)(f(W)) = g(W)$. We deduce from \eqref{equ:principal term gamma-1} that for any $n\geq 0$,
       \[x_{n+1} = \gamma^{-1}(y_n)+(\gamma^{-1}-1)(x_n).\]
       By iteration, this is equivalent to that for any $n\geq 1$,
       \[x_{n} = (\gamma^{-1}-1)^n(x_0)+\sum_{i=0}^{n-1}(\gamma^{-1}-1)^{n-1-i}\left(\gamma^{-1}(y_i)\right).\]
       By letting $x_0 = 0$, we then have that for any $n\geq 1$,
       \[x_n = \sum_{i=0}^{n-1}(\gamma^{-1}-1)^{n-1-i}\left(\gamma^{-1}(y_i)\right).\]
       So for any $\lambda>1$, we have
       \begin{equation*}
           \begin{split}
               |\frac{x_n}{n!u^n}|\lambda^n \leq &|u^{-1}||\sum_{i=0}^{n-1}(\frac{\gamma^{-1}-1}{u})^{n-1-i}(\frac{\gamma^{-1}(y_i)}{u^ii!})\frac{i!}{n!}|\lambda^n\\
               \leq & |u^{-1}|\lambda_0p^{\frac{1}{p-1}}\cdot\sup_{0\leq i\leq n-1}((\lambda_0p^{\frac{1}{p-1}})^{-(n-i)}|\frac{y_i}{u^ii!}|\cdot|\frac{i!}{n!}|)\lambda^n\\
               \leq & |u^{-1}|\lambda_0\cdot\sup_{0\leq i\leq n-1}((\frac{\lambda}{\lambda_0})^{n-i}\cdot|\frac{y_i}{u^ii!}|\lambda^i\cdot|\frac{i!(n-i)!}{n!}|).
           \end{split}
       \end{equation*}
       Recall that for any $0\leq i\leq n$, we have
       \[|\binom{n}{i}| = p^{-\sum_{\alpha\geq 1}([\frac{n}{p^\alpha}]-[\frac{i}{p^\alpha}]-[\frac{n-i}{p^\alpha}])}\geq p^{-\log_pn} = \frac{1}{n}.\]
       So by letting $\lambda = \lambda_2\in (1,\lambda_1)$ above, we have
       \[|\frac{x_n}{n!u^n}|\lambda_2^n\leq n|u^{-1}|\lambda_0\cdot\sup_{0\leq i\leq n-1}((\frac{\lambda_2}{\lambda_0})^{n-i}\cdot|\frac{y_i}{u^ii!}|\lambda_1^i\cdot(\frac{\lambda_2}{\lambda_1})^{i}).\]
       It follows directly that the right hand side above tends to zero when $n\to +\infty$ by our choice of $\lambda_2$ and $\lambda_1$, yielding that 
       \[\lim_{n\to+\infty}|\frac{x_n}{n!u^n}|\lambda_2^n\to 0.\]
       In other words, $f(W):=\sum_{n\geq 0}\sum_{i=0}^{n-1}(\gamma^{-1}-1)^{n-1-i}\left(\gamma^{-1}(y_i)\right)\binom{W}{n}_u$ is a well-defined element in $M\{W/\lambda_2\}\subset M\{W\}^{\dagger}$ satisfying $(\gamma-1)(f(W)) = g(W)$ as expected. This completes the proof of Item (2).
    \end{proof}
    \begin{lem}\label{lem:Zp cohomology-error}
        Keep assumption and notation as in Lemma \ref{lem:Zp cohomology-principal}. Let $v\in B^{\times}$ such that $1-v\in B^{\times}$ satisfying $|1-v|<1$ and $|\frac{u}{1-v}|\leq 1$. Let $B^+_v$ be the free $B^+$-module of rank $1$ with the basis $e$ and equip $B_v$ with a continuous $G$-action determined by $\gamma(e)=ve$ and let $M^+_v:= M^+\otimes_{B^+}B^+_v$. Then for any $1\leq \lambda<\lambda_0$, we have 
        \begin{enumerate}
            \item[(1)] $\rH^0(G,M^+_v\widehat \otimes_{B^+}B^+\{W/\lambda\}) = 0$, and 

            \item[(2)] for any $n\geq 1$, $\rH^n(G,M^+_v\widehat \otimes_{B^+}B^+\{W/\lambda\})$ is killed by $v-1$ and thus by $u$.
        \end{enumerate} 
    \end{lem}
    \begin{proof}
        By assumption, we have $|v| = |v^{-1}| = 1$.
        Note that for any 
        \[f(W)=\sum_{n\geq 0}x_n\binom{W}{n}_u\in M^+\widehat \otimes_{B^+}B^+\{W/\lambda\},\]
        we have
        \begin{equation}\label{equ:Zp cohomology-error}
            (\gamma-1)(f(W)e) = \sum_{n\geq 0}(\gamma(x_n)v-x_n+\gamma(x_{n+1})v)\binom{W}{n}_ue.
        \end{equation}
        Also, for any $x\in M^+$, we have 
        \[\frac{v\gamma-1}{v-1}(x) = x+v\frac{u}{v-1}\frac{\gamma-1}{u}(x).\]
        By \eqref{equ:Zp cohomology-principal}, it follows from $|\frac{u}{v-1}|\leq 1$ that 
        \[|\frac{v\gamma-1}{v-1}(x)-x|\leq \lambda_0^{-1}p^{-\frac{1}{p-1}}|x|.\]
        As a consequence, we see that $\frac{v\gamma-1}{v-1}$ acts invertibly on $M^+$ such that for any $x\in M$, 
        \begin{equation}\label{equ:Zp cohomology-error-II}
            |\frac{v\gamma-1}{v-1}(x)| = |x|.
        \end{equation}

        For item (1): Suppose that $f(W)e\in \rH^0(G,M_v^+\widehat \otimes_{B^+}B^+\{W/\lambda\})$ for some $1\leq \lambda<\lambda_0$. By \eqref{equ:Zp cohomology-error}, for any $n\geq 0$, we have
        \[v^{n+1}x_{n+1} = -\gamma^{-1}(v\gamma-1)(v^nx_n),\]
        yielding that for any $n\geq 0$, 
        \[v^nx_n = (-\gamma^{-1})^n(v\gamma-1)^n(x_0).\]
        By \eqref{equ:Zp cohomology-error-II}, we see that 
        \[|x_0| = |\frac{(-\gamma^{-1})^n(v\gamma-1)^n(x_0)}{(v-1)^n}| = |\frac{v^nx_n}{(v-1)^n}|\leq |\frac{x_n}{u^n}| = |\frac{x_n}{u^nn!}|\lambda^n \cdot \frac{|n!|}{\lambda^n}.\]
        Letting $n\to +\infty$, as $f(W)\in M^+\widehat \otimes_{B^+}B^+\{W/\lambda\}$ which forces that $\lim_{n\to +\infty}|\frac{x_n}{u^nn!}|\lambda^n = 0$, we have $|x_0| = 0$ and thus $x_0 = 0$. As a consequence, for any $n\geq 0$, 
        \[x_n = v^{-n}(-\gamma^{-1})^n(v\gamma-1)^n(x_0) = 0,\]
        yielding that $f(W) = 0$. This completes the proof of Item (1).

        For item (2): It suffices to show that $\rH^1(G,M^+_v\widehat \otimes_{B^+}B^+\{W/\lambda\})$ is killed by $v-1$. To do so, it is enough to show that for any $g(W) = \sum_{n\geq 0}y_n\binom{W}{n}_u\in M^+\widehat \otimes_{B^+}B^+\{W/\lambda\}$, there exists an $f(W)\in M^+\widehat \otimes_{B^+}B^+\{W/\lambda\}$ such that
        \[(\gamma-1)(f(W)e) = (v-1)g(W)e.\]
        
        Suppose that $g(W) = \sum_{n\geq 0}x_n\binom{W}{n}_u$ is a solution of the above equation.
        By \eqref{equ:Zp cohomology-error}, we have that for any $n\geq 0$,
        \[(v-1)y_n = (v\gamma-1)(x_n)+\gamma(x_{n+1}).\]
        In other word, if we put $\delta:= \frac{v\gamma-1}{v-1}$, then for any $n\geq 0$,
        \[x_n = \delta^{-1}(y_n)-(v-1)^{-1}\delta^{-1}\gamma(x_{n+1}).\]
        Note that for any $i\geq 0$, by \eqref{equ:Zp cohomology-error-II}, we have
        \[|(v-1)^{-i}\delta^{-(i+1)}\gamma^i(y_{n+i})| = |\frac{u^{n+i}}{(v-1)^i}\frac{y_{n+i}}{u^{n+i}}|\leq |\frac{u}{v-1}|^i\cdot|u^n|\cdot|\frac{y_{n+i}}{u^{n+i}(n+i)!}|\lambda^{n+i}\cdot\frac{|(n+i)!|}{\lambda^{n+i}}\]
        and that
        \[|\frac{(v-1)^{-i}\delta^{-(i+1)}\gamma^i(y_{n+i})}{u^nn!}|\leq\cdot |\frac{y_{n+i}}{u^{n+i}(n+i)!}|\lambda^{n+i}\cdot\frac{|i!\binom{n+i}{i}|}{\lambda^{n+i}}.\]
        As $g(W) = \sum_{n\geq 0}y_n\binom{W}{n}_u\in M^+\widehat \otimes_{B^+}B^+\{W/\lambda_1\}$, we see that for any $n\geq 0$,
        \[|\frac{y_n}{u^nn!}|\leq 1\text{ such that }\lim_{n\to+\infty}|\frac{y_n}{u^nn!}| = 0,\]
        yielding that for any $i\geq 0$
        \[|\frac{(v-1)^{-i}\delta^{-(i+1)}\gamma^i(y_{n+i})}{u^nn!}|\leq 1\text{ such that }\lim_{i\to+\infty}|(v-1)^{-i}\delta^{-(i+1)}\gamma^i(y_{n+i})| = 0.\]
        By Lemma \ref{lem:good form for calculation}, for any $n\geq 0$, 
        \[x_n:=\sum_{i\geq 0}(-1)^i(v-1)^{-i}\delta^{-(i+1)}\gamma^i(y_{n+i})\]
        is a well-defined element in $M^+$ such that
        \[|\frac{x_n}{u^nn!}|\lambda^n\leq \sup_{i\geq 0}\{|\frac{y_{n+i}}{u^{n+i}(n+i)!}|\lambda^{n+i}\cdot\frac{|i!\binom{n+i}{i}|}{\lambda^i}\}\leq \cdot\sup_{i\geq 0}\{|\frac{y_{n+i}}{u^{n+i}(n+i)!}|\lambda^{n+i}\}.\]
        So $f(W):=\sum_{n\geq 0}x_n\binom{W}{n}_u$ is a well-defined element in $M^+\widehat \otimes_{B^+}B^+\{W/\lambda\}$ for $x_n$'s given above. By construction, we have that
        \[(\gamma-1)(f(W)e) = (v-1)g(W)e,\]
        yielding that $\rH^1(G,M_v\widehat \otimes_BB\{W/\lambda\})$ is killed by $v-1$ as desired.
    \end{proof}

    \begin{prop}\label{prop:From representation to connection}
        Let $I = [r,s]$ be a closed sub-interval of $(\frac{1}{p-1},p)$ and $u = [\epsilon^{\frac{1}{p}}]-1$ as before.
        Fix a $\lambda_0\in (p^{\frac{1}{p-1}},|\xi|_I^{-1})$. For any $M_{\infty}\in \Rep^{\lambda_0}_{\Gamma}(\CCppsi)$ of rank $r$, let $M\in \Rep^{\lambda_0}_{\Gamma}(\Cppsi)$ be the representation such that $M_{\infty} = M\otimes_{\Cppsi}\CCppsi$ in the sense of Theorem \ref{thm:decompletion}.
        \begin{enumerate}
            \item[(1)] For any $1\leq \lambda<\lambda_0p^{-\frac{1}{p-1}}$, the morphism
            \[\rR\Gamma(\Gamma,M\otimes_{\Cppsi}\Cppsi\{\underline W/\lambda\}_{r,d,P})\to\rR\Gamma(\Gamma,M_{\infty}\otimes_{\CCppsi}\CCppsi\{\underline W/\lambda\}_{r,d,P})\]
            induced by the natural inclusion
            \[M\otimes_{\Cppsi}\Cppsi\{\underline W/\lambda\}_{r,d,P}\subset M\otimes_{\Cppsi}\CCppsi\{\underline W/\lambda\}_{r,d,P} =  M_{\infty}\otimes_{\CCppsi}\CCppsi\{\underline W/\lambda\}_{r,d,P}\]
            identifies the former with a direct summand of the latter whose complement is concentrated in degree $\geq 1$ and killed by $u$.

            \item[(2)] For any $1\leq \lambda_2<\lambda_1<\lambda_0p^{-\frac{1}{p-1}}$, the natural morphism
            \[M\otimes_{\Cppsi}\Cppsi\{\underline W/\lambda_2\}_{r,d,P}\subset M\otimes_{\Cppsi}\Cppsi\{\underline W/\lambda_1\}_{r,d,P}\]
            induces an identification
            \[D_{\psi}^{I,+}(M_{\infty}):=\rH^0(\Gamma,M\otimes_{\Cppsi}\Cppsi\{\underline W/\lambda_2\}_{r,d,P}) = \rH^0(\Gamma,M\otimes_{\Cppsi}\Cppsi\{\underline W/\lambda_2\}_{r,d,P}).\]
            Moreover, for $1\leq \lambda<\lambda_0$, the connection $\id_M\otimes\nabla$ on $M\otimes_{\Cppsi}\Cppsi\{\underline W/\lambda\}_{r,d,P}$ induces a flat connection
            \[\nabla_{D_{\psi}^{I,+}(M)}:D_{\psi}^{I,+}(M)\to D_{\psi}^{I,+}(M)\otimes_{\wtr_Z^{I,+}}\Omega^{1,\log}_{\wtr_Z^{I,+}}\{-1\}\]
            such that $(D_{\psi}^{I,+}(M),\nabla_{D_{\psi}^{I,+}(M)})\in \MIC^{\lambda_0}(\wtr_Z^{I,+})$.

            \item[(3)] For any $n\geq 1$, the cohomological group $\rH^n(\Gamma,M\otimes_{\Cppsi}\Cppsi\{\underline W\}^{\dagger}_{r,d,P})$ is of $p^{\infty}$-torsion.
        \end{enumerate}
    \end{prop}
    \begin{proof}
        Note that, by Lemma \ref{lem:lambda(M)}, for any $x\in M$ and any $\gamma\in \Gamma$, we have
        \begin{equation}\label{equ:From representation to connection}
            (\gamma-1)(x)\in uM \text{ such that }|\frac{(\gamma-1)}{u}(x)|\leq \lambda_0^{-1}|x|.
        \end{equation}

        For item (1): Using the $\Gamma$-equivariant decomposition \eqref{equ:decomposition for M}
        \[M_{\infty} = \widehat \oplus_{\underline \alpha\in J_{r,d}}M\cdot [(\underline T^{\flat})^{\underline \alpha}],\]
        we have a $\Gamma$-equivariant decomposition
        \[M_{\infty}\otimes_{\CCppsi}\CCppsi\{\underline W/\lambda\}_{r,d,P} = M\otimes_{\Cppsi}\Cppsi\{\underline W/\lambda\}_{r,d,P}\oplus\bigoplus_{0\neq \underline \alpha\in J_{r,d}}(M\cdot [(\underline T^{\flat})^{\underline \alpha}])\otimes_{\Cppsi}\Cppsi\{\underline W/\lambda\}_{r,d,P}.\]
        It is enough to show that for any non-zero $\underline \alpha\in J_{r,d}$,
        \[\rR\Gamma(\Gamma,(M\cdot [(\underline T^{\flat})^{\underline \alpha}])\otimes_{\Cppsi}\Cppsi\{\underline W/\lambda\}_{r,d,P})\]
        is concentrated in degree $\geq 1$ and killed by $u$. To see this, without loss of generality, we may assume that $\underline \alpha=(\alpha_0,\alpha_1,\dots,\alpha_d)$ satisfies $\alpha_0=0$ and $\alpha_1\neq 0$. We will consider the $\Gamma$-equivariant isomorphism
        \begin{equation}\label{equ:isomorphism of period ring}
            (\Cppsi\{\underline W/\lambda\}_{r,d,P},\nabla)\cong (\Cppsi\{W_1/\lambda,\dots,W_d/\lambda\},-\sum_{i=1}^d\frac{\partial}{\partial W_i}\otimes\frac{u}{t}\dlog T_i)
        \end{equation}
        in terms of Example \ref{exam:Apply to Gamma torsor}
        and use the isomorphism \eqref{equ:gamma group-II}.
        Note that $v:=[\epsilon^{\alpha_1}]$ is a unit in $\CCtilde^I(Z)$ satisfying all the conditions in Lemma \ref{lem:Zp cohomology-error} and we have 
        \[\gamma_1([(\underline T^{\flat})^{\underline \alpha}]) = [\epsilon^{\alpha_1}]\cdot [(\underline T^{\flat})^{\underline \alpha}].\]
        So by Lemma \ref{lem:Zp cohomology-error}, we have that 
        \[\rR\Gamma(\Zp\gamma_1,(M\cdot [(\underline T^{\flat})^{\underline \alpha}])\otimes_{\Cppsi}\Cppsi\{W_1/\lambda,\dots,W_d/\lambda\})\]
        is concentrated in degree $\geq 1$ and killed by $u$. Then Item (1) follows from the Hochschild--Serre spectral sequence.

        For item (2): Using \eqref{equ:gamma group-II} and \eqref{equ:isomorphism of period ring}, it follows from Lemma \ref{lem:Zp cohomology-principal}(1) that for any $1\leq \lambda<\lambda_0p^{-\frac{1}{p-1}}$, 
        \begin{equation}\label{equ:explicit D(M)}
            \begin{split}
                &\rH^0(\Gamma,M\otimes_{\Cppsi}\Cppsi\{\underline W/\lambda\}_{r,d,P})\\
                = & \rH^0(\Gamma,M\otimes_{\Cppsi}\Cppsi\{W_1/\lambda,\dots,W_d/\lambda\})\\
                = & (\prod_{i=1}^d\gamma_i^{-\frac{W_i}{u}})(M)\\
                := & \{\sum_{n_1,\dots,n_d\geq 0}(\gamma_1^{-1}-1)^{n_1}\cdots (\gamma_d^{-1}-1)^{n_d}(x)\cdot\binom{W_1}{n_1}_u\cdots \binom{W_d}{n_d}_u\mid x\in M\},
            \end{split}
        \end{equation}
        which is clearly independent of the choice of $\lambda$ and of rank $r$ as desired. 

        Now, as $\nabla_{D^{I,+}_{\psi}(M_{\infty})}$ is the restriction of $-\id_M\otimes\sum_{i=1}^d\frac{\partial}{\partial W_i}\otimes\frac{u}{t}\dlog T_i$, by \eqref{equ:explicit D(M)}, we have that for any $x\in M$,
        \[\nabla_{D^{I,+}_{\psi}(M_{\infty})}\big((\prod_{i=1}^d\gamma_i^{-\frac{W_i}{u}})(x)\big) = \sum_{i=1}^d\nabla_i\big((\prod_{i=1}^d\gamma_i^{-\frac{W_i}{u}})(x)\big)\otimes\frac{u}{t}\dlog T_i,\]
        where for any $1\leq i\leq d$,
        \begin{equation}\label{equ:explicit nabla}
            \begin{split}
                \nabla_i\big((\prod_{i=1}^d\gamma_i^{-\frac{W_i}{u}})(x)\big) = & \frac{\log\gamma_i}{u}\big((\prod_{i=1}^d\gamma_i^{-\frac{W_i}{u}})(x)\big)\\
                = & u^{-1}\sum_{n\geq 1}\frac{(-1)^n}{n}(\gamma_i^{-1}-1)^n\big((\prod_{i=1}^d\gamma_i^{-\frac{W_i}{u}})(x)\big).
            \end{split}
        \end{equation}
        It follows from \eqref{equ:From representation to connection} that 
        \[| \nabla_i\big((\prod_{i=1}^d\gamma_i^{-\frac{W_i}{u}})(x)\big)|\leq \lambda_0^{-1}|x|,\]
        yielding that $(D^{I,+}_{\psi}(M_{\infty}),\nabla_{D^{I,+}_{\psi}(M_{\infty})})\in \MIC^{\lambda_0}(\wtr_Z^{I,+})$ as desired.

        For item (3): It suffices to show that
        \[\rH^n(\Gamma,M[\frac{1}{p}]\otimes_{\Cpsi}\Cpsi\{\underline W\}^{\dagger}_{r,d,P}) = \rH^n(\Gamma,\colim_{\lambda\to 1^+}M[\frac{1}{p}]\otimes_{\Cpsi}\Cpsi\{W_1/\lambda,\dots,W_d/\lambda\}) = 0,\]
        which follows from Lemma \ref{lem:Zp cohomology-principal}(2) directly.
    \end{proof}
    \begin{cor}\label{cor:From representation to connection}
        Let $I = [r,s]$ be a closed sub-interval of $(\frac{1}{p-1},p)$ and $u = [\epsilon^{\frac{1}{p}}]-1$ as before.
        Fix a $\lambda_0\in (p^{\frac{1}{p-1}},|\xi|_I^{-1})$. For any $M_{\infty}\in \Rep^{\lambda_0}_{\Gamma}(\CCpsi)$ of rank $r$, let $M\in \Rep^{\lambda_0}_{\Gamma}(\Cpsi)$ be the representation such that $M_{\infty} = M\otimes_{\Cpsi}\CCpsi$ in the sense of Theorem \ref{thm:decompletion}.
        \begin{enumerate}
            \item[(1)] For any $1\leq \lambda<\lambda_0p^{-\frac{1}{p-1}}$, the natural inclusion
            \[M\otimes_{\Cpsi}\Cpsi\{\underline W/\lambda\}_{r,d,P}\subset M\otimes_{\Cpsi}\CCpsi\{\underline W/\lambda\}_{r,d,P} =  M_{\infty}\otimes_{\CCpsi}\CCpsi\{\underline W/\lambda\}_{r,d,P}\]
            induces a quasi-isomorphism
            \[\rR\Gamma(\Gamma,M\otimes_{\Cpsi}\Cpsi\{\underline W/\lambda\}_{r,d,P})\to\rR\Gamma(\Gamma,M_{\infty}\otimes_{\CCpsi}\CCpsi\{\underline W/\lambda\}_{r,d,P})\]
            such that $\rR\Gamma(\Gamma,M\otimes_{\Cpsi}\Cpsi\{\underline W/\lambda\}_{r,d,P})$ is a direct summand of $\rR\Gamma(\Gamma,M_{\infty}\otimes_{\CCpsi}\CCpsi\{\underline W/\lambda\}_{r,d,P})$.

            \item[(2)] For any $1\leq \lambda_2<\lambda_1<\lambda_0p^{-\frac{1}{p-1}}$, the natural morphism
            \[M\otimes_{\Cpsi}\Cpsi\{\underline W/\lambda_2\}_{r,d,P}\subset M\otimes_{\Cpsi}\Cpsi\{\underline W/\lambda_1\}_{r,d,P}\]
            induces an identification
            \[D_{\psi}^I(M_{\infty}):=\rH^0(\Gamma,M\otimes_{\Cpsi}\Cpsi\{\underline W/\lambda_2\}_{r,d,P}) = \rH^0(\Gamma,M\otimes_{\Cpsi}\Cpsi\{\underline W/\lambda_2\}_{r,d,P}).\]
            Moreover, for $1\leq \lambda<\lambda_0$, the connection $\id_M\otimes\nabla$ on $M\otimes_{\Cpsi}\Cpsi\{\underline W/\lambda\}_{r,d,P}$ induces a flat connection
            \[\nabla_{D_{\psi}^I(M)}:D_{\psi}^I(M)\to D_{\psi}^I(M)\otimes_{\wtr_Z^I}\Omega^{1,\log}_{\wtr_Z^I}\{-1\}\]
            such that $(D_{\psi}^I(M),\nabla_{D_{\psi}^I(M)})\in \MIC^{\lambda_0}(\wtr_Z^I)$.

            \item[(3)] The $\Gamma$-cohomology $\rR\Gamma(\Gamma,M\otimes_{\Cpsi}\Cpsi\{\underline W\}_{r,d,P}^{\dagger})$ is concentrated in degree $0$ such that
            \[D_{\psi}^I(M_{\infty}) = \rH^0(\Gamma,M\otimes_{\Cpsi}\Cpsi\{\underline W\}_{r,d,P}^{\dagger}).\] 
            Moreover, there exists a quasi-isomorphism
            \[\rR\Gamma(\Gamma,M_{\infty})\simeq \rD\rR(D_{\psi}^I(M_{\infty}),\nabla_{D_{\psi}^I(M_{\infty})})\]
            which is functorial in $M_{\infty}$.
        \end{enumerate}
    \end{cor}
    \begin{proof}
        Note that $u$ is invertible in $\CCtilde^I(Z)$, Items (1) and (2) together with the first part of Item (3) follow from Proposition \ref{prop:From representation to connection} directly. As a consequence, the natural map
        \[D_{\psi}^I(M_{\infty}) = \rH^0(\Gamma,M_{\infty}\otimes_{\CCpsi}\CCpsi\{\underline W\}^{\dagger}_{r,d,P})\to \rR\Gamma(\Gamma,M_{\infty}\otimes_{\CCpsi}\CCpsi\{\underline W\}^{\dagger}_{r,d,P})\]
        is a quasi-isomorphism. Therefore, the natural map
        \[\rD\rR(D_{\psi}^I(M_{\infty}),\nabla_{D_{\psi}^I(M_{\infty})})\to \rR\Gamma(\Gamma,\rD\rR(M_{\infty}\otimes_{\CCpsi}\CCpsi\{\underline W\}^{\dagger}_{r,d,P},\id_{M_{\infty}}\otimes\nabla))\]
        is a quasi-isomorphism. On the other hand, by Poincare's Lemma \ref{thm:Poincare Lemma for OC^I}, the natural map
        \[M_{\infty}\to \rD\rR(M_{\infty}\otimes_{\CCpsi}\CCpsi\{\underline W\}^{\dagger}_{r,d,P},\id_{M_{\infty}}\otimes\nabla)\]
        is a quasi-isomorphism, yielding a quasi-isomorphism
        \[\rR\Gamma(\Gamma,M_{\infty})\xrightarrow{\simeq}\rR\Gamma(\Gamma,\rD\rR(M_{\infty}\otimes_{\CCpsi}\CCpsi\{\underline W\}^{\dagger}_{r,d,P},\id_{M_{\infty}}\otimes\nabla)).\]
        Finally, we get a quasi-isomorphism
        \[\rD\rR(D_{\psi}^I(M_{\infty}),\nabla_{D_{\psi}^I(M_{\infty})})\simeq\rR\Gamma(\Gamma,M_{\infty})\]
        which is functorial in $M_{\infty}$ as desired.
    \end{proof}

\subsection{de Rham cohomology of small connections}\label{ssec:dR cohomology for small of small connection}
  We begin with a general fact.
  \begin{lem}\label{lem:partial cohomology}
      Let $B$ be a Banach $\Qp$-algebra and $B^+:=\{x\in B\mid |x|\leq 1\}$. Let $u\in B^{\times}$ with $|u|\leq 1$. Let $D$ be a Banach $B$-module with $D^+=\{x\in D\mid |x|\leq 1\}$. Suppose that $D^+$ is endowed with a $B^+$-linear endomorphism $\phi:D\to D$ such that there exists a $\lambda_0>1$ such that for any $x\in D^+$,
      \[|\phi(x)|\leq \lambda_0^{-1}p^{-\frac{1}{p-1}}|x|.\]
      For any $1\leq \lambda<\lambda_0$, let $\partial = -\frac{\partial}{\partial W}$ on $B^+\{W/\lambda\}$ and define 
      \[\partial_{D^+} = \phi\otimes\id+\id_{D^+}\otimes\partial:D^+\widehat \otimes_{B^+}B^+\{W/\lambda\}\to D^+\widehat \otimes_{B^+}B^+\{W/\lambda\}.\]
      \begin{enumerate}
          \item[(1)] The map $B^+\{W/\lambda\}\xrightarrow{W\mapsto 0}B^+$ induces an isometric isomorphism of $B^+$-modules
          \[(D^+\widehat \otimes_{B^+}B^+\{W/\lambda\})^{\partial_{D^+}=0}\xrightarrow{\cong} D^+.\]
          More precisely, we have
          \[(D^+\widehat \otimes_{B^+}B^+\{W/\lambda\})^{\partial_{D^+}=0} = \exp(\phi W)(D^+):=\{\exp(\phi W)(x):=\sum_{n\geq 0}\frac{\phi^n}{n!}(x)\mid x\in D^+\}.\]

          \item[(2)] The cokernel of $\partial_{D^+}$ is of $p^{\infty}$-torsion. Equivalently, the morphism 
          \[\partial_{D}:=\phi\otimes\id+\id_D\otimes\partial: \colim_{\lambda\to 1^+}D\widehat \otimes_BB\{W/\lambda\}\to \colim_{\lambda\to 1^+}D\widehat \otimes_BB\{W/\lambda\}\]
          is surjective.
      \end{enumerate}
  \end{lem}
  \begin{proof}
      Note that for any 
      \[f(W) = \sum_{n\geq 0}x_nW^n\in D^+\{W/\lambda\}:=  D^+\widehat \otimes_{B^+}B^+\{W/\lambda\}\] 
      with $x_n\in D$, we have
      \begin{equation}\label{equ:partial cohomology}
          \partial_{D^+}(f(W)) = \sum_{n\geq 0}(\phi(x_n)-(n+1)x_{n+1})W^n.
      \end{equation}
      
      For item (1): By \eqref{equ:partial cohomology}, $f(W)\in \ker(\partial_{D^+})$ if and only if for any $n\geq 0$, 
      \[x_{n+1} = \frac{\phi(x_n)}{n},\]
      which, by iteration, amounts to that for any $n\geq 0$,
      \[x_n = \frac{\phi^n}{n!}(x_0).\]
      So the map $f(W)\mapsto f(0) = x_0$ induces an injection $\ker(\partial_{D^+})\to D^+$. On the other hand, for any $x\in D^+$, we have that
      \[|\frac{\phi^n}{n!}(x)|\lambda^n \leq (\frac{\lambda}{\lambda_0})^np^{-\frac{n}{p-1}}|\frac{1}{n!}|\cdot|x|\leq (\frac{\lambda}{\lambda_0})^n|x|.\]
      yielding that
      \[\frac{\phi^n}{n!}(x)\in D^+ \text{ and } \lim_{n\to +\infty}|\frac{\phi^n}{n!}(x)|\lambda^n = 0.\]
      As a consequence,
      \[f_x(W):=\sum_{n\geq 0}\frac{\phi^n}{n!}(x)W^n\]
      is a well-defined element in $D^+\{W/\lambda\}^{\dagger}$ with $|f_x(W)| = |x|$, which is killed by $\partial_{D^+}$. So the morphism $f(W)\mapsto f(0) = x_0$ gives rise to an isometric isomorphism $\ker(\partial_{D^+})\cong D^+$. It is clear that
      \[(D^+\widehat \otimes_{B^+}B^+\{W/\lambda\})^{\partial_{D^+}=0} = \exp(\phi W)(D^+)\]
      as expected.

      For item (2): Put $D\{W/\lambda\}:= D\widehat \otimes_BB\{W/\lambda\}$ and let $g(W) = \sum_{n\geq 0}y_nW^n\in D\{W/\lambda_1\}$ for some $1<\lambda_1<\lambda_0$. By \eqref{equ:partial cohomology}, $\partial_D(f(W)) = g(W)$ if and only if for any $n\geq 0$, 
      \[\phi(x_n)-(n+1)x_{n+1} = y_n.\]
      By iteration, this amounts to that for any $n\geq 0$, 
      \[(n+1)!x_{n+1} = -n!y_n+\phi(n!x_n) = \cdots = -\sum_{i=0}^n\phi^{n-i}(i!y_i)+\phi^{n+1}(x_0).\]
      By letting $x_0 = 0$, we then have for any $n\geq 1$,
      \[x_n = -\sum_{i=0}^{n-1}\frac{i!}{n!}\phi^{n-1-i}(y_i).\]
      Now, for any $\lambda_2\in (1,\lambda_1)$, we have
      \[\begin{split}
          |x_n|\lambda_2^n = & |\sum_{i=0}^{n-1}\frac{i!}{n!}\phi^{n-1-i}(y_i)|\lambda_2^n\\
          \leq &\lambda_2\sup_{0\leq i\leq n-1}\{|\frac{i!}{n!}|p^{-\frac{n-1-i}{p-1}}\cdot (\frac{\lambda_2}{\lambda_0})^{n-1-i}\cdot |y_i|\lambda_1^i\cdot(\frac{\lambda_2}{\lambda_1})^i\}\\
          \leq & \lambda_2(\frac{\lambda_2}{\lambda_1})^{n-1}\sup_{0\leq i\leq n-1}\{p^{\frac{i}{p-1}}|i!|\cdot|y_i|\lambda_1^i.\}
      \end{split}\]
      Note that for any $n\geq 1$, we have 
      \[|n!|\leq p^{-\frac{n-(p-1)\frac{\log(pn)}{\log p}}{p-1}} = p^{-\frac{n}{p-1}}pn.\]
      So we have that
      \[|x_n|\lambda_2^n\leq p(n-1)\lambda_2(\frac{\lambda_2}{\lambda_1})^{n-1}\sup_{0\leq i\leq n-1}\{|y_i|\lambda_1^i\}.\]
      As $g(W) = \sum_{n\geq 0}y_nW^n\in D\{W/\lambda_1\}$, we have $\lim_{n\to +\infty}|y_n|\lambda_1^n = 0$, yielding that
      \[\lim_{n\to +\infty}|x_n|\lambda_2^n = 0\]
      because $\lambda_2<\lambda_1$.
      So for $x_n$'s constructed above, $f(W) = \sum_{n\geq 0}x_nW^n$ is a well-defined element in $D\{W/\lambda_2\}$ such that $\partial_D(f(W)) = g(W)$.
      This completes the proof of Item (2).
  \end{proof}
  \begin{prop}\label{prop:From connection to representation}
      Let $I=[r,s]$ be a closed sub-interval of $(\frac{1}{p-1},p)$ and $u = [\epsilon^{\frac{1}{p}}]-1$ as before. Fix a $\lambda_0\in (\frac{1}{p-1},|\xi|_I^{-1})$. For any $(D,\nabla_D)\in \MIC^{\lambda_0}(\wtr_Z^{I,+})$ of rank $r$ and any $1\leq \lambda<\lambda_0p^{-\frac{1}{p-1}}$, define
      \[\rd_D:=\nabla_D\otimes\id+\id_D\otimes\nabla: D\otimes_{\wtr_Z^{I,+}}\CCppsi\{W/\lambda\}_{r,d,P}\to D\otimes_{\wtr_Z^{I,+}}\CCppsi\{W/\lambda\}_{r,d,P}\otimes_{\wtr_Z^{I,+}}\Omega^{1,\log}_{\wtr_Z^{I,+}}\{-1\}\]
      and then we have $\rd_D\wedge\rd_D = 0$.
      \begin{enumerate}
          \item[(1)] Put
          
          \[M_{\infty}^+(D,\nabla_D):= ( D\otimes_{\wtr_Z^{I,+}}\CCppsi\{W/\lambda\}_{r,d,P})^{\rd_D = 0}\]
          and then it is independent of the choice of $\lambda\in [1,\lambda_0p^{-\frac{1}{p-1}})$ and of rank $r$ over $\CCppsi$. Equip $M_{\infty}^+(D,\nabla_D)$ with the $\Gamma$-action induced from $D\otimes_{\wtr_Z^{I,+}}\CCppsi\{W/\lambda\}_{r,d,P}$ and then $M_{\infty}^+(D,\nabla_D)$ is an object in $\Rep_{\Gamma}^{\lambda_0}(\CCppsi)$.

          \item[(2)] For any $n\geq 1$, the cohomological group $\rH^n_{\dR}(D\otimes_{\wtr_Z^{I,+}}\CCppsi\{\underline W\}^{\dagger}_{r,d,P},\rd_D)$ is of $p^{\infty}$-torsion; that is, we have 
          \[\rH^n_{\dR}(D[\frac{1}{p}]\otimes_{\wtr_Z^I}\CCpsi\{\underline W\}^{\dagger}_{r,d,P},\rd_D) = 0.\]
      \end{enumerate}
  \end{prop}
  \begin{proof}
      We still consider the $\Gamma$-equivariant isomorphism \eqref{equ:isomorphism of period ring}
      \[(\Cppsi\{\underline W/\lambda\}_{r,d,P},\nabla)\cong (\Cppsi\{W_1/\lambda,\dots,W_d/\lambda\},-\sum_{i=1}^d\frac{\partial}{\partial W_i}\otimes\frac{u}{t}\dlog T_i)\]
      in terms of Example \ref{exam:Apply to Gamma torsor} and use the isomorphism \eqref{equ:gamma group-II} \[\Gamma\cong \oplus_{i=1}^d\Zp\gamma_i\] in the proof of Proposition \ref{prop:From representation to connection}. Also, we write 
      \[\nabla_D = \sum_{i=1}^d\nabla_i\otimes\frac{u}{t}\dlog T_i\] 
      via the isomorphism
      \[\Omega^{1,\log}_{\wtr_Z^{I,+}}\{-1\}\cong \oplus_{i=1}^d\wtr_Z^{I,+}\cdot\frac{u}{t}\dlog T_i.\] By Lemma \ref{lem:lambda(D)}, for any $x\in D$ and any $1\leq i\leq d$, we have
      \begin{equation}\label{equ:From connection to representation}
          |\nabla_i(x)|\leq\lambda_0^{-1}|x|.
      \end{equation}
      For any $1\leq i\leq d$, put $\partial_i:=\nabla_i-\frac{\partial}{\partial W_i}$ and then we have
      \[\rd_D=\sum_{i=1}^d\partial_i\otimes\frac{u}{t}\dlog T_i.\]
      As $(D,\nabla_D)$ is flat, we see that $\nabla_i$ commutes with $\nabla_j$, yielding that $\partial_i$ commutes with $\partial_j$ for any $1\leq i,j\leq d$. In other words, we have $\rd_D\wedge\rd_D = 0$ as desired. Also, for any $1\leq \lambda<\lambda_0p^{-\frac{1}{p-1}}$, the de Rham complex
      \[\rD\rR(D\otimes_{\wtr_Z^{I,+}}\CCppsi\{\underline W\}^{\dagger}_{r,d,P} = D\otimes_{\wtr_Z^{I,+}}\CCppsi\{W_1/\lambda,\dots,W_d/\lambda\},\rd_D)\]
      coincides with the Koszul complex
      \[\rm{Kos}(D\otimes_{\wtr_Z^{I,+}}\CCppsi\{W_1/\lambda,\dots,W_d/\lambda\};\partial_1,\dots,\partial_d).\]
      Let $e_1,\dots,e_r$ be a $\wtr_Z^{I,+}$-basis of $D$ so we have 
      \[D\otimes_{\wtr_Z^{I,+}}\CCppsi\{W_1/\lambda,\dots,W_d/\lambda\} = \oplus_{i=1}^d\CCppsi\{W_1/\lambda,\dots,W_d/\lambda\}\cdot e_i.\]

      For item (1): Note that for any $1\leq \lambda<\lambda_0p^{-\frac{1}{p-1}}$, we have 
      \[\begin{split}
          ( D\otimes_{\wtr_Z^{I,+}}\CCppsi\{W/\lambda\}_{r,d,P})^{\rd_D = 0} =&\rH^0(\rm{Kos}(D\otimes_{\wtr_Z^{I,+}}\CCppsi\{W_1/\lambda,\dots,W_d/\lambda\};\partial_1,\dots,\partial_d)) \\
          =&( D\otimes_{\wtr_Z^{I,+}}\CCppsi\{W/\lambda\}_{r,d,P})^{\partial_1=\cdots=\partial_d=0}.
      \end{split}\]
      By repeatedly applying Lemma \ref{lem:partial cohomology} (1), we see that
      \begin{equation}\label{equ:explicit M(D)}
          \begin{split}
              M_{\infty}^+(D,\nabla_D)
              = & \big(\prod_{i=1}^d\exp(\nabla_iW_i)\big)(D)\\
              := & \{\sum_{j=1}^da_j\big(\prod_{i=1}^d\exp(\nabla_iW_i)\big)(e_j)\mid a_1,\dots,a_d\in\CCppsi\},
          \end{split}
      \end{equation}
      where for any $1\leq j\leq d$, 
      \[\big(\prod_{i=1}^d\exp(\nabla_iW_i)\big)(e_j):=\sum_{n_1,\dots,n_d\geq 0}\frac{\nabla_1^{n_1}}{n_1!}\cdots\frac{\nabla_d^{n_d}}{n_d!}(e_j)W_1^{n_1}\cdots W_d^{n_d}.\]
      As a consequence, we see that $M_{\infty}^+(D,\nabla_D)$ is finite free of rank $r$ over $\CCppsi$ and is indenpendent of the choice of $\lambda$. Moreover, as $\gamma_i(W_j) = W_j+\delta_{ij}u$, cf. \eqref{equ:gamma acts on W-II}, we see that for any $1\leq i,j\leq d$, we have
      \begin{equation}\label{equ:explicit gamma action}
          \begin{split}
              \gamma_i(\big(\prod_{i=1}^d\exp(\nabla_iW_i)\big)(e_j)) = & \exp(u\nabla_i)(\big(\prod_{i=1}^d\exp(\nabla_iW_i)\big)(e_j))\\
              = &\sum_{n\geq 0}\frac{u^n\nabla_i^n}{n!}(\big(\prod_{i=1}^d\exp(\nabla_iW_i)\big)(e_j)).
          \end{split}
      \end{equation}
      Therefore, we have
      \[(\gamma_i-1)(\big(\prod_{i=1}^d\exp(\nabla_iW_i)\big)(e_j))\in uM_{\infty}^+(D,\nabla_D)\]
      and it follows from \eqref{equ:From connection to representation} that
      \[|\frac{(\gamma_i-1)}{u}(\big(\prod_{i=1}^d\exp(\nabla_iW_i)\big)(e_j))|\leq\lambda_0^{-1}|\big(\prod_{i=1}^d\exp(\nabla_iW_i)\big)(e_j)|.\]
      In other words, $M_{\infty}^+(D,\nabla_D)\in \Rep_{\Gamma}^{\lambda_0}(\CCppsi)$ as expected.

      For item (2): It suffices to show that for any $n\geq 1$,
      \[\begin{split}
          &\rH^n(\colim_{\lambda\to 1^+}\rm{Kos}(D[\frac{1}{p}]\otimes_{\wtr_Z^I}\CCpsi\{W_1/\lambda,\dots,W_d/\lambda\};\partial_1,\dots,\partial_d)) \\
          =&\rH^n(\rm{Kos}(\colim_{\lambda\to 1^+}D[\frac{1}{p}]\otimes_{\wtr_Z^I}\CCpsi\{W_1/\lambda,\dots,W_d/\lambda\};\partial_1,\dots,\partial_d))\\
          =& 0.
      \end{split}\]
      By Lemma \ref{lem:partial cohomology}, we get quasi-isomorphisms
      \[\begin{split}
          &\rm{Kos}(\colim_{\lambda\to 1^+}D[\frac{1}{p}]\otimes_{\wtr_Z^I}\CCpsi\{W_1/\lambda,\dots,W_d/\lambda\};\partial_1,\dots,\partial_d)\\
          \simeq &\rm{Kos}(\big(\colim_{\lambda\to 1^+}D[\frac{1}{p}]\otimes_{\wtr_Z^I}\CCpsi\{W_1/\lambda,\dots,W_d/\lambda\}\big)^{\partial_d = 0};\partial_1,\dots,\partial_{d-1})\\
          \simeq &\rm{Kos}(\colim_{\lambda\to 1^+}\big(D[\frac{1}{p}]\otimes_{\wtr_Z^I}\CCpsi\{W_1/\lambda,\dots,W_d/\lambda\}\big)^{\partial_d = 0};\partial_1,\dots,\partial_{d-1})\\
          \simeq &\rm{Kos}(\colim_{\lambda\to 1^+}D[\frac{1}{p}]\otimes_{\wtr_Z^I}\CCpsi\{W_1/\lambda,\dots,W_{d-1}/\lambda\};\partial_1,\dots,\partial_{d-1}).
      \end{split}\]
      Thus one can conclude by induction on $d$.
  \end{proof}  
  \begin{cor}\label{cor:From connection to representation}
      Let $I=[r,s]$ be a closed sub-interval of $(\frac{1}{p-1},p)$ and $u = [\epsilon^{\frac{1}{p}}]-1$ as before. Fix a $\lambda_0\in (\frac{1}{p-1},|\xi|_I^{-1})$. For any $(D,\nabla_D)\in \MIC^{\lambda_0}(\wtr_Z^I)$ of rank $r$ and any $1\leq \lambda<\lambda_0p^{-\frac{1}{p-1}}$, define
      \[\rd_D:=\nabla_D\otimes\id+\id_D\otimes\nabla: D\otimes_{\wtr_Z^I}\CCpsi\{W/\lambda\}_{r,d,P}\to D\otimes_{\wtr_Z^I}\CCpsi\{W/\lambda\}_{r,d,P}\otimes_{\wtr_Z^I}\Omega^{1,\log}_{\wtr_Z^I}\{-1\}\]
      and then we have $\rd_D\wedge\rd_D = 0$.
      \begin{enumerate}
          \item[(1)] Put \[M_{\infty}(D,\nabla_D):= ( D\otimes_{\wtr_Z^I}\CCpsi\{W/\lambda\}_{r,d,P})^{\rd_D = 0}.\]
          Then $M_{\infty}(D,\nabla_D)$ is independent of the choice of $\lambda\in [1,\lambda_0p^{-\frac{1}{p-1}})$, which is finite free of rank $r$ over $\CCpsi$. Moreover, the $\Gamma$-action on $D\otimes_{\wtr_Z^I}\CCpsi\{W/\lambda\}_{r,d,P}$ induces a $\Gamma$-action on $M_{\infty}(D,\nabla_D)$ by restriction, which makes $M_{\infty}(D,\nabla_D)$ an object in $\Rep_{\Gamma}^{\lambda_0}(\CCpsi)$.

          \item[(2)] For any $n\geq 1$, we have
          \[\rH^n_{\dR}(D\otimes_{\wtr_Z^I}\CCpsi\{\underline W\}^{\dagger}_{r,d,P},\rd_D) = 0\]
          and thus the following sequence is exact
          \[0\to M_{\infty}(D,\nabla_D)\to \rD\rR(D\otimes_{\wtr_Z^I}\CCpsi\{\underline W\}^{\dagger}_{r,d,P},\rd_D).\]
          Moreover, we have a quasi-isomorphism
          \[\rR\Gamma(\Gamma,M_{\infty}(D,\nabla_D))\simeq \rD\rR(D,\nabla_D)\]
          which is functorial in $(D,\nabla_D)$.
      \end{enumerate}
  \end{cor}
  \begin{proof}
      Everything follows from Proposition \ref{prop:From connection to representation} except the ``moreover'' part of Item (2). For this, by construction, there is a natural morphism of de Rham complexes
      \[\rD\rR(D,\nabla_D)\to \rD\rR(D\otimes_{\wtr_Z^I}\CCpsi\{\underline W\}^{\dagger}_{r,d,P},\rd_D).\]
      As $\Gamma$ acts on $D$ trivially, the above morphism induces a morphism
      \[\rD\rR(D,\nabla_D)\to \rR\Gamma(\Gamma,\rD\rR(D\otimes_{\wtr_Z^I}\CCpsi\{\underline W\}^{\dagger}_{r,d,P},\rd_D)),\]
      which is a quasi-isomorphism because
      \[\rR\Gamma(\Gamma,\CCpsi\{\underline W\}^{\dagger}_{r,d,P})\cong \wtr_Z^I,\]
      which is a special case of Corollary \ref{cor:From representation to connection}(3).
      Using the exact sequence in Item (2), we get a quasi-isomorphism
      \[\rR\Gamma(\Gamma,M_{\infty}(D,\nabla_D))\xrightarrow{\simeq}\rR\Gamma(\Gamma,\rD\rR(D\otimes_{\wtr_Z^I}\CCpsi\{\underline W\}^{\dagger}_{r,d,P},\rd_D)).\]
      So we finally get a quasi-isomorphism
      \[\rD\rR(D,\nabla_D)\simeq \rR\Gamma(\Gamma,M_{\infty}(D,\nabla_D))\]
      which is clearly functorial in $(D,\nabla_D)$.
  \end{proof}

\subsection{Proof of Theorem \ref{thm:local correspondence over Robba ring}}\label{ssec:proof of local RH}
  \begin{prop}\label{prop:integral equivalence}
      Let $I = [r,s]\subset(\frac{1}{p-1},p)$ be a closed sub-interval and $\lambda_0\in (p^{\frac{1}{p-1}},|\xi|_I^{-1})$ be a real number as before. Then the functor
      \[D_{\psi}^{I,+}:\Rep_{\Gamma}^{\lambda_0}(\CCppsi)\to\MIC^{\lambda_0}(\wtr_Z^{I,+})\]
      in Proposition \ref{prop:From representation to connection} is an equivalence of categories whose inverse is given by the functor
      \[M_{\infty}^+:\MIC^{\lambda_0}(\wtr_Z^{I,+})\to \Rep_{\Gamma}^{\lambda_0}(\CCppsi).\]
  \end{prop}
  \begin{proof}
      Let $(P^+_{\infty},\nabla)$ be either $(\CCppsi\{\underline W/\lambda\}_{r,d,P},\nabla)$ or $(\CCppsi\{\underline W\}_{r,d,P}^{\dagger},\nabla)$, where $\lambda\in [1,\lambda_0p^{-\frac{1}{p-1}}$ for simplicity.

      For any $M_{\infty}\in \Rep_{\Gamma}^{\lambda_0}(\CCppsi)$, put $(D,\nabla_D):=(D_{\psi}^{I,+}(M_{\infty}),\nabla_{D_{\psi}^{I,+}(M_{\infty})})$ for short. We first construct a canonical morphism
      \[\iota_{M_{\infty}}: M_{\infty}(D,\nabla_D)\to M_{\infty}\]
      as follows.

      Consider the composite of the following canonical morphisms of $\Gamma$-equivairant connections:
      \begin{equation*}
          \begin{split}
              &(D\otimes_{\wtr_Z^{I,+}}P^+_{\infty},\rd_D = \nabla_D\otimes\id+\id\otimes\nabla)\\
              \xrightarrow{=}& ((M_{\infty}\otimes_{\CCppsi}P^+_{\infty})^{\Gamma}\otimes_{\wtr_z^{I,+}}P^+_{\infty},\id\otimes\nabla\otimes\id+\id\otimes\id\otimes\nabla)\\
              \hookrightarrow&(M_{\infty}\otimes_{\CCppsi}P^+_{\infty}\otimes_{\CCppsi}P^+_{\infty},\id\otimes\nabla\otimes\id+\id\otimes\id\otimes\nabla)\\
              \to&(M_{\infty}\otimes_{\CCppsi}P^+_{\infty},\id\otimes\nabla).
          \end{split}
      \end{equation*}
      Here, the first arrow follows from the definition of $(D,\nabla_D)$ as $\nabla_D$ is the restriction of $\id\otimes\nabla$ to 
      \[(M_{\infty}\otimes_{\CCppsi}P^+_{\infty})^{\Gamma}\subset M_{\infty}\otimes_{\CCppsi}P^+_{\infty},\]
      the second arrow is induced by the natural inclusion, and the last arrow follows by noting that the multiplication on $P^+_{\infty}$ is compatible with the connection $\nabla$. By taking horizontal sections, we get the desired canonical map $\iota_{M_{\infty}}$.

      For any $(D,\nabla_D)\in \MIC^{\lambda_0}(\wtr_Z^{I,+}))$, put $M_{\infty}:=M_{\infty}(D,\nabla_D)$ for short. We can similarly construct a canonical morphism
      \[\iota_{(D,\nabla_D)}:(D_{\psi}^{I,+}(M_{\infty}),\nabla_{D_{\psi}^{I,+}(M_{\infty})})\to (D,\nabla_D)\]
      as follows. 

      Consider the composite of the following canonical morphisms of $\Gamma$-equivariant connections:
      \begin{equation*}
          \begin{split}
              &(M_{\infty}\otimes_{\CCppsi}P^+_{\infty},\id\otimes\nabla)\\
              \xrightarrow{=}&((D\otimes_{\wtr_Z^{I,+}}P^+_{\infty})^{\rd_D:=\nabla_D\otimes\id+\id\otimes\nabla= 0}\otimes_{\CCppsi}P^+_{\infty},\id\otimes\nabla = \rd_D\otimes\id+\id\otimes\nabla)\\
              \hookrightarrow&(D\otimes_{\wtr_Z^{I,+}}P^+_{\infty}\otimes_{\CCppsi}P^+_{\infty},\nabla_D\otimes\id\otimes\id+\id\otimes\nabla\otimes\id+\id\otimes\id\otimes\nabla)\\
              \to&(D\otimes_{\wtr_Z^{I,+}}P^+_{\infty},\nabla_D\otimes\id+\id\otimes\nabla).
          \end{split}
      \end{equation*}
      Here, the first arrow follows from the definition of $M_{\infty}$, the second arrow is the natural inclusion, and the last arrow is again induced by the multiplication on $P^+_{\infty}$. Then one can obtain the desired $\iota_{(D,\nabla)}$ by taking $\Gamma$-invariants.

      Now, to conclude the result, it suffices to show that the above $\iota_{M_{\infty}}$ and $\iota_{(D,\nabla_D)}$ are both isomorphisms. But this follows from comparing \eqref{equ:explicit D(M)} and \eqref{equ:explicit nabla} with \eqref{equ:explicit M(D)} and \eqref{equ:explicit gamma action} directly.
  \end{proof}

  \begin{proof}[\textbf{Proof of Theorem \ref{thm:local correspondence over Robba ring}}]
      Note that 
      \[\Rep_{\Gamma}^{\dagger}(\CCpsi) = \cup_{\lambda_0>p^{\frac{1}{p-1}}}\Rep_{\Gamma}^{\lambda_0}(\CCpsi) \text{ and }\MIC^{\dagger}(\wtr_Z^I) = \cup_{\lambda_0>p^{\frac{1}{p-1}}}\MIC^{\lambda_0}(\wtr_Z^I).\]
      The assertion follows from Propositions \ref{prop:From representation to connection}, \ref{prop:From connection to representation} and \ref{prop:integral equivalence} directly.
  \end{proof}

  \begin{cor}\label{cor:shrink I-local equivalence}
      Let $J\subset I$ be a closed sub-interval. Then the equivalence $D_{\psi}^I$ is compatible with $D_{\psi}^J$ under base-change; that is, the following diagram is commutative
      \[\xymatrix@C=0.5cm{
        \Rep_{\Gamma}^{\dagger}(\CCpsi)\ar[d]^{D^I_{\psi}}\ar[rrrr]^{-\otimes_{\CCpsi}\widetilde \bfC^J_{\widetilde \psi,\infty}}&&&&\Rep_{\Gamma}^{\dagger}(\widetilde \bfC^J_{\widetilde \psi,\infty})\ar[d]^{D^J_{\psi}}\\
        \MIC^{\dagger}(\wtr_Z^I)\ar[rrrr]^{-\otimes_{\wtr_Z^I}\wtr_Z^J}&&&&\MIC^{\dagger}(\wtr_Z^J),
      }\]
      which is compatible with cohomology in the sense that for any $M_{\infty}\in \Rep_{\Gamma}^{\dagger}(\CCpsi)$, we have quasi-isomorphisms
      \[\begin{split}\rR\Gamma(\Gamma,M_{\infty}\otimes_{\CCpsi}\widetilde \bfC_{\widetilde \psi,\infty}^J)\simeq&\rD\rR(D_{\psi}^J(M_{\infty}\otimes_{\CCpsi}\widetilde \bfC_{\widetilde \psi,\infty}^J),\nabla_{D_{\psi}^J(M_{\infty}\otimes_{\CCpsi}\widetilde \bfC_{\widetilde \psi,\infty}^J)})\\
      \simeq &\rD\rR(D_{\psi}^I(M_{\infty}),\nabla_{D_{\psi}^I(M_{\infty})})\widehat \otimes_{\CCtilde^I(Z)}\CCtilde^J(Z)\\
      \simeq&\rR\Gamma(\Gamma,M_{\infty})\widehat \otimes_{\CCtilde^I(Z)}\CCtilde^J(Z).\end{split}\]
  \end{cor}
  \begin{proof}
      It follows from Lemmas \ref{lem:shrink I-local representation} and \ref{lem:shrink I-local connection} that the base-change functors corresponding to the horizontal arrows in the above diagram are both well-defined. For any $M_{\infty}\in \Rep_{\Gamma}^{\dagger}(\CCpsi)$ with the associated $M\in \Rep_{\Gamma}^{\dagger}(\Cpsi)$ via the equivalence in Theorem \ref{thm:decompletion}, as the natural morphism
      \[\CCpsi\{\underline W\}^{\dagger}_{r,d,P}\to \widetilde \bfC^J_{\widetilde \psi,\infty}\{\underline W\}^{\dagger}_{r,d,P}\]
      is compatible with $\Gamma$-actions and connections, there is a natural morphism of connections
      \[D_{\psi}^I(M_{\infty})\to D_{\psi}^J(M_{\infty}\otimes_{\CCpsi}\widetilde \bfC^J_{\widetilde \psi,\infty})\]
      and thus a morphism of connections
      \[\iota_{M_{\infty}}:D_{\psi}^I(M_{\infty})\otimes_{\wtr_Z^I}\wtr_Z^J\to D_{\psi}^J(M_{\infty}\otimes_{\CCpsi}\widetilde \bfC^J_{\widetilde \psi,\infty})\]
      by scalar extension. We have to show that $\iota_{M_{\infty}}$ is an isomorphism. 

      Suppose that $M_{\infty}\in \Rep_{\Gamma}^{\lambda}(\CCpsi)$ for some $\lambda>p^{\frac{1}{p-1}}$ and let $M\in \Rep_{\Gamma}^{\lambda}(\Cpsi)$ such that 
      \[M\otimes_{\Cpsi}\CCpsi \cong M_{\infty}.\]
      Let $e_1,\dots,e_r$ be a basis of $M$ satisfying the condition \eqref{equ:small representation}.
      It follows from \eqref{equ:explicit D(M)} that $D_{\psi}^I(M_{\infty})$ (resp. $D_{\psi}^J(M_{\infty}\otimes_{\CCpsi}\widetilde \bfC^J_{\widetilde \psi,\infty})$) is a finite free module over $\wtr_Z^I$ (resp. $\wtr_Z^J$) with basis
      \[(\prod_{i=1}^d\gamma_i^{-\frac{W_i}{u}})(e_1),\dots,(\prod_{i=1}^d\gamma_i^{-\frac{W_i}{u}})(e_r).\]
      Then one can conclude that $\iota_{M_{\infty}}$ is an isomorphism by noting that, by its construction, for any $1\leq j\leq r$,
      \[\iota_{M_{\infty}}((\prod_{i=1}^d\gamma_i^{-\frac{W_i}{u}})(e_j)) = (\prod_{i=1}^d\gamma_i^{-\frac{W_i}{u}})(e_j).\]
      The desired cohomological comparison follows from the cohomological comparison part of Lemma \ref{lem:shrink I-local connection} and Theorem \ref{thm:local correspondence over Robba ring}(3). This completes the proof.
  \end{proof}
  
\subsection{Functoriality and pull-backs}\label{ssec:local functorial in pull-back}
  In this part, we assume that $f:\frakX_Z^{\prime} = \Spf(\calR_Z^{\prime})\to \frakX_Z = \Spf(\calR_Z)$ is smooth such that \eqref{equ:relative chart-I} holds true. Let $\widetilde f:\wtx_Z^{\prime}\to \wtx_Z$ be a lifting of $f$ such that \eqref{equ:relative chart-I} holds true. Denote by $X_{\infty,Z}^{\prime}\to X_{\infty,Z}$ the morphism of affinoid perfectoid spaces appearing in Example \ref{exam:relative OC} and keep notations therein.
  Denote by
  \[(\widetilde \bfC_{\widetilde \psi^{\prime},\infty}^I,\widetilde \bfC_{\widetilde \psi^{\prime},\infty}^{I,+}) \text{ and }(\widetilde \bfC_{\widetilde \psi^{\prime}}^I,\widetilde \bfC_{\widetilde \psi^{\prime}}^{I,+})\]
  the analogue of $(\CCpsi,\CCppsi)$ and $(\Cpsi,\Cppsi)$ by replacing $\wtr_Z$ by $\wtr_Z^{\prime}$ instead. For example, we have a $\Gamma^{\prime}$-equivariant decomposition 
  \[\widetilde \bfC_{\widetilde \psi^{\prime},\infty}^{I,+} = \widehat \bigoplus_{\underline \alpha\in J_{r,d^{\prime}}}\widetilde \bfC_{\widetilde \psi^{\prime}}^{I,+}\cdot [(\underline T^{\flat})^{\underline \alpha}].\]
  Then we have natural morphisms
  \[\CCppsi\to \widetilde \bfC_{\widetilde \psi^{\prime},\infty}^{I,+} \text{ and }\Cppsi\to \widetilde \bfC_{\widetilde \psi^{\prime}}^{I,+}\]
  which are $\Gamma^{\prime}$-equivariant. 
  Here, the $\Gamma^{\prime}$-action on $\CCppsi$ (resp. $\Cppsi$) is induced from the quotient map $\Gamma^{\prime}\to \Gamma$.
  The morphism $\widetilde f:\wtr_Z\to\wtr_Z^{\prime}$ also induces a morphism $\wtr_Z^{I,+}\to \wtr_Z^{\prime I,+}$.
  \begin{lem}\label{lem:local pull-back of small}
      Keep the notation as above and let $\lambda>1$.
      \begin{enumerate}
          \item[(1)] The base-change $M_{\infty}\mapsto M_{\infty}\otimes_{\CCppsi}\widetilde \bfC_{\widetilde \psi^{\prime},\infty}^{I,+}$ (resp. $M\mapsto M\otimes_{\Cppsi}\widetilde \bfC_{\widetilde \psi^{\prime}}^{I,+}$) induces a well-defined functor
          \[\text{$\Rep_{\Gamma}^{\lambda}(\CCppsi)\to \Rep_{\Gamma^{\prime}}^{\lambda}(\widetilde \bfC_{\widetilde \psi^{\prime},\infty}^{I,+})$ (resp. $\Rep_{\Gamma}^{\lambda}(\Cppsi)\to \Rep_{\Gamma^{\prime}}^{\lambda}(\widetilde \bfC_{\widetilde \psi^{\prime}}^{I,+})$)}\]
          such that the following diagram commute
          \[\xymatrix@C=0.5cm{
            \Rep_{\Gamma}^{\lambda}(\Cppsi)\ar[d]\ar[rr]&& \Rep_{\Gamma^{\prime}}^{\lambda}(\widetilde \bfC_{\widetilde \psi^{\prime}}^{I,+})\ar[d]\\
            \Rep_{\Gamma}^{\lambda}(\CCppsi)\ar[rr]&& \Rep_{\Gamma^{\prime}}^{\lambda}(\widetilde \bfC_{\widetilde \psi^{\prime},\infty}^{I,+}),
          }\]
          where each arrow stands for the functor induced by base-change.

          \item[(2)] The base-change $(D,\nabla_D)\mapsto (D\otimes_{\wtr_Z^{I,+}}\wtr_Z^{\prime I,+},\nabla_D\otimes\id+\id\otimes\nabla)$ induces a well-defined functor
          \[\MIC^{\lambda}(\wtr_Z^{I,+})\to \MIC^{\lambda}(\wtr_Z^{\prime I,+}).\]
      \end{enumerate}
  \end{lem}
  \begin{proof}
      By Remarks \ref{rmk:isometry on base} and \ref{rmk:isometry for R^I}, for any $x\in \CCppsi$ (resp. $\Cppsi, \wtr_Z^{I,+}$), if we denote its image in $\widetilde \bfC_{\widetilde \psi^{\prime},\infty}^{I,+}$ (resp. $\widetilde \bfC_{\widetilde \psi^{\prime}}^{I,+}, \wtr_Z^{\prime I,+}$) by $x^{\prime}$, then we always have
      \[|x^{\prime}|_I\leq |x|_I.\]
      Then the result follows from definitions.
  \end{proof}
  \begin{prop}\label{prop:local functoriality via pull-back over Robba ring}
      Keep the notation as above. Then the equivalence in Theorem \ref{thm:local correspondence over Robba ring} is compatible with pull-back in the sense that the following diagram is commutative
      \[\xymatrix@C=0.5cm{
         \Rep_{\Gamma}^{\dagger}(\CCpsi)\ar[d]^{D_{\psi}^I}\ar[rrrr]^{-\otimes_{\CCpsi}\widetilde \bfC_{\widetilde \psi^{\prime},\infty}^I}&&&&\Rep_{\Gamma}^{\dagger}(\widetilde \bfC_{\widetilde \psi^{\prime},\infty}^I)\ar[d]^{D_{\psi^{\prime}}^I}\\
         \MIC^{\dagger}(\wtr_Z^I)\ar[rrrr]^{-\otimes_{\wtr_Z^I}\wtr_Z^{\prime I}}&&&&\MIC^{\dagger}(\wtr_Z^{\prime I}),
      }\]
      where the horizontal arrows are the functors induced by base-change, which are well-defined by Lemma \ref{lem:local pull-back of small}, and the vertical arrows are the equivalences established in Theorem \ref{thm:local correspondence over Robba ring}.
  \end{prop}
  \begin{rmk}
      Since $\Gamma_f = \ker(\Gamma^{\prime}\to \Gamma)$ acts trivially on the underlying lattice $M\subset M_{\infty}$, the extra torus directions in $\Gamma^{\prime}$ do not twist the basis vectors $e_j$ used in \eqref{equ:explicit D(M)}. Thus $D_{\psi^{\prime}}^I(M_{\infty}\otimes_{\CCpsi}\widetilde \bfC_{\widetilde \psi^{\prime},\infty}^I)$ is obtained from $D_{\psi}^I(M_{\infty})$ by scalar extension to $\wtr_Z^{\prime I}$, compatibly with the relative semi-stable charts.
  \end{rmk}
  \begin{proof}
      Recall we have $\Gamma^{\prime}$-equivariant morphisms 
      \[\CCpsi\{\underline W\}_{r,d,P}^{\dagger}\to \widetilde \bfC_{\widetilde \psi^{\prime},\infty}^I\{\underline W\}_{r,d,P}^{\dagger}\to \widetilde \bfC_{\widetilde \psi^{\prime},\infty}^I\{\underline W\}_{r,d^{\prime},P}^{\dagger}\]
      which are compatible with the connections, cf. Example \ref{exam:relative OC}. Now, for any $M_{\infty}\in \Rep_{\Gamma}^{\dagger}(\CCpsi)$ with corresponding $M\in \Rep_{\Gamma}^{\dagger}(\Cpsi)$ in the sense of Theorem \ref{thm:decompletion}, by construction of $D_{\psi}^I$ (and $D_{\psi^{\prime}}^I$), we get a natural map
      \[D_{\psi}^I(M_{\infty})\to D_{\psi^{\prime}}^I(M_{\infty}\otimes_{\CCpsi}\widetilde \bfC_{\widetilde \psi^{\prime},\infty}^I)\]
      which is compatible with the connections, yielding a natural map
      \[\iota_{M_{\infty}}: D_{\psi}^I(M_{\infty})\otimes_{\wtr_Z^I}\wtr_Z^{\prime I}\to D_{\psi^{\prime}}^I(M_{\infty}\otimes_{\CCpsi}\widetilde \bfC_{\widetilde \psi^{\prime},\infty}^I)\]
      which is compatible with the connections as well. It remains to show that $\iota_{M_{\infty}}$ is an isomorphism.

      Now, let $e_1,\dots,e_r$ be a basis of $M$ satisfying the condition \eqref{equ:small representation}. Then by \eqref{equ:explicit D(M)}, $D_{\psi}^I(M_{\infty})$ is a finite free $\wtr_Z^I$-module generated by 
      \[(\prod_{i=1}^d\gamma_i^{-\frac{W_i}{u}})(e_1),\dots,(\prod_{i=1}^d\gamma_i^{-\frac{W_i}{u}})(e_r).\]
      Here and in what follows, we use the isomorphisms
      \[\Gamma \cong \oplus_{i=1}^d\Zp\gamma_i \text{ and }\Gamma^{\prime}\cong \oplus_{i=1}^{d^{\prime}}\Zp\gamma_i\]
      such that the kernel $\Gamma_f:=\ker(\Gamma^{\prime}\to \Gamma)\cong \oplus_{i=d+1}^{d^{\prime}}\Zp\gamma_i$. Similarly, $D_{\psi^{\prime}}^I(M_{\infty}\otimes_{\CCpsi}\widetilde \bfC_{\widetilde \psi^{\prime},\infty}^I)$ is a finite free $\wtr_Z^{\prime I}$-module generated by 
      \[(\prod_{i=1}^{d^{\prime}}\gamma_i^{-\frac{W_i}{u}})(e_1),\dots,(\prod_{i=1}^{d^{\prime}}\gamma_i^{-\frac{W_i}{u}})(e_r).\]
      We have to show that for any $1\leq j\leq r$, 
      \[\iota_{M_{\infty}}\big((\prod_{i=1}^{d}\gamma_i^{-\frac{W_i}{u}})(e_j)\big) = (\prod_{i=1}^{d^{\prime}}\gamma_i^{-\frac{W_i}{u}})(e_j).\]
      But this follows from the fact that for any $1\leq j\leq d$,
      \[(\prod_{i=1}^{d^{\prime}}\gamma_i^{-\frac{W_i}{u}})(e_j) = (\prod_{i=1}^{d}\gamma_i^{-\frac{W_i}{u}})(e_j),\]
      because $\Gamma_f$ acts on $M$ trivially. 
  \end{proof}

\section{Local Riemann--Hilbert Correspondence: The de Rham Case}\label{sec:local RH over BdR}
   Throughout this section, let $Z = \Spa(A,A^+)\in \Perfd$ and let $\frakX_Z = \Spf(\calR_Z)$ be a small semi-stable formal scheme with a chart $\psi$ and a fixed prismatic lifting $\wtx_Z = \Spf(\wtr_Z)$ over $\bA_{\inf}(Z)$. Put $R_Z = \calR_Z[\frac{1}{p}]$, and let $\wtR_Z$ denote the lifting of $R_Z$ over $\BBdRp(Z)$ associated to $\wtr_Z$ (Example \ref{exam:small affine is liftable}). Let $X_{\infty,Z} = \Spa(\widehat R_{\infty,Z},\widehat R_{\infty,Z}^+)\to X_Z$ be the pro-finite-\'etale $\Gamma$-torsor from Example \ref{exam:Gamma torsor}. As in \S\ref{sec:local RH}, the morphism 
  \[\widetilde \psi:\wtr_Z\to \bA_{\inf}(X_{\infty,Z})\]
  sending each $T_i$ to $[T_i^{\flat}]$ induces a morphism of $\BBdRp(Z)$-algebras
  \begin{equation*}
      \widetilde \psi_{\dR}:\wtR_Z\to\BBdRp(X_{\infty,Z})
  \end{equation*}
  and $\Gamma$-equivariant decompositions of the associated completed period rings (cf.\ \eqref{equ:decomposition over dR} below).
  \begin{dfn}
      Define $\BBpsi:=\BBdRp(X_{\infty,Z})$ 
      and $\Bpsi = \widetilde \psi_{\dR}(\wtR_Z)\subset\Bppsi$, so that the decomposition of $\BBdRp(X_{\infty,Z})$ reads
      \begin{equation}\label{equ:decomposition over dR}
          \BBpsi = \widehat \bigoplus_{\underline \alpha\in J_{r,d}}\Bpsi\cdot[(\underline T^{\flat})^{\underline \alpha}].
      \end{equation}
      Note that for any $\bfB\in\{\BBpsi,\Bpsi\}$, we have
      \[\bfB \cong \varprojlim_{n\to+\infty}\bfB/t^n\]
      and thus we also denote $\bfB$ by $\bfB/t^{\infty}$ sometimes.
  \end{dfn}
  The purpose of this section is to give an analogue of Theorem \ref{thm:local correspondence over Robba ring} for $\BBpsi$-representations of $\Gamma$ and flat connections over $\wtR_Z$ by a standard d\'evissage argument.
  
  Recall that we have $\BBpsi/t = \widehat R_{\infty,Z}$ and $\Bpsi/t = R_Z$. The decomposition \eqref{equ:decomposition over dR} reduces to the decomposition \eqref{equ:decomposition original} (after inverting $p$) modulo $t$
  \[\widehat R_{\infty,Z} = \widehat \bigoplus_{\underline \alpha\in J_{r,d}}R_Z\cdot \underline T^{\underline \alpha}.\]

\subsection{Review on local Simpson correspondence}\label{ssec:local Simpson}
  Now, we give a quick review of local Simpson correspondence in \cite{Wan23}.

  \begin{dfn}\label{dfn:small local generalised representations}
      Let $(B,B^+)\in \{(R_Z,R_Z^+=\calR_Z),(\widehat R_{Z,\infty},\widehat R_{\infty,Z}^+)\}$.
      \begin{enumerate}
          \item Let $\lambda>1$. A representation $M$ of rank $r$ of $\Gamma$ over $B^+$ is called \emph{$\lambda$-trivial}, if there exists a $B^+$-basis $e_1,\dots,e_r$ of $M$ such that if we equip $M$ with the supremum norm, then for any $\gamma\in\Gamma$ and any $1\leq i\leq r$, we have 
         \begin{equation*}
             (\gamma-1)(e_i)\subset (\zeta_p-1)M \text{ such that } |\frac{(\gamma-1)}{\zeta_p-1}(e_i)|\leq \lambda^{-1}.
         \end{equation*}
         A representation $M$ of $\Gamma$ over $B$ is called \emph{$\lambda$-trivial} if it is of the form $M = M^+[\frac{1}{p}]$ for some $\lambda$-trivial $\Gamma$-representation $M^+$ over $B^+$. Denote by $\Rep_{\Gamma}^{\lambda}(B^+)$ (resp. $\Rep_{\Gamma}^{\lambda}(B)$) the category of $\lambda$-trivial representations of $\Gamma$ over $B^+$ (resp. $B$).

         \item Let $\lambda_0\geq 1$. A representation $M$ of $\Gamma$ over $B^+$ (resp. $B$) is called \emph{$\lambda_0$-small} if it is $(\lambda\frac{1}{p-1})$-trivial for some $\lambda>\lambda_0$.
         Denote by $\Rep_{\Gamma}^{\lambda_0,\dagger}(B^+)$ (resp. $\Rep_{\Gamma}^{\lambda_0,\dagger}(B)$) the category of $\lambda_0$-small representations of $\Gamma$ over $B^+$ (resp. $B$). We also say an $M$ over $B^+$ (resp. $B$) is \emph{small} if it is $1$-small and denote by $\Rep_{\Gamma}^{\dagger}(B^+)$ (resp. $\Rep_{\Gamma}^{\dagger}(B)$) the corresponding category.
      \end{enumerate}
  \end{dfn}
  \begin{dfn}\label{dfn:small local Higgs module}
      \begin{enumerate}
          \item By a \emph{Higgs module} of rank $r$ over $\calR_Z$, we mean a finite free $\calR_Z$-module $H$ of rank $r$ together with an $\calR_Z$-linear map, called \emph{Higgs field}, 
          \[\theta_H:H\to H\otimes_{\calR_Z}\Omega^{1,\log}_{\calR_Z}\{-1\}\]
          satisfying $\theta_H\wedge\theta_H = 0$. Denote by $\HIG(\calR_Z)$ the category of Higgs modules over $\calR_Z$.

          \item Let $\lambda>1$. A Higgs module $(H,\theta_H)$ over $\calR_Z$ is called \emph{$\lambda$-trivial} if for any $x\in H$, we have
          \[|\theta_H(x)|\leq \lambda^{-1}|x|,\]
          where $|\cdot|$ denotes the supremum norm on $H$.
          Denote by $\HIG^{\lambda}(\calR_Z)$ the full-subcategory of $\HIG(\calR_Z)$ consisting of $\lambda$-trivial objects.

          \item Let $\lambda_0\geq 1$. A Higgs module $(H,\theta_Z)$ is called \emph{$\lambda_0$-small} if it is $(\lambda p^{\frac{1}{p-1}})$-trivial for some $\lambda>\lambda_0$. Denote by $\HIG^{\lambda_0,\dagger}(\calR_Z)$ the full-subcategory of $\HIG(\calR_Z)$ consisting of $\lambda_0$-small objects. In particular, we call a Higgs module \emph{small} if it is $1$-small.

          \item One can similarly introduce the category $\HIG(R_Z)$ of Higgs module over $R_Z$ as in Item (1). For such a Higgs module $(H,\theta_H)\in \HIG(R_Z)$, we call it $\lambda$-trivial (resp. $\lambda_0$-small, small) if it admits a $\lambda$-trivial (resp. $\lambda_0$-small, small) lattice $(H^+,\theta_{H^+})$ over $\calR_Z$, and denote by $\HIG^{\lambda}(R_Z)$ (resp. $\HIG^{\lambda_0,\dagger}(R_Z)$, $\HIG^{\dagger}(R_Z)$) the corresponding category.

          \item For any Higgs module $(H,\theta_Z)$ over $B\in \{\calR_Z,R_Z\}$, denote by $\rD\rR(H,\theta_Z)$ the associated de Rham complex:
          \[H\xrightarrow{\theta_Z}H\otimes_{B_Z}\Omega^{1,\log}_B\{-1\}\xrightarrow{\theta_H}\cdots\xrightarrow{\theta_H}\Omega^{d,\log}_B\{-d\}. \]
          For any $n\geq 0$, define
          \[\rH^n_{\dR}(H,\theta_H):=\rH^n(\rD\rR(H,\theta_H)).\]
      \end{enumerate}
  \end{dfn}
  \begin{rmk}\label{rmk:compare local smallness with Wang}
      \begin{enumerate}
          \item Keep notations as in Definition \ref{dfn:small local generalised representations}. For any $a>\frac{1}{p-1}$ such that $\lambda = p^{a-\frac{1}{p-1}}$, it is easy to see that a $B^+$-representation of $\Gamma$ is $\lambda$-trivial if and only if it is $a$-trivial in the sense of \cite[Def. 3.3]{Wan23}. However, the $\lambda_0$-smallness of a representation $M$ of $\Gamma$ over $B^+$ is slightly different from \cite[Def. 4.1]{Wan23}. The $M$ is $\lambda_0$-small if and only if there exists some $a>\frac{1}{p-1}$ with $\lambda = p^{a-\frac{1}{p-1}}>\lambda_0$ such that $M$ is $a$-small in the sense of \cite[Def. 4.1]{Wan23}

          \item Keep notations as in Definition \ref{dfn:small local Higgs module}. For any $\lambda = p^{a-\frac{1}{p-1}}>1$, a Higgs module $(H,\theta_H)$ over $\calR_Z$ is $(\lambda p^{\frac{1}{p-1}})$-trivial if and only if it is $a$-small in the sense of \cite[Def. 4.2]{Wan23}: Note that in \emph{loc.cit.}, we regard the Higgs field $\theta$ as a map
          \[\theta_H:H\to H\otimes_{\calR_Z}\Omega^{1,\log}_{\calR_Z}(-1)\]
          taking values in the Tate twist $\Omega^{1,\log}_{\calR_Z}(-1)$ rather than the Breuil--Kisin twist $\Omega^{1,\log}_{\calR_Z}\{-1\}$. The difference is that $\Omega^{1,\log}_{\calR_Z}\{-1\} = (\zeta_p-1)\Omega^{1,\log}_{\calR_Z}(-1)$ and thus $(H,\theta_H)$ is $(\lambda p^{\frac{1}{p-1}})$-trivial if and only if 
          \[\theta_H(H)\subset p^{-(a+\frac{1}{p-1})}H\otimes_{\calR_Z}\Omega^{1,\log}_{\calR_Z}(-1).\]
      \end{enumerate}
  \end{rmk}
 \begin{lem}\label{lem:(zeta_p-1)-error}
    Fix a $\lambda> 1$ and let $M$ be a $\lambda$-trivial representation of $\Gamma$ over $\widehat R_{\infty,Z}^+$. For any $n\geq 0$, the cohomological group
    \[\rH^n(\Gamma,M\otimes_{\calR_Z}(\widehat R_{\infty,Z}^+/\calR_Z))\]
    is killed by $\zeta_p-1$ and vanishes when $n = 0$.
 \end{lem}
 \begin{proof}
     This follows from the same argument in the proof of Proposition \ref{prop:u-error}.
 \end{proof}
 The following corollary then follows immediately.
 \begin{cor}\label{cor:(zeta_p-1)-error}
     Let $\lambda>1$. Then the functor
     \[\Rep_{\Gamma}^{\lambda}(\calR_Z^+)\to\Rep_{\Gamma}^{\lambda}(\widehat R_{\infty,Z}^+)\]
     sending each $M$ to $M_{\infty}:=M\otimes_{\calR_Z}\widehat R_{\infty,Z}^+$ is fully faithful such that the morphism
     \[\rR\Gamma(\Gamma,M)\to\rR\Gamma(\Gamma,M_{\infty})\]
     identifies $\rR\Gamma(\Gamma,M)$ with a direct summand of $\rR\Gamma(\Gamma,M_{\infty})$, whose complement is concentrated in degree $\geq 1$ and killed by $\zeta_p-1$.
 \end{cor}

 \begin{thm}\label{thm:decompletion modulo t}
     Let $\lambda>p^{\frac{1}{p-1}}$. Then the functor
     \[\Rep_{\Gamma}^{\lambda}(\calR_Z^+)\to\Rep_{\Gamma}^{\lambda}(\widehat R_{\infty,Z}^+)\]
     sending each $M$ to $M_{\infty}:=M\otimes_{\calR_Z}\widehat R_{\infty,Z}^+$ is an equivalence of categories, which preserves ranks, tensor products and duals, such that the morphism
     \[\rR\Gamma(\Gamma,M)\to\rR\Gamma(\Gamma,M_{\infty})\]
     identifies $\rR\Gamma(\Gamma,M)$ with a direct summand of $\rR\Gamma(\Gamma,M_{\infty})$, whose complement is concentrated in degree $\geq 1$ and killed by $\zeta_p-1$.
 \end{thm}
 \begin{proof}
     See \cite[Th. 3.4]{Wan23} (and the proof therein).
     Alternatively, one can deduce the full faithfulness from Corollary \ref{cor:(zeta_p-1)-error} and the essential surjectivity from the same argument as in the proof of Proposition \ref{prop:essentially surjective}.
 \end{proof}
 \begin{cor}\label{cor:decompletion modulo t}
     Then the functor $M\mapsto M_{\infty}:=M\otimes_{R_Z}\widehat R_{\infty,Z}$ induces an equivalence of categories
     \[\Rep_{\Gamma}^{\dagger}(R_Z)\to\Rep_{\Gamma}^{\dagger}(\widehat R_{\infty,Z})\]
     which preserves ranks, tensor products and duals, such that the natural map $M\to M_{\infty}$ induces a quasi-isomorphism
     \[\rR\Gamma(\Gamma,M)\to\rR\Gamma(\Gamma,M_{\infty}).\]
 \end{cor}
 \begin{proof}
     Note that $\Rep_{\Gamma}^{\dagger}(R_Z) = \cup_{\lambda>p^{\frac{1}{p-1}}}\Rep_{\Gamma}^{\lambda}(R_Z)$ and $\Rep_{\Gamma}^{\dagger}(\widehat R_{\infty,Z})=\cup_{\lambda>p^{\frac{1}{p-1}}}\Rep_{\Gamma}^{\lambda}(\widehat R_{\infty,Z})$. One can conclude by Theorem \ref{thm:decompletion modulo t} after inverting $p$.
 \end{proof}

 Now, for any $1\leq \lambda<p^{\frac{1}{p-1}}$, put
 \[(S_{\infty}^{\lambda},\nabla:S_{\infty}^{\lambda}\otimes_{\calR_Z}\Omega^{1,\log}_{\calR_Z}\{-1\}):=(\calO\widehat \bC^{\lambda},\nabla)(X_{\infty,Z}).\]
 By Corollary \ref{cor:local description of OBdRp^lambda}, for $P =W_0+\cdots+W_r$, we have 
 \[(S_{\infty}^{\lambda},\nabla)\cong (\widehat R_{\infty,Z}\{\underline W/\lambda\}_{r,d,P},\sum_{i=0}^d\frac{\partial}{\partial W_i}\otimes\frac{u}{t}\dlog T_i)\]
 via the isomorphism 
 \[\Omega^{1,\log}_{\calR_Z}\{-1\} \cong \oplus_{i=0}^d\calR_Z\cdot \frac{u}{t}\dlog T_i/\calR_Z\cdot\sum_{i=0}^r\frac{u}{t}\dlog T_i.\]
 Moreover, via the isomorphism
 \[\widehat R_{\infty,Z}\{\underline W/\lambda\}_{r,d,P} = \widehat R_{\infty,Z}\{W_1/\lambda,\dots,W_d/\lambda\} \text{ and }\Omega^{1,\log}_{\calR_Z}\{-1\} \cong \oplus_{i=1}^d\calR_Z\cdot \frac{u}{t}\dlog T_i,\]
 we have
 \[(S_{\infty}^{\lambda},\nabla) \cong (\widehat R_{\infty,Z}\{W_1/\lambda,\dots,W_d/\lambda\},\sum_{i=1}^d\frac{\partial}{\partial W_i}\otimes\frac{u}{t}\dlog T_i).\]
 Note that via the isomorphism \eqref{equ:gamma group-II}, for any $1\leq i,j\leq d$, we have 
 \[\gamma_i(W_j) = W_j+\delta_{ij}(\zeta_p-1)\]
 by noting that $u\in \Ainf$ becomes $\zeta_p-1$ modulo $\xi$.
 Put
 \[S_{\infty}^{\dagger}:=\colim_{\lambda\to 1^+}S_{\infty}^{\lambda}.\]
   \begin{thm}[Local Simpson correspondence]\label{thm:local Simpson correspondence}
      \begin{enumerate}
          \item For any $M_{\infty}\in \Rep^{\dagger}_{\Gamma}(\widehat R_{\infty,Z})$ of rank $r$, we have
            \[
              \rH^n(\Gamma,M_{\infty}\otimes_{\widehat R_{\infty,Z}}S_{\infty}^{\dagger}) = 
                 \left\{
                    \begin{array}{rcl}
                      H(M_{\infty}), & n = 0 \\
                      0, & n\geq 1,
                    \end{array}
                 \right.
            \]
        where $H(M_{\infty})$ is a finite free $R_Z$-module of rank $r$. Moreover, the restriction of 
        \[\id_{M_{\infty}}\otimes\nabla=:\nabla_{M_{\infty}}:M_{\infty}\otimes_{\widehat R_{\infty,Z}}S_{\infty}^{\dagger}\to M_{\infty}\otimes_{\widehat R_{\infty,Z}}S_{\infty}^{\dagger}\otimes_{\calR_Z}\Omega^{1,\log}_{\calR_Z}\{-1\}\]
        to $H(M_{\infty})$ defines a Higgs field $\theta_{H(M_{\infty})}$ on $H(M_{\infty})$ such that $(H(M_{\infty}),\theta_{H(M_{\infty}))}$ is an object in $\HIG^{\dagger}(R_Z)$.

        \item For any $(H,\theta_H)\in \HIG^{\dagger}(R)$ of rank $r$, define
        \[\nabla_{H} = \theta_H\otimes\id+\id\otimes\nabla:H\otimes_{R_Z}S_{\infty}^{\dagger}\to H\otimes_{R_Z}S_{\infty}^{\dagger}\otimes_{\calR_Z}\Omega^{1,\log}_{\calR_Z}\{-1\}.\]
        Then it satisfies $\nabla_{H}\wedge\nabla_{H} = 0$ and we have
        \[
          \rH^n_{\dR}(H\otimes_{R_Z}S_{\infty}^{\dagger},\nabla_{H}) = 
            \left\{
              \begin{array}{rcl}
                  M_{\infty}(H,\theta_H), & n=0 \\
                  0, & n\geq 1,
              \end{array}
            \right.
        \]
        where $M_{\infty}(H,\theta_H)$ is a finite free $\widehat R_{\infty,Z}$-module of rank $r$. Moreover, the $\Gamma$-action on $P_{\infty}^{\dagger}$ induces a $\Gamma$-action on $M_{\infty}(H,\theta_H)$ such that $M_{\infty}(H,\theta_H)$ is an object in $\Rep_{\Gamma}^{\dagger}(\widehat R_{\infty,Z})$.

        \item The functors $M_{\infty}\mapsto (H(M_{\infty}),\theta_{H(M_{\infty})})$ and $(H,\theta_H)\mapsto M_{\infty}(H,\theta_H)$ induce an equivalence of categories
        \[\Rep_{\Gamma}^{\dagger}(\widehat R_{\infty,Z})\simeq \HIG^{\dagger}(R_Z)\]
        preserving ranks, tensor products and dualities. Moreover for any $M_{\infty}\in \Rep_{\Gamma}^{\dagger}(\widehat R_{\infty,Z})$ corresponding to $(H,\theta_H)\in \HIG^{\dagger}(R_Z)$, we have a $\Gamma$-equivariant isomorphism
        \[(M_{\infty}\otimes_{\widehat R_{\infty,Z}}S_{\infty}^{\dagger},\nabla_{M_{\infty}})\cong (H\otimes_{R_Z}S_{\infty}^{\dagger},\nabla_H)\]
        and a quasi-isomorphism
        \[\RGamma(\Gamma,M_{\infty})\simeq \rD\rR(H,\theta_H).\]
      \end{enumerate}
  \end{thm}
  \begin{proof}
      This is \cite[Th. 4.3]{Wan23} (by the same proof.) Alternatively, one can just apply the same argument as in the proof of Theorem \ref{thm:local correspondence over Robba ring}.
  \end{proof}

\subsection{Small representations and flat connections}\label{ssec:small BdR rep and connection}
  In this subsection, we put $\wtR_Z/t^{\infty} := \wtR_Z= \varprojlim_{n\to\infty}\wtR_Z/t^n$.
  We extend the notion of smallness to the $t$-adic setting by reduction modulo $t$: a representation or connection is small if its reduction modulo $t$ is small in the sense of Definitions \ref{dfn:small local generalised representations} and \ref{dfn:small local Higgs module}.

  \begin{dfn}\label{dfn:small local BdRp-representation}
      Let $n\in \bN\cup\{\infty\}$ and let $B\in \{\Bppsi/t^n,\Bpsi/t^n\}$. A representation $M$ of $\Gamma$ over $B$ of rank $r$ is called \emph{small} if its reduction $M/t$ is small of rank $r$ in the sense of Definition \ref{dfn:small local generalised representations}(2). Denote by $\Rep_{\Gamma}^{\dagger}(B)$ the category of small representations of $\Gamma$ over $B$.
  \end{dfn}
  
  \begin{dfn}\label{dfn:small local flat connection over BdRp}
      Fix an $n\in \bN\cup\{\infty\}$.
      \begin{enumerate}
          \item By a \emph{flat connection} of rank $r$ over $\wtR_Z/t^n$, we mean a finite free $\wtR_Z/t^n$-module $D$ of rank $r$ together with a $\BdRp$-linear map
          \[\nabla_D:D\to D\otimes_{\wtR_Z}\Omega^1_{\wtR_Z}\{-1\}\]
          satisfying $\nabla_D\wedge\nabla_D = 0$ and the Leibniz rule with respect to the derivation 
          \[\nabla:\wtR_Z\to \Omega^1_{\wtR_Z}\{-1\}.\] 
          Denote by $\MIC(\wtR_Z/t^n)$ the category of flat connections over $\wtR_Z/t^n$. Such a connection $(D,\nabla_D)$ is called \emph{small} if its reduction modulo $t$ 
          \[(H,\theta_H):=(D,\nabla_D)/t\]
          is a small Higgs module in the sense of Definition \ref{dfn:small local Higgs module}(3). Denote by $\MIC^{\dagger}(\wtR_Z/t^n)$ the full subcategory of $\MIC(\wtR_Z/t^n)$ consisting of small objects.
          
          \item For any flat connection $(D,\nabla_D)$ over $\wtR_Z/t^n$, denote by $\rD\rR(D,\nabla_D)$ the associated de Rham complex:
          \[D\xrightarrow{\nabla_D}D\otimes_{\wtR_Z}\Omega^1_{\wtR_Z}\{-1\}\xrightarrow{\nabla_D}\cdots\xrightarrow{\nabla_D}D\otimes_{\wtR_Z}\Omega^d_{\wtR_Z}\{-d\}. \]
          For any $n\geq 0$, define
          \[\rH_{\dR}^n(D,\nabla_D):=\rH^n(\rD\rR(D,\nabla_D)).\]
      \end{enumerate}
  \end{dfn}
  It follows directly that we have
  \[\Rep_{\Gamma}^{\dagger}(\BBpsi/t) = \Rep_{\Gamma}^{\dagger}(\widehat R_{\infty,Z}), \Rep_{\Gamma}^{\dagger}(\Bpsi/t) = \Rep_{\Gamma}^{\dagger}(R_Z) \text{ and }\MIC^{\dagger}(\wtR_Z/t) = \HIG^{\dagger}(R_Z).\]
  Now, we put
  \[(P_{\infty}^{\dagger},\nabla:P_{\infty}^{\dagger}\to P_{\infty}^{\dagger}\otimes_{\wtR_Z}\Omega^1_{\wtR_Z}\{-1\}):=(\calO\bB_{\dR}^{\dagger,+},\nabla)(X_{\infty,Z}), \text{ cf. Definition \ref{dfn:Overconvergent de Rham period sheaf}}.\]
  Then $(P_{\infty}^{\dagger},\nabla)$ carries a $\Gamma$-action such that 
  \[(P_{\infty}^{\dagger},\nabla)/t = (S_{\infty}^{\dagger},\nabla).\]

  \subsection{The local correspondence}\label{ssec:proof of local RH over BdR}
  The following theorem is the $\BBdRp$-analogue of Theorem \ref{thm:local correspondence over Robba ring}; the coordinates $W_i$ from \S\ref{sec:local RH} enter through the comparison with $(P_{\infty}^{\dagger},\nabla)$ and with \eqref{equ:explicit D(M)}.
  \begin{thm}[Local Riemann--Hilbert correspondence over $\BdRp$]\label{thm:local RH over BdRp}
  Fix $1\leq n\leq \infty$. 
      \begin{enumerate}
          \item For any $M_{\infty}\in \Rep^{\dagger}_{\Gamma}(\Bppsi/t^n)$ of rank $r$, we have
            \[
              \rH^n(\Gamma,M_{\infty}\otimes_{\Bppsi}P_{\infty}^{\dagger}) = 
                 \left\{
                    \begin{array}{rcl}
                      D_{{\rm RH},\psi}(M_{\infty}), & n = 0 \\
                      0, & n\geq 1,
                    \end{array}
                 \right.
            \]
        where $D_{{\rm RH},\psi}(M_{\infty})$ is a finite free $\wtR_Z/t^n$-module of rank $r$. Moreover, the restriction of 
        \[\rd_{M_{\infty}}:=\id\otimes\nabla: M_{\infty}\otimes_{\BBpsi}P_{\infty}^{\dagger}\to M_{\infty}\otimes_{\BBpsi}P_{\infty}^{\dagger}\otimes_{\wtR_Z}\Omega^1_{\wtR_Z}\{-1\}\]
        to $D_{{\rm RH},\psi}(M_{\infty})$ induces a flat connection $\nabla_{D_{{\rm RH},\psi}(M_{\infty})}$ on $D_{{\rm RH},\psi}(M_{\infty})$ such that
        \[(D_{{\rm RH},\psi}(M_{\infty}),\nabla_{D_{{\rm RH},\psi}(M_{\infty})})\in\MIC^{\dagger}(\wtR_Z/t^n).\]

        \item For any $(D,\nabla_D)\in \MIC^{\dagger}(\wtR_Z/t^n)$ of rank $r$, define
        \[\rd_{D} = \nabla_D\otimes\id+\id\otimes\nabla:D\otimes_{\wtR_Z}P_{\infty}^{\dagger}\to D\otimes_{\wtR_Z}P_{\infty}^{\dagger}\otimes_{\wtR_Z}\Omega^1_{\wtR_Z}\{-1\}.\]
        Then we have $\nabla_D\wedge\nabla_D = 0$ and
        \[
          \rH^n_{\dR}(D\otimes_{\wtR_Z}P_{\infty}^{\dagger},\rd_{D}) = 
            \left\{
              \begin{array}{rcl}
                  M_{\infty}(D,\nabla_D), & n=0 \\
                  0, & n\geq 1,
              \end{array}
            \right.
        \]
        where $M_{\infty}(D,\nabla_D)$ is a finite free $\Bppsi/t^n$-module of rank $r$. Moreover, the $\Gamma$-action on $P_{\infty}^{\dagger}$ induces a $\Gamma$-action on $M_{\infty}(D,\nabla_D)$ such that 
        \[M_{\infty}(D,\nabla_D)\in \Rep_{\Gamma}^{\dagger}(\Bppsi).\]

        \item The functors $M_{\infty}\mapsto (D_{{\rm RH},\psi}(M_{\infty}),\nabla_{D_{{\rm RH},\psi}(M_{\infty})})$ and $(D,\nabla_D)\mapsto M_{\infty}(D,\nabla_D)$ induce an equivalence of categories
        \[\Rep_{\Gamma}^{\dagger}(\Bppsi/t^n)\simeq \MIC^{\dagger}(\wtR_Z/t^n)\]
        which preserve ranks, tensor products and duals. Moreover for any $M_{\infty}\in \Rep_{\Gamma}^{\dagger}(\Bppsi)$ corresponding to $(D,\nabla_D)$, we have a quasi-isomorphism
        \[\RGamma(\Gamma,M_{\infty})\simeq \rD\rR(D,\nabla_D).\]
      \end{enumerate}
  \end{thm}
  \begin{proof}
      This follows from the same d\'evissage argument for the proof of \cite[Prop. 4.10]{LMNQW} by using Theorem \ref{thm:local Simpson correspondence} and $(P_{\infty}^{\dagger},\nabla)$ instead of \cite[Prop. 4.9]{LMNQW} and $(\BBpsi\za\underline W\ya_{\pd},\nabla) = (\calO\bB_{\dR,\pd}^+,\nabla)(X_{\infty,Z})$ respectively.
  \end{proof}
  \begin{rmk}\label{rmk:local compare with LMNQW}
      For any $1\leq n\leq \infty$, let $\Rep_{\Gamma}^{\text{H-sm}}(\BBpsi/t^n)$ and $\MIC^{\text{H-sm}}(\wtR_Z/t^n)$ be the categories of \emph{Hitchin-small} representations of $\Gamma$ over $\BBpsi/t^n$ and flat connections over $\wtR_Z$ defined in \cite{LMNQW}. Then we have
      \[\Rep_{\Gamma}^{\dagger}(\BBpsi/t^n)\subset \Rep_{\Gamma}^{\text{H-sm}}(\BBpsi/t^n)\text{ and }\MIC^{\dagger}(\wtR_Z/t^n)\subset\MIC^{\text{H-sm}}(\wtR_Z/t^n).\]
      By Proposition \ref{prop:compare with LMNQW}, the equivalence in Theorem \ref{thm:local RH over BdRp} is compatible with that established in \cite[Prop. 4.10]{LMNQW}.
  \end{rmk}
  To complete the theory, we state the following decompletion theorem for representations of $\Gamma$ over $\BBpsi/t^n$. As it will not be used anywhere in this paper, we leave the detailed proof to the reader. 
  \begin{cor}\label{cor:decompletion over BdRp}
      Fix $1\leq n\leq \infty$.
      The base-change $M\mapsto M_{\infty}:M\otimes_{\Bpsi}\BBpsi$ induces an equivalence of categories
      \[\Rep_{\Gamma}^{\dagger}(\Bpsi/t^n)\to \Rep_{\Gamma}^{\dagger}(\BBpsi/t^n)\]
      which preserve ranks, tensor products and duals. Moreover, for any $M\in \Rep_{\Gamma}^{\dagger}(\Bpsi/t^n)$, the natural inclusion $M\subset M_{\infty}$ induces a quasi-isomorphism
      \[\rR\Gamma(\Gamma,M)\xrightarrow{\simeq} \rR\Gamma(\Gamma,M_{\infty}).\]
  \end{cor}
  \begin{proof}
      This follows from \cite[Cor. 4.13]{LMNQW} and Remark \ref{rmk:local compare with LMNQW}. 
  \end{proof}

\subsection{Compatibility with Robba rings}\label{ssec:compare local RH}
  We compare Theorem \ref{thm:local RH over BdRp} with Theorem \ref{thm:local correspondence over Robba ring} from \S\ref{sec:local RH}. 
  \begin{lem}\label{lem:small modulo t is small}
      Let $I\subset (\frac{1}{p-1},p)$ be a closed sub-interval as before and suppose that $1\in I$. Let $\lambda>1$ be a real number.
      \begin{enumerate}
          \item[(1)] The base-change along $\CCppsi\to \widehat R_{\infty,Z}^+$ (resp. $\Cppsi\to \calR_Z$) gives rise to a functor
          \[\text{$\Rep_{\Gamma}^{\lambda}(\CCppsi)\to\Rep_{\Gamma}^{\lambda}(\widehat R_{\infty,Z}^+)$ (resp. $\Rep_{\Gamma}^{\lambda}(\Cppsi)\to\Rep_{\Gamma}^{\lambda}(\calR_Z)$)}\]
          such that the following diagram commutes
          \[\xymatrix@C=0.5cm{
            \Rep_{\Gamma}^{\lambda}(\Cppsi)\ar[d]\ar[rrrr]&&&&\Rep_{\Gamma}^{\lambda}(\calR_Z)\ar[d]\\
            \Rep_{\Gamma}^{\lambda}(\CCppsi)\ar[rrrr]&&&&\Rep_{\Gamma}^{\lambda}(\widehat R_{\infty,Z}^+),
          }\]
          where all functors involved are induced by base-changes.

          \item[(2)] The base-change along $\wtr_Z^{I,+}\to \calR_Z$ gives rise to a functor
          \[\MIC^{\lambda}(\wtr_Z^{I,+})\to \HIG^{\lambda}(\calR_Z).\]
      \end{enumerate}
  \end{lem}
  \begin{proof}
      This follows by unwinding definitions.
  \end{proof}
    \begin{lem}\label{lem:completion of small is small}
      Let $I\subset (\frac{1}{p-1},p)$ be a closed sub-interval as before and suppose that $1\in I$. Fix an $1\leq n\leq \infty$.
      \begin{enumerate}
          \item[(1)] The base-change along $\CCpsi\to \BBpsi/t^n$ gives rise to a functor
          \[\Rep_{\Gamma}^{\dagger}(\CCpsi)\to\Rep_{\Gamma}^{\dagger}(\BBpsi/t^n).\]

          \item[(2)] The base-change along $\wtr_Z^I\to \wtR_Z/t^n$ gives rise to a functor
          \[\MIC^{\dagger}(\wtr_Z^I)\to \MIC^{\dagger}(\wtR_Z/t^n),\]
          which is compatible with cohomology in the sense that for any $(D,\nabla_D)\in \MIC^{\lambda}(\wtr_Z^{I,+})$, there is a quasi-isomorphism
          \[\rD\rR(D,\nabla_D)\widehat \otimes_{\CCtilde^I(Z)}\wtR_Z/t^n\simeq \rD\rR(D\otimes_{\wtr_Z^I}R_Z,\nabla_{D\otimes_{\wtr_Z^I}R_Z}).\]
      \end{enumerate}
  \end{lem}
  \begin{proof}
     Note that $M_{\infty}\in \Rep_{\Gamma}^{\dagger}(\BBpsi/t^n)$ (resp. $(D,\nabla_D)\in \MIC^{\dagger}(\wtR_Z/t^n)$) if and only if $M_{\infty}/t\in \Rep_{\Gamma}^{\dagger}(\widehat R_{\infty,Z})$ (resp. $(D,\nabla_D)/t\in \HIG^{\dagger}(R_Z)$). The functors in both Items (1) and (2) are well-defined. The desired cohomological comparison follows because we have isomorphisms of connections
     \[ D\otimes_{\wtr_Z^I}\wtR_Z/t^n = D\widehat \otimes_{\wtr_Z^I}\wtR_Z/t^n \cong D\widehat \otimes_{\CCtilde^I(Z)}\BBdRp(Z)/t^n.\]
  \end{proof}
  \begin{prop}\label{prop:local Robba via dR}
      Let $I\subset (\frac{1}{p-1},p)$ be a closed sub-interval as before and suppose that $1\in I$. Then for any $1\leq n\leq \infty$, the equivalence $D_{\psi}^I$ is compatible with $D_{{\rm RH},\psi}$ under base-change; equivalently, the following diagram commutes:
      \[\xymatrix@C=0.5cm{
        \Rep_{\Gamma}^{\dagger}(\CCpsi)\ar[d]^{D^I_{\psi}}\ar[rrrr]^{-\otimes_{\CCpsi}\BBpsi/t^n}&&&&\Rep_{\Gamma}^{\dagger}(\BBpsi/t^n)\ar[d]^{D_{{\rm RH},\psi}}\\
        \MIC^{\dagger}(\wtr_Z^I)\ar[rrrr]^{-\otimes_{\wtr_Z^I}\wtR_Z/t^n}&&&&\MIC^{\dagger}(\wtR_Z/t^n).
      }\] 
      Moreover, for any $M_{\infty}\in \Rep_{\Gamma}^{\dagger}(\CCpsi)$, we have quasi-isomorphisms
      \[\begin{split}
          \rR\Gamma(\Gamma,M_{\infty})\widehat \otimes_{\CCtilde^I(Z)}\BBdRp(Z)/t^n\simeq&\rD\rR(D_{\psi}^I(M_{\infty}),\nabla_{D_{\psi}^I(M_{\infty})})\widehat \otimes_{\CCtilde^I(Z)}\BBdRp(Z)/t^n\\
          \simeq&\rD\rR(D_{{\rm RH},\psi}(M_{\infty}\otimes_{\CCpsi}\BBpsi/t^n),\nabla_{D_{{\rm RH},\psi}(M_{\infty}\otimes_{\CCpsi}\BBpsi/t^n)})\\
          \simeq&\rR\Gamma(\Gamma,M_{\infty}\otimes_{\CCpsi}\BBpsi/t^n).
      \end{split}\]
  \end{prop}
  \begin{rmk}
      The morphism $\CCpsi\{\underline W\}^{\dagger}_{r,d,P}\to P_{\infty}^{\dagger}$ intertwines $\Gamma$ and $\nabla$. Comparing the explicit generators $(\prod_i\gamma_i^{-W_i/u})(e_j)$ on the Robba and $\BBdRp^+$ sides (cf.\ \eqref{equ:explicit D(M)} and \cite[Prop. 4.9(1)]{Wan23}) shows that scalar extension of $D_{\psi}^I(M_{\infty})$ recovers $D_{{\rm RH},\psi}(M_{\infty}\otimes_{\CCpsi}\BBpsi/t^n)$.
  \end{rmk}
  \begin{proof}
      Note that the natural morphism $\CCpsi\{\underline W\}^{\dagger}_{r,d,P}\to P_{\infty}^{\dagger}$ is compatible with $\Gamma$-actions and connections, cf. Proposition \ref{prop:independent of I}. It follows from definitions of $D_{\psi}^I$ and $D_{{\rm RH},\psi}$ that for any $M_{\infty}\in\Rep_{\Gamma}^{\dagger}(\CCpsi)$, there is a natural morphism of connections
      \[D_{\psi}^I(M_{\infty})\to D_{{\rm RH},\psi}(M_{\infty}\otimes_{\CCpsi}\BBpsi/t^n)\]
      which induces a natural morphism of connections
      \[\iota_{M_{\infty}}: D_{\psi}^I(M_{\infty})\otimes_{\wtr_Z^I}\wtR_Z/t^n\to D_{{\rm RH},\psi}(M_{\infty}\otimes_{\CCpsi}\BBpsi/t^n)\]
      by scalar extension. 
      
      It remains to show that $\iota_{M_{\infty}}$ is an isomorphism. For this purpose, by Nakayama's Lemma, we are reduced to the case for $n=1$. Put $M\in \Rep_{\Gamma}^{\lambda}(\Cpsi)$ such that $M\otimes_{\Cpsi}\CCpsi \cong M_{\infty}$ in the sense of Theorem \ref{thm:decompletion} and then $M/t = M\otimes_{\Cpsi}R_Z\in \Rep^{\lambda}_{\Gamma}(R_Z)$ is the representation corresponding to $M_{\infty}/t = M_{\infty}\otimes_{\CCpsi}\widehat R_{\infty,Z}$ in the sense of Theorem \ref{thm:decompletion modulo t}, by Lemma \ref{lem:small modulo t is small}. Let $e_1,\dots,e_r$ be a basis of $M$ satisfying the condition \eqref{equ:small representation} and by abuse of notation, we also denote by $e_i$'s their reduction in $M/t$. By \eqref{equ:explicit D(M)} and \cite[Prop. 4.9(1) and Rem. 4.10]{Wan23}, both $D_{\psi}^I(M_{\infty})\otimes_{\wtr_Z^I}R_Z$ and $ D_{{\rm RH},\psi}(M_{\infty}\otimes_{\CCpsi}\widehat R_{\infty,Z})$ are finite free $R_Z$-module of rank $r$ generated by
      \[(\prod_{i=1}^d\gamma_i^{-\frac{W_i}{u}})(e_1),\dots,(\prod_{i=1}^d\gamma_i^{-\frac{W_i}{u}})(e_r).\]
      Again, one can conclude that $\iota_{M_{\infty}}$ is an isomorphism by noting that for any $1\leq j\leq r$,
      \[\iota_{M_{\infty}}((\prod_{i=1}^d\gamma_i^{-\frac{W_i}{u}})(e_j)) = (\prod_{i=1}^d\gamma_i^{-\frac{W_i}{u}})(e_j).\]

      Finally, the desired cohomological comparison follows from Theorem \ref{thm:local RH over BdRp}(3) and Lemma \ref{lem:completion of small is small}(2).
  \end{proof}

\subsection{Functoriality in pull-back for $D_{{\rm RH},\psi}$}\label{ssec:local functorial in pull-back over dR}
  Similar to \S\ref{ssec:local functorial in pull-back}, it is still true that the equivalence in Theorem \ref{thm:local RH over BdRp} is compatible with pull-back.
  We again assume that we are in the setting of \S\ref{ssec:local functorial in pull-back} and keep notations therein. Moreover, let $X_Z^{\prime} = \Spa(R_Z^{\prime},R_Z^{\prime +})\to X_Z = \Spa(R_Z,R_Z^+)$ be the generic fiber of $f$ with the lifting $\wtX_Z^{\prime} = \Spf(\wtR_Z^{\prime})\to \wtX_Z = \Spf(\wtR_Z)$. Denote respectively by $\widetilde \bfB_{\widetilde \psi^{\prime},\infty}$ and $\widetilde \bfB_{\widetilde \psi^{\prime}}$ the analogue of $\BBpsi$ and $\Bpsi$ for $\wtx_Z^{\prime}$ instead of $\wtx_Z$. Then the morphism $X_{\infty,Z}^{\prime}\to X_{\infty,Z}$ induces morphisms
  \[\Bpsi\to\widetilde \bfB_{\widetilde \psi^{\prime}} \text{ and }\BBpsi\to \widetilde \bfB_{\widetilde \psi^{\prime},\infty}\]
  which are both $\Gamma^{\prime}$-equivariant.
  \begin{lem}\label{lem:local pull-back of small dR}
      Keep the notation as above and let $\lambda>1$.
      \begin{enumerate}
          \item[(1)] For any $1\leq n\leq \infty$, the base-change $M_{\infty}\mapsto M_{\infty}\otimes_{\BBpsi}\widetilde \bfB_{\widetilde \psi^{\prime},\infty}$ induces a well-defined functor
          \[\Rep_{\Gamma}^{\lambda}(\BBpsi/t^n)\to \Rep_{\Gamma^{\prime}}^{\lambda}(\widetilde \bfB_{\widetilde \psi^{\prime},\infty}/t^n).\]

          \item[(2)] For any $1\leq n\leq \infty$, the base-change $(D,\nabla_D)\mapsto (D\otimes_{\wtR_Z}\wtR_Z^{\prime},\nabla_D\otimes\id+\id\otimes\nabla)$ induces a well-defined functor
          \[\MIC^{\lambda}(\wtR_Z/t^n)\to \MIC^{\lambda}(\wtr_Z^{\prime}/t^n).\]
      \end{enumerate}
  \end{lem}
  \begin{proof}
      It suffices to deal with the case for $n=1$. In this case, it is enough to show that the base-change $M_{\infty}\mapsto M_{\infty}\otimes_{\widehat R_{\infty,Z}^+}\widehat R_{\infty,Z}^{\prime +}$ (resp. $(H,\theta_H)\mapsto (H\otimes_{\calR_Z}\calR_Z^{\prime},\theta_H\otimes\id)$) induces a well-defined functor
      \[\text{$\Rep_{\Gamma}^{\lambda}(\widehat R_{\infty,Z}^+)\to \Rep_{\Gamma}^{\lambda}(\widehat R_{\infty,Z}^{\prime+})$(resp. $\HIG^{\lambda}(\calR_Z)\to \HIG^{\lambda}(\calR_Z^{\prime})$)}.\]
      But this follows by unwinding definition.
  \end{proof}
  \begin{prop}\label{prop:local functoriality via pull-back over BdRp}
      Keep the notation as above and let $1\leq n\leq \infty$. Then the equivalence in Theorem \ref{thm:local RH over BdRp} is compatible with pull-back in the sense that the following diagram is commutative
      \[\xymatrix@C=0.5cm{
         \Rep_{\Gamma}^{\dagger}(\BBpsi/t^n)\ar[d]^{D_{{\rm RH},\psi}}\ar[rrrr]^{-\otimes_{\BBpsi}\widetilde \bfB_{\widetilde \psi^{\prime},\infty}}&&&&\Rep_{\Gamma}^{\dagger}(\widetilde \bfB_{\widetilde \psi^{\prime},\infty}/t^n)\ar[d]^{D_{{\rm RH},\psi^{\prime}}}\\
         \MIC^{\dagger}(\wtR_Z/t^n)\ar[rrrr]^{-\otimes_{\wtR_Z}\wtR_Z^{\prime}}&&&&\MIC^{\dagger}(\wtR_Z/t^n),
      }\]
      where the horizontal arrows are the functors induced by base-change, which are well-defined by Lemma \ref{lem:local pull-back of small dR}, and the vertical arrows are the equivalences established in Theorem \ref{thm:local RH over BdRp}.
  \end{prop}
  \begin{proof}
      Similar to Proposition \ref{prop:local functoriality via pull-back over Robba ring}, we can assign to any $M_{\infty}\in \Rep_{\Gamma}^{\dagger}(\BBpsi/t^n)$ a morphism 
      \[\iota_{M_{\infty}}:D_{{\rm RH},\psi}(M_{\infty})\otimes_{\wtR_Z}\wtR_Z^{\prime}\to D_{{\rm RH},\psi^{\prime}}(M_{\infty}\otimes_{\BBpsi}\widetilde \bfB_{\widetilde \psi^{\prime},\infty})\]
      which is compatible with the connections. It is enough to show this comparison morphism $\iota_{M_{\infty}}$ is an isomorphism.

      By Nakayama's Lemma, we are reduced to the case for $n=1$. Let $M\in \Rep_{\Gamma}^{\dagger}(R_Z)$ be the small representation such that $M\otimes_{R_Z}\widehat R_{\infty,Z}\cong M_{\infty}$. Then one can conclude by the same argument in the proof of Proposition \ref{prop:local functoriality via pull-back over Robba ring} after replacing \eqref{equ:explicit D(M)} by \cite[Prop. 4.9(1)]{Wan23} instead.
  \end{proof}

\section{The $p$-adic Riemann--Hilbert correspondence}\label{sec:global RH}
  This section is devoted to proving our main results (cf. Theorems \ref{thm:intro-C} and \ref{thm:intro-DR}) by reducing the problem to the local case that we have handled in \S\ref{sec:local RH} and \S\ref{sec:local RH over BdR}.
  
  Fix a $Z = \Spa(A,A^+)\in \Perfd$ and let $\frakX_Z$ be a semi-stable formal scheme of relative dimension $d$ over $A^+$ with generic fiber $X_Z$ as before. Let $I = [r,s]\subset (\frac{1}{p-1},p)$ be a closed sub-interval with $r,s\in\bQ$. We assume that $\frakX_Z$ admits a lifting $\wtx_Z$ over $\bA_{\inf}(Z)$, with associated lifting $\wtX_Z$ of $X_Z$ over $\BBdRp(Z)$, and we work with the ringed space $\wtx_Z^I$. 

\subsection{Small local systems and small connections}\label{ssec:small objects global}
  We start with some definitions.

  \begin{dfn}\label{dfn:small local system}
      Let $\bP\in\{\CCtilde^I,\BBdRp/t^n\}$ for $1\leq n\leq \infty$.
      \begin{enumerate}
          \item[(1)] A \emph{$\bP$-local system} of rank $r$ on $X_{Z,v}$ is a sheaf $\bM$ of locally finite free $\bP$-modules of rank $r$. Morphisms are morphisms of sheaves of $\bP$-modules.
          Denote the category of $\bP$-local systems by $\rL\rS(X_Z,\bP)$. Note that we make no use of integral model here.

          \item[(2)] A $\bP$-local system $\bM$ is \emph{small} if there exists an \'etale covering $\{\frakX_{Z,i} = \Spf(\calR_{Z,i})\to \frakX_Z\}$ by small affines such that for each $i$, with $X_{Z,i}$ the generic fiber of $\frakX_{Z,i}$ and $X_{\infty,Z,i}\to X_{Z,i}$ the pro-finite-\'etale $\Gamma$-torsor from Example \ref{exam:Gamma torsor}, the global section $\bM(X_{\infty,Z,i})$ lies in $\Rep_{\Gamma}^{\dagger}(\bP(X_{\infty,Z,i}))$ (smallness as in Definition \ref{dfn:small representation}). Denote by $\rL\rS^{\dagger}(X_Z,\bP)\subset \rL\rS(X_Z,\bP)$ the full subcategory of small $\bP$-local systems.
      \end{enumerate}
  \end{dfn}

  \begin{dfn}\label{dfn:small connection}
      Let $\calS\in\{\wtx_Z^I,\wtX_Z/t^n\}$ for $1\leq n\leq \infty$. Put $\bP = \CCtilde^I$ when $\calS = \wtx_Z^I$ and $\bP = \BBdRp/t^n$ when $\calS = \wtX_Z/t^n$.
      \begin{enumerate}
          \item[(1)] A \emph{flat connection} $(\calD,\nabla_{\calD})$ on $\calS$ is a sheaf $\calD$ of locally finite free $\calO_{\calS}$-modules together with a $\bP(Z)$-linear map
          \[\nabla_{\calD}:\calD\to\calD\otimes_{\calO_{\wtx_Z}}\Omega^{1,\log}_{\wtx_Z}\{-1\}\]
          satisfying the Leibniz rule with respect to
          \[\nabla:\calO_{\calS}\to \calO_{\calS}\otimes_{\calO_{\wtx_Z}}\Omega^{1,\log}_{\wtx_Z}\{-1\}\]
          and the integrability condition $\nabla_{\calD}\wedge\nabla_{\calD} = 0$. Morphisms are $\calO_{\calS}$-linear and compatible with the connections. Denote by $\MIC(\calS)$ the category of flat connections on $\calS$. For $(\calD,\nabla_{\calD})$, let $\rD\rR(\calD,\nabla_{\calD})$ be the associated de Rham complex
          \[\calD\xrightarrow{\nabla_{\calD}}\calD\otimes_{\calO_{\wtx_Z}}\Omega^{1,\log}_{\wtx_Z}\{-1\}\to\cdots\to \calD\otimes_{\calO_{\wtx_Z}}\Omega^{d,\log}_{\wtx_Z}\{-d\}.\]

          \item[(2)] A flat connection $(\calD,\nabla_{\calD})$ is \emph{small} if it is \'etale-locally a small flat connection in the sense of Definition \ref{dfn:small local connection}; concretely, there exists an \'etale covering $\{\frakX_{Z,i} = \Spf(\calR_{Z,i})\to \frakX_Z\}$ by small affines such that $(\calD,\nabla_{\calD})(\frakX_{Z,i})\in \MIC^{\dagger}(\calO_{\calS}(\frakX_{Z,i}))$. Denote by $\MIC^{\dagger}(\calS)\subset\MIC(\calS)$ the full subcategory of small flat connections.
      \end{enumerate}
  \end{dfn}
  Before stating the main theorems, we collect several functoriality lemmas used throughout the section.
  \begin{lem}\label{lem:functorial in pull-back}
      Let $f:\frakX_Z^{\prime}\to\frakX_Z$ be a smooth morphism of semi-stable formal schemes over $A^+$. Suppose that both $\frakX_Z$ and $\frakX_Z^{\prime}$ admit liftings $\wtx_Z$ and $\wtx_Z^{\prime}$ over $\bA_{\inf}(Z)$ respectively such that $f$ extends to a morphism $\widetilde f:\wtx_Z^{\prime} \to \wtx_Z$, which induces morphisms of ringed spaces $\wtx_Z^{\prime I}\to\wtx_Z^I \text{ and }\wtX_Z^{\prime}\to \wtX_Z$.

      \begin{enumerate}
          \item[(1)] Fix a $\bP\in \{\CCtilde^I,\BBdRp/t^n\}$ where $1\leq n\leq \infty$. Then taking restriction along $X^{\prime}_{Z,v}\to X_{Z,v}$ induces a functor
          \[f^*:\rL\rS^{\dagger}(X_Z,\bP)\to \rL\rS^{\dagger}(X_Z^{\prime},\bP).\]

          \item[(2)] Let $\calS^{\prime}\to \calS$ be either $\wtx_Z^{\prime I}\to\wtx_Z^I$ or $\wtX_Z^{\prime}\to \wtX_Z$. Then the base-change 
          \[(\calD,\nabla_{\calD})\mapsto (\calD\otimes_{\calO_{\calS}}\calO_{\calS^{\prime}},\nabla_{\calO}\otimes\id+\id\otimes\nabla)\]
          along $\calO_{\calS}\to \calO_{\calS^{\prime}}$ induces a functor
          \[\widetilde f^*:\MIC^{\dagger}(\calS)\to\MIC^{\dagger}(\calS^{\prime}).\]
      \end{enumerate}
  \end{lem}
  \begin{proof}
    Here, we only show that the functor
    \[f^*:\rL\rS^{\dagger}(X_Z,\CCtilde^I)\to \rL\rS^{\dagger}(X_Z^{\prime},\CCtilde^I)\]
    is well-defined because the proof in other situation is similar. Now, fix a $\bM\in \rL\rS^{\dagger}(X_Z,\CCtilde^I)$

    Since one can check the smallness \'etale locally on $\frakX_{Z,\et}^{\prime}$, we may assume that $\frakX_Z = \Spf(\calR_Z)$ is itself small affine such that $\bM(X_{\infty,Z})\in \Rep_{\Gamma}^{\dagger}(\CCtilde^I(X_{\infty,Z}))$. Then by the smoothness of $f$, one can find an \'etale covering 
    \[\{\frakX_{Z,i}^{\prime} = \Spf(\calR_{Z,i}^{\prime})\to \frakX_Z^{\prime}\}\]
    by small affines such that for each $i$, the restriction
    \[f_i:=f_{\mid_{\frakX_{Z,i}^{\prime}}}:\frakX_{Z,i}^{\prime}\to \frakX_Z\]
    is of the form \eqref{equ:relative chart-I}. (Note that the given lifting $\widetilde f$ induces a lifting $\widetilde f_i$ of $f_i$ satisfying \eqref{equ:relative chart-II}.) Denote by $X_{Z,i}$ the generic fiber of $\frakX_{Z,i}$ and let $X_{Z,\infty,i}^{\prime}\to X_{Z,i}$ be the corresponding pro-finite-\'etale $\Gamma^{\prime}$-torsor as in Example \ref{exam:relative OC}.
    Then $f_i$ induces a morphism of affinoid perfectoid spaces $X_{Z,\infty,i}^{\prime}\to X_{Z,\infty}$ such that as a $\Gamma^{\prime}$-representation
    \[(f^*\bM)(X_{\infty,Z,i}^{\prime}) = \bM(X_{\infty,Z})\otimes_{\CCtilde^I(X_{\infty,Z})}\CCtilde^I(X^{\prime}_{\infty,Z,i}).\]
    Then one can conclude by Lemma \ref{lem:local pull-back of small}(1).
  \end{proof}
  \begin{lem}\label{lem:shrink I}
      Suppose that $J\subset I$ is a closed sub-interval.
      \begin{enumerate}
          \item[(1)] The base-change along $\CCtilde^I\to \CCtilde^J$ induces a functor
          \[\rL\rS^{\dagger}(X_Z,\CCtilde^I)\to\rL\rS^{\dagger}(X_Z,\CCtilde^J).\]

          \item[(2)] The base-change along $\calO_{\wtx_Z^I}\to\calO_{\wtx_Z^J}$
          \[(\calD,\nabla_{\calD})\mapsto (\calD\otimes_{\calO_{\wtx_Z^I}}\calO_{\wtx_Z^J},\nabla_{\calD}\otimes\id+\id\otimes\nabla)\]
          induces a functor
          \[\MIC^{\dagger}(\wtx_Z^I)\to\MIC^{\dagger}(\wtx_Z^J).\]
      \end{enumerate}
  \end{lem}
  \begin{proof}
      As in the proof of Lemma \ref{lem:functorial in pull-back}, we are reduced to the case for $\frakX_Z = \Spf(\calR_Z)$ being small semi-stable and then conclude by Lemmas \ref{lem:shrink I-local representation} and \ref{lem:shrink I-local connection}.
  \end{proof}
  \begin{lem}\label{lem:CCtilde^I to BdRp}
      Suppose that $1\in I$ and let $1\leq n\leq \infty$.
      \begin{enumerate}
          \item[(1)] The base-change along $\CCtilde^I\to \BBdRp/t^n$ induces a functor
          \[\rL\rS^{\dagger}(X_Z,\CCtilde^I)\to\rL\rS^{\dagger}(X_Z,\BBdRp/t^n).\]

          \item[(2)] The base-change along $\calO_{\wtx_Z^I}\to\calO_{\wtX_Z/t^n} = \calO_{\wtX_Z}/t^n$
          \[(\calD,\nabla_{\calD})\mapsto (\calD\otimes_{\calO_{\wtx_Z^I}}\calO_{\wtx_Z/t^n},\nabla_{\calD}\otimes\id+\id\otimes\nabla)\]
          induces a functor
          \[\MIC^{\dagger}(\wtx_Z^I)\to\MIC^{\dagger}(\wtX_Z/t^n).\]
      \end{enumerate}
  \end{lem}
  \begin{proof}
      As in the proof of Lemma \ref{lem:functorial in pull-back}, we are reduced to the case for $\frakX_Z = \Spf(\calR_Z)$ being small semi-stable and then conclude by Lemma \ref{lem:completion of small is small}.
  \end{proof}

  Finally, we give a local description of $\CCtilde^I$-local systems (resp. $\BBdRp/t^n$-local systems).

  \begin{lem}\label{lem:vanishing on perfectoid}
      Let $\bP\in \{\CCtilde^I,\BBdRp/t^n\}$ where $1\leq n\leq \infty$. For any affinoid perfectoid $U = \Spa(S,S^+)\in X_{Z,v}$, any $\bM\in\rL\rS(X_Z,\bP)$ and any $n\geq 1$, we have
      \begin{equation*}
          \rH^n(U,\bM) = \left\{
            \begin{array}{rcl}
                \bM(U), & n=0 \\
                0, & n\geq 1,
            \end{array}
          \right.
      \end{equation*}
      where $\bM(U)$ is a finite projective $\bP(U)$-module such that the natural morphism
      \[\bM(U)\otimes_{\bP(U)}\bP_{\mid_U}\to\bM_{\mid_U}\]
      is an isomorphism of sheaves of $\bP_{\mid_U}$-modules.
  \end{lem}
  \begin{proof}
      For $\bP = \CCtilde^I$, the result follows from \cite[Th. 4.3.1 and Th. 4.3.3]{KL16}. For $\bP = \BBdRp/t^n$, the result follows from \cite[Lem. 8.2]{GMW24}.
  \end{proof}
  \begin{lem}\label{lem:local description of local systems}
      Suppose that $\frakX_Z = \Spf(\calR_Z)$ is semi-stable with generic fiber $X_Z = \Spa(R_Z,R_Z^+)$ and the pro-finite-\'etale $\Gamma$-torsor $X_{\infty,Z} = \Spa(\widehat R_{\infty,Z},\widehat R_{\infty,Z}^+)\to X_Z$ as in Example \ref{exam:Gamma torsor}.
      Let $\bP\in \{\CCtilde^I,\BBdRp/t^n\}$ where $1\leq n\leq \infty$. Then taking global sections $\bM\mapsto \bM(X_{\infty,Z})$ induces an equivalence of categories
      \[\rL\rS(X_Z,\bP)\to \Rep_{\Gamma}(\bP(X_{\infty,Z}))\]
      such that there is a quasi-isomorphism
      \[\rR\Gamma(X_{Z,v},\bM)\simeq \rR\Gamma(\Gamma,\bM(X_{\infty,Z})).\]
  \end{lem}
  \begin{proof}
      Denote by $X_{\infty,Z}^{\bullet}\to X_Z$ the \v Cech nerve associated to the covering $X_{\infty,Z}\to X_Z$. Then for any $\bullet\geq 0$, we have
      \[X_{\infty,Z}^{\bullet} = \Spa(\rC(\Gamma^{\bullet},\widehat R_{\infty,Z}),\rC(\Gamma^{\bullet},\widehat R_{\infty,Z}^+)),\]
      where and in what follows, for any topological ring $S$, $\rC(\Gamma^{\bullet},S)$ denotes the ring of continuous functions from $\Gamma^{\bullet}$ to $S$, yielding that
      \[\bP(X_{\infty,Z}^{\bullet}) = \rC(\Gamma^{\bullet},\bP(X_{\infty,Z})).\]
      By Lemma \ref{lem:vanishing on perfectoid}, giving a $\bP$-local system $\bM$ amounts to giving a locally finite projective $\bP(X_{\infty,Z})$-module together with a descent data corresponding to the cosimplicial ring $\bP(X_{\infty,Z}^{\bullet}) = \rC(\Gamma^{\bullet},\bP(X_{\infty,Z}))$. So the equivalence
      \[\rL\rS(X_Z,\bP)\to \Rep_{\Gamma}(\bP(X_{\infty,Z}))\]
      follows by Galois descent. Using Lemma \ref{lem:vanishing on perfectoid} again, the desired quasi-isomorphism then follows from \v Cech-to-derived spectral sequence.
  \end{proof}

\subsection{The main correspondence theorems}\label{ssec:statement of main results}
  Recall that the lifting $\wtx_Z$ of $\frakX_Z$ induces period sheaves with connection $(\calO\CCtilde^{I,\dagger},\nabla)$ and $(\calO\bB_{\dR}^{\dagger,+},\nabla)$.
  \begin{thm}[The Riemann--Hilbert correspondence over Robba ring]\label{thm:RH correspondence}
      Let $I = [r,s]\subset(\frac{1}{p-1},p)$ be a closed sub-interval. Denote by $\nu:X_{Z,v}\to \frakX_{Z,\et}$ the morphism of sites. 
      \begin{enumerate}
          \item[(1)] For any $\bM$ in $\rL\rS^{\dagger}(X_Z,\CCtilde^I)$ of rank $r$, put 
          \[\begin{split}
              \nabla_{\bM} = \id_{\bM}\otimes\nabla: \bM\otimes_{\CCtilde^I}\calO\CCtilde^{I,\dagger}\to \bM\otimes_{\CCtilde^I}\calO\CCtilde^{I,\dagger}\otimes_{\calO_{\wtx_Z}}\Omega^{1,\log}_{\wtx_Z}\{-1\}
          \end{split}\]
          Then for any $i\geq 0$, we have 
          \[\begin{split}
              \rR^i\nu_*(\bM\otimes_{\CCtilde^I}\calO\CCtilde^{I,\dagger}) = \left\{
                \begin{array}{rcl}
                    \calD^I(\bM), & i=0 \\
                    0, & i\geq 1
                \end{array}
              \right.
          \end{split}\]
          where $\calD^I(\bM)$ is a locally finite free $\calO_{\wtx_Z^I}$-module of rank $r$ and the restriction of $\nabla_{\calD^I(\bM)}$ to $\calD^I(\bM)$ defines a flat connection such that 
          \[(\calD^I(\bM),\nabla_{\calD^I(\bM)})\in \MIC^{\dagger}(\wtx_Z^I).\]

          \item[(2)] For any $(\calD,\nabla_{\calD})\in \MIC^{\dagger}(\wtx_Z^I)$ of rank $r$, put 
          \[\nabla_{(\calD,\nabla_{\calD})}:=\nabla_{\calD}\otimes\id+\id\otimes\nabla:\calD\otimes_{\calO_{\wtx_Z^I}}\calO\CCtilde^{I,\dagger}\to \calD\otimes_{\calO_{\wtx_Z^I}}\calO\CCtilde^{I,\dagger}\otimes_{\calO_{\wtx_Z}}\Omega^{1,\log}_{\wtx_Z}\{-1\}.\]
          Then $\nabla_{(\calD,\nabla_{\calD)}}\wedge\nabla_{(\calD,\nabla_{\calD})} = 0$ and for any $i\geq 0$, we have
          \begin{equation*}
              \rH^i_{\dR}(\calD\otimes_{\calO_{\wtx_Z^I}}\calO\CCtilde^{I,\dagger},\nabla_{(\calD,\nabla_{\calD})}) = \left\{
                \begin{array}{rcl}
                    \bM(\calD,\nabla_{\calD}), & i=0 \\
                    0, & i\geq 1
                \end{array}
              \right.
          \end{equation*}
          where $\bM(\calD,\nabla_{\calD})\in \rL\rS^{\dagger}(X_Z,\CCtilde^I)$ is a small $\CCtilde^I$-local system of rank $r$.

          \item[(3)] The functor $(\bM\mapsto \calD^I(\bM),\nabla_{\calD^I(\bM)})$ in Item (1) defines an equivalence of categories
          \[\calD^I:\rL\rS^{\dagger}(X_Z,\CCtilde^I)\to\MIC^{\dagger}(\wtx_Z^I)\]
          whose quasi-inverse in given by the functor $(\calD,\nabla_{\calD})\mapsto \bM(\calD,\nabla_{\calD})$ given in Item (2). Moreover, for any $\bM\in \rL\rS^{\dagger}(X_Z,\CCtilde^I)$ with associated $(\calD,\nabla_{\calD})\in \MIC^{\dagger}(\wtx_Z^I)$, there exists a quasi-isomorphism
          \[\rR\nu_*(\bM)\simeq \rD\rR(\calD,\nabla_{\calD}),\]
          yielding a quasi-isomorphism
          \[\rR\Gamma(X_{Z,v},\bM)\simeq \rR\Gamma(\frakX_{Z,\et},\rD\rR(\calD,\nabla_{\calD})).\]

          \item[(4)] The functor $\calD^I$ is functorial in smooth pull-back in the following sense: Let $f:\frakX_Z^{\prime}\to \frakX_Z$ be a smooth morphism such that $\frakX_Z^{\prime}$ admits a lifting $\wtx_Z^{\prime}$ over $\bA_{\inf}(Z)$ so that $f$ extends to a morphism $\widetilde f:\wtx_Z^{\prime}\to \wtx_Z$. Denote by $\calD^I_{\wtx_Z^{\prime}}$ and $\calD^I_{\wtx_Z}$ the equivalences obtained above for $\frakX_Z^{\prime}$ and $\frakX_Z$ respectively. Then the following diagram commutes
          \begin{equation*}
              \xymatrix@C=0.5cm{
                \rL\rS^{\dagger}(X_{Z,v},\CCtilde^I)\ar[d]^{\calD^I_{\wtx_Z}}\ar[rrrr]^{f^*}&&&&\rL\rS^{\dagger}(X^{\prime}_{Z,v},\CCtilde^I)\ar[d]^{\calD^I_{\wtx_Z^{\prime}}}\\
                \MIC^{\dagger}(\wtx_Z^I)\ar[rrrr]^{\widetilde f^*}&&&&\MIC^{\dagger}(\wtx_Z^{\prime I})
              },
          \end{equation*}
          where the horizontal arrows stand for the functors defined in Lemma \ref{lem:functorial in pull-back}.
      \end{enumerate}
  \end{thm}
  A similar result also holds true for working with coefficients over $\BdRp$. More precisely, we have the following result.
  \begin{thm}[The Riemann--Hilbert correspondence over $\BdRp$]\label{thm:RH correspondence over dR}
      Let $I = [r,s]\subset(\frac{1}{p-1},p)$ be a closed sub-interval and let $1\leq n\leq \infty$. Denote by $\nu:X_{Z,v}\to \frakX_{Z,\et}$ the morphism of sites (cf. Remark \ref{rmk:x_Z via X_Z}(1)). 
      \begin{enumerate}
          \item[(1)] For any $\bM$ in $\rL\rS^{\dagger}(X_Z,\BBdRp/t^n)$ of rank $r$, put 
          \[\begin{split}
              \nabla_{\bM} = \id_{\bM}\otimes\nabla: \bM\otimes_{\BBdRp}\calO\bB_{\dR}^{\dagger,+}\to \bM\otimes_{\BBdRp}\calO\bB_{\dR}^{\dagger,+}\otimes_{\calO_{\wtx_Z}}\Omega^{1,\log}_{\wtx_Z}\{-1\}
          \end{split}\]
          Then for any $i\geq 0$, we have 
          \[\begin{split}
              \rR^i\nu_*(\bM\otimes_{\BBdRp}\calO\bB_{\dR}^{\dagger,+}) = \left\{
                \begin{array}{rcl}
                    \widehat \calD(\bM), & i=0 \\
                    0, & i\geq 1
                \end{array}
              \right.
          \end{split}\]
          where $\widehat \calD(\bM)$ is a locally finite free module over $\calO_{\wtX_Z/t^n} = \calO_{\wtX_Z}/t^n$ of rank $r$ and the restriction of $\nabla_{\widehat \calD(\bM)}$ to $\widehat \calD(\bM)$ defines a flat connection such that 
          \[(\widehat \calD(\bM),\nabla_{\widehat \calD(\bM)})\in \MIC^{\dagger}(\wtX_Z/t^n).\]

          \item[(2)] For any $(\calD,\nabla_{\calD})\in \MIC^{\dagger}(\wtX_Z/t^n)$ of rank $r$, put 
          \[\nabla_{(\calD,\nabla_{\calD})}:=\nabla_{\calD}\otimes\id+\id\otimes\nabla:\calD\otimes_{\calO_{\wtX_Z}}\calO\bB_{\dR}^{\dagger,+}\to \calD\otimes_{\calO_{\wtX_Z}}\calO\bB_{\dR}^{\dagger,+}\otimes_{\calO_{\wtx_Z}}\Omega^{1,\log}_{\wtx_Z}\{-1\}.\]
          Then $\nabla_{(\calD,\nabla_{\calD)}}\wedge\nabla_{(\calD,\nabla_{\calD})} = 0$ and for any $i\geq 0$, we have
          \begin{equation*}
              \rH^i_{\dR}(\calD\otimes_{\calO_{\wtX_Z}}\calO\bB_{\dR}^{\dagger,+},\nabla_{(\calD,\nabla_{\calD})}) = \left\{
                \begin{array}{rcl}
                    \bM(\calD,\nabla_{\calD}), & i=0 \\
                    0, & i\geq 1
                \end{array}
              \right.
          \end{equation*}
          where $\bM(\calD,\nabla_{\calD})\in \rL\rS^{\dagger}(X_Z,\BBdRp/t^n)$ is a small $\CCtilde^I$-local system of rank $r$.

          \item[(3)] The functor $(\bM\mapsto \widehat \calD(\bM),\nabla_{\widehat \calD(\bM)})$ in Item (1) defines an equivalence of categories
          \[\widehat \calD:\rL\rS^{\dagger}(X_Z,\BBdRp/t^n)\to\MIC^{\dagger}(\wtX_Z/t^n)\]
          whose quasi-inverse in given by the functor $(\calD,\nabla_{\calD})\mapsto \bM(\calD,\nabla_{\calD})$ given in Item (2). Moreover, for any $\bM\in \rL\rS^{\dagger}(X_Z,\BBdRp/t^n)$ with associated $(\calD,\nabla_{\calD})\in \MIC^{\dagger}(\wtX_Z)$, there exists a quasi-isomorphism
          \[\rR\nu_*(\bM)\simeq \rD\rR(\calD,\nabla_{\calD}),\]
          which, as a consequence, yields a quasi-isomorphism
          \[\rR\Gamma(X_{Z,v},\bM)\simeq \rR\Gamma(\frakX_{Z,\et},\rD\rR(\calD,\nabla_{\calD})).\]

          \item[(4)] The functor $\widehat \calD$ is functorial in smooth pull-back in the following sense: Let $f:\frakX_Z^{\prime}\to \frakX_Z$ be a smooth morphism such that $\frakX_Z^{\prime}$ admits a lifting $\wtx_Z^{\prime}$ over $\bA_{\inf}(Z)$ so that $f$ extends to a morphism $\widetilde f:\wtx_Z^{\prime}\to \wtx_Z$. Denote by $\widehat \calD_{\wtx_Z^{\prime}}$ and $\widehat \calD_{\wtx_Z}$ the equivalences obtained above for $\frakX_Z^{\prime}$ and $\frakX_Z$ respectively. Then the following diagram commutes
          \begin{equation*}
              \xymatrix@C=0.5cm{
                \rL\rS^{\dagger}(X_{Z,v},\BBdRp/t^n)\ar[d]^{\widehat \calD_{\wtx_Z^{\prime}}}\ar[rrrr]^{f^*}&&&&\rL\rS^{\dagger}(X^{\prime}_{Z,v},\BBdRp/t^n)\ar[d]^{\widehat \calD_{\wtx_Z}}\\
                \MIC^{\dagger}(\wtX_Z^{\prime}/t^n)\ar[rrrr]^{\widetilde f^*}&&&&\MIC^{\dagger}(\wtX_Z/t^n)
              },
          \end{equation*}
          where the horizontal arrows stand for the functors defined in Lemma \ref{lem:functorial in pull-back}.
      \end{enumerate}
  \end{thm}

  Moreover, the above two theorems are compatible in the following sense.
  \begin{thm}[Functoriality in coefficients]\label{thm:functoriality in coefficient}
      Suppose that $I = [r,s]\subset(\frac{1}{p-1},p)$ is a closed sub-interval as before.
      \begin{enumerate}
          \item[(1)] Let $J\subset I$ be a closed sub-interval. Then the following diagram commutes
          \begin{equation*}
              \xymatrix@C=0.5cm{
                 \rL\rS^{\dagger}(X_Z,\CCtilde^I)\ar[d]^{\calD^I}\ar[rrrr]^{-\otimes_{\CCtilde^I}\CCtilde^J}&&&&\rL\rS^{\dagger}(X_Z,\CCtilde^J)\ar[d]^{\calD^J}\\
                 \MIC^{\dagger}(\wtx_Z^I)\ar[rrrr]^{-\otimes_{\calO_{\wtx_Z^I}}\calO_{\wtx_Z^J}}&&&&\MIC^{\dagger}(\wtx_Z^J),
              }
          \end{equation*}
          where the horizontal arrows stand for the base-change functors given in Lemma \ref{lem:shrink I}.

          \item[(2)] Suppose that $1\in I$ and let $1\leq n\leq \infty$. Then the functor $\calD^I$ is compatible with $\widehat \calD$ in the sense that the following diagram commutes
          \begin{equation*}
              \xymatrix@C=0.5cm{
                 \rL\rS^{\dagger}(X_Z,\CCtilde^I)\ar[d]^{\calD^I}\ar[rrrr]^{-\otimes_{\CCtilde^I}\BBdRp/t^n}&&&&\rL\rS^{\dagger}(X_Z,\BBdRp/t^n)\ar[d]^{\widehat \calD}\\
                 \MIC^{\dagger}(\wtx_Z^I)\ar[rrrr]^{-\otimes_{\calO_{\wtx_Z^I}}\calO_{\wtX_Z}/t^n}&&&&\MIC^{\dagger}(\wtX_Z/t^n),
              }
          \end{equation*}
          where the horizontal arrows stand for the base-change functors given in Lemma \ref{lem:CCtilde^I to BdRp}.
      \end{enumerate}
  \end{thm}
  \begin{proof}
      For item (1): As $J\subset I$, there is a morphism of period sheaves with connection
      \[(\calO\CCtilde^{I,\dagger},\nabla)\to (\calO\CCtilde^{J,\dagger},\nabla).\]
      By construction in Theorem \ref{thm:RH correspondence}(1), this gives a natural morphism 
      \[\calD^I(\bM)\to\calD^J(\bM),\]
      which is compatible with the connections, and thus a morphism
      \[\iota_{\bM}:\calD^I(\bM)\otimes_{\calO_{\wtx_Z^I}}\calO_{\wtx_Z^J}\to \calD^J(\bM),\]
      which is also compatible with the connections. It suffices to show that $\iota_{\bM}$ is an isomorphism. Since the problem is local on $\frakX_{Z,\et}$, we may assume that $\frakX_Z = \Spf(\calR_Z)$ is small semi-stable with generic fiber $X_Z$ and the pro-finite-\'etale $\Gamma$-torsor $X_{\infty,Z}\to X_Z$ such that $\bM(X_{\infty,Z})\in \Rep_{\Gamma}^{\dagger}(\CCtilde^I(X_{\infty,Z}))$. Then one can conclude by Lemma \ref{lem:shrink I-local representation}.

      For item (2): It follows from the same argument above by considering the morphism of period sheaves with connections
      \[(\calO\CCtilde^{I,\dagger},\nabla)\to (\calO\bB_{\dR}^{\dagger,+},\nabla)\]
      and applying Proposition \ref{prop:local Robba via dR} instead of Lemma \ref{lem:shrink I-local representation} in the local case.
  \end{proof}

  We first record a few remarks and the complete arguments are in \S\ref{ssec:proof}.

  \begin{rmk}\label{rmk:x_Z via X_Z}
   Denote by $\nu^{\prime}:X_{Z,v}\to X_{Z,\et}$ and $i:X_{Z,\et}\to\frakX_{Z,\et}$ the projections of sites so we have $\nu = i\circ \nu^{\prime}:X_{Z,\et}\to \frakX_{Z,\et}$.
          
      \begin{enumerate}
         
          \item[(1)] For any $1\leq n\leq \infty$, denote by $\wtx_Z/\xi^n[\frac{1}{p}]$ the ringed space $(\frakX_{Z,\et},\calO_{\wtx_Z[\frac{1}{p}]}/\xi^n)$, where for $n=\infty$, 
          \[\calO_{\wtx_Z[\frac{1}{p}]}/\xi^{\infty}:= \varprojlim_m\calO_{\wtx_Z[\frac{1}{p}]}/\xi^m.\]
          Replacing $\wtX_Z/t^n$ by $\wtx_Z[\frac{1}{p}]/\xi^n$ in Definition \ref{dfn:small connection}, one can define the category $\MIC^{\dagger}(\wtx_Z[\frac{1}{p}]/\xi^n)$ of small flat connections on $\wtx_Z[\frac{1}{p}]/\xi^n$.

          Note that for any vector bundle $\calE$ of rank $r$ on $\wtX_Z/t^n$ and for any $m\geq 0$, the push-forward $\rR^mi_*\calE$ is the sheafification of the presheaf
          \begin{equation*}
              \frakY_Z = \Spf(\calR_Z)\in \frakX_{Z,\et}\mapsto \rH^m(Y_{Z,\et},\calE) = \left\{
                \begin{array}{rcl}
                    \calE(Y_Z), & m = 0 \\
                    0, & m\geq 1
                \end{array}
              \right.
          \end{equation*}
          where $Y_Z$ is the generic fiber of $\frakY_Z$. So the complex $\rR i_*\calE$ is concentrated in degree $0$ such that $i_*\calE$ is a vector bundle of rank $r$ on $\wtx_Z[\frac{1}{p}]/\xi^n$. On the other hand, for any vector bundle $\calD$ of rank $r$ on $\wtx_Z[\frac{1}{p}]/\xi^n$, the pull-back $i^*\calD$ is a vector bundle on $\wtX_Z/t^n$ such that the adjunction morphism 
          \[\calD\to i_*i^*\calD\]
          is an isomorphism by above discussion. Thus, $i_*$ together with $i^*$ induces an equivalences between the categories of vector bundles on $\wtX_Z/t^n$ and $\wtx_Z[\frac{1}{p}]/\xi^n$. As a consequence, $i_*$ induces an equivalence of categories
          \[\MIC^{\dagger}(\wtX_Z/t^n)\simeq \MIC^{\dagger}(\wtx_Z[\frac{1}{p}]/\xi^n)\]
          whose quasi-inverse of $i^*$. Via this equivalence, when studying connections (resp. vector bundles, Higgs bundles, and etc.), there is no harm in identifying $\wtX_Z$ with $\wtx_Z[\frac{1}{p}]/\xi^n$. This applies to Theorem \ref{thm:RH correspondence over dR}.

          \item[(2)] On the other hand, in Theorem \ref{thm:RH correspondence}, we have to work with $\wtx_Z^I$ whose underlying site is exactly $\frakX_{Z,\et}$. Recall that $\calO_{\wtx_Z^I} = \nu_*\calO\CCtilde^{I,\dagger}$. Nevertheless, we do not know what the push-forward $\nu_*(\calO\CCtilde^{I,\dagger})$ is. 
      \end{enumerate}
  \end{rmk}

  \begin{rmk}\label{rmk:Zariski case}
      Assume that $\frakX_Z = \frakX\times_{\Spf(\calO_C)}A^+$ with $\wtx_Z = \wtx\times_{\Spf(\bfA_{\inf})}\Spf(\AAinf(Z))$ for some smooth formal scheme $\frakX$ over $\calO_C$ together with a lifting $\wtx$ over $\bfA_{\inf}$. 
      \begin{enumerate}
          \item By a result of Kedlaya (cf. \cite[Lem. 4.9]{Bha18}), there exists a Zariski covering $\{\frakX_{\alpha}\to\frakX\}$ of $\frakX$ such that each $\frakX_{\alpha}$ is small. So in this case, all results above also holds true if we replace $\frakX_{Z,\et}$ by $\frakX_{Z,\rm{Zar}}$.

          \item Put $\calO_K = \rW(\kappa)$. Suppose that $\frakX = \frakX_0\times_{\Spf(\calO_K)}\Spf(\calO_C)$ for some smooth formal scheme $\frakX_0$ over $\calO_K$ and $\wtx = \frakX_0\times_{\Spf(\calO_K)}\Spf(\bfA_{\inf})$. By the same argument in the proof of \cite[Lem. 4.9]{Bha18}), $\frakX_0$ admits a Zariski covering by small affine schemes $\frakX_{0,\alpha} = \Spf(\calR_{0,\alpha})$ over $\calO_K$. Let $X_0$ be the generic fiber of $\frakX_0$. Define 
          \[\wtX_Z^I:= \frakX_0\times_{\Spf(\calO_K)}\Spa(\CCtilde^I(Z),\CCtilde^{I,+}(Z))\]
          then we get a morphism $\iota_*:\wtX^I_{Z,\an}\to \frakX_{Z,\Zar}$ such that
          \[\iota^*\calO_{\wtx_Z^I} = \calO_{\wtX_Z^I} \text{ and }\rR \iota_*\calO_{\wtX_Z^I} = \calO_{\wtx_Z^I}.\]
          Thus, we get a fully faithful functor
          \[\iota^*:\MIC(\wtx_Z^I)\to\MIC(\wtX_Z^I)\]
          which is compatible with cohomology.
          In this case, denote by $\MIC^{\dagger}(\wtX_Z^I)$ the essential image of $\MIC^{\dagger}(\wtx_Z^I)$ under $\iota^*$.
      \end{enumerate}
  \end{rmk}

  \begin{rmk}\label{rmk:global comparison with Wang and LMNQW}
      By Proposition \ref{prop:compare with Wang} and Lemma \ref{lem:CCtilde^I to BdRp} for $n=1$, one can check that Theorems \ref{thm:RH correspondence} and Theorem \ref{thm:RH correspondence over dR} are compatible with the $p$-adic Simpson correspondence established in \cite{Wan23}. Denote by $\rL\rS^{\text{H-sm}}(X_Z,\BBdRp/t^n)$ and $\MIC^{\text{H-sm}}(\wtX_Z/t^n)$ the category of \emph{Hitchin-small} $\BBdRp/t^n$-local systems on $X_{Z,v}$ and \emph{Hitchin-small} flat connections on $\wtX_Z/t^n$, respectively (cf. \cite{LMNQW}). One can check by definition that 
      \[\rL\rS^{\dagger}(X_Z,\BBdRp/t^n)\subset \rL\rS^{\text{H-sm}}(X_Z,\BBdRp/t^n) \text{ and }\MIC^{\dagger}(\wtX_Z/t^n)\subset \MIC^{\text{H-sm}}(\wtX_Z/t^n).\]
      Using Proposition \ref{prop:compare with LMNQW}, one can check that Theorems \ref{thm:RH correspondence} and Theorem \ref{thm:RH correspondence over dR} are compatible with the $p$-adic Riemann--Hilbert correspondence established in \cite{LMNQW}.
  \end{rmk}

\subsection{Proof of main results}\label{ssec:proof}
  This subsection is devoted to proving Theorems \ref{thm:RH correspondence} and Theorem \ref{thm:RH correspondence over dR}. We begin with some lemmas.

  \begin{lem}\label{lem:almost vanishing for period sheaf}
      Suppose that $\frakX_Z$ is small semi-stable. Then for any affinoid perfectoid $U\in X_{Z,v}$, any $\bM\in \rL\rS(X_Z,\CCtilde^I)$ (resp. $\rL\rS(X_Z,\BBdRp/t^n)$), and any $i\geq 1$, we have 
      \[\text{$\rH^i(U,\bM\otimes_{\CCtilde^I}\calO\CCtilde^{I,\dagger})=0$ (resp. $\rH^i(U,\bM\otimes_{\BBdRp}\calO\bB_{\dR}^{\dagger,+})=0$)}.\]
  \end{lem}
  \begin{proof}
      We first deal with the case for $\CCtilde^I$-local systems.
      By Lemma \ref{lem:vanishing on perfectoid}, $\bM_{\mid U}$ is a direct summand of $(\CCtilde^I)^{\oplus N}$ for some $N\geq 1$. So it suffices to prove the result for $\bM = \CCtilde^I$. By Corollary \ref{cor:local description}, it is enough to show that for any $i\geq 1$,
      \[\rH^i(U,\CCtilde^I\{V_1.\dots,V_d\}^{\dagger}) = 0.\]
      By \v Cech-to-derived spectral sequence, it suffices to show that for any covering $U^0\to U$ in $X_{Z,v}$ with associated \v Cech nerve $U^{\bullet}\to U$, the corresponding \v Cech complex
      \[\cdots\to \CCtilde^I(U^{\bullet})\{V_1,\dots,V_d\}^{\dagger}\xrightarrow{\rd^{\bullet}}\CCtilde^I(U^{\bullet+1})\{V_1,\dots,V_d\}^{\dagger}\to\cdots\]
      is acyclic in degree $\geq 1$.

      Now, fix an $n\geq 0$ and let $f(\underline V) = \sum_{J\in\bN^d}f_J\underline V^J\in \CCtilde^I(U^n)\{V_1,\dots,V_d/\}^{\dagger}$ such that
      \[\rd^n(f(\underline V)) = \sum_{J\in\bN^d}\rd^n(f_J)\underline V^J = 0.\]
      It is enough to find a $g(\underline V)\in \CCtilde^I(U^n)\{V_1,\dots,V_d/\}^{\dagger}$ such tha $\rd^n(g(\underline V)) = f(\underline V)$.
      
      As there is a $\lambda>1$ such that 
      \[\lim_{|J|\to +\infty}|f_J|_I\lambda^{|J|} = 0,\]
      up to a multiplication by $p^N$ for some $N\geq 0$, we may further assume that $f_J\in \CCtilde^{I,+}(U^n)$ for any $J\in\bN^d$. Define
      \[N_J:=[\frac{|J|\log \lambda}{\log p}]\]
      and then for any $J\in\bN^d$, we have
      \[p^{N_J}\leq \lambda^{|J|}<p^{N_J+1}.\]

      As $|\frac{f_J}{p^{N_J}}|_J = |f_J|p^{N_J}$, we see that for any $J\in\bN^d$,
      \[\frac{f_J}{p^{N_J}}\in \CCtilde^{I,+}(U^n)\]
      For any $\varepsilon\in\bQ_{>0}$, by Proposition \ref{prop:C^I is a sheaf}, for any $J\in \bN^d$, there exists a $g_{J,\varepsilon}\in \CCtilde^{I,+}(U^n)$ such that 
      \[\rd^n(g_{J,\varepsilon}) = [\varpi^{\varepsilon}]\frac{f_J}{p^{N_J}}.\]
      Then for any $1<\rho<\lambda$, we have 
      \[|\frac{g_{\varepsilon,J}p^{N_J}}{[\varpi^{\varepsilon}]}|_I\rho^{|J|}\leq p^{-N_J}|[\varpi^{-\varepsilon}]|_I\rho^{|J|} = p^{\varepsilon s}\frac{\lambda^{|J|}}{p^{N_J}}\big(\frac{\rho}{\lambda}\big)^{|J|}\leq p^{1+\varepsilon s}\big(\frac{\rho}{\lambda}\big)^{|J|}.\]
      Put $g_{\varepsilon}(\underline Y):=\sum_{J\in\bN^d}\frac{g_{\varepsilon,J}p^{N_J}}{[\varpi^{\varepsilon}]}\underline V^J$ and then it is a well-defined element in $\CCtilde^I(U^n)\{V_1/\rho,\dots,V_d/\rho\}$ such that
      \[\rd^n(g_{\varepsilon}(\underline V)) = f(\underline V)\]
      as required. This completes the proof for $\CCtilde^I$-local systems.

      For $\BBdRp/t^n$-local systems, by derived Nakayama's Lemma, we are reduced to the $n=1$ case and then, by Corollary \ref{cor:local description of OBdRp^lambda}, are reduced to showing that for any $i\geq 1$,
      \[\rH^i(U,\calO\bC^{\dagger} = \OX\{V_1,\dots,V_d\}^{\dagger}) = 0.\]
      Using \cite[Prop. 8.5(iii)]{Sch-Diamond}, one can conclude by the same argument as above. (See also \cite[Lem. 5.11]{Wan23} for direct proof.) 
  \end{proof}
  \begin{lem}\label{lem:v-cohomology and group cohomology}
      Suppose that $\frakX_Z$ is small semi-stable with generic fiber $X_{\infty,Z}\to X_Z$. Denote by $X_{\infty,Z}\to X_Z$ the pro-finite-\'elate $\Gamma$-torsor in Example \ref{exam:Gamma torsor}. Then for any local system $\bM$ in $\rL\rS(X_Z,\CCtilde^I)$ (resp. $\rL\rS(X_Z,\BBdRp/t^n)$), there exists a quasi-isomorphism
      \[\begin{split}
          &\rR\Gamma(\Gamma,(\bM\otimes_{\CCtilde^I}\calO\CCtilde^{I,\dagger})(X_{\infty,Z}))\simeq \rR\Gamma(X_{Z,v},\bM\otimes_{\CCtilde^I}\calO\CCtilde^{I,\dagger})\\
          \text{(resp. }&\rR\Gamma(\Gamma,(\bM\otimes_{\BBdRp}\calO\bB_{\dR}^{\dagger,+})(X_{\infty,Z}))\simeq \rR\Gamma(X_{Z,v},\bM\otimes_{\BBdRp}\calO\bB_{\dR}^{\dagger,+})\text{)}.
      \end{split}\]
  \end{lem}
  \begin{proof}
      Using Lemma \ref{lem:almost vanishing for period sheaf} instead of Lemma \ref{lem:vanishing on perfectoid}, one can conclude by the same proof of Lemma \ref{lem:local description of local systems}.
  \end{proof}
  Now we are able to prove Theorems \ref{thm:RH correspondence} and \ref{thm:RH correspondence over dR}.
  \begin{proof}[\textbf{Proof of Theorem \ref{thm:RH correspondence}}]
      For item (1): Fix a $\bM\in \rL\rS^{\dagger}(X_Z,\CCtilde^I)$ of rank $r$. Let $\{\frakX_{Z,i} \to \frakX_Z\}$ be the \'etale covering by small affines as in Definition \ref{dfn:small local system}. Let $\Sigma_{\bM}$ be the sub-category of $\frakX_{Z,\et}$ consisting of small affines $\frakY_Z$ whose structure morphism $\frakY_Z\to \frakX_Z$ factors through $\frakX_{Z,i}$ for some $i$. Then $\Sigma_{\bM}$ forms a basis for $\frakX_{Z,\et}$.
      
      Recall that we have to show that the complex $\rR\nu_*(\bM\otimes_{\CCtilde^I}\calO\CCtilde^{I,\dagger})$ is concentrated in degree $0$ such that 
      \[(\calD^I(\bM),\nabla_{\calD^I(\bM)}) := \nu_*(\bM\otimes_{\CCtilde^I}\calO\CCtilde^{I,\dagger},\nabla_{\bM})\]
      is an object in $\MIC^{\dagger}(\wtx_Z^I)$ of rank 
      $r$. Note that for any $n\geq 0$, $\rR^n\nu_*(\bM\otimes_{\CCtilde^I}\calO\CCtilde^{I,\dagger})$ is the sheafification of the presheaf
      \[\frakY_Z=\Spf(\calR_Z)\in \Sigma_{\bM}\mapsto \rH^n(Y_{Z,v},\bM\otimes_{\CCtilde^I}\calO\CCtilde^{I,\dagger}),\]
      where $Y_Z = \Spa(R_Z,R_Z^+)$ admits the generic fiber of $\frakY_Z$.
      So it suffices to show that for any $\frakY_Z\in\Sigma_{\bM}$, the complex
      \[\rR\Gamma(Y_{Z,v},\bM\otimes_{\CCtilde^I}\calO\CCtilde^{I,\dagger})\]
      is concentrated in degree $0$ such that 
      \[D:=\rH^0(Y_{Z,v},\bM\otimes_{\CCtilde^I}\calO\CCtilde^{I,\dagger})\]
      is a finite free module of rank $r$ over $\wtr_Z^I = \calO_{\wtx_Z^I}(\frakY_Z)$ such that the connection on $D$
      \[\nabla_D:D\to D\otimes_{\wtr_Z}\Omega^{1,\log}_{\wtr_Z}\{-1\}\]
      induced by the restriction of $\nabla_{\bM}$ defines an object $(D,\nabla_D)\in \MIC^{\dagger}(\wtr_Z^I)$.

      Let $Y_{\infty,Z}\to Y_Z$ be the pro-finite-\'etale $\Gamma$-torsor as in Example \ref{exam:Gamma torsor}. By Lemma \ref{lem:v-cohomology and group cohomology}, we are reduced to showing that the complex
      \[\rR\Gamma(\Gamma,(\bM\otimes_{\CCtilde^I}\calO\CCtilde^{I,\dagger})(Y_{\infty,Z}))\]
      is concentrated in degree $0$ such that 
      \[D:=\rH^0(\Gamma,\bM(Y_{\infty,Z})\otimes_{\CCtilde^I(Y_{\infty,Z})}\calO\CCtilde^{I,\dagger}(Y_{\infty,Z}))\]
      with the connection $\nabla_D$ induced by $\nabla_{\bM}$ defines an object $(D,\nabla_D)\in \MIC^{\dagger}(\wtr_Z^I)$ of rank $r$. As $\frakY_Z\in\Sigma_{\bM}$, which implies that $\bM(Y_{\infty,Z})\in \Rep_{\Gamma}^{\dagger}(\CCtilde^I(Y_{\infty,Z}))$, one can conclude by Theorem \ref{thm:local correspondence over Robba ring}(1). 

      For item (2): Fix a $(\calD,\nabla_{\calD})\in \MIC^{\dagger}(\wtx_Z^I)$ of rank $r$.  Let $\{\frakX_{Z,i} \to \frakX_Z\}$ be the \'etale covering by small affines as in Definition \ref{dfn:small connection}. Let $\Sigma_{(\calD,\nabla_{\calD})}$ be the sub-category of $\frakX_{Z,\et}$ consisting of small affines $\frakY_Z$ whose structure morphism $\frakY_Z\to \frakX_Z$ factors through $\frakX_{Z,i}$ for some $i$ as in Item (1). Then $\Sigma_{(\calD,\nabla_{\calD})}$ forms a basis for $\frakX_{Z,\et}$ as well. 

      It is easy to see that $\nabla_{(\calD,\nabla_{\calD})}\wedge \nabla_{(\calD,\nabla_{\calD})}$ because $\nabla_{\calD}\wedge\nabla_{\calD} = 0$. We have to show that the de Rham complex $\rD\rR(\calD\otimes_{\calO_{\wtx_Z^I}}\calO\CCtilde^{I,\dagger},\nabla_{(\calD,\nabla_{\calD})})$
      \[\calD\otimes_{\calO_{\wtx_Z^I}}\calO\CCtilde^{I,\dagger}\xrightarrow{\nabla_{(\calD,\nabla_{\calD})}}\calD\otimes_{\calO_{\wtx_Z^I}}\calO\CCtilde^{I,\dagger}\otimes_{\calO_{\wtx_Z}}\Omega^{1,\log}_{\wtx_Z}\{-1\}\to \cdots\]
      is acyclic in degree $\geq 1$ and 
      \[\bM:=(\calD\otimes_{\calO_{\wtx_Z^I}}\calO\CCtilde^{I,\dagger})^{\nabla_{(\calD,\nabla_{\calD})} = 0}\]
      is an object in $\rL\rS^{\dagger}(X_Z,\CCtilde^I)$ of rank $r$.

      Since the problem is local on $\frakX_{Z,\et}$, we may work with $\frakY_Z\in \Sigma_{(\calD,\nabla_{\calD})}$. Since the problem is also local on $Y_{Z,v}$, we may work locally on $Y_{\infty,Z,v}$, where $Y_{\infty,Z}\to Y_Z$ is the pro-finite-\'etale $\Gamma$-torsor as in Example \ref{exam:Gamma torsor}. By Lemma \ref{lem:vanishing on perfectoid} and Lemma \ref{lem:almost vanishing for period sheaf}, it suffices to show that the complex
      \[\begin{split}
          &\rD\rR(\calD(\frakY_Z)\otimes_{\calO_{\wtx_Z}(\frakY_Z)}\calO\CCtilde^{I,\dagger}(Y_{\infty,Z}),\rd_{\calD(\calY_Z)}:=\nabla_{\calD(\frakY_Z)}\otimes\id+\id\otimes\nabla) \\
          =&\Gamma\big(Y_{\infty,Z},\rD\rR(\calD\otimes_{\calO_{\wtx_Z^I}}\calO\CCtilde^{I,\dagger},\nabla_{(\calD,\nabla_{\calD})})\big)
      \end{split}\]
      is concentrated in degree $0$ such that 
      \[M_{\infty} = \rH^0_{\dR}(\calD(\frakY_Z)\otimes_{\calO_{\wtx_Z}(\frakY_Z)}\calO\CCtilde^{I,\dagger}(Y_{\infty,Z}),\rd_{\calD(\frakY_Z)})\]
      is an object in $\Rep_{\Gamma}^{\dagger}(\CCtilde^I(Y_{\infty,Z}))$ of rank $r$. But this follows from Theorem \ref{thm:local correspondence over Robba ring}(2). 

      For item (3): We first show that $\calD^I$ does induce an equivalence of categories
      \[\rL\rS^{\dagger}(X_Z,\CCtilde^I)\xrightarrow{\simeq}\MIC^{\dagger}(\wtx_Z^I) \]
      with the desired quasi-inverse. 

      For any $\bM\in \rL\rS^{\dagger}(X_Z,\CCtilde^I)$ (resp. $(\calD,\nabla_{\calD})\in \MIC^{\dagger}(\wtx_Z^I)$), by the same argument as in the proof of Proposition \ref{prop:integral equivalence}, we have a natural morphism 
      \begin{equation*}
          \begin{split}
              &\iota_{\bM}:\bM(\calD^I(\bM),\nabla_{\calD^I(\bM)})\to \bM\\
              \text{(resp. }&\iota_{(\calD,\nabla_{\calD})}:(\calD^I(\bM(\calD,\nabla_{\calD})),\nabla_{\calD^I(\bM(\calD,\nabla_{\calD}))})\to (\calD,\nabla_{\calD})\text{)}
          \end{split}
      \end{equation*}
      of $\CCtilde^I$-local systems (resp. flat connections on $\wtx_Z^I$). It suffices to show that both $\iota_{\bM}$ and $\iota_{(\calD,\nabla_{\calD})}$ are isomorphisms. But this is a local problem on $\frakX_{Z,\et}$ (and on $X_{Z,v}$ for $\iota_{\bM}$). So as above, we are reduced to the small affine case and then conclude by Theorem \ref{thm:local correspondence over Robba ring}(3) (and the proof of Proposition \ref{prop:integral equivalence}). 
      
      It remains to show that 
      \[\rR\nu_*(\bM)\simeq \rD\rR(\calD^I(\bM),\nabla_{\calD^I(\bM)}).\]
      By Poincar\'e's Lemma (cf. Theorem \ref{thm:Poincare Lemma for OC^I}), there is a quasi-isomorphism
      \[\rR\nu_*(\bM)\xrightarrow{\simeq}\rR\nu_*(\rD\rR(\bM\otimes_{\CCtilde^I}\calO\CCtilde^{I,\dagger},\nabla_{\bM})).\]
      By what we have seen in Item (1), the natural map
      \[(\calD^I(\bM),\nabla_{\calD^I(\bM)}) = \nu_*(\rD\rR(\bM\otimes_{\CCtilde^I}\calO\CCtilde^{I,\dagger},\nabla_{\bM}))\to \rR\nu_*(\rD\rR(\bM\otimes_{\CCtilde^I}\calO\CCtilde^{I,\dagger},\nabla_{\bM}))\]
      is a quasi-isomorphism. Combining these two quasi-isomorphisms together, we get the desired quasi-isomorphism
      \[\rR\nu_*(\bM)\simeq (\calD^I(\bM),\nabla_{\calD^I(\bM)}).\]

      For item (4): Fix an $\bM\in \rL\rS^{\dagger}(X_Z,\CCtilde^I)$. Using the natural morphism of period sheaves with connection \eqref{equ:relative morphism on OC^Idagger}, by the same argument as in the proof of Proposition \ref{prop:local functoriality via pull-back over Robba ring}, we get a natural morphism
      \[\iota_{\bM}:\widetilde f^*(\calD^I_{\wtx_Z}(\bM))\to \calD_{\wtx_Z^{\prime}}(f^*(\bM))\]
      which is compatible with the connections. It is enough to show that $\iota_{\bM}$ is an isomorphism. Note that this is a local problem on $\frakX_{Z,\et}$. So we are reduced to the same local case as in the proof of Item (1) and then conclude by Proposition \ref{prop:local functoriality via pull-back over Robba ring}.

      The proof is complete.
  \end{proof}
  \begin{proof}[\textbf{Proof of Theorem \ref{thm:RH correspondence over dR}}]
      Using Theorem \ref{thm:local RH over BdRp}, Proposition \ref{prop:local functoriality via pull-back over BdRp} and the morphism of period sheaves with connection \eqref{equ:relative morphism on OBdRdagger} instead of Theorem \ref{thm:local correspondence over Robba ring}, Proposition \ref{prop:local functoriality via pull-back over Robba ring} and \eqref{equ:relative morphism on OC^Idagger}, respectively, we can conclude the result by the same argument in the proof of Theorem \ref{thm:RH correspondence}.
  \end{proof}

\section{The local Frobenius action}\label{sec:local Frobenius action}
  \subsection{Prismatic liftings}\label{ssec:prismatic liftings}
  Throughout this section, fix $Z=\Spa(A,A^+)\in \Perfd$ and let $\frakX_Z$ be a semi-stable formal scheme over $A^+$ as before. We explain how Frobenius on a prismatic lift interacts with period sheaves and prove that the Riemann--Hilbert correspondence is compatible with these Frobenius actions.
  \begin{dfn}\label{dfn:prismatic lifting}
      A flat lifting $\wtx$ of $\frakX_Z$ over $\AAinf(Z)$ is called a \emph{prismatic lifting} if there exists a Frobenius endomorphism
      \[\varphi_{\wtx_Z}:\widetilde \frakX_Z\to \widetilde \frakX_Z\]
      compatible with the Frobenius on $\AAinf(Z)$, whose reduction modulo $p$ agrees with the absolute Frobenius of $\frakX_Z^{\flat} = \wtx_Z\times_{\Spf\AAinf(Z)}\Spf(A^{\flat})$, and which is compatible with the log-structure on $\wtx_Z$.
  \end{dfn}
  \begin{rmk}\label{rmk:local prismatic lifting}
      Existence of a prismatic lifting is a non-trivial geometric condition: it always holds for small affine formal schemes (Example \ref{exam:prismatic lifting}), while in general it imposes global constraints on the formal model.

      Assume $\frakX_Z = \Spf(\calR_Z)$ is small semi-stable with lifting $\wtx_Z = \Spf(\wtr_Z)$. Since $\wtr_Z$ is $p$-torsion free, specifying a prismatic lifting is equivalent to equipping $\wtr_Z$ with a $\delta_{\log}$-structure in the sense of \cite{Kos20}. To see this, let $\varphi:\wtr_Z\to \wtr_Z$ be the morphism lifting the absolute Frobenius map
      \[\wtr_Z/p=:\wtr_Z^{\flat}\xrightarrow{x\mapsto x^p}\wtr_Z^{\flat}\]
      associated to $\varphi_{\wtx_Z}$ and denote by $f:M\to M$ the morphism of monoids corresponding to $\varphi_{\wtx_Z}$, where $M$ denotes the log-structure on $\wtr_Z$. Then for any $m\in M$, we have $\iota(m) \equiv m^p\mod \wtr_Z^{\times}$. Write $\iota(m) = m^pu(m)$ for some $u(m)\in \wtr_Z^{\times}\subset M$. As $\varphi_{\wtx_Z}$ lifts the absolute Frobenius of $\wtx_Z^{\flat}$, we have $u(m)\equiv 1\mod p$, yielding that $u(m) =1+p\delta_{\log}(m)$ for some $\delta_{\log}(m)\in \wtr_Z$. For any element $x\in \wtr_Z$ such that $\iota(m) = m^p(1+px)$, as $M$ is integral, we have
      \[1+p\delta_{\log}(m) = 1+px,\]
      yielding that $\delta_{\log}(m) = x$ as $\wtr_Z$ is $p$-torsion free. So we get a well-defined morphism
      \[\delta_{\log}:M\to \wtr_Z.\]
      Using $\wtr_Z$ is $p$-torsion free again, it is easy to see that $(\wtr_Z,(\xi),M,\delta_{\log})$ defines a log-prism in the sense of \cite{Kos20}. If $T_0,\dots, T_r, T_{r+1},\dots, T_d$ denote the coordinates on $\wtr_Z$, then we have $\varphi(T_i) = T_i^p(1+p\delta_{\log}(T_i))$ for any $0\leq i\leq d$. Recall that in the chart we have $T_0\cdots T_r = [\varpi^a]$ (or $p^a$ in the non-lifted setting); this relation forces the $\delta_{\log}$-structure on the product to satisfy $(1+p\delta_{\log}(T_0))\cdots(1+p\delta_{\log}(T_r)) = 1$, as used in Example \ref{exam:Frobenius on OC^I}.
  \end{rmk}
  \begin{exam}\label{exam:prismatic lifting}
      Every small semi-stable affine $\frakX_Z = \Spf(\calR_Z)$ with lifting $\wtx_Z = \Spf(\wtr_Z)$ admits a prismatic lifting (cf.\ Example \ref{exam:Gamma torsor}). Equip
      \[\AAinf(Z)\za T_0,\dots,T_r,T_{r+1}^{\pm 1},\dots, T_d^{\pm 1}\ya/\za T_0\cdots T_r-[\varpi^a]\ya\]
      with Frobenius $\varphi$ compatible with $\AAinf(Z)$ and satisfying $\varphi(T_i) = T_i^p$ for $0\leq i\leq d$. The chart on $\frakX_Z$ induces a $(p,\xi)$-completely \'etale morphism from this algebra to $\wtr_Z$, so $\varphi$ extends uniquely to $\varphi:\wtr_Z\to\wtr_Z$. This is the $\delta_{\log}$-structure with $\delta_{\log}(T_i) = 0$ for $0\leq i\leq d$.
  \end{exam}
  
  \subsection{The Frobenius action on $\calO\CCtilde^I$}\label{ssec:phi on OC^I}
  For any prismatic lifting $\wtx_Z$ of $\frakX_Z$, the Frobenius $\varphi:\calO_{\wtx_Z}\to\calO_{\wtx_Z}$ induces $\varphi:\Omega^1_{\wtx_Z}\to \Omega^1_{\wtx_Z}$ by sending $\rd f$ to $\rd\varphi(f)$ for local sections $f\in \calO_{\wtx_Z}$.
  For closed sub-intervals $I,J\subset (0,+\infty)$ with $J\subset I/p$, $\varphi$ further induces $\varphi:\calO_{\wtx_Z^I}\to \calO_{\wtx_Z^J}$ compatible with $\CCtilde^I(Z)\xrightarrow{\phi}\CCtilde^{I/p}(Z)\to\CCtilde^J(Z)$, and $\varphi:\Omega^1_{\wtx_Z^I}\to\Omega^1_{\wtx_Z^J}$ making the following diagram commutative:
  \begin{equation*}
      \begin{tikzcd}
          \calO_{\wtx_Z^I} \arrow[d, "\varphi"] \arrow[r, "\rd"] & \Omega^1_{\wtx_Z^I} \arrow[d, "\varphi"] \\
          \calO_{\wtx_Z^J} \arrow[r, "\rd"] & \Omega^1_{\wtx_Z^J}. 
      \end{tikzcd}
  \end{equation*}

  The goal here is to show that for a suitable closed interval $I\subset (0,+\infty)$ and a prismatic lifting $\wtx_Z$, the Frobenius $\varphi_{\wtx_Z}$ induces $\varphi:\calO\bC^{I,\dagger}\to\calO\bC^{I/p,\dagger}$ compatible with $\phi:\CCtilde^I\to\CCtilde^{I/p}$ and with the connection, so that Theorem \ref{thm:RH correspondence} is compatible with Frobenius whenever $\frakX_Z$ admits a prismatic lifting.

  We construct the Frobenius morphism $\varphi:\calO\bC^{I,\dagger}\to\calO\bC^{I/p,\dagger}$ associated to a prismatic lifting $\wtx_Z$.

  We first treat the case $\frakX_Z = \frakY_Z = \Spf(\calR_Z)$, so that $\wtx_Z = \wty_Z = \Spf(\wtr_Z)$. As discussed in Remark \ref{rmk:local prismatic lifting}, the Frobenius $\varphi_{\wtx_Z}$ defines a log-prism $(\wtr_Z,(\xi))$ in the sense of \cite{Kos20}. Let $U=\Spa(S,S^+)\in X_v$ be affinoid perfectoid as before. Then there exists a unique $\delta_{\log}$ structure on $\AAinf(U)$ (which is equipped with the canonical log-structure) by letting $\delta_{\log}([x]) = 0$ for any $x\in S^{\flat\times}\cap S^{\flat+}$. It follows from \cite{Kos20} that the $\delta_{\log}$-structure on $\wtr_Z\widehat \otimes_{\AAinf(Z)}\AAinf(U)$ uniquely extends to a $\delta_{\log}$-structure on 
  \[\wtr_Z\widehat \otimes_{\AAinf(Z)}\AAinf(U)[(\frac{T_0}{t_0})^{\pm 1},\dots,(\frac{T_r}{t_r})^{\pm 1}]^{\wedge_{(p,\xi)}}.\]
  More precisely, if we write 
  \[\varphi(T_i) = T_i^p(1+p\delta_{\log}(T_i)) \text{ and }\varphi(t_i) = t_i^p(1+p\delta_{\log}(t_i)),\]
  then $\delta_{\log}(\frac{T_i}{t_i})$ is uniquely determined by
  \[1+p\delta_{\log}(\frac{T_i}{t_i}) = \frac{1+p\delta_{\log}(T_i)}{1+p\delta_{\log}(t_i)}\]
  by noting that
  \[(\frac{T_i}{t_i})^p(1+p\delta_{\log}(\frac{T_i}{t_i})) = \varphi(\frac{T_i}{t_i}) = \frac{\varphi(T_i)}{\varphi(t_i)} = \frac{T_i^p(1+p\delta_{\log}(T_i))}{t_i^p(1+p\delta_{\log}(t_i))} = (\frac{T_i}{t_i})^p\frac{1+p\delta_{\log}(T_i)}{1+p\delta_{\log}(t_i)}.\]
  Recall that $A_{\frakY_Z} = A_{\frakY_Z}^+[\frac{1}{p}]$ with $A_{\frakY_Z}^+:=\wtr_Z\widehat \otimes_{\AAinf(Z)}\AAinf(U)[(\frac{T_0}{t_0})^{\pm 1},\dots,(\frac{T_r}{t_r})^{\pm 1}]^{\wedge_{(p,\xi)}}$ and that for any closed interval $I\subset(0,+\infty)$, $\Spa(A_{\frakY_Z}^I,A_{\frakY_Z}^{I,+})$ is the rational localization of $\Spa(A_{\frakY_Z},A_{\frakY_Z}^{\circ})$ with respect to the condition $|[\varpi]|\in |p|^I$.
  \begin{lem}\label{lem:phi is isometry}
      \begin{enumerate}
          \item[(1)] For any $x\in A_{\frakY_Z}$, we have $|\varphi(x)|_{I/p} = |x|_I$.

          \item[(2)] For any $x\in A_{\frakY_Z}^I$, we have $|\varphi(x)|_{I/p} = |x|_I$.
      \end{enumerate}
  \end{lem}
  \begin{proof}
      Item~(1) is an isometry on $A_{\frakY_Z}$ for the $I$- and $I/p$-norms; Item~(2) follows because $A_{\frakY_Z}$ is dense in $A_{\frakY_Z}^I$ and both norms are continuous.
      For Item (1), for any $0\leq i\leq d$ we have
      \[|\varphi(\frac{T_i}{t_i})|_{I/p} = 1 = |\frac{T_i}{t_i}|_I.\]
      Then the result follows from Lemma \ref{lem:second description of C^I}(3) and Lemma \ref{lem:property on I-norm}(5).
  \end{proof}
  
  \begin{lem}\label{lem:extends to blow-up algebra}
      Suppose $J = [r,s]\subset (\frac{1}{p},+\infty)\cap \frac{I}{p}$ is a non-empty closed interval (this condition ensures that overconvergence is preserved under $\varphi$, i.e., the radius of convergence does not shrink too rapidly). Then the morphism $\varphi: A_{\frakY_Z}^I\to A_{\frakY_Z}^J$ given by the composite
      \[A_{\frakY_Z}^I\xrightarrow{\varphi}A_{\frakY_Z}^{I/p}\to A_{\frakY_Z}^J\]
      uniquely extends to a morphism
      \[\varphi:A_{\frakY_Z}^I[\frac{\ker(\theta_{\frakY_Z}^{\prime})}{\xi}]\to A_{\frakY_Z}^J[\frac{\ker(\theta_{\frakY_Z}^{\prime})}{\xi}],\]
      where $\ker(\theta_{\frakY_Z}^{\prime})$ denotes the kernel of the surjection $\theta_{\frakY_Z}^{\prime}:A_{\frakY_Z}^+\to S^+$ as before.
  \end{lem}
  \begin{proof}
      \emph{Idea.} For $x\in\ker(\theta_{\frakY_Z}^{\prime})$ one has $\varphi(x)=x^p+p\delta(x)$; in the blow-up algebra this rewrites as \eqref{equ:extends to blow-up algebra}, so extension is well-defined once $p$ lies in the ideal $(\phi(\xi))$ (Lemma \ref{lem:phi(xi)/p is a unit}).

      It suffices to show that the morphism
      \[A_{\frakY_Z}^I\xrightarrow{\varphi}A_{\frakY_Z}^J\to A_{\frakY_Z}^{I/p}[\frac{\ker(\theta_{\frakY_Z}^{\prime})}{\xi}]\]
      sends $\ker(\theta_{\frakY_Z}^{\prime})$ into the ideal generated by $\varphi(\xi)$. For any $x\in \ker(\theta_{\frakY_Z}^{\prime})$, write 
      $\varphi(x) = x^p+p\delta(x)$. So in $A_{\frakY_Z}^J[\frac{\ker(\theta_{\frakY_Z}^{\prime})}{\xi}]$, we have
      \begin{equation}\label{equ:extends to blow-up algebra}
          \varphi(x) = (\frac{x}{\xi})^p\xi^p+p\delta(x) = (\frac{x}{\xi})^p\phi(\xi)+p\delta(x)-p\delta(\xi) (\frac{x}{\xi})^p.
      \end{equation}
      So it suffices to show that $p$ belongs to the ideal generated by $\phi(\xi)$. This follows from the lemma below.
  \end{proof}
  \begin{lem}\label{lem:phi(xi)/p is a unit}
      Let $J = [r,s]\subset(\frac{1}{p},+\infty)$ be a closed interval. Then $\frac{\phi(\xi)}{p}$ is a unit in $\Ctilde^{J,+}$.
  \end{lem}
  \begin{proof}
      This follows from Example \ref{exam:t in C^I}(3)(ii) together with Lemma \ref{lem:second description of C^I}(3).
  \end{proof}
  \begin{rmk}
      Lemma \ref{lem:phi(xi)/p is a unit} is the key point for extending $\varphi$ to $A_{\frakY_Z}^I[\frac{\ker(\theta_{\frakY_Z}^{\prime})}{\xi}]$: having $\frac{\phi(\xi)}{p}$ a unit means one can ``divide by $\xi$'' after Frobenius without leaving the overconvergent tower, so the induced map on the blow-up is well-defined.
  \end{rmk}
  \begin{lem}\label{lem:extends to blow-up algebra-norm}
      Suppose that $J = [r,s]\subset (\frac{1}{p},+\infty)\cap \frac{I}{p}$ is a closed interval. Then for any $x\in \ker(\theta_{\frakY_Z}^{\prime})$, we have
      \[|\varphi(\frac{x}{\xi})|_J\leq |\frac{x}{\xi}|_I\]
      and the equality holds when $I = pJ$.
  \end{lem}
  \begin{proof}
      As $pJ\subset I$, we have $|\cdot|_{pJ}\leq|\cdot|_I$. Thus, it suffices to deal with the case for $pJ = I$. In this case, we have 
      \[|\varphi(\frac{x}{\xi})|_J = |\frac{\varphi(x)}{p}|_J = |\varphi(\frac{x}{p})|_J = |\frac{x}{p}|_{pJ} = |\frac{x}{\xi}\cdot\frac{\xi}{p}|_{pJ} = |\frac{x}{\xi}|_{pJ}\]
      as expected. Here, the last equality follows as $\frac{\xi}{p}$ is a unit in $\Ctilde^{pJ,+}$ by Example \ref{exam:t in C^I}(3)(ii).
  \end{proof}
  \begin{cor}\label{cor:extends to blow-up algebra-norm}
      Let $I,J$ be closed sub-intervals of $(0,+\infty)$ such that $J\subset(\frac{1}{p},+\infty)\cap I/p$. Then for any $x\in A_{\frakY_Z}^I[\frac{\ker(\theta_{\frakY_Z}^{\prime})}{\xi}]$, we have $|\varphi(x)|_J\leq |x|_I$, where the equality holds when $I = pJ$. As a consequence, we get a well-defined morphism
      \[\varphi:(A_{\frakY_Z}^I[\frac{\ker(\theta_{\frakY_Z}^{\prime})}{\xi}],A_{\frakY_Z}^I[\frac{\ker(\theta_{\frakY_Z}^{\prime})}{\xi}]^{\circ})\to (A_{\frakY_Z}^J[\frac{\ker(\theta_{\frakY_Z}^{\prime})}{\xi}],A_{\frakY_Z}^J[\frac{\ker(\theta_{\frakY_Z}^{\prime})}{\xi}]^{\circ})\]
      which is compatible with the Frobenius map $\varphi:\wtr_Z\to\wtr_Z$ and the composite $\CCtilde^I\xrightarrow{\phi}\CCtilde^{I/p}\to\CCtilde^J$.
  \end{cor}
  \begin{proof}
       This follows directly from Lemma \ref{lem:phi is isometry}(2) and Lemma \ref{lem:extends to blow-up algebra-norm}.
  \end{proof}
  Recall that for any $1\leq \lambda<|\xi|_I^{-1}$, the rational localization 
  \[\Spa(A_{\frakY_Z}^{I,\lambda},A_{\frakY_Z}^{I,\lambda,+})\to \Spa(A_{\frakY_Z}^I[\frac{\ker(\theta_{\frakY_Z}^{\prime})}{\xi}],A_{\frakY_Z}^I[\frac{\ker(\theta_{\frakY_Z}^{\prime})}{\xi}]^{\circ})\]
  is defined by requiring that $|\frac{\ker(\theta_{\frakY_Z}^{\prime})}{\xi}|\leq \lambda$.
  \begin{lem}\label{lem:extend to period ring}
      Let $I,J$ be two closed sub-intervals of $(0,+\infty)$ such that $J\subset(\frac{1}{p},+\infty)\cap I/p$. Then for any $1\leq \lambda<|\xi|_I^{-1}$, the morphism
      \[\varphi:(A_{\frakY_Z}^I[\frac{\ker(\theta_{\frakY_Z}^{\prime})}{\xi}],A_{\frakY_Z}^I[\frac{\ker(\theta_{\frakY_Z}^{\prime})}{\xi}]^{\circ})\to (A_{\frakY_Z}^J[\frac{\ker(\theta_{\frakY_Z}^{\prime})}{\xi}],A_{\frakY_Z}^J[\frac{\ker(\theta_{\frakY_Z}^{\prime})}{\xi}]^{\circ})\]
      extends to a morphism
      \[\varphi:(A_{\frakY_Z}^{I,\lambda^p},A_{\frakY_Z}^{I,\lambda^p,+})\to (A_{\frakY_Z}^{J,\lambda},A_{\frakY_Z}^{J,\lambda,+})\]
      which is compatible with the Frobenius map $\varphi:\wtr_Z\to\wtr_Z$ and the composite $\CCtilde^I\xrightarrow{\phi}\CCtilde^{I/p}\to\CCtilde^J$.
  \end{lem}
  \begin{proof}
      \emph{Idea.} Rational localization asks that $|\frac{y}{\xi}|\leq\lambda$ imply $|\varphi(\frac{y}{\xi})|\leq\lambda^p$; this follows from \eqref{equ:extends to blow-up algebra} and the fact that $\frac{p}{\phi(\xi)}$ is a unit (Lemma \ref{lem:phi(xi)/p is a unit}), so $\varphi(\frac{y}{\xi})$ is controlled by $(\frac{y}{\xi})^p$ up to terms multiplied by $p/\phi(\xi)$.

      It is enough to show that the morphism
      \[\Spa(A_{\frakY_Z}^{J,\lambda},A_{\frakY_Z}^{J,\lambda,+})\to \Spa(A_{\frakY_Z}^J[\frac{\ker(\theta_{\frakY_Z}^{\prime})}{\xi}],A_{\frakY_Z}^J[\frac{\ker(\theta_{\frakY_Z}^{\prime})}{\xi}]^{\circ})\to \Spa(A_{\frakY_Z}^I[\frac{\ker(\theta_{\frakY_Z}^{\prime})}{\xi}],A_{\frakY_Z}^I[\frac{\ker(\theta_{\frakY_Z}^{\prime})}{\xi}]^{\circ})\]
      factors through $\Spa(A_{\frakY_Z}^{I,\lambda^p},A_{\frakY_Z}^{I,\lambda^p,+})$. In other words, we need to show that for any point 
      \[|\cdot|\in \Spa(A_{\frakY_Z}^J[\frac{\ker(\theta_{\frakY_Z}^{\prime})}{\xi}],A_{\frakY_Z}^J[\frac{\ker(\theta_{\frakY_Z}^{\prime})}{\xi}]^{\circ})\]
      satisfying $|\frac{x}{\xi}|\leq \lambda$ for any $x\in \ker(\theta_{\frakY_Z}^{\prime})$, it must be true that for any $y\in \ker(\theta_{\frakY_Z}^{\prime})$,
      \[|\varphi(\frac{y}{\xi})| \leq \lambda^p.\]
      But this follows from Lemma \ref{lem:phi(xi)/p is a unit} and \eqref{equ:extends to blow-up algebra} because we have
      \[|\varphi(\frac{y}{\xi})| = |(\frac{y}{\xi})^p+\frac{p}{\varphi(\xi)}\delta(y)-\frac{p}{\varphi(\xi)}\delta(\xi)(\frac{y}{\xi})^p|\leq \sup\{|(\frac{y}{\xi})^p|,|\delta(y)|\}\leq \sup\{\lambda^p,1\} = \lambda^p.\]
  \end{proof}
  \begin{construction}\label{construction:Frobenius on OC^I}
      Suppose that $\frakX_Z$ admits a prismatic lifting $\wtx_Z$.
      Let $I,J$ be two closed sub-intervals of $(0,+\infty)$ such that $J\subset(\frac{1}{p},+\infty)\cap I/p$ and denote by $\calO\CCtilde^I$ and $\calO\CCtilde^J$ the corresponding period sheaves, respectively. For each $1\leq \lambda<|\xi|_I^{-1}$, we define a morphism 
      \[\varphi:\calO\CCtilde^{I,\lambda^p}\to\calO\CCtilde^{J,\lambda}\]
      compatible with $\varphi:\calO_{\wtx_Z}\to\calO_{\wtx_Z}$ and with $\CCtilde^I\xrightarrow{\phi}\CCtilde^{I/p}\to\CCtilde^J$ by sheafifying the morphisms from Lemma \ref{lem:extend to period ring} on $X_{Z,v}^{\Box}$, as follows: 

      \begin{enumerate}
          \item Suppose that $\frakX_Z = \Spf(\calR_Z)$ is small affine. For any $U\in X_{Z,v}$, the morphism
          \[\colim_{\frakY_Z = \Spf(\calR_Z^{\prime})\in\Sigma_U}\varphi_{\frakY_Z}:\colim_{\frakY_Z = \Spf(\calR_Z^{\prime})\in\Sigma_U}A_{\frakZ}^{I,\lambda^p}\to \colim_{\frakY_Z = \Spf(\calR_Z^{\prime})\in\Sigma_U}A_{\frakY_Z}^{J,\lambda},\]
          where $\varphi_{\frakY_Z}$ is the morphism defined in Lemma \ref{lem:extend to period ring} associated to $\frakY_Z$, induces a morphism
          \[\varphi_U:\calO\CCtilde^{I,\lambda^p}_{\mid_U}\to \calO\CCtilde^{J,\lambda}_{\mid_U}\]
          by Definition \ref{dfn:local OC^I}, and thus defines the desired
          \[\varphi:\calO\CCtilde^{I,\lambda^p}\to \calO\CCtilde^{J,\lambda}\]
          in this case.

          \item In general, as in Definition \ref{dfn:OC^I}, we may define 
          \[\varphi:\calO\CCtilde^{I,\lambda^p}\to \calO\CCtilde^{J,\lambda}\]
          by letting $\varphi(U) = \varphi_U:\calO\CCtilde^{I,\lambda^p}(U)\to \calO\CCtilde^{J,\lambda}(U)$ for any $U\in X_{Z,v}^{\Box}$.
      \end{enumerate}
      By construction, this morphism is also compatible with the morphism
      \[\calO_{\wtx_Z^I}\xrightarrow{\varphi}\calO_{\wtx_Z^{I\cap I/p}}\to\calO_{\wtx_Z^J}.\]
      Now, letting $\lambda\to 1^+$, we get a morphism 
      \[\varphi:\calO\CCtilde^{I,\dagger}\to\calO\CCtilde^{J,\dagger},\]
      which is also compatible with the Frobenius map $\varphi:\calO_{\wtx_Z}\to\calO_{\wtx_Z}$ (and thus the morphism $\calO_{\wtx_Z^I}\xrightarrow{\varphi}\calO_{\wtx_Z^{I\cap I/p}}\to\calO_{\wtx_Z^J}$) and the composite of the morphism $\CCtilde^I\xrightarrow{\phi}\CCtilde^{I/p}\to\CCtilde^J$. Note that the Frobenius morphism $\varphi:\calO_{\wtx_Z}\to\calO_{\wtx_Z}$ also induces a morphism
      \[\varphi:\Omega^1_{\wtx_Z}\to \Omega^1_{\wtx_Z}\]
      by sending each local section $\rd f$ to $\rd\varphi(f)$. It follows that the morphism $\varphi:\calO\CCtilde^{I,\lambda^p}\to \calO\CCtilde^{J,\lambda}$ is compatible with the connections; that is, the following diagram commutes:
      \begin{equation}\label{equ: construction of Frobenius on OC^I}
          \begin{tikzcd}
             0 \arrow[r] & \CCtilde^I \arrow[d, "\phi"] \arrow[r] &{\calO\CCtilde^{I,\lambda^p}} \arrow[d, "\varphi"] \arrow[r] &{\calO\CCtilde^{I,\lambda^p}\otimes\Omega^1_{\wtx_Z}\{-1\}} \arrow[d, "\varphi"] \arrow[r] & \cdots \\
             0 \arrow[r] & \CCtilde^J \arrow[r] & {\calO\CCtilde^{J,\lambda}} \arrow[r] & {\calO\CCtilde^{I,\lambda^p}\otimes\Omega^1_{\wtx_Z}\{-1\}} \arrow[r] & \cdots,
          \end{tikzcd}
      \end{equation}
      where the Frobenius action on the Breuil--Kisin twist $\bullet\{-1\}$ is induced by $\phi(\frac{1}{\xi}) = \frac{\xi}{p}\frac{p}{\phi(\xi)}\frac{1}{\xi}$ (cf. Lemma \ref{lem:phi(xi)/p is a unit}); the same holds for $\varphi:\calO\CCtilde^{I,\dagger}\to\calO\CCtilde^{J,\dagger}$.
  \end{construction}

  \begin{exam}\label{exam:Frobenius on OC^I}
      Let $I,J$ be intervals with $I\subset(\frac{1}{p-1},p)$ and $J\subset (\frac{1}{p-1},1)\cap I/p$, and let $\frakX_Z = \Spf(\calR_Z)$ be small affine with generic fiber $X_Z = \Spa(R_Z,R_Z^+)$. Let $\wtx_Z = \Spf(\wtr_Z)$ be the prismatic lifting of $\frakX_Z$ so that $\varphi(T_i) = T_i^p(1+p\delta_{\log}(T_i))$ for any $0\leq i\leq d$. As $T_0\cdots T_r = [\varpi^a]$, we have
      \[(1+p\delta_{\log}(T_0))\cdots(1+p\delta_{\log}(T_r))=1.\]
      Let $X_{\infty,Z}\to X_Z$ be the pro-finite-\'etale $\Gamma$-torsor as in Example \ref{exam:Gamma torsor} so for any $U = \Spa(S,S^+)\in X_{Z,\infty,v}$, we can choose $t_i = [T_i^{\flat}]$ for any $0\leq i\leq d$. In particular, we have $\delta_{\log}([T_i^{\flat}]) = 0$. Now, for any $\lambda\in[1,\min\{|\xi|_I^{-1}p^{\frac{1}{p-1}},|\xi|_J^{-\frac{1}{p}}p^{\frac{1}{p(p-1)}}\})$ (which implies that $\lambda^p\in [1,|\xi|_J^{-1}p^{\frac{1}{p-1}})$), by Example \ref{exam:Apply to Gamma torsor}, for $P = W_0+\dots+W_r$, we have 
      \[\begin{split}\calO\CCtilde^{J,
      \lambda^p} = \CCtilde^J\{\underline W/\lambda^p\}_{r,d,P} = \CCtilde^J\{W_0/\lambda^p,\dots,W_d/\lambda^p\}/(W_0+\cdots+W_r)\end{split}\]
      and 
      \[\begin{split}\calO\CCtilde^{I,
      \lambda} = \CCtilde^I\{\underline W/\lambda\}_{r,d,P} = \CCtilde^I\{W_0/\lambda,\dots,W_d/\lambda\}/(W_0+\cdots+W_r)\end{split}\]
      where $W_i = \frac{u}{t}\log\frac{[T_i^{\flat}]}{T_i}$ for any $0\leq i\leq d$. In particular, we have 
      \[\varphi(W_i) = \frac{\phi(u)}{\phi(t)}\log\big(\frac{[T_i^{\flat}]^p}{T_i^p}(1+p\delta_{\log}(T_i))^{-1}\big) = \xi W_i - \frac{\phi(u)}{t}\frac{\log(1+p\delta_{\log}(T_i))}{p}\]
      for any $0\leq i\leq d$ and thus 
      \[\varphi(P) = \xi P-\sum_{i=0}^r\frac{\phi(u)}{t}\frac{\log(1+p\delta_{\log}(T_i))}{p} = \xi P-\frac{\phi(u)}{t}\frac{\log\big(\prod_{i=0}^r(1+p\delta_{\log}(T_i))\big)}{p} =  \xi P.\]
      Recall that $\frac{\phi(u)}{t}$ is a unit in $\CCtilde^J(Z)$ by Example \ref{exam:t in C^I}. Via the isomorphisms
      \[\begin{split}\calO\CCtilde^{J,
      \lambda^p} 
      \cong \CCtilde^J\{W_1/\lambda^p,\dots,W_d/\lambda^{p}\}\text{ and }\calO\CCtilde^{I,
      \lambda} \cong  \CCtilde^I\{W_1/\lambda,\dots,W_d/\lambda\}\end{split},\]
      the morphism $\varphi:\calO\CCtilde^{I,
      \lambda} \to \calO\CCtilde^{J,
      \lambda^p}$ coincides with the map 
      \[\phi:\CCtilde^I\{W_1/\lambda,\dots,W_d/\lambda\}\to \CCtilde^J\{W_1/\lambda^p,\dots,W_d/\lambda^{p}\}\]
      which acts on $\CCtilde^I$ via the composite $\CCtilde^I\xrightarrow{\phi}\CCtilde^{I/p}\to\CCtilde^J$ and sends $W_i$ to $\xi W_i - \frac{\phi(u)}{t}\frac{\log(1+p\delta_{\log}(T_i))}{p}$ for any $1\leq i\leq d$. Letting $\lambda\to 1^+$, we get a similar description of the morphism
      \[\varphi:\calO\CCtilde^{I,
      \dagger} \to \calO\CCtilde^{J,
      \dagger}\]
      via the isomorphisms
      \[\begin{split}\calO\CCtilde^{J,\dagger}\cong \colim_{\lambda\to 1^+} \CCtilde^J\{W_1/\lambda^p,\dots,W_d/\lambda^{p}\}\text{ and }\calO\CCtilde^{I,
      \dagger} \cong \colim_{\lambda\to 1^+} \CCtilde^I\{W_1/\lambda,\dots,W_d/\lambda\}\end{split}.\]
      In particular, if the $\delta_{\log}$-structure on $\wtr_Z$ is given in Example \ref{exam:prismatic lifting} (that is, $\delta_{\log}(T_i) = 0$ for any $0\leq i\leq d$), we have
      \[\varphi(W_i) = \xi W_i\]
      for any $0\leq i\leq d$. This identity is the key link between the ``connection'' side (the operators $\nabla_i$ and the $W_i$) and the ``arithmetic'' side (the Frobenius $\varphi$): $\varphi$ scales the logarithmic coordinates $W_i$ by $\xi$, which is exactly what is needed for the compatibility $\nabla_i\circ\varphi_D = \xi\varphi_D\circ\nabla_i$ in the local Riemann--Hilbert correspondence.
  \end{exam}

\subsection{$\varphi$-modules and $\varphi$-connections}\label{ssec:phi-module and phi-connection}
  In this subsection we assume $I\subset(\frac{1}{p-1},p)$ is a closed interval with $I\cap\frac{I}{p}\neq \emptyset$, so that the Breuil--Kisin twist $\Ainf\{-1\}$ is trivialized by $\Ainf\cdot\frac{u}{t}$. Such an interval exists if and only if $\frac{1}{p-1}<1$, i.e., $p\geq 3$. We also fix a prismatic lifting $\wtx_Z$ of $\frakX_Z$.

  \begin{dfn}\label{dfn:phi-module}
      By a \emph{$\varphi$-module} over $\CCtilde^I$ of rank $r$, we mean a $\CCtilde^I$-local system $\bM\in \rL\rS(X_Z,\CCtilde^I)$ of rank $r$ together with a $\CCtilde^{I\cap I/p}$-linear isomorphism
      \[\varphi_{\bM}:\bM\otimes_{\CCtilde^I,\phi}\CCtilde^{I\cap I/p}\xrightarrow{\cong}\bM\otimes_{\CCtilde^I}\CCtilde^{I\cap I/p}.\]
      The cohomology $\rR\Gamma(X_Z,\bM)^{\varphi=1}$ of a $\varphi$-module $(\bM,\varphi_{\bM})$ is defined by
      \[\rR\Gamma(X_Z,\bM)^{\varphi=1}:=\rR\Gamma(X_Z,\fib(\bM\xrightarrow{\varphi-1}\bM\otimes_{\CCtilde^I}\CCtilde^{I\cap I/p})) = \fib(\rR\Gamma(X_Z,\bM)\xrightarrow{\varphi-1}\rR\Gamma(X_Z,\bM)),\]
      where $\varphi$ denotes the morphism induced by the following composite
      \[\bM\xrightarrow{x\mapsto x\otimes1}\bM\otimes_{\CCtilde^I,\varphi}\CCtilde^{I\cap I/p}\xrightarrow{\varphi_{\bM}}\bM\otimes_{\CCtilde^I}\CCtilde^{I\cap I/p}\]
      while $1$ denotes the morphism induced by
      \[\bM\xrightarrow{x\mapsto x\otimes 1}\bM\otimes_{\CCtilde^I}\CCtilde^{I\cap I/p}.\]
      A $\varphi$-module $(\bM,\varphi_{\bM})$ over $\CCtilde^I$ is called \emph{small} if $\bM\in \rL\rS^{\dagger}(X_Z,\CCtilde^I)$. Denote by $\rL\rS^{\varphi}(X_Z,\CCtilde^I)$ the category of $\varphi$-module over $\CCtilde^I$ and by $\rL\rS^{\dagger,\varphi}(X_Z,\CCtilde^I)$ the full sub-category of small objects.
  \end{dfn}

  \begin{rmk}\label{rmk:Qp-local system as phi-module}
      By a theorem of Kedlaya--Liu \cite{SW-Berkeley}, the functor 
      \[\rL\mapsto (\bM(\rL):=\rL\otimes_{\Qp}\CCtilde^I,\varphi_{\bM(\rL)}:=\id_{\rL}\otimes\phi)\]
      defines a fully faithful functor
      \[\rL\rS(X_{Z,\proet},\Qp)\to \rL\rS^{\varphi}(X_Z,\CCtilde^I)\]
      from the category $\rL\rS(X_{Z,\proet},\Qp)$ of pro-\'etale $\Qp$-local systems on $X_Z$ to the category $\rL\rS^{\varphi}(X_Z,\CCtilde^I)$ of $\varphi$-modules over $\CCtilde^I$ on $X_Z$ whose essential image consists of the objects that are pointwise pure of slope zero. This functor preserves ranks, tensor products and dualities, and is compatible with cohomology in the sense that the inclusion $\rL\to \bM(\rL)$ induces a quasi-isomorphism
      \[\rR\Gamma(X_{Z,\proet},\rL)\xrightarrow{\simeq}\rR\Gamma(X_Z,\bM(\rL))^{\varphi=1}.\]
      In particular, we have a quasi-isomorphism
      \[\rR\Gamma(X_{Z,\proet},\Qp)\xrightarrow{\simeq}\rR\Gamma(X_Z,\CCtilde^I)^{\varphi=1}.\]
  \end{rmk}

  The Frobenius $\varphi:\calO_{\wtx_Z}\to\calO_{\wtx_Z}$ induces $\varphi:\calO_{\wtx_Z^I}\to\calO_{\wtx_Z^{I\cap I/p}}$ compatible with the usual connection $\rd$ and with $\nabla$, in the sense that the following diagram commutes:
  \begin{equation*}
      \begin{tikzcd}
          \calO_{\wtx_Z^I} \arrow[r,"\nabla"] \arrow[d, "\varphi"] & \Omega^1_{\wtx_Z^I}\{-1\} \arrow[d, "\varphi"] \\
          \calO_{\wtx_Z^{I\cap I/p}} \arrow[r,"\nabla"] & \Omega^1_{\wtx_Z^{I\cap I/p}}\{-1\},       
      \end{tikzcd}
  \end{equation*}
  where the Frobenius action on the Breuil--Kisin twist $\bullet\{-1\}$ is given by 
  \[\varphi(\frac{u}{t}) = \frac{\phi(u)}{pt} = \frac{\xi}{p}\frac{u}{t}.\]
  \begin{dfn}\label{dfn:phi-connection}
      Fix a prismatic lifting $\wtx_Z$ of $\frakX_Z$. By a \emph{$\varphi$-connection} of rank $r$ on $\wtx_Z^I$, we mean a flat connection $(\calD,\nabla_{\calD})$ of rank $r$ on $\wtx_Z^I$ together with an isomorphism
      \[\varphi_{\calD}:\calD\otimes_{\calO_{\wtx_Z^I},\varphi}\calO_{\wtx_Z^{I\cap I/p}}\to \calD\otimes_{\calO_{\wtx_Z^I}}\calO_{\wtx_Z^{I\cap I/p}}\]
      which is compatible with the connections. For a $\varphi$-connection $(\calD,\nabla_{\calD},\varphi_{\calD})$, denote by $\rD\rR(\calD,\nabla_{\calD})^{\varphi=1}$ the complex
      \[\fib(\rD\rR(\calD,\nabla_{\calD})\xrightarrow{\varphi-1}\rD\rR(\calD\otimes_{\calO_{\wtx_Z^I}}\calO_{\wtx_Z^{I\cap I/p}},\nabla_{\calD\otimes_{\calO_{\wtx_Z^I}}\calO_{\wtx_Z^{I\cap I/p}}})),\]
      which stands for the total complex of
      \begin{equation*}
          \begin{tikzcd}
              \calD \arrow[d, "\varphi-1"] \arrow[r, "\nabla_{\calD}"] & \calD\otimes_{\calO_{\wtx_Z^I}}\Omega^1_{\wtx_Z^I}\{-1\} \arrow[d, "\varphi-1"] \arrow[r, "\nabla_{\calD}"] & \calD\otimes_{\calO_{\wtx_Z^I}}\Omega^2_{\wtx_Z^I}\{-2\} \arrow[d, "\varphi-1"] \arrow[r] & \cdots \\
              \calD\otimes_{\calO_{\wtx_Z^I}}\calO_{\wtx_Z^{I\cap I/p}} \arrow[r, "\nabla_{\calD\otimes_{\calO_{\wtx_Z^I}}\calO_{\wtx_Z^{I\cap I/p}}}"] & \calD\otimes_{\calO_{\wtx_Z^I}}\Omega^1_{\wtx_Z^{I\cap I/p}}\{-1\} \arrow[r, "\nabla_{\calD\otimes_{\calO_{\wtx_Z^I}}\calO_{\wtx_Z^{I\cap I/p}}}"] & \calD\otimes_{\calO_{\wtx_Z^I}}\Omega^2_{\wtx_Z^{I\cap I/p}}\{-2\} \arrow[r] & \cdots.
          \end{tikzcd}
      \end{equation*}
      Here, $\varphi$ denotes the morphism induced by the composite
      \[\calD\otimes_{\calO_{\wtx_Z^I}}\Omega^{\bullet}_{\wtx_Z^I}\{-\bullet\}\xrightarrow{x\otimes\omega\mapsto x\otimes\omega\otimes 1}\calD\otimes_{\calO_{\wtx_Z^I}}\Omega^{\bullet}_{\wtx_Z^I}\{-\bullet\}\otimes_{\calO_{\wtx_Z^I},\varphi}\calO_{\wtx_Z^{I\cap I/p}}\xrightarrow{\varphi_{\calD}\otimes\varphi\otimes\id}\calD\otimes_{\calO_{\wtx_Z^I}}\Omega^{\bullet}_{\wtx_Z^{I\cap I/p}}\{-\bullet\}\]
      and $1$ denotes the morphism induced by
      \[\calD\otimes_{\calO_{\wtx_Z^I}}\Omega^{\bullet}_{\wtx_Z^I}\{-\bullet\}\xrightarrow{x\otimes\omega\mapsto x\otimes\omega\otimes 1}\calD\otimes_{\calO_{\wtx_Z^I}}\Omega^{\bullet}_{\wtx_Z^I}\{-\bullet\}\otimes_{\calO_{\wtx_Z^I}}\calO_{\wtx_Z^{I\cap I/p}} \cong \calD\otimes_{\calO_{\wtx_Z^I}}\Omega^{\bullet}_{\wtx_Z^{I\cap I/p}}\{-\bullet\}.\]
      The cohomology $\rR\Gamma_{\dR}(\frakX_Z,\rD\rR(\calD,\nabla_{\calD}))^{\varphi=1}$ of a $\varphi$-connection $(\calD,\nabla_{\calD},\varphi_{\calD})$ is defined by
      \[\rR\Gamma(\frakX_Z,\rD\rR(\calD,\nabla_{\calD})^{\varphi=1}) = \fib(\rR\Gamma(\frakX_Z,\rD\rR(\calD,\nabla_{\calD}))\xrightarrow{\varphi-1}\rR\Gamma(\frakX_Z,\rD\rR(\calD\otimes_{\calO_{\wtx_Z^I}}\calO_{\wtx_Z^{I\cap I/p}},\nabla_{\calD\otimes_{\calO_{\wtx_Z^I}}\calO_{\wtx_Z^{I\cap I/p}}}))).\]
      A $\varphi$-connection $(\calD,\nabla_{\calD},\varphi_{\calD})$ is called \emph{small} if $(\calD,\nabla_{\calD})\in \MIC^{\dagger}(\wtx_Z^I)$. Denote by $\MIC^{\varphi}(\wtx_Z^I)$ the category of $\varphi$-connections on $\wtx_Z^I$ and by $\MIC^{\dagger,\varphi}(\wtx_Z^I)$ the full sub-category of small objects.
  \end{dfn}
  \begin{thm}\label{thm:prismatic RH}
      Fix a prismatic lifting $\wtx_Z$ of $\frakX_Z$ and let $\varphi:\calO\CCtilde^{I,\dagger}\to\calO\CCtilde^{I\cap I/p,\dagger}$ be the associated Frobenius action on the period sheaves with connection. Then the Riemann--Hilbert correspondence \ref{thm:RH correspondence} is compatible with Frobenius action in the sense that it induces an equivalence of categories
      \[\rL\rS^{\dagger,\varphi}(X_Z,\CCtilde^I)\simeq\MIC^{\dagger,\varphi}(\wtx_Z^I),\]
      which is compatible with cohomology. In particular, the induced Frobenius action on the cohomology groups of both sides is preserved under this equivalence; this is the foundational step for the Frobenius-equivariance of the rigid cohomology comparison in \S\ref{sec:rigid cohomology}.
  \end{thm}
  \begin{proof}
      We keep the notation of Theorem \ref{thm:RH correspondence}.

      \smallskip\noindent\emph{From $\varphi$-modules to $\varphi$-connections.}
      Let $(\bM,\varphi_{\bM})$ be a small $\varphi$-module over $\CCtilde^I$. Then $\varphi_{\bM}$ induces an isomorphism
      \[\varphi_{\bM}\otimes\id:\bM\otimes_{\CCtilde^I,\varphi}\CCtilde^{I\cap I/p}\otimes_{\CCtilde^{I\cap I/p}}\calO\CCtilde^{I\cap I/p,\dagger}\to \bM\otimes_{\CCtilde^I}\CCtilde^{I\cap I/p}\otimes_{\CCtilde^{I\cap I/p}}\calO\CCtilde^{I\cap I/p,\dagger},\]
      which amounts to the isomorphism
      \[(\bM\otimes_{\CCtilde^I}\calO\CCtilde^{I,\dagger})\otimes_{\calO\CCtilde^{I,\dagger},\varphi}\calO\CCtilde^{I\cap I/p,\dagger}\to (\bM\otimes_{\CCtilde^I}\CCtilde^{I\cap I/p})\otimes_{\CCtilde^{I\cap I/p}}\calO\CCtilde^{I\cap I/p,\dagger}.\]
      Applying $\nu_*$, we get an isomorphism
      \[\varphi_{\calD(\bM)}:\calD(\bM)\otimes_{\calO_{\wtx_Z^I},\varphi}\calO_{\wtx_Z^{I\cap I/p}}\to \calD(\bM\otimes_{\CCtilde^I}\CCtilde^{I\cap I/p})\cong \calD(\bM)\otimes_{\calO_{\wtx_Z^I}}\calO_{\wtx_Z^{I\cap I/p}},\]
      where the last isomorphism follows from Theorem \ref{thm:functoriality in coefficient}. As $\varphi:\calO\CCtilde^{I,\dagger}\to\calO\CCtilde^{I\cap I/p,\dagger}$ is compatible with the connections, we see that
      $\varphi_{\calD(\bM)}$ is also compatible with the connections and thus $(\calD(\bM),\nabla_{\calD(\bM)},\varphi_{\calD(\bM)})$ is a small $\varphi$-connection.

      \smallskip\noindent\emph{From $\varphi$-connections to $\varphi$-modules.}
      Conversely, let $(\calD,\nabla_{\calD},\varphi_{\calD})$ be a small $\varphi$-connection. Then $\varphi_{\calD}$ induces an isomorphism compatible with the connection
      \[\varphi_{\calD}\otimes\id:\calD\otimes_{\calO_{\wtx_Z^I},\varphi}\calO_{\wtx_Z^{I\cap I/p}}\otimes_{\calO_{\wtx_Z^{I\cap I/p}}}\calO\CCtilde^{I\cap I/p,\dagger}\to \calD\otimes_{\calO_{\wtx_Z^I}}\calO_{\wtx_Z^{I\cap I/p}}\otimes_{\calO_{\wtx_Z^{I\cap I/p}}}\calO\CCtilde^{I\cap I/p,\dagger}\]
      which amounts to the isomorphism
      \[\calD\otimes_{\calO_{\wtx_Z^I}}\calO\CCtilde^{I,\dagger}\otimes_{\calO\CCtilde^{I,\dagger},\varphi}\calO\CCtilde^{I\cap I/p,\dagger}\to \calD\otimes_{\calO_{\wtx_Z^I}}\calO_{\wtx_Z^{I\cap I/p}}\otimes_{\calO_{\wtx_Z^{I\cap I/p}}}\calO\CCtilde^{I\cap I/p,\dagger}.\]
      By taking horizontal sections, this induces an isomorphism
      \[\varphi_{\bM(\calD,\nabla_{\calD})}:\bM(\calD,\nabla_{\calD})\otimes_{\CCtilde^I,\varphi}\CCtilde^{I\cap I/p}\to \bM(\calD\otimes_{\calO_{\wtx_Z^I}}\calO_{\wtx_Z^{I\cap I/p}},\nabla_{\calD\otimes_{\calO_{\wtx_Z^I}}\calO_{\wtx_Z^{I\cap I/p}}})\cong \bM(\calD,\nabla_{\calD})\otimes_{\CCtilde^I}\CCtilde^{I\cap I/p},\]
      where for the last isomorphism, we again apply Theorem \ref{thm:functoriality in coefficient}.       Hence $(\bM(\calD,\nabla_{\calD}),\varphi_{\bM(\calD,\nabla_{\calD})})$ is a small $\varphi$-module over $\CCtilde^I$.

      \smallskip\noindent\emph{Equivalence and cohomology.}
      For any small $\bM$ and small $(\calD,\nabla_{\calD})$, the comparison isomorphisms $\iota_{\bM}$ and $\iota_{(\calD,\nabla_{\calD})}$ from Theorem \ref{thm:RH correspondence} intertwine the $\varphi$-actions constructed above, which gives the equivalence of categories. Cohomological compatibility follows by the same functoriality.
  \end{proof}

\subsection{Local description of Theorem \ref{thm:prismatic RH}}\label{ssec:local prismatic RH}
  The passage from local connections to $(\varphi,\Gamma)$-modules below uses the decompletion theorem of \S\ref{sec:decompletion} (Lemma \ref{lem:prismatic decompletion}), which gives the equivalence between completed and decompleted representations.

  Let $I = [r,s]\subset (\frac{1}{p-1},p)$ be a closed interval as in the previous subsection. We assume $\frakX_Z = \Spf(\calR_Z)$ is small semi-stable with prismatic lifting $\wtx_Z = \Spf(\wtr_Z)$ satisfying $\delta_{\log}(T_i) = 0$ for $0\leq i\leq d$ (Example \ref{exam:prismatic lifting}), and let $X_{\infty,Z}\to X_Z$ be the pro-finite-\'etale $\Gamma$-torsor of the rigid generic fiber $X_Z$ (Example \ref{exam:Gamma torsor}). The morphism $\widetilde \psi: \wtr_Z\to \AAinf(X_{\infty,Z})$ from \S\ref{sec:decompletion} is $\varphi$-equivariant. As a consequence, the inclusion $\Cppsi\subset\CCppsi$ is also $\varphi$-equivariant in the sense that for any $J\subset I\cap I/p$, the following diagram
  \begin{equation*}
      \begin{tikzcd}
          {\widetilde \bfC_{\widetilde \psi}^{I,+}} \arrow[d, "\cap"] \arrow[r, "\varphi"] & {\widetilde \bfC_{\widetilde \psi}^{I/p,+}} \arrow[d, "\cap"] \arrow[r] & {\widetilde \bfC_{\widetilde \psi}^{J,+}} \arrow[d, "\cap"] \\
          {\widetilde \bfC_{\widetilde \psi,\infty}^{I,+}} \arrow[r, "\varphi"] & {\widetilde \bfC_{\widetilde \psi,\infty}^{I/p,+}} \arrow[r] & {\widetilde \bfC_{\widetilde \psi,\infty}^{J,+}} .          
       \end{tikzcd}
  \end{equation*}
  The same also holds for the inclusion $\Cpsi\subset\CCpsi$. As the $\varphi$-action commutes with the $\Gamma$-action on $\CCpsi$, the $\varphi$-action on $\Cpsi$ also commutes with the $\Gamma$-action.
  \begin{dfn}\label{dfn:phi-Gamma module over C^I}
      By a \emph{$(\varphi,\Gamma)$-module} of rank $r$ over $\widetilde \bfC_{\widetilde \psi,\infty}^I$ (resp. $\widetilde \bfC_{\widetilde \psi}^I$), we mean a $\Gamma$-representation $M_{\infty}\in \Rep_{\Gamma}(\widetilde \bfC_{\widetilde \psi,\infty}^I)$ (resp. $M\in \Rep_{\Gamma}(\widetilde \bfC_{\widetilde \psi}^I)$) of rank $r$ together with a $\Gamma$-equivariant $\widetilde \bfC_{\widetilde \psi,\infty}^{I\cap I/p}$-linear (resp. $\widetilde \bfC_{\widetilde \psi}^{I\cap I/p}$-linear) isomorphism
      \[\begin{split}&\varphi_{M_{\infty}}:M_{\infty}\otimes_{\widetilde \bfC_{\widetilde \psi,\infty}^I,\varphi}\widetilde \bfC_{\widetilde \psi,\infty}^{I\cap I/p}\to M_{\infty}\otimes_{\widetilde \bfC_{\widetilde \psi,\infty}^I}\widetilde \bfC_{\widetilde \psi,\infty}^{I\cap I/p}\\
      &\text{(resp. }\varphi_{M}:M\otimes_{\widetilde \bfC_{\widetilde \psi}^I,\varphi}\widetilde \bfC_{\widetilde \psi}^{I\cap I/p}\to M\otimes_{\widetilde \bfC_{\widetilde \psi}^I}\widetilde \bfC_{\widetilde \psi}^{I\cap I/p}\text{)}.\end{split}\]
      Such a $(\varphi,\Gamma)$-module $M_{\infty}$ (resp. $M$) is called \emph{small}, if its underlying $\Gamma$-representation is small in the sense of Definition \ref{dfn:small representation}. Denote by $\Rep_{\Gamma}^{\varphi}(\CCpsi)$ (resp. $\Rep_{\Gamma}^{\varphi}(\Cpsi)$) the category of $(\varphi,\Gamma)$-modules over $\CCpsi$ (resp. $\Cpsi$) and by $\Rep_{\Gamma}^{\dagger,\varphi}(\CCpsi)$ (resp. $\Rep_{\Gamma}^{\dagger,\varphi}(\Cpsi)$) the full sub-category of small objects.
  \end{dfn}
  The following lemma is clear.
  \begin{lem}\label{lem:local description of phi-modules}
      taking global sections $\bM\mapsto \bM(X_{\infty,Z})$ induces an equivalence of categories
      \[\rL\rS^{\varphi}(X_Z,\CCtilde^I)\xrightarrow{\simeq}\Rep_{\Gamma}^{\varphi}(\CCpsi)\]
      and thus an equivalence of categories
      \[\rL\rS^{\dagger,\varphi}(X_Z,\CCtilde^I)\xrightarrow{\simeq} \Rep_{\Gamma}^{\dagger,\varphi}(\CCpsi).\]
      These equivalences are compatible with cohomology.
  \end{lem}
  \begin{proof}
      This directly follows from Lemma \ref{lem:local description of local systems}.
  \end{proof}
  The decompletion functor is also compatible with Frobenius action.
  \begin{lem}\label{lem:prismatic decompletion}
      The functor $M\mapsto M_{\infty}:=M\otimes_{\Cpsi}\CCpsi$ in Theorem \ref{thm:decompletion} induces an equivalence of categories
      \[\Rep_{\Gamma}^{\dagger,\varphi}(\Cpsi)\xrightarrow{\simeq}\Rep_{\Gamma}^{\dagger,\varphi}(\CCpsi)\]
      which is compatible with cohomology.
  \end{lem}
  \begin{proof}
      We first establish the desired equivalence of categories.
      \emph{Idea.} Starting from a $\CCppsi$-basis of a small lattice, apply $\varphi$ and use $\Gamma$-equivariance; the estimates $|\frac{(\gamma-1)}{u}(f_i)|<p^{-1/(p-1)}$ show that Frobenius preserves smallness, and the induced $\varphi_M$ on the completed side is the unique map whose scalar extension recovers $\varphi_{M_{\infty}}$.

      For any $M_{\infty}\in\Rep_{\Gamma}^{\dagger,\varphi}(\CCpsi)$ with associated $M\in \Rep_{\Gamma}^{\dagger}(\Cpsi)$, it is enough to show that the isomorphism $\varphi_{M_{\infty}}$ induces an isomorphism
      \[\varphi_M:M\otimes_{\Cpsi,\varphi}\widetilde \bfC_{\widetilde \psi}^{I\cap I/p}\to M\otimes_{\Cpsi}\widetilde \bfC_{\widetilde \psi}^{I\cap I/p}\]
      whose scalar extension along $ \widetilde \bfC_{\widetilde \psi}^{I\cap I/p}\to \widetilde \bfC_{\widetilde \psi,\infty}^{I\cap I/p}$ is exactly $\varphi_{M_{\infty}}$
      
      Let $e_1,\dots,e_r$ be a $\CCppsi$-basis of $M_{\infty}^+$, a $\Gamma$-stable $\CCppsi$-lattice of $M_{\infty}$, satisfying the condition \eqref{equ:small representation} for some $\lambda>p^{\frac{1}{p-1}}$. Put $f_i = \varphi(e_i)$ for any $1\leq i\leq r$. Then $M_{\infty}^{\prime,+}:=\oplus_{i=1}^r\widetilde \bfC_{\widetilde \psi,\infty}^{I\cap I/p,+}\cdot f_i$ is a $\Gamma$-stable $\widetilde \bfC_{\widetilde \psi,\infty}^{I\cap I/p,+}$-lattice of $M_{\infty}\otimes_{\CCpsi,\varphi}\widetilde \bfC_{\widetilde \psi,\infty}^{I\cap I/p}$.
      Then we have $\varphi_{M_{\infty}}(M_{\infty}^+)\subset M_{\infty}^{\prime,+}$ such that for any $\gamma\in \Gamma$, we have
      \[(\gamma-1)(f_i) = \varphi_{M_{\infty}}((\gamma-1)(e_i))\subset \varphi_{M_{\infty}}(uM_{\infty}^+)\subset \phi(u)M_{\infty}^{\prime,+}.\]
      In particular, we have
      \[|\frac{(\gamma-1)}{u}(f_i)|\leq |\xi|_{I\cap I/p}\leq\max\{|p|,\norm\varpi\norm^r\}<p^{-\frac{1}{p-1}},\]
      yielding that $M_{\infty}\otimes_{\CCpsi,\varphi}\widetilde \bfC_{\widetilde \psi,\infty}^{I\cap I/p}\in \Rep_{\Gamma}^{\dagger}(\widetilde \bfC_{\widetilde \psi,\infty}^{I\cap I/p})$. As
      \[(M\otimes_{\Cpsi,\varphi}\widetilde \bfC_{\widetilde \psi}^{I\cap I/p})\otimes_{\widetilde \bfC_{\widetilde \psi}^{I\cap I/p}}\widetilde \bfC_{\widetilde \psi,\infty}^{I\cap I/p} \cong (M\otimes_{\Cpsi}\CCpsi)\otimes_{\CCpsi,\varphi}\widetilde \bfC_{\widetilde \psi,\infty}^{I\cap I/p} = M_{\infty}\otimes_{\CCpsi,\varphi}\widetilde \bfC_{\widetilde \psi,\infty}^{I\cap I/p}.\]
      We see that $M\otimes_{\Cpsi,\varphi}\widetilde \bfC_{\widetilde \psi}^{I\cap I/p}\in \Rep_{\Gamma}^{\dagger}(\widetilde \bfC_{\widetilde \psi}^{I\cap I/p})$ is associated to $M_{\infty}\otimes_{\CCpsi,\varphi}\widetilde \bfC_{\widetilde \psi,\infty}^{I\cap I/p}$ via the equivalence in Corollary \ref{cor:decompletion} (2). So the isomorphism
      \[\varphi_{M_{\infty}}:M_{\infty}\otimes_{\CCpsi,\varphi}\widetilde \bfC_{\widetilde \psi,\infty}^{I\cap I/p}\to M_{\infty}\otimes_{\CCpsi}\widetilde \bfC_{\widetilde \psi}^{I\cap I/p}\]
      induces an isomorphism
      \[\varphi_M:M\otimes_{\Cpsi,\varphi}\widetilde \bfC_{\widetilde \psi}^{I\cap I/p}\to M\otimes_{\Cpsi}\widetilde \bfC_{\widetilde \psi}^{I\cap I/p}\]
      which satisfies all the expected conditions.

      To see this equivalence is compatible with cohomology, it suffices to show that for any $M\in \Rep_{\Gamma}^{\dagger,\varphi}(\Cpsi)$, the morphism $M\mapsto M_{\infty}:=M\otimes_{\Cpsi}\CCpsi$ induces a quasi-isomorphism
      \[\fib(\rR\Gamma(\Gamma,M)\xrightarrow{\varphi-1}\rR\Gamma(\Gamma,M\otimes_{\Cpsi}\widetilde \bfC_{\widetilde \psi}^{I\cap I/p}))\to \fib(\rR\Gamma(\Gamma,M_{\infty})\xrightarrow{\varphi-1}\rR\Gamma(\Gamma,M_{\infty}\otimes_{\CCpsi}\widetilde \bfC_{\widetilde \psi,\infty}^{I\cap I/p})).\]
      As we have quasi-isomorphisms
      \[\rR\Gamma(\Gamma,M)\simeq \rR\Gamma(\Gamma,M_{\infty})\text{ and }\rR\Gamma(\Gamma,M\otimes_{\Cpsi}\widetilde \bfC_{\widetilde \psi}^{I\cap I/p}))\simeq \rR\Gamma(\Gamma,M_{\infty}\otimes_{\CCpsi}\widetilde \bfC_{\widetilde \psi,\infty}^{I\cap I/p})),\]
      the desired quasi-isomorphism follows from a standard spectral sequence argument as in the proof \cite[Th. 3.11]{MW25}.
  \end{proof}

  \begin{dfn}\label{dfn:local phi-connection}
      By a \emph{$\varphi$-connection over $\wtr_Z^I$} of rank $r$, we mean a flat connection $(D,\nabla_D)\in \MIC(\wtr_Z)$ together with an $\wtr_Z^{I\cap I/p}$-linear isomorphism
      \[\varphi_D:D\otimes_{\wtr_Z^I,\varphi}\wtr_Z^{I\cap I/p}\to D\otimes_{\wtr_Z^I}\wtr_Z^{I\cap I/p}\]
      which is compatible with the connection. Here the Breuil--Kisin twist $\{-1\}$ and its trivialization by $\frac{u}{t}$ are the same as in \S\ref{sec:local RH}, Application I (Cohomology of $\bA^d_Z$), and \S\ref{sec:rigid cohomology}. Note that if we write 
      \[\nabla_D = \sum_{i=1}^d\nabla_i\otimes\frac{u}{t}\dlog T_i,\]
      then $\varphi_D$ is compatible with the connection if and only if for any $1\leq i\leq d$
      \[\nabla_i\circ\varphi_D = \xi\varphi_D\circ\nabla_i.\]
      A $\varphi$-connection $(D,\nabla_D,\varphi_D)$ is called \emph{small} if its underlying connection $(D,\nabla_D)$ is small in the sense of Definition \ref{dfn:small local connection}. Denote by $\MIC^{\varphi}(\wtr_Z^I)$ the category of $\varphi$-connections over $\wtr_Z^I$ and by $\MIC^{\dagger,\varphi}(\wtr_Z^I)$ the full sub-category of small objects.
  \end{dfn}
  taking global sectionss on $\frakX_Z = \Spf(\calR_Z)$ gives equivalences $\MIC^{\varphi}(\wtx_Z^I)\simeq \MIC^{\varphi}(\wtr_Z^I)$ and $\MIC^{\dagger,\varphi}(\wtx_Z^I)\simeq \MIC^{\dagger,\varphi}(\wtr_Z^I)$.
  \begin{prop}\label{prop:local prismatic RH}
      Keep the notation of Proposition \ref{prop:From representation to connection} and Proposition \ref{prop:From connection to representation}. The local Riemann--Hilbert correspondence (Theorem \ref{thm:local correspondence over Robba ring}) induces equivalences of categories
      \[\Rep_{\Gamma}^{\dagger,\varphi}(\CCpsi)\simeq \Rep_{\Gamma}^{\dagger,\varphi}(\Cpsi)\simeq \MIC^{\dagger,\varphi}(\wtr_Z^I)\]
      that are compatible with cohomology.
      More precisely, we have the following explicit description of the above equivalence:
      \begin{enumerate}
          \item[(1)] For any $(M_{\infty},\varphi_{M_{\infty}})\in \Rep_{\Gamma}^{\dagger,\varphi}(\CCpsi)$ with associated $(M,\varphi_M)\in \Rep_{\Gamma}^{\dagger,\varphi}(\Cpsi)$, the corresponding $(D,\nabla_D,\varphi_D)$ is given by
          \[(D = (\prod_{i=1}^d\gamma_i^{-\frac{W_i}{u}})(M),\nabla_D = \sum_{i=1}^d\frac{\log\gamma_i}{u}\otimes\frac{u}{t}\dlog T_i,\varphi_D = \varphi_M)\]
          in the sense of \eqref{equ:explicit D(M)} and \eqref{equ:explicit nabla} such that for any $x\in M$, 
          \begin{equation}\label{equ:explicit phi_D}
              \varphi_D(\big(\prod_{i=1}^d\gamma_i^{-\frac{W_i}{u}}\big)(x)) = \big(\prod_{i=1}^d\gamma_i^{-\frac{W_i}{u}}\big)(\varphi_M(x)).
          \end{equation}

          \item For any $(D,\nabla_D = \sum_{i=1}^d\nabla_i\otimes\frac{u}{t}\dlog T_i,\varphi_D)\in \MIC^{\dagger,\varphi}(\wtr_Z^I)$, the corresponding $(M,\varphi_M)\in \Rep_{\Gamma}^{\dagger,\varphi}(\Cpsi)$ is given by
          \[(M_{\infty} = \exp(\sum_{i=1}^d\nabla_iW_i)(M),\varphi_{M_{\infty}} = \varphi_D)\]
          in the sense of \eqref{equ:explicit M(D)} and \eqref{equ:explicit gamma action} such that for any $x\in D$,
          \begin{equation}\label{equ:explicit phi_M}
              \varphi_{M_{\infty}}\big(\exp(\sum_{i=1}^d\nabla_iW_i)(x)\big) = \exp(\sum_{i=1}^d\nabla_iW_i)\big(\varphi_D(x)\big).
          \end{equation}
      \end{enumerate}
  \end{prop}
  \begin{proof}
      The desired equivalences follow from Theorem \ref{thm:prismatic RH} and Lemma \ref{lem:prismatic decompletion}. It remains to prove Items (1) and (2).
      
      For item (1): Everything is clear by Proposition \ref{prop:From representation to connection} except the explicit description \eqref{equ:explicit phi_D} of $\varphi_D$. To prove \eqref{equ:explicit phi_D}, as $\varphi(W_i) = \xi W_i$ (cf. Example \ref{exam:Frobenius on OC^I}), it remains only to observe that 
      \[\varphi_D(\big(\prod_{i=1}^d\gamma_i^{-\frac{W_i}{u}}\big)(x)) = \big(\prod_{i=1}^d\gamma_i^{-\frac{\varphi(W_i)}{\phi(u)}}\big)(\varphi_M(x)) = \big(\prod_{i=1}^d\gamma_i^{-\frac{\xi W_i}{\phi(u)}}\big)(\varphi_M(x)) = \big(\prod_{i=1}^d\gamma_i^{-\frac{W_i}{u}}\big)(\varphi_M(x))\]
      as desired. Here, we implicitly use that $\varphi_M$ commutes with $\gamma_i$ for any $1\leq i\leq d$.

      For item (2): Everything is clear by Proposition \ref{prop:From connection to representation} except the explicit description \eqref{equ:explicit phi_M}. To prove \eqref{equ:explicit phi_M}, as $\varphi(W_i) = \xi W_i$ (cf. Example \ref{exam:Frobenius on OC^I}), one concludes by noting that
      \begin{equation*}
          \begin{split}
              \varphi_{M_{\infty}}\big(\exp(\sum_{i=1}^d\nabla_iW_i)(x)\big) = &\varphi_{M_{\infty}}\big(\sum_{n_1,\dots,n_d\geq 0}\frac{\nabla_1^{n_1}W_1^{n_1}}{n_1!}\cdots \frac{\nabla_d^{n_d}W_d^{n_d}}{n_d!}(x)\big)\\
              = & \sum_{n_1,\dots,n_d\geq 0}\frac{\varphi(W_1)^{n_1}}{n_1!}\cdots \frac{\varphi(W_d)^{n_d}}{n_d!}\varphi_D\circ(\nabla_1^{n_1}\cdots\nabla_d^{n_d})(x)\\
              = & \sum_{n_1,\dots,n_d\geq 0}\frac{W_1^{n_1}}{n_1!}\cdots \frac{W_d^{n_d}}{n_d!}\xi^{n_1+\cdots+n_d}\varphi_D\circ(\nabla_1^{n_1}\cdots\nabla_d^{n_d})(x)\\
              = & \sum_{n_1,\dots,n_d\geq 0}\frac{W_1^{n_1}}{n_1!}\cdots \frac{W_d^{n_d}}{n_d!}\nabla_1^{n_1}\cdots\nabla_d^{n_d}(\varphi_D(x))\\
              = & \exp(\sum_{i=1}^d\nabla_iW_i)\big(\varphi_D(x)\big).
          \end{split}
      \end{equation*}
      Here, we use $\nabla_i\circ\varphi_D = \xi\varphi_D\circ\nabla_i$ for any $1\leq i\leq d$ to conclude that 
      \[\xi^{n_1+\cdots+n_d}\varphi_D\circ(\nabla_1^{n_1}\cdots\nabla_d^{n_d}) = (\nabla_1^{n_1}\cdots\nabla_d^{n_d})\circ\varphi_D\]
      for any $n_1,\dots,n_d\geq 0$.
      Thus the identity $\varphi(W_i)=\xi W_i$ (Example \ref{exam:Frobenius on OC^I}) and the flatness relation $\nabla_i\circ\varphi_D = \xi\varphi_D\circ\nabla_i$ exactly cancel the extra $\xi$-powers, so $\varphi_{M_{\infty}}$ and $\exp(\sum_i\nabla_iW_i)$ commute as in \eqref{equ:explicit phi_M}.
  \end{proof}

\section{Application I: Cohomology of $\bA^d_Z$}\label{sec:pro-etale cohomology of Affine space}
  Let $I\subset(\frac{1}{p-1},p)$ be a closed interval.
  Fix an affinoid perfectoid $Z = \Spa(S,S^+)\in \Perfd$ as before.
  Using Theorem \ref{thm:prismatic RH}, we give a new proof of the following result of Le Bras \cite[Th. 3.2]{LB18}.
  \begin{thm}\label{thm:pro-etale cohomology of Affine space}
      Let $\bA^d_Z$ be the affine space of relative dimension $d$ over $Z$.
      \begin{enumerate}
          \item[(1)] We have $\rH^0_v(\bA^d_Z,\Qp) = \Qp$ and $\rH^0_v(\bA^d_Z,\BBdRp) = \BBdRp(Z)$.
          \item[(2)] For any $n\geq 1$, we have $\rH^n_v(\bA^d_Z,\Qp) \cong \Omega^n(\bA_Z^d)^{\rd = 0}\{-n\}\cong \rH^n_v(\bA^d_Z,\BBdRp)$.
      \end{enumerate}
      Here, $\rd:\Omega^n(\bA_Z^d)\to \Omega^{n+1}(\bA_Z^d)$ denotes the differential in the usual de Rham complex of $\bA_Z^d$:
      \[0\to \calO(\bA_Z^d)\xrightarrow{\rd}\Omega^1(\bA_Z^d)\to\cdots\to\Omega^n(\bA_Z^d)\xrightarrow{\rd}\Omega^{n+1}(\bA_Z^d)\to\cdots.\]
  \end{thm}

  Let $\bB^d_Z = \Spa(R_Z:=S\za T_1,\dots,T_d\ya,\calR_Z:=S^+\za T_1,\dots,T_d\ya)$ be the closed unit ball over $Z$ of relative dimension $d$.
  It admits $\frakB_Z^d:=\Spf(\calR_Z)$ as a smooth formal model, which has a prismatic lifting $\widetilde \frakB^d_Z:=\Spf(\wtr_Z:=\AAinf(Z)\za T_1,\dots,T_d\ya)$ such that $\delta(T_i) = 0$ for any $1\leq i\leq d$.
  So we have 
  \[\wtr_Z^I = \CCtilde^I(Z)\za T_1,\dots,T_d\ya \text{ and }\wtR_Z = \BBdRp(Z)\za T_1,\dots,T_d\ya.\]

  \begin{lem}\label{lem:cohomology for unit ball}
      For any $n\geq 0$, we have
      \begin{enumerate}
          \item[(1)] $\rH^n_v(\bB_Z^d,\CCtilde^I) = \frac{\Omega^n(\widetilde \frakB_Z^{d,I})^{\rd = 0}}{\xi \rd(\Omega^{n-1}(\widetilde \frakB_Z^{d,I}))}\{-n\}$,
          where $\widetilde \frakB_Z^{d,I}:=\Spa(\wtr_Z^I,\wtr_Z^{I,+})$, and 

          \item[(2)] $\rH^n_v(\bB_Z^d,\BBdRp) = \frac{\Omega^n(\widetilde \bB_Z^d)^{\rd = 0}}{\xi\rd(\Omega^{n-1}(\widetilde \bB_Z^d))}\{-n\}$,
          where $\widetilde \bB_Z^d:=\Spa(\wtR_Z)$.
      \end{enumerate}
      In particular, we have
      \[\rH^0_v(\bB_Z^d,\CCtilde^I) = \CCtilde^I(Z) \text{ and }\rH^0_v(\bB_Z^d,\BBdRp) = \BBdRp(Z).\]
  \end{lem}
  \begin{proof}
      The Breuil--Kisin twists $\{-n\}$ in the statement are those of \S\ref{sec:local RH}: $\CCtilde^{I,+}\{1\}$ is trivialized by $\frac{t}{u}$, with $t$ a unit-multiple of $\xi u$ in $\CCtilde^I$ (Example \ref{exam:t in C^I}(1)). Thus the grading $\Omega^{\bullet}\{-\bullet\}$ in the de Rham complex below agrees with the one from the Riemann--Hilbert correspondence.
      Let $\nu: \bB_{Z,v}^d\to \frakB_{Z,\et}^d$ be the morphism of sites.
      Then by Theorem \ref{thm:RH correspondence}, we have 
      \[\rR\Gamma(\bB_{Z,v}^d,\CCtilde^I) = \rR\Gamma(\frakB_{Z,\et}^d,\rR\nu_*\CCtilde^I)\simeq \rR\Gamma(\frakB_{Z,\et}^d,\rD\rR(\calO_{\widetilde \frakB_Z^{d,I}},\nabla))\]
      where $\rD\rR(\calO_{\widetilde \frakB_Z^{d,I}},\nabla)$ denotes the complex
      \[\calO_{\widetilde \frakB_Z^{d,I}}\xrightarrow{\nabla}\Omega^1_{\widetilde \frakB_Z^{d,I}}\{-1\}\xrightarrow{\nabla}\Omega^2_{\widetilde \frakB_Z^{d,I}}\{-2\}\to\cdots\]
      and $\nabla = \sum_{i=1}^d\xi\frac{\partial}{\partial T_i}\otimes\frac{\rd T_i}{\xi}= \rd$.
      As $\frakB_Z^d$ is affine, $\rR\Gamma(\frakB_{Z,\et}^d,\rD\rR(\calO_{\widetilde \frakB_Z^{d,I}},\nabla))$ is exactly the complex 
      \[\calO(\widetilde \frakB_Z^{d,I})\xrightarrow{\rd}\Omega^1(\widetilde \frakB_Z^{d,I})\{-1\} \xrightarrow{\rd}\Omega^2(\widetilde \frakB_Z^{d,I})\{-2\}\to \cdots.\]
      Since $\rd$ is $\CCtilde^I(Z)$-linear, we have 
      \[(\Omega^n(\widetilde \frakB_Z^{d,I})\{-n\})^{\rd = 0} = \Omega^n(\widetilde \frakB_Z^{d,I})^{\rd = 0}\{-n\}\]
      and
      \[\rd(\Omega^{n-1}(\widetilde \frakB_Z^{d,I})\{-(n-1)\}) = \rd(\xi\Omega^{n-1}(\widetilde \frakB_Z^{d,I})\{-n\}) = \xi\rd(\Omega^{n-1}(\widetilde \frakB_Z^{d,I}))\{-n\},\]
      which yields Item (1).
      Using Theorem \ref{thm:RH correspondence over dR} instead of Theorem \ref{thm:RH correspondence}, we obtain Item (2) in a similar way.
      The desired description of the corresponding $\rH^0$ follows from a direct calculation.
  \end{proof}
  \begin{cor}\label{cor:cohomology for unit ball}
      For any $n\geq 1$, the following assertions hold:
      \begin{enumerate}
          \item[(1)] When $1\in I$, there exists a short exact sequence
          \[0\to \rH^n_{\dR}(\widetilde \frakB_Z^{d,I})\{-(n-1)\}\to \rH_v^n(\bB_Z^d,\CCtilde^I)\to \Omega^n(\bB_Z^d)^{\rd = 0}\{-n\}\to 0,\]
          and when $1\notin I$, there is an isomorphism
          \[\rH^n_{\dR}(\widetilde \frakB_Z^{d,I})\{-n\}\cong \rH^n_v(\bB_Z^d,\CCtilde^I).\]

          \item[(2)] There exists a short exact sequence
          \[0\to \rH^n_{\dR}(\widetilde \bB_Z^d)\{-(n-1)\}\to \rH_v^n(\bB_Z^d,\BBdRp)\to \Omega^n(\bB_Z^d)^{\rd = 0}\{-n\}\to 0.\]
      \end{enumerate}
  \end{cor}
  \begin{proof}
      The Breuil--Kisin twists $\{-(n-1)\}$ and $\{-n\}$ here are the same as in Lemma \ref{lem:cohomology for unit ball}, hence consistent with the trivialization $\CCtilde^{I,+}\{1\} \cong \CCtilde^{I,+}\cdot\frac{t}{u}$ of \S\ref{sec:local RH}.
      For item (1): Note that we have a short exact sequence
      \[0\to \frac{\xi\Omega^n(\widetilde \frakB_Z^{d,I})^{\rd = 0}}{\xi \rd(\Omega^{n-1}(\widetilde \frakB_Z^{d,I}))}\{-n\}\to \frac{\Omega^n(\widetilde \frakB_Z^{d,I})^{\rd = 0}}{\xi \rd(\Omega^{n-1}(\widetilde \frakB_Z^{d,I}))}\{-n\}\to \frac{\Omega^n(\widetilde \frakB_Z^{d,I})^{\rd = 0}}{\xi \Omega^n(\widetilde \frakB_Z^{d,I})^{\rd = 0}}\{-n\}\to 0.\]
      Combining Lemma \ref{lem:cohomology for unit ball}(1) with the isomorphism
      \[\frac{\xi\Omega^n(\widetilde \frakB_Z^{d,I})^{\rd = 0}}{\xi \rd(\Omega^{n-1}(\widetilde \frakB_Z^{d,I}))}\{-n\} = \rH^n_{\dR}(\widetilde \frakB_Z^{d,I})\{-(n-1)\},\]
      we have that
      \[0\to \rH^n_{\dR}(\widetilde \frakB_Z^{d,I})\{-(n-1)\}\to \rH^n_v(\bB_Z^d,\CCtilde^I)\to \frac{\Omega^n(\widetilde \frakB_Z^{d,I})^{\rd = 0}}{\xi \Omega^n(\widetilde \frakB_Z^{d,I})^{\rd = 0}}\{-n\}\to 0.\]
      Since $\xi\in\CCtilde^I(Z)^{\times}$ is a unit when $1\notin I$ (cf. Example \ref{exam:t in C^I}(3)) and $\CCtilde^I(Z)/\xi = S$, we can conclude by noting that
      \[\frac{\Omega^n(\widetilde \frakB_Z^{d,I})^{\rd = 0}}{\xi \Omega^n(\widetilde \frakB_Z^{d,I})^{\rd = 0}}\{-n\} = \left\{\begin{array}{cc}
          \Omega^n(\bB_Z^d)^{\rd = 0}\{-n\}, & \text{if $1\in I$} \\
          0, & \text{if $1\notin I$.}
      \end{array}\right.\]

      Item (2) follows from Lemma \ref{lem:cohomology for unit ball} in a similar way.
  \end{proof}

  Let $\iota_p:\frakB_Z^d\to\frakB_Z^d$ be the morphism corresponding to the morphism
  \[S^+\za T_1,\dots,T_d\ya\to S^+\za T_1,\dots,T_d\ya\]
  sending each $T_i$ to $pT_i$.
  Denote the associated morphism $\bB_Z^d\to \bB_Z^d$ also by $\iota_p$. Via the identification $\bB_Z^d = \bB_Z^d(1)\subset \bB_Z^d(p)$, where $\bB_Z^d(r)$ denotes the closed ball of radius $r$ over $Z$, $\iota_p$ is the rational localization; thus $\bA_Z^d = \colim_{\iota_p}\bB_Z^d$.
  While $\bB_Z^d$ is the unit ball (radius $1$), this exhaustion by $\iota_p$ (transition from radius $1$ to radius $p$) correctly captures the overconvergent nature of the affine space that is needed for comparisons with rigid cohomology later.

  The following result is well-known.
  \begin{lem}\label{lem:usual dR cohomology of affine space}
      Let $(A,|\cdot|)$ be a Banach $\Qp$-algebra.
      Let $\bB_A^d$ be the closed unit ball of relative dimension $d$ over $A$ and $\iota_p:\bB_A^d\to\bB_A^d$ the map corresponding to the morphism
      \[A\za T_1,\dots,T_d\ya\to A\za T_1,\dots,T_d\ya\]
      sending each $T_i$ to $pT_i$, so that $\bA_A^d = \colim_{\iota_p}\bB_A^d$.
      Then
      \[\rH^n_{\dR}(\bA_A^d) = \left\{\begin{array}{rcl}
          A, & \text{if $n = 0$} \\
          0, & \text{if $n\geq 1$}.
      \end{array}\right.\]
      More precisely, for any $n\geq 1$ we have $\varprojlim_{\iota_p}\rH^n_{\dR}(\bB_A^d) = 0$, and for any $n\geq 0$ we have $\rR^1\varprojlim_{\iota_p}\rH^n_{\dR}(\bB_A^d) = 0$.
  \end{lem}
  \begin{proof}
      We provide a proof here for completeness.
      Note that we have a short exact sequence
      \begin{equation}\label{equ:usual dR cohomology of affine space}
          0\to\rR^1\varprojlim_{\iota_p}\rH^{n-1}_{\dR}(\bB_A^d)\to \rH^n_{\dR}(\bA_A^d)\to \varprojlim_{\iota_p}\rH_{\dR}^n(\bB_A^d)\to 0.
      \end{equation}
      Taking $n = 0$, we have 
      \[\rH^0_{\dR}(\bA_A^d)\cong \varprojlim_{\iota_p}\rH_{\dR}^0(\bB_A^d) = \varprojlim_{\iota_p}A = A,\]
      where we use the fact, which can be deduced by direct calculation, that $\rH_{\dR}^0(\bB_A^d) = A$ and $\iota_p$ acts as $\id_A$ on $A\subset A\za T_1,\dots,T_d\ya$.
      Also, this implies that $\rR^1\varprojlim_{\iota_p}\rH^0_{\dR}(\bB_A^d) = 0$.
      To conclude, by \eqref{equ:usual dR cohomology of affine space} it suffices to show that $\rH_{\dR}^n(\bA_A^d) = 0$ for any $n\geq 1$.
      For any $1\leq i\leq d$, put 
      \[\partial_i:=\frac{\partial}{\partial T_i}: A\za T_1,\dots,T_d\ya\to A\za T_1,\dots,T_d\ya\]
      and let
      \[\iota_i:A\za T_1,\dots,T_d\ya\to A\za T_1,\dots,T_d\ya\]
      be the morphism sending $T_j$ to $p^{\delta_{ij}}T_j$, where $\delta_{ij}$ denotes the Kronecker delta symbol.
      Then the usual de Rham complex of $\bB_Z^d$ is given by the Koszul complex $\rK(A\za T_1,\dots,T_d\ya;\partial_1,\dots,\partial_d)$:
      \[0\to A\za T_1,\dots,T_d\ya\xrightarrow{\partial_1,\dots,\partial_d}A\za T_1,\dots,T_d\ya^{\oplus d}\to\cdots\]
      and the usual de Rham complex of $\bA_A^d$ is given by 
      \[\varprojlim_{\iota_p}\rK(A\za T_1,\dots,T_d\ya;\partial_1,\dots,\partial_d) \simeq \varprojlim_{\iota_1}\big(\varprojlim_{\iota_2}(\cdots\varprojlim_{\iota_d}\rK(A\za T_1,\dots,T_d\ya;\partial_1,\dots,\partial_d))\big).\]
      We first deal with the $d = 1$ case.
      In other words, we have to show that the natural morphism
      \[A\to \varprojlim_{\iota_p}\rK(A\za T\ya;\partial:=\frac{\partial}{\partial T})\]
      is a quasi-isomorphism.
      To do this, it suffices to show that for any $f(T)\rd T$, there exists a $g(T)\in A\za T\ya$ such that 
      \[\iota_p(f(T)\rd T) = pf(pT)\rd T = (\partial g)(T)\rd T.\]
      This is clear: For any $f(T) = \sum_{n\geq 0}a_nT^n$, we can check
      \[g(T) = p\sum_{n\geq 0}\frac{a_np^n}{n+1}T^{n+1}\]
      is well-defined in $A\za T\ya$ and satisfies the desired condition.
      This completes the proof in the $d = 1$ case.
      In general, as $\rK(A\za T_1,\dots,T_d\ya;\partial_1,\dots,\partial_d)$ is the total complex of 
      \[\rK(A\za T_1,\dots,T_d\ya;\partial_1,\dots,\partial_{d-1})\xrightarrow{\partial_d}\rK(A\za T_1,\dots,T_d\ya;\partial_1,\dots,\partial_{d-1}),\]
      by what we have proved in the $d = 1$ case, it follows that
      \begin{equation*}
          \begin{split}
              \varprojlim_{\iota_1}\big(\varprojlim_{\iota_2}(\cdots\varprojlim_{\iota_d}\rK(A\za T_1,\dots,T_d\ya;\partial_1,\dots,\partial_d))\big)&\simeq \varprojlim_{\iota_1}\big(\varprojlim_{\iota_2}(\cdots\varprojlim_{\iota_{d-1}}\rK(A\za T_1,\dots,T_{d-1}\ya;\partial_1,\dots,\partial_{d-1}))\big)\\
              &\simeq \cdots\\
              &\simeq \varprojlim_{\iota_1}\rK(A\za T_1\ya;\partial_1)\\
              &\simeq A.
          \end{split}
      \end{equation*}
      This completes the proof.
  \end{proof}
  The following lemma will also be used.
  \begin{lem}\label{lem:Mittag-Leffler}
      Let $\iota:\{A_n\}_{n\geq 0}\to\{B_n\}_{n\geq 0}$ be a morphism of inverse systems, and denote by $f_n:A_{n+1}\to A_n$ and $g_n:B_{n+1}\to B_n$ the transition maps.
      Suppose that:
      \begin{enumerate}
          \item[(1)] For each $n\geq 0$, $B_n$ is a Banach algebra and $g_n$ is bounded: there exists $M_n>0$ such that $\norm g_n(x)\norm\leq M_n\norm x\norm$ for all $x\in B_{n+1}$.
          \item[(2)] The inverse system $\{A_n\}_{n\geq 0}$ satisfies the Mittag--Leffler condition.
          \item[(3)] For each $n\geq 0$, $\iota(A_n)\subset B_n$ is dense.
      \end{enumerate}
      Then $\rR^1\varprojlim B_n = 0$.
  \end{lem}
  \begin{proof}
      Recall that for any inverse system $\{C_n\}_{n\geq 0}$ of abelian groups with transition maps $h_n:C_{n+1}\to C_n$, there is an exact sequence
      \[0\to \varprojlim C_n\to \prod_{n\geq 0}C_n\xrightarrow{\delta}\prod_{n\geq 0}C_n\to \rR^1\varprojlim C_n\to 0,\]
      where $\delta$ sends any $(c_n)_{n\geq 0}\in \prod_{n\geq 0}C_n$ to $(c_n-h_n(c_{n+1}))_{n\geq 0}$.
      Thus, it suffices to show that for any $(b_n)_{n\geq 0}\in \prod_{n\geq 0}B_n$, there exists a sequence $(x_n)_{n\geq 0}\in \prod_{n\geq 0}B_n$ such that for any $n\geq 0$,
      \[x_n-g_n(x_{n+1}) = b_n.\]

      By Item (3), there exists a sequence $(a_n)_{n\geq 0}\in \prod_{n\geq 0}A_n$ such that for any $n\geq 0$, 
      \[\norm b_n-\iota(a_n)\norm\leq \frac{1}{2^n\prod_{i=0}^{n-1}M_i}.\]
      So, for any $r\geq 0$, by Item (1), we have
      \[\norm (g_{n}\circ\cdots\circ g_{n+r})(b_{n+r+1}-\iota(a_{n+r+1}))\norm\leq (\prod_{i=0}^rM_{n+i})\cdot\norm b_{n+r+1}-\iota(a_{n+r+1})\norm\leq \frac{1}{2^n\cdot(\prod_{i=0}^{n-1}M_i)\cdot2^{r+1}}.\]
      Since $B_n$ is a Banach space,
      \[z_n:=b_n-\iota(a_n)+\sum_{r\geq 0}(g_{n}\circ\cdots\circ g_{n+r})(b_{n+r+1}-\iota(a_{n+r+1}))\]
      is a well-defined element in $B_n$ such that
      \[z_n-g_n(z_{n+1}) = b_n-\iota(a_n)\]
      for any $n\geq 0$.
      On the other hand, by Item (2) we have $\rR^1\varprojlim A_n = 0$, which, in particular, implies that there exists a sequence $(y_n)_{n\geq 0}\in\prod_{n\geq 0}A_n$ such that
      \[y_n-f_n(y_{n+1}) = a_n\]
      for any $n\geq 0$.
      Now, for any $n\geq 0$, define 
      \[x_n:=z_n+\iota(y_n).\]
      Then we have
      \[x_n-g_n(x_{n+1}) = z_n-g_n(z_{n+1})+\iota(y_n)-g_n(\iota(y_{n+1})) = b_n-\iota(a_n)+\iota(y_n-f_n(y_{n+1})) = b_n\]
      as desired, which completes the proof.
  \end{proof}
  \begin{cor}\label{cor:usual dR cohomology of affine space}
      Keep notations as in Lemma \ref{lem:usual dR cohomology of affine space}.
      For any $n\geq 1$, we have
      \[\rR^1\varprojlim_{\iota_p}(\rd(\Omega^{n-1}(\bB_A^d))) = \rR^1\varprojlim_{\iota_p}(\Omega^n(\bB_A^d)^{\rd = 0}) = \rR^1\varprojlim_{\iota_p}\Omega^n(\bB_A^d) = 0.\]
      As a consequence, for any $n\geq 0$, we have
      \[\varprojlim_{\iota_p}(\Omega^n(\bB_A^d)^{\rd = 0}) = \Omega^n(\bA_A^d)^{\rd = 0}\]
      and 
      \[\varprojlim_{\iota_p}(\rd(\Omega^n(\bB_A^d))) = \rd(\Omega^n(\bA_A^d)).\]
  \end{cor}
  \begin{proof}
      We start with the following two short exact sequences:
      \begin{equation*}
          0\to \Omega^n(\bB_A^d)^{\rd = 0}\xrightarrow{\subset}\Omega^n(\bB_A^d)\xrightarrow{\rd}\rd(\Omega^n(\bB_A^d))\to 0
      \end{equation*}
      and 
      \begin{equation*}
          0\to \rd(\Omega^{n-1}(\bB_A^d))\to \Omega^n(\bB_A^d)^{\rd = 0}\to\rH_{\dR}^n(\bB_A^d)\to 0.
      \end{equation*}
      Taking $\rR\varprojlim_{\iota_p}$, we obtain the long exact sequences
      \begin{equation}\label{equ:usual dR cohomology of affine space II}
      \begin{split}
          0&\to \varprojlim_{\iota_p}\Omega^n(\bB_A^d)^{\rd = 0}\to \varprojlim_{\iota_p}\Omega^n(\bB_A^d)\xrightarrow{\rd}\varprojlim_{\iota_p}\rd(\Omega^n(\bB_A^d))\\
          &\to \rR^1\varprojlim_{\iota_p}\Omega^n(\bB_A^d)^{\rd = 0}\to \rR^1\varprojlim_{\iota_p}\Omega^n(\bB_A^d)\to\rR^1\varprojlim_{\iota_p}\rd(\Omega^n(\bB_A^d))\to 0
       \end{split}
      \end{equation}
      and 
      \begin{equation}\label{equ:usual dR cohomology of affine space III} 
      \begin{split}
          0&\to \varprojlim_{\iota_p}\rd(\Omega^{n-1}(\bB_A^d))\to \varprojlim_{\iota_p}\Omega^n(\bB_A^d)^{\rd = 0}\to\varprojlim_{\iota_p}\rH_{\dR}^n(\bB_A^d)\\
          &\to \rR^1\varprojlim_{\iota_p}\rd(\Omega^{n-1}(\bB_A^d))\to \rR^1\varprojlim_{\iota_p}\Omega^n(\bB_A^d)^{\rd = 0}\to\rR^1\varprojlim_{\iota_p}\rH_{\dR}^n(\bB_A^d)\to 0.
      \end{split}
      \end{equation}
      Applying Lemma \ref{lem:usual dR cohomology of affine space} to \eqref{equ:usual dR cohomology of affine space III}, for any $n\geq 1$, we have
      \[\varprojlim_{\iota_p}\rd(\Omega^{n-1}(\bB_A^d))\cong \varprojlim_{\iota_p}\Omega^n(\bB_A^d)^{\rd = 0} \text{ and }\rR^1\varprojlim_{\iota_p}\rd(\Omega^{n-1}(\bB_A^d))\cong \rR^1\varprojlim_{\iota_p}\Omega^n(\bB_A^d)^{\rd = 0}.\]

      Assuming $\rR^1\varprojlim_{\iota_p}\Omega^m(\bB_A^d) = 0$ for all $m\geq 0$, \eqref{equ:usual dR cohomology of affine space II} implies that for any $n\geq 1$,
      \[\rR^1\varprojlim_{\iota_p}\rd(\Omega^{n-1}(\bB_A^d)) = 0\]
      and thus
      \[\rR^1\varprojlim_{\iota_p}\Omega^n(\bB_A^d)^{\rd = 0}  = 0.\]
      Then, by \eqref{equ:usual dR cohomology of affine space II} again, we have
      \[\varprojlim_{\iota_p}\Omega^n(\bB_A^d)^{\rd = 0}\cong(\varprojlim_{\iota_p}\Omega^n(\bB_A^d))^{\rd = 0} = \Omega^n(\bA_A^d)^{\rd = 0}\]
      and that
      \[\varprojlim_{\iota_p}\rd(\Omega^n(\bB_A^d)) = \rd(\Omega^n(\bA_A^d)),\]
      as desired.
      It remains to prove that $\rR^1\varprojlim_{\iota_p}\Omega^m(\bB_A^d) = 0$ for any $m\geq 0$.
      Equip $A\za T_1,\dots,T_d\ya$ with the Gauss norm associated to $A$.
      This gives rise to a norm on $\Omega^m(\bB_A^d)$ that makes it a Banach $A$-module such that for any $\omega\in \Omega^m(\bB_A^d)$, we have
      \[\norm\iota_p(\omega)\norm\leq p^{-m}\norm\omega\norm.\]
      The key to the Mittag--Leffler argument is that the \emph{algebraic} $m$-forms are dense in the Banach space of forms on the unit ball, and that $\iota_p$ is surjective on them.
      Let $\Omega^m_{\bA_A^{d,{\rm alg}}}$ be the module of $m$-forms on the algebraic affine space $\bA_A^{d,{\rm alg}} = \Spec(A[T_1,\dots,T_d])$.
      Then $\Omega^m_{\bA_A^{d,{\rm alg}}}\subset \Omega^m(\bB_A^d)$ is dense, and the restriction $\iota_p:\Omega^m_{\bA_A^{d,{\rm alg}}}\to \Omega^m_{\bA_A^{d,{\rm alg}}}$ is surjective; thus the inverse system 
      \[\cdots\to\Omega^m_{\bA_A^{d,{\rm alg}}}\xrightarrow{\iota_p}\Omega^m_{\bA_A^{d,{\rm alg}}}\xrightarrow{\iota_p}\cdots\xrightarrow{\iota_p}\Omega^m_{\bA_A^{d,{\rm alg}}}\]
      satisfies the Mittag--Leffler condition.
      Lemma \ref{lem:Mittag-Leffler} then yields $\rR^1\varprojlim_{\iota_p}\Omega^m(\bB_A^d) = 0$: density of algebraic forms provides the approximation (condition (3)), and surjectivity of $\iota_p$ on them provides the Mittag--Leffler condition (2) for the subsystem.
  \end{proof}

  \begin{prop}\label{prop:cohomology for affine space}
      The following assertions are true.
      \begin{enumerate}
          \item[(1)] For $n = 0$, we have $\rH^0_v(\bA_Z^d,\CCtilde^I) = \CCtilde^I(Z)$ and $\rH^0_v(\bA_Z^d,\BBdRp) = \BBdRp(Z)$.
          \item[(2)] For $n \geq 1$, we have
          \[\rH^n_v(\bA_Z^d,\CCtilde^I) = \left\{\begin{array}{rcl}
                  \Omega^n(\bA_Z^d)^{\rd = 0}\{-n\}, & \text{if $1\in I$} \\
                  0, & \text{if $1\notin I$}
              \end{array}\right.\]
          and   
          \[\rH^n_v(\bA_Z^d,\BBdRp) = \Omega^n(\bA_Z^d)^{\rd = 0}\{-n\}.\]
      \end{enumerate}
  \end{prop}
  \begin{proof}
      The $n=0$ case is trivial according to Lemma \ref{lem:cohomology for unit ball}.
      For any $n\geq 1$, consider the short exact sequence
      \begin{equation}\label{equ:cohomology for affine space}
          0\to\rR^1\varprojlim_{\iota_p}\rH^{n-1}_v(\bB_Z^d,\CCtilde^I)\to \rH^n_v(\bA_Z^d,\CCtilde^I)\to\varprojlim_{\iota_p}\rH^n_v(\bB_Z^d,\CCtilde^I)\to 0.
      \end{equation}
      When $1\notin I$, by Corollary \ref{cor:cohomology for unit ball}(1), we have 
      \[\rH^n_v(\bA_Z^d,\CCtilde^I)\cong \rH^n_{\dR}(\widetilde \frakB_Z^{d,I}),\]
      where $\widetilde \frakB_Z^{d,I} = \bB_{\CCtilde^I(Z)}^d$.
      Applying Lemma \ref{lem:usual dR cohomology of affine space} (for $A = \CCtilde^I(Z)$) to \eqref{equ:cohomology for affine space}, we conclude that
      \[\rH^n_v(\bA_Z^d,\CCtilde^I) = 0\]
      as desired.       When $1\in I$, by Corollary \ref{cor:cohomology for unit ball}(1) again, we have the long exact sequence
      \[\begin{split}0&\to \varprojlim_{\iota_p}\rH^n_{\dR}(\widetilde \frakB_Z^{d,I})\{-(n-1)\}\to \varprojlim_{\iota_p}\rH_v^n(\bB_Z^d,\CCtilde^I)\to \varprojlim_{\iota_p}\Omega^n(\bB_Z^d)^{\rd = 0}\{-n\}\\
      &\to\rR^1\varprojlim_{\iota_p}\rH^n_{\dR}(\widetilde \frakB_Z^{d,I})\{-(n-1)\}\to \rR^1\varprojlim_{\iota_p}\rH_v^n(\bB_Z^d,\CCtilde^I)\to \rR^1\varprojlim_{\iota_p}\Omega^n(\bB_Z^d)^{\rd = 0}\{-n\}\to 0.\end{split}\]
      By Lemma \ref{lem:usual dR cohomology of affine space}, we have
      \[\varprojlim_{\iota_p}\rH^n_{\dR}(\widetilde \frakB_Z^{d,I})\{-(n-1)\} = 0.\]
      By Corollary \ref{cor:usual dR cohomology of affine space}, we have
      \[\rR^1\varprojlim_{\iota_p}\rH^n_{\dR}(\widetilde \frakB_Z^{d,I})\{-(n-1)\} = 0 = \rR^1\varprojlim_{\iota_p}\Omega^n(\bB_Z^d)^{\rd = 0}\{-n\}.\]
      
      So for any $n\geq 1$ we obtain
      \[\varprojlim_{\iota_p}\rH_v^n(\bB_Z^d,\CCtilde^I)\cong \varprojlim_{\iota_p}\Omega^n(\bB_Z^d)^{\rd = 0}\{-n\} =\Omega^n(\bA_Z^d)^{\rd = 0}\{-n\}\]
      and
      \[\rR^1\varprojlim_{\iota_p}\rH_v^n(\bB_Z^d,\CCtilde^I) = 0.\]
      Since $\rR^1\varprojlim_{\iota_p}\rH_v^0(\bB_Z^d,\CCtilde^I) = \rR^1\varprojlim_{\iota_p}\CCtilde^I(Z) = 0$, \eqref{equ:cohomology for affine space} gives
      \[\rH^n_v(\bA_Z^d,\CCtilde^I)\cong\varprojlim_{\iota_p}\rH^n_v(\bB_Z^d,\CCtilde^I)\cong \Omega^n(\bA_Z^d)^{\rd = 0}\{-n\}\]
      as desired.
      This completes the proof of Item (2) for $\CCtilde^I$. 

      By the same argument as in the $1\in I$ case, we get
      \[\rH^n_v(\bA_Z^d,\BBdRp) = \Omega^n(\bA_Z^d)^{\rd = 0}\{-n\}\]
      for any $n\geq 1$ as desired.
  \end{proof}
  Now, we are able to prove Theorem \ref{thm:pro-etale cohomology of Affine space}.
  \begin{proof}[\textbf{Proof of Theorem \ref{thm:pro-etale cohomology of Affine space}}]
      The $\rH^n_v(\bA_Z^d,\BBdRp)$ has been calculated in Proposition \ref{prop:cohomology for affine space}.
      It remains to compute $\rH^n_v(\bA_Z^d,\Qp)$.
      It suffices to show that 
      \[\rH^n_v(\bA_Z^d,\Qp(1)) = \left\{
      \begin{array}{rcl}
          \Qp(1), & n=0 \\
          \Omega^n(\bA_Z^d)^{\rd = 0}\{-(n-1)\}, & n\geq 1,
      \end{array}\right.\]
      where $\Qp(1)$ denotes the level-one Tate twist of $\Qp$.
      By Example \ref{exam:t in C^I}, we always have $\Ctilde^J\{1\} = \Ctilde^J(1)$ when $J$ is a closed sub-interval of $(\frac{1}{p-1},p)$.
      By \cite[Prop. 6.2.2]{KL15}, we have the fundamental exact sequence
      \[0\to\Qp(1)\to \CCtilde^I(1)\xrightarrow{\varphi-1}\CCtilde^{I\cap I/p}(1)\to 0.\]
      Under the Riemann--Hilbert correspondence, the pro-\'etale $\Qp$-local system $\Qp(1)$ corresponds to the $\varphi$-module $\CCtilde^I(1)$ over $\CCtilde^I$ (Remark \ref{rmk:Qp-local system as phi-module}); by Theorem \ref{thm:prismatic RH}, this identification is compatible with the equivalence between small local systems and small connections. Hence
      \begin{equation}\label{equ:pro-etale cohomology of Affine space}
          \rR\Gamma(\bA_{Z,v}^d,\Qp(1)) = [\rR\Gamma(\bA_{Z,v}^d,\CCtilde^I(1))\xrightarrow{\varphi-1}\rR\Gamma(\bA_{Z,v}^d,\CCtilde^{I\cap I/p}(1))],
      \end{equation}
      where we take $I\subset(\frac{1}{p-1},p)$ with $I\cap I/p\neq \emptyset$.
      This condition is critical: for such $I$ one has $1\in I$ and $1\notin I/p$ (so $\xi$ is a unit in $\CCtilde^I(Z)$ but not in $\CCtilde^{I\cap I/p}(Z)$), which will imply that the $E_2^{n1}$ terms vanish for $n\geq 1$ in the spectral sequence below.
      The complex \eqref{equ:pro-etale cohomology of Affine space} is the mapping cone of $\varphi-1$, so it gives rise to a convergent spectral sequence
      \[\rE_2^{ij}\Rightarrow \rH^n_v(\bA_Z^d,\Qp(1)),\]
      with $\rE_2^{ij} \neq 0$ only if $0\leq i\leq d$ and $j\in\{0,1\}$, and
      \[0\to \rE_2^{n0}\xrightarrow{\subset}\rH^n_v(\bA^d_Z,\CCtilde^I(1))\xrightarrow{\varphi-1}\rH^n_v(\bA_Z^d,\CCtilde^{I\cap I/p}(1))\to \rE_2^{n1}\to 0.\]
      By Proposition \ref{prop:cohomology for affine space}, we have
      \[\rE_2^{00}=\CCtilde^I(1)^{\varphi = 1} = \Qp(1) \text{ and }\rE_2^{01} = \CCtilde^{I\cap I/p}(1)/(\varphi-1)(\CCtilde^I(1)) = 0,\]
      and for any $n\geq 1$,
      \[\rE_2^{n0} = \Omega^n(\bA_Z^d)^{\rd = 0}\{-(n-1)\}\text{ and }\rE_2^{n1} = 0.\]
      The vanishing $\rE_2^{n1}=0$ for $n\geq 1$ is precisely what we deduced from the condition $1\in I$ and $1\notin I/p$.
      Thus
      \[\rH^n_v(\bA_Z^d,\Qp(1)) = \rE_2^{n0} = \left\{
      \begin{array}{rcl}
          \Qp(1), & n=0 \\
          \Omega^n(\bA_Z^d)^{\rd = 0}\{-(n-1)\}, & n\geq 1,
      \end{array}\right.\]
      as desired.
  \end{proof}
  \begin{rmk}
      There are several ways to compute $\rR\Gamma(\bA_{Z,v}^d,\Qp)$.
      In \cite{LB18}, Le Bras used the fundamental exact sequence
      \[0\to \Qp\to \CCtilde^I[\frac{1}{t}]^{\varphi=1}\to\BBdR/\BBdRp\to 0\]
      and then reduced the problem to computing 
      \[\rR\Gamma(\bA_{Z,v}^d,\CCtilde^I[\frac{1}{t}])^{\varphi=1}\text{ and }\rR\Gamma(\bA_{Z,v}^d,\BBdRp).\]
      In \cite{CN20}, Colmez and Nizioł computed $\rR\Gamma(\bA_{Z,v}^d,\Qp)$ using Hyodo--Kato cohomology.
      Compared with their methods, ours is more direct: it avoids the complexity of Hyodo--Kato cohomology by relying on the ``smallness'' property inherent in the Riemann--Hilbert correspondence established in \S\ref{sec:global RH}; after applying that correspondence, the result follows from the well-known de Rham cohomology of $\bA_Z^d$.
  \end{rmk}

\section{Application II: Comparison with Rigid Cohomology}\label{sec:rigid cohomology}

  Fix a prime $p\geq 3$ and let $I \subset (\frac{1}{p-1},p)$ be a closed sub-interval. Let $Z = \Spa(A,A^+)\in \Perfd$ be an affinoid perfectoid space. Throughout this section, we set $\calO_K = \rW(\kappa)$. We assume that $\frakX = \frakX_0\times_{\Spf(\calO_K)}\Spf(\calO_C)$ for some quasi-compact, quasi-separated smooth formal scheme $\frakX_0$ over $\calO_K$, and we fix the lifting $\wtx = \frakX_0\times_{\Spf(\calO_K)}\Spf(\bfA_{\inf})$ of $\frakX$. As Remark \ref{rmk:Zariski case} applies, we retain the notations established therein. For either period sheaf $\bP\in \{\CCtilde^I,\BBdRp\}$, we construct a perfect complex of $\bP(Z)$-modules $\rR\Gamma^{\dagger}(X_Z,\bP)$ and prove the following theorem.

\begin{thm}\label{thm:P-cohomology}
    Keep the notation as above and denote by $\calX$ the special fiber of $\frakX_0$. Suppose that $I = [r,s]\subset (\frac{1}{p-1},p)$ and let $\bP\in \{\CCtilde^I,\BBdRp\}$.
    \begin{enumerate}
        \item The Frobenius morphism $\phi:\CCtilde^r\to \CCtilde^{rp}$ induces a quasi-isomorphism
        \[\rR\Gamma^{\dagger}(X_Z,\CCtilde^r)\xrightarrow{\simeq}\rR\Gamma^{\dagger}(X_Z,\CCtilde^{rp}).\]

        \item If $1\in I$, then the natural morphism $\CCtilde^I\to\BBdRp$ induces a quasi-isomorphism
        \[\rR\Gamma^{\dagger}(X_Z,\CCtilde^I)\otimes_{\CCtilde^I(Z)}\BBdRp(Z)\xrightarrow{\simeq}\rR\Gamma^{\dagger}(X_Z,\BBdRp).\]

        \item Denote by $\rR\Gamma_{\rig}(\calX/K)$ the rigid cohomology of $\calX$. Then there exists a canonical quasi-isomorphism
        \[\iota_{\bP}:\rR\Gamma_{\rig}(\calX/K)\otimes_K\bP(Z)\xrightarrow{\simeq}\rR\Gamma^{\dagger}(X_Z,\bP).\]
        This quasi-isomorphism is $\phi$-equivariant: the following diagram commutes, where the left vertical $\phi$ is the Frobenius on rigid cohomology $\rR\Gamma_{\rig}(\calX/K)$ (via the coefficient map $\CCtilde^r(Z)\to \CCtilde^{rp}(Z)$) and the right vertical $\phi$ is the Frobenius on period cohomology constructed in \S\ref{sec:local Frobenius action}.
        \begin{equation}\label{equ:Frobenius-equivariant}
            \begin{tikzcd}
                {\rR\Gamma_{\rig}(\calX/K)\otimes_K\CCtilde^r(Z)} \arrow[rr, "\iota_{\CCtilde^r}"] \arrow[d,"\phi"] &  & {\rR\Gamma^{\dagger}(X_Z,\CCtilde^r)} \arrow[d,"\phi"] \\
                \rR\Gamma_{\rig}(\calX/K)\otimes_K\CCtilde^{rp}(Z) \arrow[rr, "\iota_{\CCtilde^{rp}}"] &  & {\rR\Gamma^{\dagger}(X_Z,\CCtilde^{rp})}
            \end{tikzcd}
        \end{equation}
    \end{enumerate}
\end{thm}
\subsection{The construction of $\rR\Gamma^{\dagger}(X_Z,\bP)$}

Let $\frakX_0$ be a smooth formal scheme over $\calO_K$. Assume we have embeddings
\[\frakX_0\to\frakY_0\hookrightarrow\overline \frakY_0\]
with $\overline \frakY_0$ proper smooth over $\calO_K$ and $\frakY_0\subset \overline \frakY_0$ a closed sub-scheme.
Let $\calX\to\calY\hookrightarrow\overline \calY$ and $X_0\to Y_0\hookrightarrow \overline Y_0$ be the special and generic fibers of these embeddings, $\sp:\overline Y_0\to \overline \calY$ the specialization map, and $]\calX[_{\overline Y_0}\subset \overline Y_0$ the preimage of $\calX\subset \overline \calY$. Base-changing the open embedding $]\calX[_{\overline Y_0}\subset \overline Y_0$ along $K\to \CCtilde^I(Z)$ gives $]\widetilde X_Z^I[_{{\widetilde{\overline Y}}_Z^I}\subset \widetilde{\overline Y}_Z^I$.

Define $j^{\dagger}:=\colim_Vj_{\widetilde V_Z^I,*}j_{\widetilde V_Z^I}^{-1}$ on abelian sheaves on $]\wtY_Z^I[_{\widetilde{\overline Y}_Z^I}$, where $j_V:V\hookrightarrow ]\calY[_{\overline Y_0}$ runs over all strict neighborhoods of $]\calX[_{\overline Y_0}$ in $]\calY[_{\overline Y_0}$ (cf. \cite[Def. 3.1.1 and Prop. 5.1.12]{LeStum}) and $\widetilde V_Z^I$ is the base-change of $V$ along $K\to\CCtilde^I(Z)$. By \cite[\S1.2]{Ber97}, $j^{\dagger}$ is exact. For any Zariski covering $\{\frakX_{0,\alpha}\to\frakX_0\}_{\alpha\in\Lambda}$ of $\frakX_0$ with $\Lambda$ finite and $\frakX_{0,\Sigma}:=\bigcap_{\alpha\in \Sigma}\frakX_{0,\alpha}$ for $\Sigma\subset \Lambda$, and for any abelian sheaf $E$ on $]\wtY_Z^I[_{\widetilde{\overline Y}_Z^I}$, we have an exact sequence
\begin{equation}\label{equ:simplicial sequence}
    0\to j^{\dagger}E\to \prod_{\alpha\in\Lambda}j_{\alpha}^{\dagger}E\to \prod_{\{\alpha_1,\alpha_2\}\subset\Lambda}j_{\alpha_1,\alpha_2}^{\dagger}E\to\cdots\to j_{\Lambda}^{\dagger}E\to 0.
\end{equation}

\begin{dfn}\label{dfn:P-cohomology}
    Let $\frakX_0$ be a quasi-compact smooth formal scheme over $\calO_K$ as above, and let $\bP\in\{\CCtilde^I,\BBdRp\}$.
    \begin{enumerate}
        \item Suppose there is a closed embedding $\frakY_0\hookrightarrow\overline \frakY_0$ of formal schemes such that we have embeddings $\frakX_0\to\frakY_0\hookrightarrow\overline \frakY_0$ as above. Let $\nu:\overline Y_{Z,v}\to \overline \frakY_{Z,\Zar}$ be the morphism from the $v$-site of the generic fiber $\overline Y_Z$ to the Zariski site of the formal scheme $\overline \frakY_Z$. Define
        \[\rR\Gamma^{\dagger}(X_Z,\bP):=\rR\Gamma(]\widetilde X_Z^I[_{\widetilde{\overline Y}_Z^I},j^{\dagger}\iota^*\rL\eta_t\rR\nu_*\bP).\]
        Here, $\iota^*: \overline \frakY_{Z,\Zar}\to \widetilde{\overline Y}_{Z,an}^{I}$ is the functor described in Remark \ref{rmk:Zariski case}, and $j^{\dagger} := \colim_{V}j_{\widetilde V_Z^I,\ast}j_{\widetilde V_Z^I}^{-1}$ is defined as above. As shown in \cite[Prop. 6.5.4]{LeStum}, the complex $\rR\Gamma^{\dagger}(X_Z,\bP)$ is independent of the choice of the embeddings $\frakX_0\to\frakY_0\hookrightarrow\overline \frakY_0$.

        \item In general, let $\{\frakX_{0,\alpha}\to\frakX_0\}_{\alpha\in \Lambda}$ be a Zariski covering of $\frakX_0$ with $\Lambda$ finite such that for each $\alpha\in \Lambda$, there are embeddings $\frakX_{0,\alpha}\to \frakY_{0,\alpha}\hookrightarrow\overline \frakY_{0,\alpha}$. Define
        \[\rR\Gamma^{\dagger}(X_Z,\bP):=\rR\varprojlim\big( \prod_{\alpha\in\Lambda}\rR\Gamma^{\dagger}(X_{\alpha,Z},\bP)\to \prod_{\{\alpha_1,\alpha_2\}\subset\Lambda}\rR\Gamma^{\dagger}(X_{\{\alpha_1,\alpha_2\},Z},\bP)\to\cdots\to \rR\Gamma^{\dagger}(X_{\Lambda},\bP)\big).\]
    \end{enumerate}
\end{dfn}

\begin{rmk}\label{rmk:rigid cohomology}
    Recall that if there are embeddings $\frakX_0\to\frakY_0\hookrightarrow\overline \frakY_0$ as in Definition \ref{dfn:P-cohomology}(1), then the rigid cohomology of $\calX$ over $K$ is defined by
    \[\rR\Gamma_{\rig}(\calX/K):=\rR\Gamma(]\calY[_{\overline Y_0},j^{\dagger}\rD\rR(\calO_{\overline Y_0},\rd)),\]
    where $j^{\dagger} = \colim_Vj_{V,\ast}j_V^{-1}$ with $j_V:V\hookrightarrow]\calY[_{\overline Y_0}$ running over all strict neighborhoods of $]\calX[_{\overline Y_0}$ in $]\calY[_{\overline Y_0}$, and $\rD\rR(\calO_{\overline Y_0},\rd)$ denotes the usual de Rham complex associated to the standard trivial connection $\calO_{\overline Y_0}\xrightarrow{\rd}\Omega^1_{\overline Y_0}$ on $\overline Y_0$. This definition is independent of the choice of embeddings $\frakX_0\to\frakY_0\hookrightarrow\overline \frakY_0$.
    
    In general, let $\{\frakX_{0,\alpha}\to\frakX_0\}_{\alpha\in\Lambda}$ be the Zariski covering of $\frakX_0$ as in Definition \ref{dfn:P-cohomology}(2). Then we have 
    \begin{equation}\label{equ:simplicial for rigid cohomology}
        \rR\Gamma_{\rig}(\calX/K) = \rR\varprojlim\big(\prod_{\alpha\in\Lambda}\rR\Gamma_{\rig}(\calX_{\alpha}/K)\to\prod_{\{\alpha_1,\alpha_2\}\subset\Lambda}\rR\Gamma_{\rig}(\calX_{\{\alpha_1,\alpha_2\}}/K)\to\cdots\to\rR\Gamma_{\rig}(\calX_{\Lambda}/K)\big).
    \end{equation}
    This definition is independent of the choice of covering $\{\frakX_{0,\alpha}\to\frakX_0\}_{\alpha\in\Lambda}$.
\end{rmk}

We now construct the comparison morphism $\iota_{\bP}$ for $\bP\in \{\CCtilde^I,\BBdRp\}$.

\begin{construction}\label{construction:comparison map}
    We keep the above notation and the site morphism $\nu:\overline Y_{Z,v}\to \overline \frakY_{Z,\Zar}$ from Definition \ref{dfn:P-cohomology}(1).

    Suppose first that we have embeddings $\frakX_0\to\frakY_0\hookrightarrow\overline \frakY_0$. For $\bP = \CCtilde^I$, Theorem \ref{thm:RH correspondence} gives a quasi-isomorphism
    \[\rR\nu_*\CCtilde^I\simeq \rD\rR(\calO_{\widetilde{\overline \frakY}_Z^I},\nabla),\]
    where $\nabla:\calO_{\widetilde{\overline \frakY}_Z^I}\to \Omega^1_{\widetilde{\overline \frakY}_Z^I}\{-1\}$ is the trivial connection in the sense of Definition \ref{dfn:small connection}. The $\rL\eta_t$-functor (cf. \cite[\S 6]{BMS18}) ``shifts'' the cohomology of this complex so that it matches the integral structures required for rigid cohomology: since $t$ is a unit-multiple of $\xi$ in $\CCtilde^I$ (Example \ref{exam:t in C^I}), we obtain a quasi-isomorphism
    \[\rL\eta_t\rR\nu_*\CCtilde^I\simeq \rD\rR(\calO_{\widetilde{\overline \frakY}_Z^I},\rd),\]
    where $\rd:\calO_{\widetilde{\overline \frakY}_Z^I}\to \Omega^1_{\widetilde{\overline \frakY}_Z^I}$ is the usual trivial connection on $\wtx_Z^I$. By Remark \ref{rmk:Zariski case}, we get a quasi-isomorphism
    \[j^{\dagger}\iota^*\rL\eta_t\rR\nu_*\CCtilde^I\simeq j^{\dagger}\rD\rR(\calO_{\widetilde{\overline Y}_Z^I},\rd),\]
    where $\rd:\calO_{\widetilde{\overline Y}_Z^I}\to \Omega^1_{\widetilde{\overline Y}_Z^I}$ is the usual trivial connection on $\widetilde{\overline Y}_Z^I$. Thus, we have 
    \[\rR\Gamma^{\dagger}(X_Z,\CCtilde^I)\simeq \rR\Gamma(]\widetilde X_Z^I[_{\widetilde{\overline Y}_Z^I},j^{\dagger}\rD\rR(\calO_{\widetilde{\overline Y}_Z^I},\rd)).\]
    Comparison with Remark \ref{rmk:rigid cohomology} shows that scalar extension along $K\to\CCtilde^I(Z)$ induces
    \[\iota_{\CCtilde^I}:\rR\Gamma_{\rig}(\calX/K)\otimes_K\CCtilde^I(Z)\to \rR\Gamma^{\dagger}(X_Z,\CCtilde^I).\]
    For general $\frakX_0$, a Zariski covering as in Definition \ref{dfn:P-cohomology}(2) yields a natural morphism $\iota_{\CCtilde^I}$ by the same base-change. This completes the construction for $\bP = \CCtilde^I$.

    For $\bP = \BBdRp$, Theorem \ref{thm:RH correspondence over dR} and base-change along $K\to\BBdRp(Z)$ give
    \[\iota_{\BBdRp}:\rR\Gamma_{\rig}(\calX/K)\otimes_K\BBdRp(Z)\to \rR\Gamma^{\dagger}(X_Z,\BBdRp).\]
\end{construction}

\begin{construction}\label{construction:Frobenius}
    Assume $I = [r,rp]$ and that we have embeddings $\frakX_0\to\frakY_0\hookrightarrow\overline \frakY_0$ as in Definition \ref{dfn:P-cohomology}. Let $\nu: \overline Y_{Z,v}\to \overline \frakY_{Z,\Zar}$ be the morphism of sites. The Frobenius $\phi:\CCtilde^r\to\CCtilde^{rp}$ induces a quasi-isomorphism
    \[\rL\eta_t\rR\nu_*\CCtilde^r\xrightarrow{\phi}\rL\eta_{\phi(t)}\rR\nu_*\CCtilde^{rp} = \rL\eta_{pt}\rR\nu_*\CCtilde^{rp}\simeq \rL\eta_t\rR\nu_*\CCtilde^{rp}.\]
    (For any complex of sheaves of $\Qp$-vector spaces $\bK$, \cite[Lem. 6.11]{BMS18} gives $\rL\eta_{pt}\bK\simeq \rL\eta_t\rL\eta_p\bK\simeq\rL\eta_t\bK$.) Hence we obtain a quasi-isomorphism
    \[\phi:j^{\dagger}\iota^*\rL\eta_t\rR\nu_*\CCtilde^r\xrightarrow{\simeq} j^{\dagger}\iota^*\rL\eta_t\rR\nu_*\CCtilde^{rp}\]
    and thus a quasi-isomorphism
    \[\phi:\rR\Gamma^{\dagger}(]\wtX_Z^I[_{\widetilde{\overline Y}_Z^I},j^{\dagger}\iota^*\rL\eta_t\rR\nu_*\CCtilde^r)\xrightarrow{\simeq}\rR\Gamma^{\dagger}(]\wtX_Z^I[_{\widetilde{\overline Y}_Z^I},j^{\dagger}\iota^*\rL\eta_t\rR\nu_*\CCtilde^{rp}).\]
    
    Suppose that we also have embeddings $\frakX_0\to\frakY_0^{\prime}\to\overline \frakY_0^{\prime}$. Then we get the diagonal embedding
    \[\frakX_0\to \frakY_0\times\frakY_0^{\prime}\to \overline \frakY_0\times\overline \frakY_0^{\prime}.\]
    The projection $\overline \frakY_0\times\overline \frakY_0^{\prime}\to \overline \frakY_0$ induces a natural morphism
    \[\rR\Gamma^{\dagger}(]\wtX_Z^I[_{\widetilde{\overline Y}_Z^I},j^{\dagger}\iota^*\rL\eta_t\rR\nu_*\CCtilde^I)\to \rR\Gamma^{\dagger}(]\wtX_Z^I[_{\widetilde{\overline Y}_Z^I\times\widetilde{\overline Y}_Z^{\prime I}},j^{\dagger}\iota^*\rL\eta_t\rR\nu_*\CCtilde^I)\]
    which is a quasi-isomorphism by the argument of \cite[Prop. 6.5.4]{LeStum}. It follows directly that this morphism is compatible with Frobenius. Considering the projection $\overline \frakY_0\times\overline \frakY_0^{\prime}\to \overline \frakY_0^{\prime}$, we then obtain a quasi-isomorphism
    \[\rR\Gamma^{\dagger}(]\wtX_Z^I[_{\widetilde{\overline Y}_Z^I},j^{\dagger}\iota^*\rL\eta_t\rR\nu_*\CCtilde^I)\to \rR\Gamma^{\dagger}(]\wtX_Z^I[_{\widetilde{\overline Y}_Z^{\prime I}},j^{\dagger}\iota^*\rL\eta_t\rR\nu_*\CCtilde^I)\]
    which is compatible with Frobenius. This defines the Frobenius morphism
    \[\phi:\rR\Gamma^{\dagger}(X_Z,\CCtilde^r)\xrightarrow{\simeq}\rR\Gamma^{\dagger}(X_Z,\CCtilde^{rp}),\]
    so $\phi$ is a quasi-isomorphism in this case. For general $\frakX_0$, the Frobenius morphism $\phi:\rR\Gamma^{\dagger}(X_Z,\CCtilde^r)\xrightarrow{\simeq}\rR\Gamma^{\dagger}(X_Z,\CCtilde^{rp})$ is defined via Definition \ref{dfn:P-cohomology}(2) and is again a quasi-isomorphism.
\end{construction}

\begin{prop}\label{prop:Frobenius compatibility}
    The comparison morphism $\iota_{\CCtilde^I}$ is Frobenius-equivariant in the sense that \eqref{equ:Frobenius-equivariant} commutes.
\end{prop}
\begin{proof}
    By Constructions \ref{construction:comparison map} and \ref{construction:Frobenius} and \eqref{equ:simplicial sequence}, the assertion is Zariski-local on $\frakX_{0,\Zar}$. We may therefore assume embeddings $\frakX_0\to \frakY_0\to \overline \frakY_0$ with $\overline \frakY_0$ proper smooth and a morphism $\phi:\overline \frakY_0\to \overline \frakY_0$ over the Frobenius $\phi:\calO_K\to \calO_K$ lifting the absolute Frobenius $\phi:\overline \calY\to\overline \calY$. Then $\widetilde{\overline \frakY}_Z:=\overline \frakY_0\times_{\Spf(\calO_K)}\Spf(\AAinf(Z))$ is a prismatic lifting of $\overline \frakY_Z$, and $\phi$ induces a Frobenius morphism on the de Rham complex $\rD\rR(\calO_{\widetilde{\overline \frakY}_Z},\rd)$ for the usual differential $\rd:\calO_{\widetilde{\overline \frakY}_Z}\to\Omega^1_{\widetilde{\overline \frakY}_Z}$.
    
    Let $\nu:\overline Y_{Z,v}\to \overline \frakY_{Z,\Zar}$ be the morphism of sites. By Theorem \ref{thm:prismatic RH}, as discussed at the beginning of Construction \ref{construction:comparison map}, we get a quasi-isomorphism
    \[\rL\eta_t\rR\nu_*\CCtilde^I\simeq \rD\rR(\calO_{\widetilde{\overline \frakY}_Z^I},\rd)\]
    which is compatible with the Frobenius morphism on both sides. Thus, we have a Frobenius-equivariant quasi-isomorphism
    \begin{equation}\label{equ:Frobenius compatibility-I}
      \rR\Gamma^{\dagger}(X_Z,\CCtilde^I) = \rR\Gamma(]\widetilde X_Z^I[_{\widetilde{\overline Y}_Z^I},j^{\dagger}\iota^*\rL\eta_t\rR\nu_*\CCtilde^I)\simeq \rR\Gamma(]\widetilde X_Z^I[_{\widetilde{\overline Y}_Z^I},j^{\dagger}\iota^*\rD\rR(\calO_{\widetilde{\overline \frakY}_Z^I},\rd)).
    \end{equation}
    However, in this case, by \cite[Prop. 8.3.12]{LeStum}, we have a Frobenius-equivariant quasi-isomorphism
    \begin{equation}\label{equ:Frobenius compatibility-II}
      \rR\Gamma_{\rig}(\calX/K)\simeq \rR\Gamma(]X_0[_{\overline Y_0},j^{\dagger}\rD\rR(\calO_{\overline Y_0},\rd)).
    \end{equation}
    Since $\iota_{\CCtilde^I}$ is given by scalar extension, \eqref{equ:Frobenius compatibility-I} and \eqref{equ:Frobenius compatibility-II} imply that it is Frobenius-equivariant.
\end{proof}

\subsection{Proof of Theorem \ref{thm:P-cohomology}}

We prove that the comparison morphism $\iota_{\bP}:\rR\Gamma_{\rig}(\calX/K)\otimes_K\bP(Z)\to\rR\Gamma^{\dagger}(X_Z,\bP)$ is a quasi-isomorphism. By construction we may work Zariski-locally on $\frakX_{0,\Zar}$, so we assume embeddings $\frakX_0\to\frakY_0\hookrightarrow\overline \frakY_0$ as in Definition \ref{dfn:P-cohomology}. Then
\[\iota^*\rL\eta_t\rR\nu_*\CCtilde^I\simeq \rD\rR(\calO_{\wtX_Z^I},\rd)\]
and 
\[\iota^*\rL\eta_t\rR\nu_*\BBdRp\simeq \rD\rR(\calO_{\wtX_Z},\rd).\]
As in the Introduction (Remark on dagger structures), the only ``integral dagger variety'' we treat is that of this section: $]\calX[_{\overline Y_0}$ is a dagger space over $K$ \cite{GK00}, realized as the intersection of strict neighborhoods in $]\calY[_{\overline Y_0}$. After further localization we may assume $]\calX[_{\overline Y_0}$ is affine with dagger structure $]\calX[_{\overline Y_0} = \lim V_n$ for smooth affinoids $V_n=\Spa(A_n,A_n^+)\subset \overline Y_0$; this exhaustion by strict neighborhoods $V_n$ yields the colimit descriptions below. We then have quasi-isomorphisms
\begin{equation*}
    \begin{split}
        & \rR\Gamma_{\rig}(\calX/K)\simeq \colim_{n\to\infty}\rR\Gamma(V_n,\rD\rR(\calO_{V_n},\rd)) = \colim_{n\to\infty}\rD\rR(A_n,\rd),\\
        & \rR\Gamma^{\dagger}(X_Z,\CCtilde^I) \simeq \colim_{n\to\infty}\rR\Gamma(\widetilde V_{n,Z}^I,\rD\rR(\calO_{\widetilde V_Z^I},\rd)) = \colim_{n\to\infty}\rD\rR(A_n\widehat \otimes_K\CCtilde^I(Z),\rd),\\
        &\rR\Gamma^{\dagger}(X_Z,\BBdRp) \simeq \colim_{n\to\infty}\rR\Gamma(\widetilde V_{n,Z},\rD\rR(\calO_{\widetilde V_{n,Z}},\rd)) = \colim_{n\to\infty}\rD\rR(A_n\widehat \otimes_K\BBdRp(Z),\rd).
    \end{split}
\end{equation*}

Thus it suffices to prove the following.
\begin{prop}\label{prop:iota_P is quasi-isomorphism}
    For any $k\geq 0$, the natural map
    \[(\colim_{n\to \infty}\rH^k(\rD\rR(A_n,\rd)))\otimes_K\bP(Z)\to\colim_{n\to\infty}(\rH^k(\rD\rR(A_n\widehat \otimes_K\bP(Z),\rd)))\]
    is an isomorphism.
\end{prop}
\begin{proof}
    Note that $\{\rD\rR(A_n,\rd)\}_{n\geq 0}$ is a direct system of complexes of Banach $K$-modules of degrees in $[0,d]$ such that the translation map $\rD\rR(A_n,\rd)\to \rD\rR(A_{n+1},\rd)$ is continuous and that
    \[\colim_{n\to \infty}\rH^k(\rD\rR(A_n,\rd)) = \rH^k(\colim_{n\to \infty}\rD\rR(A_n,\rd)).\]
    By \cite{Ked06}, rigid cohomology is finite-dimensional over $K$, so $\rH^k(\colim_{n\to \infty}\rD\rR(A_n,\rd))$ is finite-dimensional. The result is then a special case of Proposition \ref{prop:technique proposition}.
\end{proof}

\begin{prop}\label{prop:technique proposition}
    Let $\{(C_n^{\bullet},\partial_n)\}_{n\geq 0}$ be a direct system of uniformly bounded complexes of Banach $K$-modules and let $(C^{\bullet},\partial):=\colim_n(C_n^{\bullet},\partial_n)$. Suppose that:
    \begin{enumerate}
        \item[(1)] The differential $\partial_n^i:C_n^{i-1}\to C_n^i$ is continuous for all $n$ and all $i$,
        \item[(2)] The map $C^{\bullet}_n\to C_{n+1}^{\bullet}$ is continuous for all $n$, and
        \item[(3)] The cohomology group $\rH^{k}(C^{\bullet})$ is finite-dimensional over $K$ for all $k$.
    \end{enumerate}
    For any Banach $K$-module $V$, define 
    \[(C_n^{\bullet}(V),\partial_n):= (C_n^{\bullet}\widehat \otimes_KV,\partial\otimes\id) \text{ and }(C^{\bullet}(V),\partial):=\colim_n(C^{\bullet}_n(V),\partial_n).\]
    Then the natural morphism 
    \[\rH^k(C^{\bullet})\otimes_KV\to \rH^k(C^{\bullet}(V))\]
    is an isomorphism for any $k$.
\end{prop}
\begin{proof}
    Any continuous morphism $U\to F$ with $F$ finite-dimensional over $K$ admits a continuous section, so $F$ is a direct summand of $U$. The key step is the following lemma.
    \begin{lem}\label{lem:key for technique proposition}
        Let $\{\partial_n: C_n\to D_n\}_{n\geq 0}$ be a direct system of continuous morphisms of Banach $K$-modules with continuous translation maps. Suppose that the direct limit 
        \[\partial:=\colim_n\partial_n:C:=\colim_nC_n\to \colim_nD_n=:D\]
        has a finite-dimensional cokernel $H$. For any Banach $K$-module $V$, define 
        \[C_n(V):=C_n\widehat \otimes_KV,~D_n(V):=D_n\widehat \otimes_KV,~C(V):=\colim_nC_n(V),\text{ and }D(V):=\colim_nD_n(V).\]
        Then the natural map 
        \[H\otimes_KV\to H(V):=\Coker(C(V)\to D(V))\]
        is an isomorphism.
    \end{lem}
    We first derive the proposition from the lemma, and then prove the lemma. After shifting the complexes, we may assume that for each $n\geq 0$, $C_n^{\bullet}$ is a complex concentrated in degrees $[0,d]$ for some $d\geq 0$; that is, we require $C_n^i= 0$ for $i\notin[0,d]$ for all $n$. Applying the lemma to the direct system $\{\partial_n:C_n^{d-1}\to C_n^d\}$, we conclude the desired isomorphism:
    \[\rH^d(C^{\bullet})\otimes_KV\xrightarrow{\cong}\rH^d(C^{\bullet}(V)).\]
    Now, let $Z_n^{d-1} = \Ker(\partial:C_n^{d-1}\to C_n^d)$. This is a Banach $K$-module (since $\partial$ is continuous) such that $Z:=\colim_nZ_n$ is the kernel of $C^{d-1}\to C^d$. Replace the complexes $C_n^{\bullet}$ and $C$ with $D_n^{\bullet}$, where 
    \[D^{\bullet} = \left\{
       \begin{array}{rcl}
           C^{\bullet}, & \bullet\leq d-2 \\
           Z^{d-1}, & \bullet = d-1
       \end{array}
    \right.\]
    and $D:=\colim_nD_n$ (equipped with the induced differential maps), respectively. Then $D_n^{\bullet}$ is a complex in degrees $[0,d-1]$ such that for any Banach $K$-module and any $k\in[0,d-1]$,
    \[\rH^k(D^{\bullet}(V)) = \rH^k(C^{\bullet}(V)).\]
    Applying Lemma \ref{lem:key for technique proposition} again we obtain
    \[\rH^{d-1}(C^{\bullet})\otimes_KV = \rH^{d-1}(D^{\bullet})\otimes_KV \xrightarrow{\cong}\rH^{d-1}(D^{\bullet}(V)) = \rH^{d-1}(C^{\bullet}(V)).\]
    Iterating this procedure proves the proposition.
\end{proof}

\begin{proof}[Proof of Lemma \ref{lem:key for technique proposition}]
    Fix an $\calO_K$-lattice $V^+$ of $V$. For any $n\geq 0$, let $K_n:=\Ker(D_n\to H)$. Then $K_n$ is a Banach $K$-module with $D_n = K_n\oplus H$ such that $\colim_nK_n = K:=\Ker(D\to H)$. By replacing $D_n$ with $K_n$, we may assume without loss of generality that $H=0$, thereby reducing the problem to showing that $H(V)=0$.

    For each $n\geq 0$, let $j_n:D_n\to D$ be the natural map. We have
    \[\cup_{k\in\bZ}p^kD_n^+ = D_n = j_n^{-1}(D) = j_n^{-1}(\partial C) = j_n^{-1}(\cup_{m\geq n}\partial_mC_m)\]
    because $H = 0$ and, a fortiori, $\partial$ is surjective. For each $n$, fix $\calO_K$-lattices $C_n^+\subset C_n$ and $D_n^+\subset D_n$ such that $\partial_n(C_n^+)\subset D_n^+$. Then we have
    \[D_n = j_n^{-1}(\cup_{m\geq n}\cup_{l\in\bZ}\partial_m(p^lC_m^+)) = \cup_{m\geq n}\cup_{l\in\bZ}j_n^{-1}(\partial_m(p^lC_m^+)) = \cup_{m\geq n}\cup_{l\in\bZ}\overline{j_n^{-1}(\partial_m(p^lC_m^+))},\]
    where $\overline{j_n^{-1}(\partial_m(p^lC_m^+))}$ denotes the closure of $j_n^{-1}(\partial_m(p^lC_m^+))$ in $D_n$. To lift elements from the image of $\partial$ in the direct limit back to a specific stage: by Baire's category theorem, there exist $N\geq n$ and $k,l\in\bZ$ such that
    \[p^kD_n^+\subset \overline{j_n^{-1}(\partial_N(p^lC_N^+))}.\]
    Replacing $C_N^+, D_N^+, C_n^+$, and $D_n^+$ with $p^lC_N^+, p^lD_N^+, p^kC_n^+$, and $p^kD_n^+$ respectively, we may assume 
    \[D_n^+\subset \overline{j_n^{-1}(\partial_NC_N^+)}.\]
       
    We claim $D_n^+\subset j_n^{-1}(\partial_NC_N^+)$. Indeed, for any $x\in D_n^+$, we shall construct sequences $x_i\in D_n^+$ and $y_i\in C_N^+$ such that for any $j\geq 0$, we have
    \[x - (\sum_{i=0}^jp^ix_i)\in p^{j+1}D_n^+ \text{ and }j_n(x_i) = \partial_N(y_i).\]
    Using this, $x = \sum_{i\geq 0}p^ix_i$ and $y:=\sum_{i\geq 0}p^iy_i$ are well-defined elements in $D_n^+$ and $C_N^+$, respectively, such that
    \[j_n(x) = \sum_{i\geq 0}p^ij_n(x_i) = \partial_N(\sum_{i\geq 0}p^iy_i) = \partial_N(y),\]
    yielding that $x\in j_n^{-1}(\partial_N(C_N^+))$ as desired.

    We now construct $(x_k,y_k)$ by induction on $k$. For $k = 0$, since $x$ belongs to the closure of $j_n^{-1}(\partial_N(C_N^+))$, there exists some $(x_0,y_0)\in C_n^+\times C_N^+$ such that $x - x_0\in pD_n^+$ and $j_n(x_0) = \partial_N(y_0)$, as required. Suppose that we have obtained the required $(x_i,y_i)$ for $0\leq i\leq k-1$. Then, by assumption, there exists some $z\in D_n^+$ such that $x - \sum_{i=0}^{k-1}p^ix_i = p^kz$. Similar to the construction of $(x_0,y_0)$, there exists some $(x_k,y_k)\in D_n^+\times C_N^+$ such that 
    \[z-x_k\in pD_n^+ \text{ and }j_n(x_k) = \partial_N(y_k).\]
    We can then finish the induction argument (and complete the proof of the claim) by noting that
    \[x - \sum_{i=0}^{k}p^ix_i = p^k(z-x_k)\in p^{k+1}D_n^+.\]

    By the claim above, $D_n^+\subset j_n^{-1}(\partial_NC_N^+)$. Since $D_n^+$ is a topologically free $\calO_K$-module, there is a continuous map $f:D_n^+\to C_N^+$ such that $j_n(x) = \partial_N(f(x))$ for all $x\in D_n^+$, i.e., the following diagram commutes:
    \begin{equation*}   
        \xymatrix@C=0.5cm{C_n^+\ar[rr]\ar[d]^{\partial_n}&&C_N^+\ar[d]^{\partial_N}\\
           D_n^+\ar[rr]\ar@{.>}[rru]^f&&D_N^+.}
    \end{equation*}
    Applying $-\widehat \otimes_{\calO_K}V^+$ and inverting $p$, we obtain a commutative diagram:
     \begin{equation*} 
        \xymatrix@C=0.5cm{C_n(V)\ar@{=}[r]&C_n\widehat \otimes_KV\ar[rr]\ar[d]^{\partial_n}&&C_N\widehat \otimes_KV\ar[d]^{\partial_N}\ar@{=}[r]&C_N(V)\\
         D_n(V)\ar@{=}[r]&D_n\widehat \otimes_KV\ar[rr]\ar@{.>}[rru]^f&&D_N\widehat \otimes_KV\ar@{=}[r]&D_N(V).}
    \end{equation*}
    Thus, we have 
    \[j_n(D_n(V)) = j_N(\partial_N(f(D_n(V))))\subset j_N(\partial_N(C_N(V)))\subset \partial(C(V)),\]
    yielding that $D(V)\subset\partial(C(V))$ as desired.
\end{proof}

\end{document}